\documentclass[11pt]{article}

\usepackage{amsmath} 
\usepackage{amsthm} 
\usepackage{amsfonts} 
\usepackage{amssymb} 
\usepackage{stmaryrd} 
\usepackage{mathtools} 
\usepackage{thm-restate} 

\usepackage[style=alphabetic, maxnames=100, maxalphanames=100]{biblatex} 
\DeclareFieldFormat{postnote}{\mknormrange{#1}} 

\usepackage{xifthen} 
\usepackage{xstring} 

\usepackage{xcolor} 
\usepackage{bm} 
\usepackage{bbm} 

\usepackage{array} 
\usepackage{multirow} 
\usepackage[inline]{enumitem} 

\usepackage{microtype} 
\usepackage{lmodern} 
\usepackage{graphicx} 

\usepackage{pdflscape} 
\usepackage{afterpage} 

\usepackage{tikz} 
\usetikzlibrary{cd} 
\usetikzlibrary{calc} 
\usetikzlibrary{decorations.pathreplacing} 
\usepackage[%
  algoruled, 
  algosection, 
  vlined, 
  fillcomment, 
  linesnumbered, 
  resetcount
]{algorithm2e} 
\SetKw{Break}{break}
\SetKw{Continue}{continue}
\SetKw{Failure}{Failure}

\makeatletter
\SetKwIF{forwhile@dummy@if}{forwhile@dummy@elseif}{forwhile@dummy@else}{for}{do}{while}{else}{endif}
\let\uFor\uforwhile@dummy@if
\let\uWhile\uforwhile@dummy@elseif
\makeatother

\makeatletter
\SetKwIF{invisiblebegin@dummy@if}{invisiblebegin@dummy@elseif}{invisiblebegin@dummy@else}{if}{then}{else if}{}{}
\let\InvisibleBegin\invisiblebegin@dummy@else
\let\uInvisibleBegin\uinvisiblebegin@dummy@else
\makeatother

\newcounter{algosplit}

\usepackage[in]{fullpage} 

\usepackage[hypcap=true]{subcaption} 

\usepackage[%
  plainpages=false,
  pdfpagelabels,
  colorlinks,
  citecolor=black,
  linkcolor=black,
  urlcolor=black,
  filecolor=black,
  bookmarksopen=false
]{hyperref} 
\usepackage[all]{hypcap} 
\hypersetup{hypertexnames=false}

\makeatletter
\newcommand*{\newletterthm@internal}{}
\newcommand*{\newletterthm}[1]{%
  \def\newletterthm@name{#1}%
  \renewcommand*{\newletterthm@internal}[1][]{%
    \ifthenelse{\isempty{##1}}{%
      \expandafter\expandafter\expandafter\newtheorem%
      \expandafter\expandafter\expandafter{%
        \expandafter\newletterthm@name%
        \expandafter}%
      \expandafter{%
        \newletterthm@text}%
      \expandafter\renewcommand%
      \expandafter*%
      \expandafter{%
        \csname the#1\endcsname}{\Alph{#1}}%
    }{%
      \expandafter\expandafter\expandafter\newtheorem%
      \expandafter\expandafter\expandafter{%
        \expandafter\newletterthm@name%
        \expandafter}%
      \expandafter{%
        \newletterthm@text}[##1]%
    \expandafter\renewcommand%
        \expandafter*%
        \expandafter{%
          \csname the#1\endcsname}{\csname the##1\endcsname.\Alph{#1}}%
    }%
  }%
  \newletterthm@newthm%
}

\newcommand*{\newletterthm@newthm}[2][]{%
  \ifthenelse{\isempty{#1}}{%
    \def\newletterthm@text{#2}%
    \newletterthm@internal%
  }{%
    \expandafter\newtheorem\expandafter{\newletterthm@name}[#1]{#2}%
  }%
}
\makeatother

\newtheoremstyle{thmstyle}
  {\medskipamount}
  {\smallskipamount}
  {\slshape}
  {0pt}
  {\bfseries}
  {.}
  { }
  {\thmname{#1}\thmnumber{ #2}{\normalfont\thmnote{ (#3)}}}
  
\newtheoremstyle{plainstyle}
  {\medskipamount}
  {\smallskipamount}
  {\rmfamily}
  {0pt}
  {\bfseries}
  {.}
  { }
  {\thmname{#1}\thmnumber{ #2}{\normalfont\thmnote{ (#3)}}}

\theoremstyle{thmstyle}
\newtheorem{theorem}{Theorem}[section]
\newtheorem{lemma}[theorem]{Lemma}
\newtheorem{corollary}[theorem]{Corollary}
\newtheorem{proposition}[theorem]{Proposition}

\newletterthm{theoremLet}{Theorem}
\newletterthm{problemLet}{Problem}

\newtheorem{claim}{Claim}[theorem]

\theoremstyle{plainstyle}
\newtheorem{definition}[theorem]{Definition}

\newtheorem{notation}[theorem]{Notation}
\newtheorem{remark}[theorem]{Remark}
\newtheorem{discussion}[theorem]{Discussion}

\newenvironment{proofof}[1]{\begin{proof}[Proof of #1.]}{\end{proof}}

\makeatletter
\def\refdescformat#1{%
  \phantomsection%
  \let\oldlabel\label%
  \let\label\@gobble%
  \edef\@currentlabel{#1}
  \let\label\oldlabel%
  #1:
}
\makeatother

\newlist{refdesc}{description}{1}
\setlist[refdesc]{format={\refdescformat}}

\newlist{enumdef}{enumerate}{1}
\setlist[enumdef]{before={\leavevmode}, label={\arabic*.}, ref={\thetheorem.\arabic*}}

\setlist[enumerate]{label={\roman*.}, ref={(\roman*)}} 

\let\emph\textit

\makeatletter
\newcommand{\notoc@internal}[2][]{}
\newcommand{\notoc}[1]{
  \renewcommand{\notoc@internal}[2][]{%
    \let\old@addtocontents\addtocontents%
    \let\addtocontents\@gobbletwo%
    \ifthenelse{\isempty{##1}}{%
      #1{##2}%
    }{%
      #1[##1]{##2}%
    }%
    \let\addtocontents\old@addtocontents%
  }%
  \notoc@internal%
}
\makeatother

\makeatletter
\newcommand{\biggg}{\bBigg@\thr@@}
\def\bigggl{\mathopen\biggg}
\def\bigggm{\mathrel\biggg}
\def\bigggr{\mathclose\biggg}
\newcommand{\Biggg}{\bBigg@{3.5}}

\newcommand{\bigggg}{\bBigg@{4}}

\newcommand{\Bigggg}{\bBigg@{4.5}}

\newcommand{\biggggg}{\bBigg@{5}}

\newcommand{\Biggggg}{\bBigg@{5.5}}

\makeatother

\numberwithin{equation}{section} 

\makeatletter
\newcommand{\combinesymbols}[3][\mathord]{#1{\mathpalette\combinesymbols@paletted{{#2}{#3}}}}
\newcommand{\combinesymbols@paletted}[2]{\combinesymbols@internal{#1}#2}
\newcommand{\combinesymbols@internal}[3]{\ooalign{\hss$\m@th #1#2$\hss\cr\hss$\m@th #1#3$\hss}}
\makeatother

\let\epsilon\varepsilon

\newcommand{\rn}{\bm}
\newcommand{\df}{\stackrel{\text{def}}{=}}
\newcommand{\place}{\mathord{-}}
\newcommand{\symdiff}{\mathbin{\triangle}}

\newcommand{\rest}{\mathord{\vert}}
\newcommand{\floor}[1]{\lfloor#1\rfloor}
\newcommand{\ceil}[1]{\lceil#1\rceil}
\newcommand{\Floor}[1]{\left\lfloor#1\right\rfloor}
\newcommand{\Ceil}[1]{\left\lceil#1\right\rceil}

\newcommand{\conc}{\mathbin{{}^\smallfrown}}

\newcommand{\assign}{\leftarrow}

\DeclareMathOperator{\im}{im}

\DeclareMathOperator{\atom}{atom}
\DeclareMathOperator{\VC}{VC}
\DeclareMathOperator{\Lit}{Lit}
\DeclareMathOperator{\OP}{OP}
\DeclareMathOperator{\len}{len}
\DeclareMathOperator{\UC}{UC}
\DeclareMathOperator{\abs}{abs}
\DeclareMathOperator{\Root}{root}
\DeclareMathOperator{\samp}{samp}
\DeclareMathOperator{\adj}{adj}
\DeclareMathOperator{\Sum}{sum}
\DeclareMathOperator{\nonequi}{nonequi}
\DeclareMathOperator{\equi}{equi}
\DeclareMathOperator{\equirand}{equirand}
\DeclareMathOperator{\equiPRG}{equiPRG}
\DeclareMathOperator{\good}{good}
\DeclareMathOperator{\num}{num}
\DeclareMathOperator{\den}{den}

\newcommand{\EE}{\mathbb{E}}
\newcommand{\FF}{\mathbb{F}}
\newcommand{\NN}{\mathbb{N}}
\newcommand{\PP}{\mathbb{P}}
\newcommand{\QQ}{\mathbb{Q}}
\newcommand{\RR}{\mathbb{R}}
\newcommand{\ZZ}{\mathbb{Z}}

\newcommand{\One}{\mathbbm{1}}

\newcommand{\cF}{\mathcal{F}}
\newcommand{\cH}{\mathcal{H}}

\newcommand{\cO}{\mathcal{O}}
\newcommand{\cP}{\mathcal{P}}
\newcommand{\cQ}{\mathcal{Q}}

\def\Szemeredi{Szemer\'{e}di}
\def\Lovasz{Lov\'{a}sz}

\def\Cervonenkis{\v{C}ervonenkis}

\def\Turan{Tur\'{a}n}
\def\Rodl{R\"{o}dl}
\def\Erdos{Erd\H{o}s}

\newcommand{\WARNING}[2][Warning]{
  \typeout{^^J#1 on line \the\inputlineno: #2^^J}
  \textbf{#1 on line \the\inputlineno: #2}
}

\title{Stable Regularity Lemmas: Efficient Algorithms\\
         and Essentially Tight Littlestone Bounds}
\author{%
  Leonardo N.~Coregliano\thanks{This study was financed, in part, by the São Paulo Research Foundation (FAPESP), Brasil.
    Process Numbers \#2025/10990-5 and \#2026/16585-8.}
  \and
  Fernando G.~Jeronimo
}
\date{\today}

\begin{document}
\maketitle

\begin{abstract}
  Regularity lemmas are fundamental in theoretical computer science and combinatorics, and have found a vast range of
  applications including property testing, approximation algorithms, and several structural results. More recently, important
  connections to model theory and differentially private learning were established. In this context, the stable regularity lemma
  of Malliaris and Shelah strengthens Szemerédi regularity lemma for stable graphs by producing partitions with no exceptional
  pairs. The more common way of characterizing stable graphs is to bound the size of half-graphs, but a classic Ramsey-like
  argument says that this is equivalent to bounding the Littlestone dimension.

  In this paper, we determine the precise asymptotics of the number of parts of stable regularity equipartitions in terms of the
  Littlestone dimension. Namely, we show that every graph $G$ of Littlestone dimension $\Lit(G)\leq\ell$ has a stable regular
  equipartition (into excellent sets) with $(1 + o_{\epsilon\to 0,\ell}(1))\cdot\epsilon^{-\ell-1}$ parts and in the other
  direction, for every $\ell\in\NN_+$, we produce an infinite family of graphs, all of Littlestone dimension $\ell$, whose
  stable regularity equipartitions (into good sets) must have size at least $(1 + o_{\epsilon\to
    0,\ell}(1))\cdot\epsilon^{-\ell-1}$.

  Dropping the equitability condition, we determine the asymptotics of (not necessarily equitable) stable regularity partitions
  up to a multiplicative logarithmic factor. Namely, we show that every graph $G$ with $\Lit(G)\leq\ell$ has a stable regular
  partition (into excellent sets) with $(1 + o_{\epsilon\to 0,\ell}(1))\cdot\epsilon^{-\ell}\cdot\ln(1/\epsilon)$ parts and in
  the other direction, for every $\ell\in\NN_+$, we produce an infinite family of graphs, all of Littlestone dimension $\ell$,
  whose stable regularity partitions (into good sets) must have size at least $(1 + o_{\epsilon\to
    0,\ell}(1))\cdot\epsilon^{-\ell}$.

  Finally, we also show that stable regularity partitions with these number of parts can be obtained algorithmically efficiently
  in an approximation scheme fashion. Namely, replacing the $o_{\epsilon\to 0,\ell}(1)$ term in the bounds above by a constant
  $c > 0$, we obtain randomized $O_{c,\epsilon,\ell}(n\cdot\log(n))$-time algorithms for partitions/equipartitions into good
  sets, a deterministic $O_{c,\epsilon,\ell}(n^2)$-time algorithm for partitions into good sets, a deterministic
  $O_{c,\epsilon,\ell}(n^6)$-time algorithm for equipartitions into good sets, a deterministic $O_{c,\ell,\epsilon}(1)\cdot
  n^{O(\ell\cdot 2^{2\cdot\ell+4})}$-time algorithm for partitions/equipartitions into excellent sets, and
  $O_{c,\epsilon,\ell}(\log(n+1))$-space algorithms for partitions/equipartitions into good/excellent sets.

  The main new ingredient for algorithmically efficiently handling excellent sets is the notion of atomicity, a finite-witness
  strengthening of excellence. In fact, all results above about excellent partitions are derived as corollaries of results
  proved for atomic partitions instead.
\end{abstract}

\clearpage

\tableofcontents

\clearpage

\section{Introduction}

Regularity lemmas are one of the main ways of turning large objects into finite-dimensional approximations. In graph-theoretic
language, a regularity lemma partitions the vertex set of a large graph into a bounded number of parts so that the edge relation
between most pairs of parts is quasirandom, or even nearly constant. From the point of view of theoretical computer science,
such a partition is a small sketch of the input graph: it is a constant-size object, for fixed accuracy. This viewpoint has been
central in the use of regularity methods in property testing, approximation algorithms for CSPs, graph limits, and counting and
removal lemmas, see, for example,~\cite{KSSE02,GGR98,AKK99,AFKK03,AS08}.

The regularity method originated in \Szemeredi's proof of his theorem on arithmetic progressions~\cite{Sz75} and was
subsequently formulated as the graph regularity lemma~\cite{Sz78}. Roughly speaking, \Szemeredi's Regularity Lemma states that
for every $\epsilon>0$, every sufficiently large graph admits an equitable partition into at most
$m_{\operatorname{Sz}}(\epsilon)$ parts such that all but an $\epsilon$-fraction of pairs of parts are $\epsilon$-regular. This
is an extremely robust approximation theorem, but it has two well-known costs. First, the number of parts is enormous: the best
possible bound is of exponential tower of height a (fixed) power of $1/\epsilon$~\cite{Gow97} (see also the subsequent lower
bound by Moshkovitz and Shapira~\cite{MS16}). Second, the lemma necessarily allows (small amount) of irregular pairs of parts
for which no regularity conclusion is guaranteed.

A major theme in the development of regularity methods has been to understand which relaxations or structural assumptions lead
to smaller partitions, better algorithms, or stronger guarantees. The Frieze--Kannan weak regularity lemma replaces \Szemeredi's
regularity by approximation in cut norm and obtains only $m_{\operatorname{FK}}(\epsilon)=2^{O(\epsilon^{-2})}$
parts~\cite{FK96,FK99a}. This weaker notion is still strong enough for many algorithmic purposes, including additive
approximation of CSPs. Moreover, weak regularity admits efficient algorithmic versions, with subsequent work giving fast
deterministic constructions~\cite{FK99b,FLZ19}. Spectral and matrix-decomposition viewpoints give related algorithmic regularity
results, including for low-threshold-rank graphs and sparse graph classes~\cite{OT15,BV22}. The exponential dependence on
$1/\epsilon$ is necessary in general~\cite{CF12}.

A route to polynomial-size partitions is to impose bounded VC dimension\footnote{Bounded VC dimension is the combinatorial
counterpart of NIP in model theory.}. For graphs whose neighborhood set system has VC dimension at most $d$, work of
Alon--Fischer--Newman~\cite{AFN07} and \Lovasz--Szegedy~\cite{LS10} gives much stronger regularity partitions with polynomial
dependence on $1/\epsilon$, namely $m_{\VC}(\epsilon) = (1/\epsilon)^{O_d(1)}$. Fox--Pach--Suk made major progress on the
Ramsey-theoretic \Erdos--Hajnal problem for bounded-VC graphs~\cite{FPS19} using a VC regularity lemma and the polynomial number
of parts was crucial. The conjecture for this class of graphs was recently resolved by Nguyen, Scott, and Seymour~\cite{NSS25}.
Nevertheless, VC-type regularity lemmas still require irregular pairs in general.

The stable regularity lemma of Malliaris and Shelah~\cite{MS14} gives a qualitatively stronger conclusion under a stronger
structural assumption. In this paper we measure stability by Littlestone dimension\footnote{For graphs, finite Littlestone
dimension is qualitatively equivalent to excluding arbitrarily long order patterns, or half-graphs; see also
Hodges~\cite[Lemma~6.7.9]{Hod93} stated as Lemma~\ref{lem:Hodges} below.}, the online-learning dimension introduced
in~\cite{Lit88}. The Malliaris--Shelah theorem says, informally, that if the edge relation is stable, then every sufficiently
large finite graph admits an equitable partition with polynomially many parts such that every pair of parts is regular, and in
fact almost homogeneous: each pair has density close either to $0$ or to $1$. Thus, the half-graph obstruction is not merely one
example showing that \Szemeredi's regularity needs irregular pairs; in a precise sense, it is the only obstruction to
eliminating them. Malliaris and Shelah gives an elegant expository account of this proof and its model-theoretic context
in~\cite{MS21}.

The stable regularity lemma achieves its strength through a pair of structural properties called goodness and excellence. Very
roughly, a set is good if every vertex has a coherent majority opinion about it, and it is excellent if it has coherent majority
behavior with respect to every good test set. Once a partition into good/excellent sets has been constructed, all pairs of parts
are automatically almost homogeneous (see Lemma~\ref{lem:good->hom} below). However, excellence in general is strictly stronger
than goodness which in turn is strictly stronger than goodness. Furthermore, for excellence, an algorithmic difficulty appears:
excellence is a non-local condition, since verifying it a priori requires quantifying over all good sets. The Malliaris--Shelah
proof is a beautiful existence argument, but it does not directly give an efficient procedure for finding an excellent partition
of a given input graph. In contrast to the other regularity lemmas discussed above for which efficient algorithms were known, no
efficient algorithm was known for it since the seminal work of~\cite{MS14}.

This leads to the main algorithmic question of the paper.

\begin{center}
  \emph{Does the stable regularity lemma admit an efficient algorithm?}
\end{center}

We answer this question in the affirmative by establishing efficient algorithms under various computational models: randomized,
deterministic, and space bounded. A new conceptual ingredient is the notion of atomicity which strengthens the corresponding
notions of good and excellent parts in the stable regularity lemmas~\cite{MS14} (see their formal definitions in
Section~\ref{sec:prelim}), which we will discuss in more detail shortly in Section~\ref{subsec:strat}.

\begin{theorem}[Algorithmic stable regularity, informal]
  For every integer $\ell\geq 1$, every $\epsilon>0$ and every $c > 0$, there is an algorithm which, given access to an
  $n$-vertex graph $G$ with $\Lit(G)\leq\ell$ in one of the computation models considered below, computes a stable regular
  partition of $G$ with no irregular pairs and with number of parts at most
  \begin{gather*}
    \begin{dcases*}
      (1+c)\cdot\epsilon^{-\ell}\cdot\ln(1/\epsilon), & for not necessarily equitable partitions,\\
      (1+c)\cdot\epsilon^{-\ell-1}, & for equitable partitions.
    \end{dcases*}
  \end{gather*}
  More precisely, for fixed $c$, $\ell$ and $\epsilon$, we obtain the following algorithmic forms:
  \begin{itemize}
  \item Randomized Monte Carlo algorithm computing a stable partition (equitable or not) into good parts that runs in
   time $O_{c,\ell,\epsilon}(n\cdot\log n)$.
  \item Deterministic algorithms computing a stable partition into good parts with run time $O_{c,\ell,\epsilon}(n^2)$ in the
    non-equitable case and $O_{c,\ell,\epsilon}(n^6)$ in the equitable case.
  \item Deterministic algorithms for a stable partition into atomic (hence excellent) parts, (equitable or not) that runs in
    time $O_{c,\ell,\epsilon}(1)\cdot n^{O(\ell\cdot 2^{2\cdot\ell+4})}$ (the exponent $O$ hides an absolute constant).
  \item A $O_{c,\ell,\epsilon}(\log(n+1))$-space algorithm for a stable partition (equitable or not) into atomic parts.
  \end{itemize}
\end{theorem}

Having obtained an algorithmic stable regularity lemma with $(1+o_{\epsilon\to 0,\ell}(1))\cdot\epsilon^{-\ell-1}$ parts in the
equitable case and $(1+o_{\epsilon\to 0,\ell}(1))\cdot\epsilon^{-\ell}\cdot\ln(1/\epsilon)$ in the non-equitable case, it is
natural to ask whether these bounds in the number of parts are tight.

\begin{center}
  \emph{What is the polynomial dependence of the number of parts on $1/\epsilon$ in stable regularity?}
\end{center}

Our second main result is that in the equitable case, the number of parts above is tight:
\begin{theorem}[Equitable stable regularity lower bound, informal]
  For every integer $\ell\geq 1$ and all sufficiently small $\epsilon>0$, there are infinitely many graphs $G$ with
  $\Lit(G)=\ell$ such that every equitable partition into $\epsilon$-good sets requires
  \begin{gather*}
    \bigl(1-o_{\lvert G\rvert\to\infty,\epsilon\to 0,\ell}(1)\bigr)\cdot\epsilon^{-\ell-1}
  \end{gather*}
  parts. In particular, the dependence on $1/\epsilon$ in the algorithmic equitable stable regularity lemma is tight up to lower
  order terms.
\end{theorem}

For the non-equitable case, we prove almost tightness, missing only a logarithmic factor:
\begin{theorem}[Non-equitable stable regularity lower bound, informal]
  For every integer $\ell\geq 1$ and all sufficiently small $\epsilon>0$, there are infinitely many graphs $G$ with
  $\Lit(G)=\ell$ such that every (not necessarily equitable) partition into $\epsilon$-good sets requires
  \begin{gather*}
    \bigl(1-o_{\lvert G\rvert\to\infty,\epsilon\to 0,\ell}(1)\bigr)\cdot\epsilon^{-\ell}
  \end{gather*}
  parts. In particular, the dependence on $1/\epsilon$ in the algorithmic (not necessarily equitable) stable regularity lemma is
  tight up to a $\ln(1/\epsilon)$ multiplicative factor.
\end{theorem}

\paragraph{Related work.} There is a long line of work on algorithmic regularity. Constructive forms of \Szemeredi's regularity
go back to Alon--Duke--Lefmann--\Rodl--Yuster and Frieze--Kannan~\cite{ADLRY94,FK99b}. Weak regularity has particularly
efficient algorithmic versions and has become a standard tool for graph and matrix approximation~\cite{FK96,FK99a,FLZ19};
spectral, low-threshold-rank, and matrix-decomposition variants give further algorithmic regularity tools~\cite{OT15,BV22}.
Regularity methods are also central in property testing and approximation algorithms~\cite{GGR98,AKK99,AFKK03,AS08}.

On the model-theoretic side, Malliaris and Shelah introduced stable regularity and the use of excellent sets in finite stable
graphs~\cite{MS14}; their later expository note~\cite{MS21} gives a helpful overview of the proof and the role of half-graphs.
More recently, Malliaris and Moran used methods from online learning to obtain new existence proofs for excellent sets, valid in
a wider range of parameters~\cite{MM24}, and developed an analogue of Shelah's unstable formula theorem~\cite{MM25}. Related
recent work studies finitary saturation phenomena for stable graphs and Littlestone classes~\cite{MMM25}. These works reveal
deep connections between stability, majority structure, and online learning. The present paper is complementary: it addresses
the computational problem of constructing, from a given finite input graph, a stable regular partition with the quantitative
bounds above.

Besides its original importance in online learning, Littlestone dimension was recently shown to govern differentially private
PAC learning in~\cite{ABLMM22}.


\section{Preliminaries}\label{sec:prelim}

\begin{notation}
  The set of non-negative integers is denoted $\NN$ and we set $\NN_+\df\NN\setminus\{0\}$. For $k\in\NN$ and a set $V$, we let
  $\cP(V)\df\{A\subseteq V\}$ be the set of all subsets of $V$ and we let $\binom{V}{k}\df\{A\subseteq V \mid \lvert A\rvert =
  k\}$ be the set of subsets of $V$ of cardinality $k$, we let $V^{< k}\df\bigcup_{i=0}^{k-1} V^i$ be the set of strings of
  length less than $k$ over the alphabet $V$.

  A graph\footnote{We will mostly be concerned with \emph{finite graphs}, but some of the results readily extend to infinite
  graphs.} is a pair $G=(V,E)$ where $V$ is a set of \emph{vertices} and $E\subseteq\binom{V}{2}$ is a set of \emph{edges}. We
  use the shorthand notation $V(G)\df V$ and $E(G)\df E$. We also let $\vec{E}(G)\df\{(u,v)\in V(G)^2 \mid \{u,v\}\in E(G)\}$ be
  the ordered version of the set of edges. We set $\vec{E}^1(G) \df \vec{E}(G)$ and $\vec{E}^0(G) \df V(G)^2\setminus\vec{E}(G)$
  and we let $\overline{G}\df(V(G),\binom{V}{2}\setminus E(G))$ be the complement graph of $G$.
  
  For a vertex $v\in V(G)$, we let $N_G(v)\df\{u\in V(G) \mid \{u,v\}\in E(G)\}$ denote the \emph{neighborhood} of $v$ in $G$
  and we let $d_G(v)\df\lvert N_G(v)\rvert$ denote the \emph{degree} of $v$ in $G$. We also define
  \begin{align*}
    N_G^1(v) & \df N_G(v), &
    N_G^0(v) & \df V(G)\setminus N_G(v), &
    d_G^1(v) & \df d_G(v), &
    d_G^0(v) & \df \lvert N_G^0(v)\rvert.
  \end{align*}

  For non-empty sets $U_1,U_2\subseteq V(G)$, the \emph{density} of the pair $(U_1,U_2)$ in $G$ is $d_G(U_1,U_2) \df
  \lvert\vec{E}(G)\cap(U_1\times U_2)\rvert/(\lvert U_1\rvert\cdot\lvert U_2\rvert)$. For $b\in\{0,1\}$, we set $d_G^b(U_1,U_2)
  \df \lvert\vec{E}^b(G)\cap(U_1\times U_2)\rvert/(\lvert U_1\rvert\cdot\lvert U_2\rvert)$.
\end{notation}

\begin{definition}[$\VC$-dimension and Littlestone trees]
  Let $\cH\subseteq\cP(X)$.
  \begin{enumdef}
  \item We say that $\cH$ \emph{shatters} $A\subseteq X$ if for the family $\cH\rest_A \df \{H\cap A \mid H\in\cH\}$ we have
    $\cH\rest_A = \cP(A)$, that is, when we restrict the family $\cH$ to $A$, we see all possible patterns on $A$. The
    \emph{Vapnik--\Cervonenkis\ ($\VC$) dimension} of $\cH$ is $\VC(\cH)\df\sup\{\lvert A\rvert \mid A\subseteq X\text{ finite
      set shattered by }\cH\}$. Note that we always have the trivial bound $\VC(\cH)\leq\min\{\lvert
    X\rvert,\log_2(\lvert\cH\rvert)\}$.
  \item For $t\in\NN$, a \emph{$t$-Littlestone tree} in $\cH$ is a pair $(H,x)$ of sequences, where $x =
    (x_\sigma)_{\sigma\in\{0,1\}^{< t}}$ is a sequence in $X$ (indexed by binary strings of length less than $t$) and
    $H=(H_\tau)_{\tau\in\{0,1\}^t}$ is a sequence in $\cH$ (indexed by binary strings of length exactly $t$) such that for every
    $\tau\in\{0,1\}^t$ and every $i\in[t]$, we have $x_{\tau\rest_{[i-1]}}\in H_\tau\iff\tau_i = 1$. See
    Figure~\ref{fig:littlestone}. The \emph{Littlestone dimension} of $\cH$ is $\Lit(\cH)\df\sup\{t\in\NN \mid \cH\text{ has a
      $t$-Littlestone tree}\}$. Note that we always have the trivial bound $\Lit(\cH)\leq\min\{\lvert
    X\rvert,\log_2(\lvert\cH\rvert)\}$.
    \begingroup
\def\depth{3}
\def\ptsize{2pt}
\def\layersep{1}
\def\absep{2}
\def\ycontrol{1}
\def\horbasesep{2}
\def\scale{1}

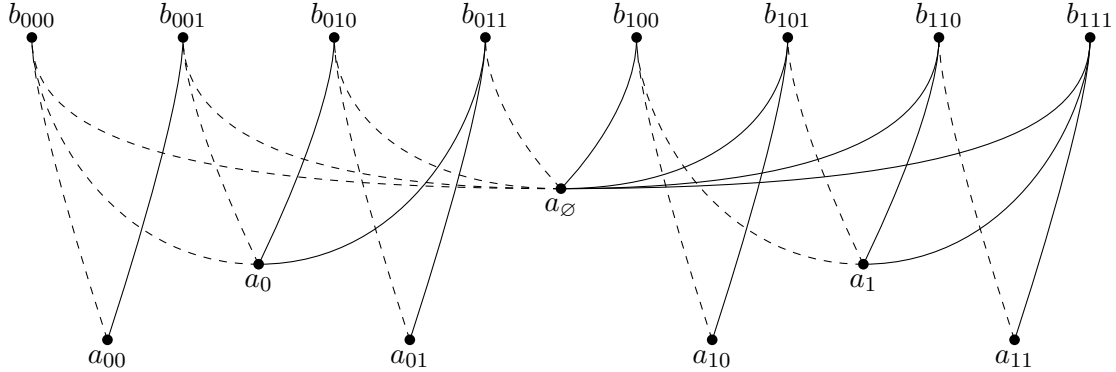
\begin{figure}[htbp]
  \centering
  \begin{tikzpicture}[scale=\scale]
    \def\strings{v,}
    \foreach \t [%
      evaluate=\t as \y using -(\t-1)*\layersep,
      evaluate=\t as \basehorstep using 2^(\depth-\t+1)%
    ] in {1,...,\depth}{
      \let\oldstrings\strings
      \def\strings{}
      \foreach \s [%
        count=\c,
        evaluate=\s as \x using (\basehorstep*(\c-1)+(\basehorstep-1)/2)*\horbasesep%
      ]
      in \oldstrings {
        \ifx\s\empty\relax
        \else
        \coordinate (\s) at (\x,\y);
        \fill (\s) circle (\ptsize);
        \StrGobbleLeft{\s}{1}[\l]
        \IfEq{\s}{v}{%
          \node[below] at (\s) {$a_\varnothing$};%
        }{%
          \node[below] at (\s) {$a_{\l}$};%
        }
        \edef\newstrings{\strings \s0,\s1,}
        \global\let\strings\newstrings
        \fi
      }
    }
    \pgfmathsetmacro{\y}{\absep}
    \foreach \s [%
      count=\c,
      evaluate=\s as \x using (\c-1)*\horbasesep%
    ] in \strings {
      \ifx\s\empty\relax
      \else
      \coordinate (\s) at (\x,\y);
      \fill (\s) circle (\ptsize);
      \StrGobbleLeft{\s}{1}[\l]
      \node[above] at (\s) {$b_{\l}$};
      \foreach \t [%
        evaluate=\t as \cy using -(\t-1)*\layersep%
      ] in {1,...,\depth} {
        \StrLeft{\s}{\t}[\p]
        \pgfmathtruncatemacro{\n}{\t+1}
        \StrChar{\s}{\n}[\char]
        \IfEq{\char}{1}{%
          \draw (\p) .. controls ($(\x,\cy) - 0.5*(\horbasesep,0)$) and +(0,-\ycontrol).. (\s);%
        }{%
          \draw[dashed] (\p) .. controls ($(\x,\cy) + 0.5*(\horbasesep,0)$) and +(0,-\ycontrol) .. (\s);%
        }
      }
      \fi
    }
  \end{tikzpicture}
  \caption{A $\depth$-Littlestone tree in a graph $G$. Solid arrows are edges, dashed arrows are non-edges. Pairs of vertices
    without any kind of line can be either adjacent or non-adjacent.}
  \label{fig:littlestone}
\end{figure}
\endgroup


  \item For a graph $G$, the \emph{$\VC$ dimension} and \emph{Littlestone dimension} of $G$ are defined as the corresponding
    concepts of the family $\cH_G\df\{N_G(v) \mid v\in V(G)\}$ of neighborhoods of vertices of $G$ (as a family of subsets of
    $V(G)$): $\VC(G)\df \VC(\cH_G)$ and $\Lit(G)\df \Lit(\cH_G)$.
  \end{enumdef}
\end{definition}

\begin{definition}[Homogeneity, goodness and excellence]
  Let $G$ be a finite graph and let $\delta,\epsilon\in[0,1]$.
  \begin{enumdef}
  \item A pair $(U_1,U_2)$ of non-empty subsets of $V(G)$ is \emph{$\epsilon$-homogeneous} in $G$ if
    $d_G(U_1,U_2)\in[0,\epsilon]\cup[1-\epsilon,1]$. (Equivalently, there exists $b_G^\epsilon(U_1,U_2)\in\{0,1\}$ such that
    $d_G^{b_G^\epsilon(U_1,U_2)}(U_1,U_2)\geq 1-\epsilon$.)
  \item For a finite set $U\subseteq V(G)$ and a vertex $v\in V(G)$, a \emph{$\delta$-majority opinion} of $U$ with respect to
    $v$ in $G$ is a value $t_G^\delta(v,U)\in\{0,1\}$ (if it exists) such that
    \begin{gather*}
      \lvert N_G^{t_G^\delta(v,U)}(v)\cap U\rvert \geq (1-\delta)\cdot\lvert U\rvert.
    \end{gather*}
    That is, if $t_G^\delta(v,U)=1$, then $v$ is adjacent to all but a $\delta$-proportion of vertices of $U$ and if
    $t_G^\delta(v,U)=0$, then $v$ is non-adjacent to all but a $\delta$-proportion of vertices of $U$. Note that if $\delta <
    1/2$ and $U$ is non-empty, then at most one $\delta$-majority opinion $t_G^\delta(v,U)$ exists.
  \item A non-empty finite set $U\subseteq V(G)$ is \emph{$\epsilon$-good} in $G$ if for every $v\in V(G)$, there exists a
    $\epsilon$-majority opinion $t_G^\epsilon(v,U)\in\{0,1\}$ of $U$ with respect to $v$ in $G$.
  \item For a finite set $W\subseteq V(G)$ and a $\delta$-good set $U$, a \emph{$(\delta,\epsilon)$-majority opinion} of $W$
    with respect $U$ in $G$ is a value $t_G^{\delta,\epsilon}(W,U)\in\{0,1\}$ (if it exists) such that
    \begin{gather*}
      \lvert\{w\in W \mid t_G^\delta(w,U) = t_G^{\delta,\epsilon}(W,U)\}\rvert
      \geq
      (1-\epsilon)\cdot\lvert W\rvert.
    \end{gather*}
    That is, for all but an $\epsilon$-proportion of the vertices $w$ in $W$, the $\delta$-majority opinion of $U$ with respect
    to $w$ (which exists since $U$ is $\delta$-good) is $t_G^{\delta,\epsilon}(W,U)$, i.e., the $\delta$-majority opinions of
    $U$ with respect to vertices of $W$ are the same up to $\epsilon$-error. Again, if $\epsilon < 1/2$ and $W$ is non-empty,
    then at most one $(\delta,\epsilon)$-majority opinion of $W$ with respect to $U$ in $G$ exists.
  \item A non-empty finite set $W\subseteq V(G)$ is \emph{$(\delta,\epsilon)$-excellent} in $G$ if for every $\delta$-good in
    $G$ set $U$, there exists an $(\delta,\epsilon)$-majority opinion $t_G^{\delta,\epsilon}(W,U)$ of $W$ with respect to $U$ in
    $G$. Note that since every singleton is $\delta$-good, it follows that $(\delta,\epsilon)$-excellence implies
    $\epsilon$-goodness.
  \end{enumdef}
\end{definition}

\begin{definition}[Atoms]
  Let $G$ be a graph, let $\delta,\epsilon > 0$, let $W\subseteq V(G)$ be a finite set, let $m\in\NN_+$ and let
  $x=(x_1,\ldots,x_m)\in V(G)^m$ be an $m$-tuple (possibly with repetitions) of vertices of $G$.
  \begin{enumdef}
  \item A \emph{$\delta$-majority opinion}\footnote{This is essentially the $\delta$-majority opinion concept when we view $x$
  as a multiset.} of a vertex $v\in V(G)$ with respect to $x$ in $G$ is a value $t_G^\delta(v,x)\in\{0,1\}$ such that
    $d_G^{t_G^\delta(v,x)}(v,x) \geq 1 - \delta$, where $d_G^b(v,x)\df\lvert\{i\in[m] \mid x_i\in N_G^b(v)\}\rvert/m$.
  \item We say that $W$ is \emph{$(\delta,\epsilon)$-split} by the tuple $x$ in $G$ if for every $b\in\{0,1\}$, we have
    \begin{gather*}
      \lvert\{w\in W \mid t_G^\delta(w,x) = b\}\rvert
      >
      \epsilon\cdot\lvert W\rvert.
    \end{gather*}
  \item A non-empty finite set $W\subseteq V(G)$ is a \emph{$(\delta,\epsilon,m)$-atom} of $G$ if for every $m$-tuple $x\in
    V(G)^m$, the set $W$ is \emph{not} $(\delta,\epsilon)$-split by $x$ in $G$. Clearly, $(\delta,\epsilon,1)$-atomicity is
    equivalent to $\epsilon$-goodness provided $\delta < 1$. Furthermore, the atomicity condition is trivially satisfied when
    $\epsilon\geq 1/2 < \delta$. We say that $W$ is a \emph{$(\delta,\epsilon)$-atom} of $G$ if it is $(\delta,\epsilon,m)$-atom
    of $G$ for every $m\in\NN_+$.
  \end{enumdef}
\end{definition}

\begin{definition}[Partitions]
  A \emph{partition} of a graph $G$ is a partition $\cP$ of $V(G)$. For $\delta,\epsilon > 0$, $\cP$ is:
  \begin{description}
  \item[Equitable (or an equipartition)] if $\lvert\lvert P\rvert - \lvert P'\rvert\rvert\leq 1$ for every $P,P'\in\cP$.
  \item[Totally $\epsilon$-homogeneous] if every\footnote{The adverb ``totally'' is to differentiate from notion of
  $\epsilon$-homogeneous partitions that allow for ``few'' pairs to fail $\epsilon$-homogeneity.} $(P,P')\in\cP^2$ is
    $\epsilon$-homogeneous in $G$.
  \item[$\epsilon$-good] if every $P\in\cP$ is $\epsilon$-good in $G$.
  \item[$(\delta,\epsilon)$-excellent] if every $P\in\cP$ is $(\delta,\epsilon)$-excellent in $G$.
  \item[$(\delta,\epsilon)$-atomic] if every $P\in\cP$ is a $(\delta,\epsilon)$-atom of $G$.
  \end{description}
\end{definition}


\subsection{Proof Overview}
\label{subsec:strat}

We first recall the mechanism behind the stable regularity lemma of Malliaris and Shelah. This will also explain why their
proof, although conceptually very clean, is not immediately algorithmic.

\subsubsection*{The Malliaris--Shelah proof: goodness, excellence, and failed Littlestone trees}

Let $G=(V,E)$ be a finite graph. The basic local notion is that of a majority opinion. If $U\subseteq V$ is non-empty and $v\in
V$, we say that $v$ has $\delta$-majority opinion $1$ about $U$ if $v$ is adjacent to all but a $\delta$-fraction of the
vertices of $U$. Similarly, $v$ has $\delta$-majority opinion $0$ about $U$ if $v$ is non-adjacent to all but a
$\delta$-fraction of $U$. When $\delta<1/2$, such a majority opinion, if it exists, is unique.

A set $U$ is $\epsilon$-good if every vertex $v\in V$ has an $\epsilon$-majority opinion about $U$. Equivalently, no single
vertex can split $U$ into two large pieces according to adjacency and non-adjacency. Thus, goodness is a one-vertex uniformity
condition: every vertex sees $U$ as almost entirely one color.

The next notion is excellence. Suppose $U$ is $\delta$-good. Then every vertex $w$ has a well-defined $\delta$-majority opinion
about $U$. We say that $W$ has a $(\delta,\epsilon)$-majority opinion about $U$ if, after discarding at most an
$\epsilon$-fraction of the vertices of $W$, all remaining vertices of $W$ have the same $\delta$-majority opinion about $U$.
Finally, $W$ is $(\delta,\epsilon)$-excellent if this happens for every $\delta$-good set $U$.

Thus, one should keep the following distinction in mind:
\begin{center}
  \begin{tabular}{rl}
    goodness: & no single vertex splits $W$,\\
    excellence: & no already-good set $U$ splits $W$ by majority opinions.
  \end{tabular}
\end{center}
Goodness is local: it tests $W$ against vertices. Excellence is non-local: it tests $W$ against good sets.

These notions are strong enough to imply the stable form of regularity. If two sets are good, then their pair is already almost
homogeneous: the edge density between them is close either to $0$ or to $1$. Excellence is a stronger coherence condition,
ensuring that majority opinions are themselves stable when tested against good sets. Consequently, a partition into excellent
pieces gives a stable regular partition with no irregular pairs. This is the qualitative difference between stable regularity
and \Szemeredi\ regularity: in the stable setting, one can eliminate exceptional pairs altogether.

We use Littlestone dimension to measure a forbidden tree pattern. For the neighborhood set system of $G$, a Littlestone tree of
height $t$ consists of vertices $(a_\sigma)_{\sigma\in\{0,1\}^{<t}}$ at the internal nodes and vertices
$(u_\tau)_{\tau\in\{0,1\}^t}$ at the leaves such that, for every leaf $\tau\in\{0,1\}^t$ and every level $i\in[t]$,
$a_{\tau\rest{[i-1]}}\in N_G(u_\tau) \iff \tau_i=1$. The Littlestone dimension $\Lit(G)$ is the largest height of such a tree.
Thus the assumption $\Lit(G)\leq \ell$ says precisely that no binary decision tree of height $\ell+1$ can be encoded by
neighborhoods in $G$.

The Malliaris--Shelah proof can be viewed as a failed Littlestone tree argument. Start with a large set. If it is not good, then
some vertex splits it into two large subsets. If it is good but not excellent, then some good set splits it into two large
subsets according to majority opinions. If this process could be continued for $\ell+1$ levels, it would produce a complete
binary tree of witnesses (provided $\delta < 1/2^{\ell+1}$), and hence a Littlestone tree of height $\ell+1$. This contradicts
$\Lit(G)\leq\ell$. Therefore, the splitting process must produce an excellent set before height $\ell+1$.

This proof is elegant, but it is not obviously algorithmic. The obstacle is the quantifier in the definition of excellence. To
certify that $W$ is excellent, one must verify coherent majority behavior with respect to every $\delta$-good set $U$. The
family of good sets is not explicitly given, and it can be enormous. Hence the original proof gives a structural certificate
that excellent sets exist, but it does not directly give a finite certificate or an efficient search procedure for finding one.

\subsubsection*{The new algorithmic ingredient: atomicity}

The main new idea of this paper is to replace excellence by a stronger intermediate notion, which we call atomicity. Atomicity
is designed so that it implies excellence but can be searched for using fixed length tuples of vertices. Instead of testing $W$
against arbitrary good sets $U$, we test it against finite tuples $x=(x_1,\ldots,x_m)\in V^m$. A vertex $w$ has a
$\delta$-majority opinion about a tuple $x$ if $w$ is adjacent to at least a $1-\delta$ fraction of the entries of $x$, or
non-adjacent to at least a $1-\delta$ fraction of the entries of $x$. We say that $W$ is $(\delta,\epsilon)$-split by $x$ if
both opinions occur on more than an $\epsilon$-fraction of $W$. A set $W$ is a $(\delta,\epsilon,m)$-atom if no $m$-tuple splits
$W$, and $W$ is a $(\delta,\epsilon)$-atom if it is a $(\delta,\epsilon,m)$-atom for every $m\geq 1$.

This definition is a finite-pattern analogue of excellence. For $m=1$, atomicity recovers goodness, provided $\delta<1$:
forbidding $1$-tuples from splitting $W$ says exactly that every vertex has a majority opinion about $W$. For larger $m$,
atomicity says that no finite pattern of vertices can separate $W$ into two large populations with different majority behavior.
Thus an atom behaves like one coherent type with respect to all finite tests.

The implication from atomicity to excellence is simple and is one of the key conceptual bridges. Let $W$ be a
$(\delta,\epsilon)$-atom, and let $U$ be any $\delta$-good set. Enumerate the elements of $U$ in a tuple
$x=(x_1,\ldots,x_{\lvert U\rvert})$. Since $U$ is $\delta$-good, every vertex $w\in W$ has a $\delta$-majority opinion about
$U$, and this is exactly the same as its $\delta$-majority opinion about the tuple $x$. If $W$ failed to have a coherent
majority opinion about $U$, then more than an $\epsilon$-fraction of $W$ would have opinion $0$, and more than an
$\epsilon$-fraction would have opinion $1$. But then the tuple $x$ would split $W$, contradicting atomicity. Therefore every
atom is excellent.

At first sight, atomicity seems even harder to verify than excellence, because it quantifies over tuples of all lengths (i.e.,
it seems that we have passed from having to check an exponential amount of conditions to having to check countably infinitely
many conditions). The second key observation is that bounded $\VC$ dimension reduces this infinite family of tests to a
bounded-length family. Since bounded Littlestone dimension implies bounded $\VC$ dimension, this is still a consequence of
stability. Suppose $\VC(G)\leq d$. If a long tuple $z$ splits $W$, view $z$ as defining a probability measure on its entries
(counting multiplicities). By uniform convergence for $\VC$ classes, there is a much shorter tuple $x$, of length $m =
O(d/(c^2\cdot\delta^2))$ (this will turn out to be $O(d\cdot 2^{2\cdot\ell+4})$) that approximates all relevant neighborhood
densities up to error $c\cdot\delta$. Consequently, any vertex that had a $\delta$-majority opinion about $z$ has the same
$(1+c)\cdot\delta$-majority opinion about $x$. Hence the short tuple still splits $W$, with a slight loss in the majority
parameter. Equivalently, checking $((1+c)\delta,\epsilon,m)$-atomicity for this bounded value of $m$ certifies full
$(\delta,\epsilon)$-atomicity. This is the parameter adjustment used by the algorithms when they search for finite tuple
witnesses.

This is the point at which the proof becomes algorithmic. We no longer search over all good sets. We search over all tuples of
one fixed length $m$, where $m$ depends only on the $\VC$ dimension and the error parameters, not on $n$. For fixed parameters,
atomicity becomes a finite condition.

The algorithmic splitting procedure is then a finite version of the Malliaris--Shelah failed Littlestone tree argument. Start
with a large set $U_\varnothing$. At a node $\sigma$, search for an $m$-tuple $x_\sigma$ that $(\delta,\epsilon)$-splits
$U_\sigma$. If no such tuple exists, then $U_\sigma$ is a $(\delta,\epsilon,m)$-atom, and the algorithm stops. If such a tuple
exists, then both majority classes are large, and we choose large children
\begin{align*}
  U_{\sigma\conc 0} & \subseteq \{u\in U_\sigma \mid t^\delta_G(u,x_\sigma)=0\},
  &
  U_{\sigma\conc 1} & \subseteq \{u\in U_\sigma \mid t^\delta_G(u,x_\sigma)=1\}.
\end{align*}
Thus, each successful split records a finite tuple and divides the current set according to majority opinion about that tuple.

The stopping argument is the same stability principle as in the original proof. If the algorithm did not stop for $\ell+1$
levels, then we would obtain a full binary tree of sets $U_\sigma$ and splitting tuples $x_\sigma$. Choose one leaf vertex
$u_\tau$ from each terminal set. For every internal node $\sigma$ and every descendant leaf $\tau$, the construction guarantees
that $u_\tau$ has the branch-prescribed majority opinion about $x_\sigma$. Since this is a majority statement over the
coordinates of $x_\sigma$, and imposing a constraint $\delta<2^{-\ell-1}$, a union bound lets us choose a single coordinate
$(x_\sigma)_{i_\sigma}$ that works simultaneously for all descendant leaves. The selected vertices
$a_\sigma=(x_\sigma)_{i_\sigma}$ together with the neighborhoods $N_G(u_\tau)$ of the leaf vertices, form a Littlestone tree of
height $\ell+1$. This contradicts $\Lit(G)\leq\ell$. Therefore, the algorithm must stop before height $\ell+1$, and the set
returned at the stopping node is atomic.

We also note that in such a procedure sets $U_\sigma$ are guaranteed to have size at least
$\epsilon^{\lvert\sigma\rvert}\cdot\lvert U_\varnothing\rvert$ simply because at each level the sets are guaranteed to have a
proportion $\epsilon$ of the vertices of that generated them. This means gives our crucial ingredient for producing atomic
partitions:
\begin{center}
  Given a set $U$, one can algorithmically efficiently extract a $(\delta,\epsilon,m)$-atom $W\subseteq U$ with $\lvert
  W\rvert\geq\epsilon^\ell\cdot\lvert U\rvert$.
\end{center}

This means that we can repeat this extraction step enough times to produce a partition of all but a small leftover set into
atomic pieces. The leftover vertices can then be distributed among the existing pieces. This uses a continuity lemma: adding or
removing a small relative fraction of vertices only mildly worsens goodness, excellence, or atomicity. Hence, after absorbing
the leftover vertices, the partition remains atomic with slightly weaker parameters. Since atoms are excellent, this gives the
desired stable regular (not necessarily equitable) partition. The choice of parameters that yields the bound
$(1+c)\cdot\epsilon^{-\ell}\cdot\ln(1/\epsilon)$ on the number of parts is to extract atomic sets with parameter
$\widetilde{\epsilon}\df (1-\epsilon^{c'})\cdot\epsilon$ (roughly), where $c' > 0$ is a constant that depends only on $c$.

\subsubsection*{Equipartitions}

An equipartition is a partition whose parts have sizes differing by at most one. Even with a slightly modified version that was
already known for excellent sets by Malliaris--Shelah~\cite{MS14}, the extraction procedure above can only extract sets in a
controlled geometric sequence $\lvert U\rvert\geq s_0,s_1,\ldots,s_\ell=\epsilon^\ell\cdot s_0$ of ratio $\epsilon$. The main
obstacle is that goodness, excellence and atomicity are not necessarily preserved under arbitrary subsets. However,
Malliaris--Shelah observed that if one instead picked a uniformly random subset of a fixed desired size, then hypergeometric
distribution bounds say that any give trace of a majority opinion will be approximately preserved except for an exponentially
small probability of failure. For goodness, exponential concentration is enough: we only need to approximately preserve $n$
possible majority opinions, so a union bound goes through. For excellence, Malliaris--Shelah use $\VC$ dimension: one can show
that since the $\VC$ dimension is bounded, the number of patterns of majority opinions of good sets in a given excellent set can
only be polynomial in $n$ (this is a consequence of the Sauer--Shelah--Perles Lemma), so an exponential concentration still
allows a polynomial sized union bound to go through.

In the case of atomicity, this step is even easier: we already used $\VC$ theory to only consider $m$-tuples of a fixed length
$m$; these clearly can only generate at most $n^m$ patterns. This gives a first (randomized) strategy for equipartitions: choose
carefully parameters $\widetilde{\epsilon}\leq\xi\leq\epsilon$ and an initial size $s_0$ and set
$s_i\df\widetilde{\epsilon}\cdot s_{i-1}$ for $i\in[\ell]$; then, while we have at least $s_0$ vertices left, by applying the
extraction procedure extract an atom of size $s_i$ for some $i\in\{0,\ldots,\ell\}$ and parameter $\widetilde{\epsilon}$; next,
subpartition all parts randomly so that they all have sizes $s_\ell$ and with high probability are atoms with parameter $\xi$;
finally, take the remaining at most $s_0$ vertices and distribute them equally among the parts. Our choice of parameters needs
to be such that this final step only deteriorates the atomicity parameter from $\xi$ to $\epsilon$. The best choice of
parameters for this strategy is roughly
\begin{align*}
  \xi & \df \frac{\ell\cdot\epsilon}{\ell+1}, &
  s_0 & \df \frac{\epsilon}{\ell+1}, &
  \widetilde{\epsilon} & \df (1-c')\cdot\xi
\end{align*}
and yields a total number of parts of $(1+c)\cdot(1 + 1/\ell)^\ell\cdot(\ell+1)\cdot\epsilon^{-\ell-1}$, where $c > 0$ is a
constant that can be made arbitrarily small by adjusting $c'$ alone. However, this is off from our promised bound by a
multiplicative factor of roughly $e\cdot(\ell+1)$.

To actually obtain $K\df (1+c)\cdot\epsilon^{-\ell-1}$ parts, we make an adaptive extraction. Namely, the $i$th extraction
procedure is ran with parameter roughly $\widetilde{\epsilon}_i\df(1-c')\cdot\xi_i$, where $\xi_i\in(0,\epsilon)$ satisfies
\begin{gather}\label{eq:strategy:xii}
  -1 + \epsilon - \xi_i + \widetilde{\epsilon}_i^\ell\cdot\frac{n_{i-1}}{r} = 0,
\end{gather}
where $n_{i-1}$ is the number of vertices we have left before the $i$th extraction and $r\df n/K$ and we later run the random
subsampling procedure to yield an atom of size $\widetilde{\epsilon}_i^\ell\cdot n_{i-1}$ and parameter $\xi_i$. The constraint
on $\xi_i$ above is exactly what allows a part of size $\widetilde{\epsilon}_i^\ell\cdot n_{i-1}$ to absorb enough
vertices to become of size $r$ while deteriorating its atomicity parameter from $\xi$ to at most $\epsilon$. The bulk of the
work is to show that such $\xi_i$ exist for all $i\in[K]$ and given our choice of $c'$ in terms of $c$ alone.

Before we talk about how to eliminate randomness from this procedure, let us point out that in the proof intuition above (both
in the non-equitable and equitable cases) we completely ignored divisibility constraints for simplicity of presentation. To take
these into account, all arguments instead need to replace $\epsilon$ by $\zeta\df(1-c_{\ZZ})\cdot\epsilon$ for some $c_{\ZZ} >
0$ that can be taken to be as small as desired; this is adjustment only affects the final constant $c > 0$.

\subsubsection*{Derandomization}

To derandomize our algorithms, we want to use the standard derandomization trick: if we can find an algorithm that uses a
logarithmic amount of bits of randomness, then by enumerating all such bits and testing what the result of the randomized
algorithm would be, we can obtain a polynomial time deterministic algorithm. The issue is that to perform our random subset
sampling, the number of random bits required is not logarithmic at all. So instead, we need to change our randomized algorithm
slightly. Instead of subsampling a set of size $\widetilde{\epsilon}_i^\ell\cdot n_{i-1}$ from our extracted set $W$, we
consider instead setting $p\df\lvert\widetilde{\epsilon}_i^\ell\cdot n_{i-1}/\lvert W\rvert\in[0,1]$ and we flip a $p$-biased
coin for each vertex of and take the random subset $\rn{W}$ of all vertices of $W$ whose coin was heads. This is not quite a
hypergeometric distribution but it behaves sufficiently similarly that the argument still goes through: with high probability
the resulting set will be an atom and will have size approximately (rather than exactly) $\widetilde{\epsilon}_i^\ell\cdot
n_{i-1}$. However, even flipping these biased coins requires too many bits of randomness; so instead we intentionally prove a
weaker concentration result, via a high-moment argument as opposed to a Chernoff bound, in this way any distribution of
$p$-biased random bits that the same moments as that of independent $p$-biased random bits up to the moment used in the
concentration bound also satisfies the same concentration. For goodness, we need only to consider the fourth moment, for
atomicity, we end up considering the $(2\cdot m+2)$th moment. Later we show that it is possible to generate a $p$-biased
distribution of $n$ bits that has the correct first $h$ moments, we need only $O_\epsilon(1) + h\cdot\log_2(n)$ bits of
randomness (and we can do this in time $O_\epsilon(h\cdot n^2\cdot\log(n+1))$). Since checking goodness and atomicity take time
$O_\epsilon(n^2)$ and $O_\epsilon(n^{m+1})$ respectively, this yields a total time of $O_{c,\ell,\epsilon}(n^6)$ for good
equipartitions and of $O_{c,\ell,\epsilon}(n^{3\cdot m + 3})$ for atomic equipartitions.

\subsubsection*{The almost tight lower bounds}

The lower-bound constructions are variations of the same $\ell$-level stable branching example. One starts with the selector
graph $G_{\ell,r}$, which is a bipartite graph whose left vertices are non-empty strings in $[r]^{<\ell+1}$ and whose right
vertices are strings in $[r]^\ell$, which are connected precisely to the left vertices that are prefixes of it. One then
considers two kinds of independent blow-ups of right vertices of $G_{\ell,r}$: in the non-equitable case, one replaces all right
vertices with independent sets of the same arbitrarily large size $n$, and in the equitable case, one does the same, except that
one of the parts is made of size $\lambda\cdot n$ for some $\lambda\geq 1$. It is straightforward to check that all these graphs
have Littlestone dimension $\ell$. Let us call each of the resulting sets on the right side a \emph{leaf} and call the one
resulting from $\One$ in the equitable construction the \emph{special leaf}.

When we take $n$ to be arbitrarily large, any partition into $\epsilon$-good sets of these graphs can be turned into a partition
into $\epsilon'$-good sets of the right side alone for an $\epsilon'$ arbitrarily close to $\epsilon$ (provided $n$ is
sufficiently large); this is a consequence of the continuity lemma for goodness. Furthermore, we can also preserve equitability
in this construction.

The key observation is that if we take $r < (1-\epsilon')/\epsilon'$ (in fact in the construction, we take the largest possible
such $r\df\ceil{(1-\epsilon')/\epsilon'} - 1$), then for every $\epsilon'$-good set containing only right vertices, there must
exists exactly one leaf containing a $1-\epsilon'$ proportion of the vertices of the $\epsilon'$-good set.

In the non-equitable case, this immediately puts an upper bound on the size of the $\epsilon'$-good part: it can be at most
$n/(1-\epsilon')$, hence this gives a lower bound of $(1-\epsilon')\cdot r^\ell = (1-o_{n\to\infty,\epsilon\to
  0,\ell}(1))\cdot\epsilon^{-\ell}$ on the number of parts.

In the equitable case, we set $\lambda\df(r^\ell-1)\cdot r$ so that if an $\epsilon'$-good equipartition of the right side has
fewer than some number $K_*$ of parts, then every $\epsilon'$-good set must have at least a $(1-\epsilon')$-proportion of its
vertices in the special leaf of size $\lambda\cdot n$ (as the parts will be too big for the size $n$ of the small leaves to be
$(1-\epsilon')$ proportion of it). This happens provided $K_* < (1-\epsilon')\cdot(r^\ell-1+\lambda)$. Summing the mass of the
parts outside the special leaf, we then get $(r^\ell-1)\cdot n\leq\epsilon'\cdot(r^\ell-1+\lambda)\cdot n$, which when
rearranged gives $\lambda \geq (1-\epsilon')\cdot (r^\ell-1)/\epsilon' > (r^\ell-1)\cdot r = \lambda$, a contradiction. Thus,
the equipartition must have had at least $(1-\epsilon')\cdot(r^\ell-1+\lambda) = (1-o_{n\to\infty,\epsilon\to
  0,\ell}(1))\cdot\epsilon^{-\ell-1}$ parts.

\subsubsection*{Space-efficient algorithms}

Let us briefly list the key observations that lead to $O_{c,\epsilon,\ell}(\log(n+1))$-space algorithms:
\begin{itemize}
\item When we extract a good/atom set, we do not need to keep track of the set itself, rather, it suffices to keep track of path
  and vertices/tuples used in the failed Littlestone tree to construct it. This takes space $O(\ell\cdot\log(n+1))$ for good
  sets and $O(\ell\cdot m\cdot\log(n+1))$ for atoms (per set).
\item Adjusting sizes of sets can be done greedily: if we need to add vertices, we simply add vertices that were not used in
  the natural order of $V(G)=[n]$; if we need to remove vertices, we remove them in reverse of the natural order of $V(G)=[n]$.
  This means that we only need to keep track of the final size of the size adjustment; this takes space $O(\log(n+1))$ per set.
\item For the subsampling via $p$-biased distribution, we need to keep track of which set of $O_\epsilon(1) + h\cdot\log_2(n)$
  bits generated the correct subsample. This obviously takes $O_\epsilon(1) + h\cdot\log_2(n)$ per set.
\item We need to argue that the generation of the $n$ pseudorandom bits from the $O_\epsilon(1) + h\cdot\log_2(n)$ actually
  random bits can be done space-efficiently. The space bound we end up getting is $h\cdot(O_\epsilon(1) + O(\log(n+1)))$ (recall
  that $h$ is the pseudorandomness quality in terms of moments).
\end{itemize}

\subsubsection*{Other algorithmic technicalities}

Let us briefly mention that there are extra technicalities regarding computation of parameters that were not covered in this
proof intuition: for example, several of the calculations (e.g., finding $\xi_i$ satisfying~\eqref{eq:strategy:xii}) cannot be
done exactly as they typically yield irrational values; instead the algorithms have to perform sufficiently good approximations
to these parameters and we have to justify that these can be done using a reasonable amount of time/space.


\section{Deferred Preliminaries}

\begin{notation}
  In algorithm inputs, rational numbers are represented by specifying their numerator and denominator (in reduced form, in base
  $2$) and for a rational $q$, we denote by $\len(q)$ the bitlength of the representation of $q$.

  Our undecorated big Oh's and related notation only hide universal constants. For decorated ones, the order of the terms
  specifies the order of the limits taken, e.g., a notation of the form $o_{n\to\infty,\epsilon\to 0,\delta}(1)$ means a function
  $f(n,\epsilon,\delta,c)$ (where $c$ collects all other parameters) such that for every function $c(n,\epsilon,\delta)$ and
  every fixed $\delta$, we have
  \begin{gather*}
    \lim_{\epsilon\to 0} \lim_{n\to\infty} f\bigl(n,\epsilon,\delta,c(n,\epsilon,\delta)\bigr) = 0,
  \end{gather*}
  or, equivalently, for every $\alpha > 0$, there exists $\epsilon_0 > 0$ such that for every
  $\epsilon\in(-\epsilon_0,\epsilon_0)$, there exists $n_0\in\NN$ such that for every $n\geq n_0$ and every $\delta$, we have
  \begin{gather*}
    \Bigl\lvert f\bigl(n,\epsilon,\delta,c(n,\epsilon,\delta)\bigr)\Bigr\rvert < \alpha.
  \end{gather*}
  Similarly, a notation of the form $O_{n\to\infty,\epsilon\to 0,\delta}(g(n,\epsilon,c))$ means a function
  $f(n,\epsilon,\delta,c))$ such that there exists a function $h(\delta)$ such that for every fixed $\delta$, we have
  \begin{gather*}
    \limsup_{\epsilon\to 0} \limsup_{n\to\infty}
    \frac{\lvert f(n,\epsilon,\delta,c)\rvert}{\lvert g(n,\epsilon,c)\rvert}
    \leq
    h(\delta).
  \end{gather*}
\end{notation}

\begin{definition}[Order property]
  Let $\cH\subseteq\cP(X)$ be a family of subsets of some set $X$.
  \begin{enumdef}
  \item For $t\in\NN$, a \emph{$t$-order property} (or \emph{$t$-half graph}) in $\cH$ is a pair $(H,x)$ of sequences, where $H
    = (H_i)_{i=1}^t$ is a sequence in $\cH$ and $x=(x_j)_{j=1}^t$ is a sequence in $X$ such that for every $i,j\in[t]$, we have
    \begin{gather*}
      x_j\in H_i \iff i\leq j.
    \end{gather*}
    The \emph{order property dimension} of $\cH$ is $\OP(\cH)\df\sup\{t\in\NN \mid \cH\text{ has a $t$-order property}\}$. Note
    that we always have the trivial bound $\OP(\cH)\leq\min\{\lvert X\rvert,\lvert\lvert\cH\rvert\}$.
  \item For a graph $G$, the \emph{order property dimension} is defined as the corresponding concepts of the family
    $\cH_G\df\{N_G(v) \mid v\in V(G)\}$ of neighborhoods of vertices of $G$ (as a family of subsets of $V(G)$):
    $\OP(G)\df\OP(\cH_G)$.
  \end{enumdef}
\end{definition}

\begin{definition}[Independent blow-ups]\label{def:indepblowup}
  Given a graph $G$ and a function $b\colon V(G)\to\NN$, the \emph{$b$-independent blow-up} of $G$ is the graph $G^b$ obtained
  from $G$ by replacing each vertex $v\in V(G)$ by $b$ copies of it (and making edges become complete bipartite graphs between
  the copies of its endpoints). Formally, we have
  \begin{align*}
    V(G^b) & \df \Bigl\{(v,i) \;\Bigm\vert\; v\in V(G)\land i\in\bigl[b(v)\bigr]\Bigr\},
    \\
    E(G^b) & \df \Bigl\{\bigl\{(v,i),(w,j)\bigr\} \;\Bigm\vert\; \{v,w\}\in E(G)\Bigr\}.
  \end{align*}
\end{definition}


\subsection{Models of computation}

In this paper we will consider the following models of computation:
\begin{description}
\item[Deterministic:] Standard model of computation. More formally, multi-tape Turing machine with a fixed finite alphabet
  (in particular, arithmetic operations are typically done representing the number in binary ).
\item[Randomized:] Multi-tape turing machine with access to as many oracles to uniform distributions on $[m]$ for some
  $m\in\NN_+$ as wanted. Each oracle consultation costs one unit of time. Our randomized algorithms will be Monte Carlo
  algorithms, meaning that they will have a bound on the time-complexity that always holds, but their answer will only be
  correct with at least some promised probability $1-\rho$. In our case, at the expense of taking more time, the probability of
  failure $\rho$ can be made arbitrarily small.
\item[Random-query:] a randomized model that further has access to a special random-query oracle $\cO_G$ about a graph $G$.
  Given a vertex $v$ and a non-empty set $U\subseteq V(G)$, consulting $\cO_G(v,U)$ returns a random bit $\rn{b}$ corresponding
  to sampling $\rn{u}$ uniformly at random in $U$ and returning $\rn{b}\df\One[\{\rn{u},v\}\in E(G)]$. Note that if one has
  (standard) oracle access to the graph $G$, the randomized model can simulate such random-query oracle in time $O(\log(n+1))$,
  which means that any random-query algorithm can be transformed into a randomized algorithm by paying an extra multiplicative
  time cost of $\log(n+1)$.
\item[Query-into-oracle:] This model of computation is used for space-efficient computations in which a typical representation
  of the input and/or the output is much larger than the space-complexity of its computation (e.g., the typical representation
  of a partition of $[n]$ into $K$ labelled parts takes space $O(n^{\log(K+1)})$). This model then contains two algorithms, one
  called the \emph{computation algorithm} and one called the \emph{oracle algorithm}. Both algorithms are given oracle access to
  their input (e.g., this might be an oracle $\cO_G$ for adjacency of a graph $G$). The computation algorithm performs a
  preprocessing that produces an intermediate object $\sigma$ of small size (e.g., the computation algorithm might produce a
  list of vertices that play a key role in determining a stable partition). The oracle algorithm receives this intermediate
  object and functions as an oracle for some big output, meaning that upon receiving some extra information $u$, it outputs some
  small amount of information $b_{\sigma,u}$ about $u$ such that the full sequence $(b_u)_u$ over all $u$ is the desired output
  (e.g., if the output is a meant to be subset $S_\sigma$ of $[n]$ implicitly encoded by $\sigma$, then upon receiving $\sigma$
  and $u\in[n]$, the oracle outputs $b_{\sigma,u}=\One[u\in S_\sigma]$).
\end{description}

Table~\ref{tab:complexities} contains a bird-eye view of the different algorithms used to compute stable partitions and their
complexities.
\begin{landscape}
  \vfill
  \begin{small}
    \begin{table}[p]
      \centering
      \begingroup
      \renewcommand{\arraystretch}{1.2}
      \begin{tabular}{l*{9}{c}}
        Alg.
        &
        \multicolumn{2}{c}{Theorem}
        &
        \multicolumn{2}{c}{Model}
        &
        \multicolumn{2}{c}{Partition}
        &
        Number of parts
        &
        Time-complexity
        &
        Space-complexity
        \\
        \hline\hline
        \ref{alg:nonequi:good}
        &
        \multirow{3}{*}{\ref{thm:nonequi}}
        &
        \ref{thm:nonequi:good}
        &
        \multicolumn{2}{c}{Random-query}
        &
        \multirow{3}{*}{%
          \begin{tabular}{c}
            Non-\\
            equitable
          \end{tabular}
        }
        &
        \multirow{2}{*}{good}
        &
        \multirow{3}{*}{$\displaystyle\frac{(1+\widetilde{c})\cdot\ln(1/\epsilon)}{\epsilon^{\ell}}$}
        &
        $O_{\epsilon,c,\ell}\bigl(n\cdot\log(n/\rho)\bigr)$
        \\
        \cline{4-5}
        \ref{alg:nonequi:gooddet}
        & &
        \ref{thm:nonequi:gooddet}
        &
        \multicolumn{2}{c}{\multirow{2}{*}{Deterministic}}
        & & & &
        $O_{\epsilon,c,\ell}(n^2)$
        \\
        \cline{7-7}
        \ref{alg:nonequi:atom}
        & &
        \ref{thm:nonequi:atom}
        & & & & 
        atomic
        & &
        $O_{\epsilon,c,\ell}(1)\cdot n^{O(d\cdot 2^{2\cdot\ell+4})}$
        \\
        \cline{1-9}
        \ref{alg:equiPRG:good}
        &
        \ref{thm:equiPRG}
        &
        \ref{thm:equiPRG:good}
        &
        \multicolumn{2}{c}{Deterministic}
        &
        \multirow{4}{*}{Equitable}
        &
        \multirow{2}{*}{good}
        &
        \multirow{4}{*}{$\displaystyle\frac{(1+\widetilde{c})}{\epsilon^{\ell}}$}
        &
        $O_{\epsilon,c,\ell}(n^6)$
        \\
        \cline{2-5}
        \ref{alg:equirand:good}
        &
        \multirow{2}{*}{\ref{thm:equirand}}
        &
        \ref{thm:equirand:goodalg}
        &
        \multicolumn{2}{c}{Random-query}
        & & & &
        $O_{\epsilon,c,\ell}(n\cdot\log(n/\rho))$
        \\
        \cline{4-5}\cline{7-7}
        \ref{alg:equirand:atom}
        & &
        \ref{thm:equirand:atomalg}
        & 
        \multicolumn{2}{c}{Randomized}
        & &
        \multirow{2}{*}{atomic}
        & &
        $O_{\epsilon,c,\ell}(1)\cdot n^{O(d\cdot 2^{2\cdot\ell+4})}$
        \\
        \cline{2-5}
        \ref{alg:equiPRG:atom}
        &
        \ref{thm:equiPRG}
        &
        \ref{thm:equiPRG:atom}
        &
        \multicolumn{2}{c}{Deterministic}
        & & & &
        $O_{\epsilon,c,\ell}(1)\cdot n^{O(d\cdot 2^{2\cdot\ell+4})}$
        \\
        \hline
        \ref{alg:nonequi:goodspace}
        &
        \multirow{4}{*}{\ref{thm:nonequi}}
        &
        \ref{thm:nonequi:goodspace}
        &
        \multirow{8}{*}{\rotatebox{90}{Query-into-oracle}}
        &
        \multirow{2}{*}{computation}
        &
        \multirow{4}{*}{%
          \begin{tabular}{c}
            Non-\\
            equitable
          \end{tabular}
        }
        &
        good
        &
        \multirow{4}{*}{$\displaystyle\frac{(1+\widetilde{c})\cdot\ln(1/\epsilon)}{\epsilon^{\ell}}$}
        &
        $n^{O_{\epsilon,c,\ell}(1)}$
        &
        \multirow{8}{*}{$O_{\epsilon,c,\ell}(\log(n+1))$}
        \\
        \ref{alg:nonequi:atomspace}
        & &
        \ref{thm:nonequi:atomspace}
        & & & &
        atomic
        & &
        $n^{O_{\epsilon,c,\ell}(1)}$
        \\
        \cline{1-1}\cline{5-5}
        \multirow{2}{*}{\ref{alg:nonequi:atomspaceoracle}}
        & &
        \ref{thm:nonequi:goodspace}
        & &
        \multirow{2}{*}{oracle}
        & &
        good
        & &
        $n^{O(\epsilon^{-\ell}\cdot\log(1/\epsilon))}$
        \\
        & &
        \ref{thm:nonequi:atomspace}
        & & & &
        atomic
        & &
        $n^{O(d\cdot 2^{2\cdot\ell+4})}$
        \\
        \cline{1-3}\cline{5-9}
        \ref{alg:equiPRG:goodspace}
        &
        \multirow{4}{*}{\ref{thm:equiPRG}}
        &
        \ref{thm:equiPRG:goodspace}
        & &
        \multirow{2}{*}{computation}
        &
        \multirow{4}{*}{Equitable}
        &
        good
        &
        \multirow{4}{*}{$\displaystyle\frac{(1+\widetilde{c})}{\epsilon^{\ell}}$}
        &
        $O_{\epsilon,c,\ell}(1)\cdot n^{O(\epsilon^{-\ell-1})}$
        \\
        \ref{alg:equiPRG:atomspace}
        & &
        \ref{thm:equiPRG:atomspace}
        & & & &
        atomic
        & &
        $O_{\epsilon,c,\ell}(1)\cdot n^{O(\epsilon^{-\ell-1}\cdot d\cdot 2^{2\cdot\ell+4})}$
        \\
        \cline{1-1}\cline{5-5}
        \multirow{2}{*}{\ref{alg:equiPRG:atomspaceoracle}}
        & &
        \ref{thm:equiPRG:goodspace}
        & &
        \multirow{2}{*}{oracle}
        & &
        good
        & &
        $O_{\epsilon,c,\ell}(1)\cdot n^{O(\epsilon^{-\ell-1})}$
        \\
        & &
        \ref{thm:equiPRG:atomspace}
        & & & &
        atomic
        & &
        $O_{\epsilon,c,\ell}(1)\cdot n^{O(\epsilon^{-\ell-1})}$
      \end{tabular}
      \endgroup
      \caption{Algorithms for stable partitions. All algorithms further have several fixed parameters that control
        the $o(1)$ term of the multiplicative complexities, which are summarized simply as $c$ in the above (and the $o(1)$ term is
        summarized as $\widetilde{c}$); for the specific complexity dependence on these parameters, see the corresponding theorem.
        In the above $\ell$ is the bound in the Littlestone dimension and $d\leq\ell$ is the bound in the $\VC$ dimension and $\rho$
        is a parameter that controls the probability of failure. We do \emph{not} compute the space-complexity of the algorithms that
        were specifically designed for time-efficiency alone.}
      \label{tab:complexities}
    \end{table}
  \end{small}
  \vfill
\end{landscape}

\section{Basic facts}

\begin{lemma}[Monotonicity in parameters]\label{lem:mono}
  Let $G$ be a graph, let $\delta,\widetilde{\delta},\epsilon,\widetilde{\epsilon}\in(0,1)$ with $\delta\geq\widetilde{\delta}$
  and $\epsilon\leq\widetilde{\epsilon}$. Let $b\in\{0,1\}$, let $m\in\NN_+$, let $U\subseteq V(G)$ be a non-empty set, let
  $v\in V(G)$ and let $x\in V(G)^m$ be an $m$-tuple of vertices. Then the following hold:
  \begin{enumerate}
  \item\label{lem:mono:good} If $U$ is $\epsilon$-good in $G$, then $U$ is $\widetilde{\epsilon}$-good in $G$.
  \item\label{lem:mono:maj} If $t_G^{\widetilde{\delta}}(v,x)=b$ then $t_G^\delta(v,w)=b$.
  \item\label{lem:mono:exc} If $U$ is $(\delta,\epsilon)$-excellent in $G$, then $U$ is
    $(\widetilde{\delta},\widetilde{\epsilon})$-excellent in $G$.
  \item\label{lem:mono:atomtuple} If $U$ is a $(\delta,\epsilon,m)$-atom of $G$, then $U$ is a
    $(\widetilde{\delta},\widetilde{\epsilon},m)$-atom of $G$.
  \item\label{lem:mono:atom} If $U$ is a $(\delta,\epsilon)$-atom of $G$, then $U$ is a
    $(\widetilde{\delta},\widetilde{\epsilon})$-atom of $G$.
  \end{enumerate}
\end{lemma}

\begin{proof}
  For item~\ref{lem:mono:good}, we know that for every $v\in V(G)$, there exists an $\epsilon$-majority opinion
  $t_G^\epsilon(v,U)$ of $U$ with respect to $v$ in $G$, so we have
  \begin{gather*}
    \lvert N_G^{t_G^\epsilon(v,U)}\cap U\rvert \geq (1-\epsilon)\cdot\lvert U\rvert\geq (1-\widetilde{\epsilon})\cdot\lvert U\rvert,
  \end{gather*}
  so $t_G^\epsilon(v,U)$ is also a $\widetilde{\epsilon}$-majority opinion of $U$ with respect to $v$ in $G$, hence $U$ is
  $\widetilde{\epsilon}$-good in $G$.

  \medskip

  For item~\ref{lem:mono:maj}, since $t_G^{\widetilde{\delta}}(v,x)=b$, we know that
  \begin{gather*}
    \lvert\{i\in[m] \mid x_i\in N_G^b(v)\}\rvert
    \geq
    (1-\widetilde{\delta})\cdot m
    \geq
    (1-\delta)\cdot m,
  \end{gather*}
  so $t_G^\delta(v,x)=b$.

  \medskip

  For item~\ref{lem:mono:exc}, let $W$ be a $\widetilde{\delta}$-good set in $G$. By item~\ref{lem:mono:good}, we know that $W$
  is $\delta$-good in $G$ and since $U$ is $(\delta,\epsilon)$-excellent in $G$, there exists a $(\delta,\epsilon)$-majority
  opinion $t_G^{\delta,\epsilon}(U,W)$ of $W$ with respect to $U$ in $G$, which means that
  \begin{gather*}
    \{u\in U \mid t_G^\delta(u,W) = t_G^{\delta,\epsilon}(U,W)\}
  \end{gather*}
  has size at least $(1-\epsilon)\cdot\lvert U\rvert$. We claim that the set above contains the set
  \begin{gather*}
    \{u\in U \mid t_G^{\widetilde{\delta}}(u,W) = t_G^{\delta,\epsilon}(U,W)\}.
  \end{gather*}
  Indeed, this follows from item~\ref{lem:mono:maj} by taking $x\in V(G)^{\lvert W\rvert}$ to be an $\lvert W\rvert$-tuple that
  enumerates the elements of $W$. Thus, we get
  \begin{gather*}
    \lvert\{u\in U \mid t_G^{\widetilde{\delta}}(u,W) = t_G^{\delta,\epsilon}(U,W)\}\rvert
    \geq
    (1-\epsilon)\cdot\lvert U\rvert
    \geq
    (1-\widetilde{\epsilon})\cdot\lvert U\rvert,
  \end{gather*}
  that is, $t_G^{\delta,\epsilon}(U,W)$ is also a $(\widetilde{\delta},\widetilde{\epsilon})$-majority opinion of $U$ with
  respect to $W$ in $G$, hence $U$ is also $(\widetilde{\delta},\widetilde{\epsilon})$-excellent in $G$.

  \medskip

  Item~\ref{lem:mono:atomtuple} has a similar proof to that of item~\ref{lem:mono:exc}: given an $m$-tuple $x\in V(G)^m$, since
  $U$ is a $(\delta,\epsilon,m)$-atom of $G$, we know that $x$ does not $(\delta,\epsilon)$-split $U$ in $G$, which means that
  there must exist $b\in\{0,1\}$ such that
  \begin{gather*}
    \{u\in U\mid t_G^\delta(u,x) = b\}
  \end{gather*}
  is strictly smaller than $\epsilon\cdot\lvert U\rvert$. By item~\ref{lem:mono:maj}, the set above contains the set
  \begin{gather*}
    \{u\in U\mid t_G^{\widetilde{\delta}}(u,x) = b\},
  \end{gather*}
  so this latter set is strictly smaller than $\epsilon\cdot\lvert U\rvert\leq\widetilde{\epsilon}\cdot\lvert U\rvert$, which
  means that $x$ also does not $(\widetilde{\delta},\widetilde{\epsilon})$-split $U$ in $G$, hence $U$ is a
  $(\widetilde{\delta},\widetilde{\epsilon},m)$-atom of $G$>

  \medskip

  Item~\ref{lem:mono:atom} follows directly from item~\ref{lem:mono:atomtuple}.
\end{proof}

\begin{lemma}[Sauer--Shelah--Perles]\label{lem:SSP}
  Let $\cH\subseteq\cP(X)$ be a family of subsets of some set $X$ with $\VC(\cH)\leq d < \infty$. Then for a finite set
  $A\subseteq X$ with $\lvert A\rvert=n\in\NN$, we have
  \begin{gather*}
    \lvert\cH\rest_A\rvert
    \leq
    \sum_{i=0}^d \binom{n}{i}.
  \end{gather*}

  In particular, if $d\leq n/2$, then
  \begin{gather*}
    \lvert\cH\rest_A\rvert \leq 2^{h_2(d/n)\cdot n},
  \end{gather*}
  where
  \begin{gather*}
    h_2(p) \df p\cdot\log_2\frac{1}{p} + (1-p)\log_2\frac{1}{1-p}
  \end{gather*}
  is the binary entropy.
\end{lemma}

\begin{lemma}[Independent blow-ups]\label{lem:indepblowup}
  Let $G$ be a graph and $b\colon V(G)\to\NN$, then
  \begin{align*}
    \OP(G^b) & \leq \OP(G), &
    \Lit(G^b) & \leq \Lit(G), &
    \VC(G^b) & \leq \VC(G).
  \end{align*}
  Furthermore, equality holds when $b(v)\geq 1$ for every $v\in V(G)$.
\end{lemma}

\begin{proof}
  If we have either an order property, a Littlestone tree or a shattered set in $G^b$ then the corresponding vertices, or
  vertices that generate neighborhoods as hypotheses can be distinguished by neighborhoods, i.e., for any two distinct vertices
  $v,w$ within the pattern, there exists $u\in V(G)\setminus\{v,w\}$ such that $u$ is adjacent to exactly one of $v$ and $w$,
  provided either both $v$ and $w$ come from the vertex side of the pattern or both come from the hypotheses side of the
  pattern.

  This means means that each vertex in the pattern must be a copy of a different vertex of the original graph $G$, hence the
  original graph $G$ must also contain the same pattern. Hence the inequalities hold.

  The equality statement follows since $b(v)\geq 1$ for every $v\in V(G)$ implies that $G$ is an induced subgraph of $G^b$,
\end{proof}

\begin{lemma}[Hodges~\protect{\cite[Lemma~6.7.9]{Hod93}}]\label{lem:Hodges}
  For $\cH\subseteq\cP(X)$, we have
  \begin{align*}
    \OP(\cH) & \leq 2^{\Lit(\cH)+1}-2, &
    \Lit(\cH) & \leq 2^{\OP(\cH)+3} - 2.
  \end{align*}
\end{lemma}

\begin{lemma}[Folklore]\label{lem:LitOPleqVC}
  For $\cH\subseteq\cP(X)$, we have
  \begin{gather*}
    \VC(\cH)\leq\max\{\Lit(\cH),\OP(\cH)\}.
  \end{gather*}
\end{lemma}

\begin{proof}
  If $\cH$ shatters a set $A$ of size $d$, then for every $U\subseteq A$, there exists $G_U\in\cH$ with $G_U\cap A = U$.

  We can then construct a $d$-order property $(H,x)$ in $\cH$ as follows: pick $x_1,\ldots,x_d$ in $A$ distinct and take $H_i\df
  G_{\{x_i,\ldots,x_d\}}$.

  Similarly, we can construct a $d$-Littlestone tree $(H,x)$ as follows: pick $a_1,\ldots,a_d$ in $A$ distinct, set
  $x_\sigma\df a_{\lvert\sigma\rvert}$ and for every $\tau\in\{0,1\}^d$, we let $H_\tau\df G_{\{a_i \mid \tau_i=1\}}$.
\end{proof}

The lemma below shows that when we consider $\delta$-majority opinions of a bounded $\VC$-dimension graph with
$\delta\in(0,1/(d+1))$, the $\VC$-dimension stays the same. We refer the reader to~\cite{MMM25} for a more thorough study of
the $\VC$ and Littlestone dimension of $\delta$-majority opinions.
\begin{lemma}\label{lem:VCgood}
  Let $G$ be a graph with $\VC(G)\leq d < \infty$, let $\delta\in(0,1/(d+1))$ and define
  \begin{gather*}
    \cH_{G,\delta}
    \df
    \{N_{G,\delta}(U) \mid U\text{ is $\delta$-good in $G$}\}
    \subseteq
    \cP\bigl(V(G)\bigr),
  \end{gather*}
  where
  \begin{gather*}
    N_{G,\delta}(U) \df \{v\in V(G) \mid t_G^\delta(v,U) = 1\}.
  \end{gather*}

  Then $\VC(\cH_{G,\delta})\leq d$.
\end{lemma}

\begin{proof}
  Let $A\subseteq V(G)$ be a set that is shattered by $\cH_{G,\delta}$ and suppose for a contradiction that $\lvert A\rvert =
  d+1$. Let us show that $A$ is also shattered by $\cH_G\df\{N_G(v) \mid v\in V(G)\}$, which will then contradict
  $\VC(G)=\VC(\cH_G)\leq d$.

  Since $A\subseteq V(G)$ is shattered by $\cH_{G,\delta}$, then for every $B\subseteq A$, there exists a $\delta$-good set $U_B$
  such that $N_{G,\delta}(U_B)\cap A = B$, which in particular implies that for every $a\in A$, we have $t_G^\delta(a,U_B) = \One_B(a)$.

  Now note
  \begin{align*}
    \left\lvert U_B\cap\bigcap_{a\in A} N_G^{\One_B(a)}(a)\right\rvert
    & \geq
    \lvert U_B\rvert - \sum_{a\in A} \lvert U_B\cap N_G^{1 - \One_B(a)}(a)\rvert
    \geq
    \lvert U_B\rvert - \sum_{a\in A}\delta\cdot\lvert U_B\rvert
    \\
    & =
    \lvert U_B\rvert\cdot\bigl(1 - (d+1)\cdot\delta\bigr)
    >
    0,
  \end{align*}
  hence there must exist a vertex $u_B\in U_B\cap\bigcap_{a\in A} N_G^{\One_B(a)}(a)$, from which it follows that $N_G(u_B)\cap
  A = B$, hence $\cH_G$ shatters $A$.
\end{proof}

\begin{lemma}[a slight improvement over Terry~\protect{\cite[Lemma~2.3]{Ter26}} using a trick from Fox--Pach--Suk~\protect{\cite[Lemma~8]{FPS19}}]\label{lem:good->hom}
  Let $G$ be a graph, let $\epsilon\in(0,1)$ and let $U_1,U_2\subseteq V(G)$ be $\epsilon$-good in $G$ sets. Then the following
  hold:
  \begin{enumerate}
  \item\label{lem:good->hom:hom} The pair $(U_1,U_2)$ is $\epsilon'$-homogeneous, where
    \begin{gather}\label{eq:good->hom:epsilon}
      \begin{aligned}
        \epsilon'
        & \df
        \begin{dcases*}
          \frac{1 - \sqrt{1 - 8\cdot\epsilon\cdot(1-\epsilon)}}{2}, & if $\epsilon < (2-\sqrt{2})/4$,\\
          1/2, & otherwise,
        \end{dcases*}
        \\
        & =
        2\cdot\epsilon + O_{\epsilon\to 0}(\epsilon^2).
      \end{aligned}
    \end{gather}
  \item\label{lem:good->hom:exc} There exists an $(\epsilon,\epsilon'')$-majority opinion $t_G^{\epsilon,\epsilon''}(U_1,U_2)$ of $U_1$
    with respect to $U_2$ in $G$, where
    \begin{gather*}
      \epsilon''
      \df
      \frac{\epsilon'}{1-\epsilon}
      =
      \frac{1 - \sqrt{1 - 8\cdot\epsilon\cdot(1-\epsilon)}}{2\cdot(1-\epsilon)}
      =
      2\cdot\epsilon + O_{\epsilon\to 0}(\epsilon^2).
    \end{gather*}
  \item\label{lem:good->hom:partition} If $\cP$ is an $\epsilon$-good partition of $G$, then $\cP$ is a totally
    $\epsilon'$-homogeneous partition of $G$ and an $(\epsilon,\epsilon'')$-excellent partition of $G$.
  \end{enumerate}
\end{lemma}

\begin{proof}
  The results are trivial when $\epsilon\geq(2-\sqrt{2})/4$ as in this case we have $\epsilon',\epsilon''\geq 1/2$. So suppose
  $\epsilon < (2-\sqrt{2})/4$.
  
  We prove item~\ref{lem:good->hom:hom} by the contra-positive: suppose $(U_1,U_2)$ is not $\epsilon'$-homogeneous, that is, we
  have $d_G(U_1,U_2)\in(\epsilon',1-\epsilon')$. Let us pick $\rn{a}_1,\rn{b}_1$ independently uniformly at random in $U_1$ and
  $\rn{a}_2,\rn{b}_2$ independently (and independently from the previous random vertices) uniformly at random in $U_2$. Consider
  the events
  \begin{align*}
    E & \df \{\rn{a}_1,\rn{a}_2\}\in E(G)\land\{\rn{b}_1,\rn{b}_2\}\notin E(G),
    \\
    E_1 & \df \{\rn{a}_1,\rn{a}_2\}\in E(G)\land\{\rn{a}_1,\rn{b}_2\}\notin E(G),
    \\
    E_2 & \df \{\rn{a}_1,\rn{b}_2\}\in E(G)\land\{\rn{b}_1,\rn{b}_2\}\notin E(G).
  \end{align*}
  Clearly $E_1$ and $E_2$ are disjoint and the event $E$ implies the disjunction of the events $E_1$ and $E_2$, so we have
  \begin{align*}
    \PP[E_1] + \PP[E_2]
    =
    \PP[E_1\lor E_2]
    \geq
    \PP[E]
    =
    d_G(U_1,U_2)\cdot (1-d_G(U_1,U_2))
    >
    \epsilon'\cdot(1-\epsilon').
  \end{align*}

  Thus, there must exist $i\in[2]$ such that $\PP[E_i]\geq\epsilon'\cdot(1-\epsilon')/2$. If
  $\PP[E_1]\geq\epsilon'\cdot(1-\epsilon')/2$, by conditioning on $\rn{a}_1$, we conclude that there must exist $a_1\in U_1$ such
  that
  \begin{gather*}
    \frac{\lvert N_G(a_1)\cap U_2\rvert}{\lvert U_2\rvert}
    \cdot
    \left(1 - \frac{\lvert N_G(a_1)\cap U_2\rvert}{\lvert U_2\rvert}\right)
    >
    \frac{\epsilon'\cdot(1-\epsilon')}{2}.
  \end{gather*}
  On the other hand, if $\PP[E_2]\geq\epsilon'\cdot(1-\epsilon')/2$, by conditioning on $\rn{b}_2$, we conclude that there exists
  $b_2\in U_2$ such that
  \begin{gather*}
    \frac{\lvert N_G(b_2)\cap U_1\rvert}{\lvert U_1\rvert}
    \cdot
    \left(1 - \frac{\lvert N_G(b_2)\cap U_1\rvert}{\lvert U_1\rvert}\right)
    >
    \frac{\epsilon'\cdot(1-\epsilon')}{2}.
  \end{gather*}
  In either case, it then follows that there exists $i\in[2]$ and $v\in U_i$ with
  \begin{gather*}
    \frac{\lvert N_G(v)\cap U_{3-i}\rvert}{\lvert U_{3-i}\rvert}
    \cdot
    \left(1 - \frac{\lvert N_G(v)\cap U_{3-i}\rvert}{\lvert U_{3-i}\rvert}\right)
    >
    \frac{\epsilon'\cdot(1-\epsilon')}{2},
  \end{gather*}
  from which we conclude that
  \begin{gather*}
    \frac{\lvert N_G(v)\cap U_{3-i}\rvert}{\lvert U_{3-i}\rvert}
    \in
    \left(
    \frac{1 - \sqrt{1 - 2\cdot\epsilon'\cdot(1-\epsilon')}}{2},
    1 - \frac{1 - \sqrt{1 - 2\cdot\epsilon'\cdot(1-\epsilon')}}{2},
    \right)
    =
    (\epsilon,1-\epsilon),
  \end{gather*}
  hence $U_{3-i}$ is not $\epsilon$-good in $G$.

  \medskip

  We now prove item~\ref{lem:good->hom:exc}. By item~\ref{lem:good->hom:hom}, we know that $(U_1,U_2)$ is
  $\epsilon'$-homogeneous (where $\epsilon'$ is given by~\eqref{eq:good->hom:epsilon}), so there exists $b =
  b_G^{\epsilon'}(U_1,U_2)\in\{0,1\}$ such that $d_G^b(U_1,U_2)\geq 1 - \epsilon'$. We will show that $b$ is a
  $(\epsilon,\epsilon'')$-majority opinion $t_G^{\epsilon,\epsilon''}$ of $U_1$ with respect to $U_2$ in $G$.

  Let us pick a vertex $\rn{v}$ uniformly at random from $U_1$ and note that
  \begin{gather*}
    \EE\left[\frac{\lvert N_G^{1-b}\rn{v}\cap U_2\rvert}{\lvert U_2\rvert}\right]
    =
    d_G^{1 - b}(U_1,U_2)
    \leq
    \epsilon'.
  \end{gather*}
  so by Markov's Inequality, we have
  \begin{gather*}
    \PP\left[\frac{\lvert N_G^{1-b}\rn{v}\cap U_2\rvert}{\lvert U_2\rvert}\geq 1-\epsilon\right]
    \leq
    \frac{\epsilon'}{1 - \epsilon}
    =
    \epsilon''.
  \end{gather*}

  Since $U_2$ is $\epsilon$-good in $G$, we have
  \begin{gather*}
    \{v\in U_1 \mid t_G^\epsilon(v) \neq b\}
    =
    \{v\in U_1 \mid \lvert N_G^{1-b}(v)\cap U_2\rvert\geq 1 - \epsilon\}
  \end{gather*}
  and from our earlier derivation, this set has size at most $\epsilon''\cdot\lvert U_1\rvert$, so we conclude that
  \begin{gather*}
    \lvert\{v\in U_1 \mid t_G^\epsilon(v) = b\}\rvert\geq(1-\epsilon'')\cdot\lvert U_1\rvert,
  \end{gather*}
  hence $b$ is a $(\epsilon,\epsilon'')$-majority opinion $t_G^{\epsilon,\epsilon''}$ of $U_1$ with respect to $U_2$ in $G$.

  \medskip

  Finally, for item~\ref{lem:good->hom:partition}, the assertion that $\cP$ is totally $\epsilon'$-homogeneous follows directly
  from item~\ref{lem:good->hom:hom}. For the excellence assertion, note that if $P\in\cP$ and $U\subseteq V(G)$ is
  $\epsilon$-good in $G$, then by item~\ref{lem:good->hom:exc}, there exists an $(\epsilon,\epsilon'')$-majority opinion
  $t_G^{\epsilon,\epsilon''}(P,U)$ of $P$ with respect to $U$ in $G$. Since $U$ is arbitrary, it follows that $P$ is
  $(\epsilon,\epsilon'')$-excellent in $G$.
\end{proof}

\begin{remark}
  The bounds of Lemma~\ref{lem:good->hom}\ref{lem:good->hom:hom} are asymptotically tight. Namely, consider a bipartite
  graph $G$ given by
  \begin{gather*}
    V(G) \df V_{11}\cup V_{12}\cup V_{21}\cup V_{22},\\
    E(G) \df \{\{u,v\} \mid u\in V_{11}, v\in V_{21}\}\cup\{\{u,v\} \mid u\in V_{12}, v\in V_{22}\},
  \end{gather*}
  where $V_{11},V_{12},V_{21},V_{22}$ are pairwise disjoint sets of sizes
  \begin{align*}
    \lvert V_{11}\rvert = \lvert V_{21}\rvert & = (1-\xi)\cdot n,
    &
    \lvert V_{12}\rvert = \lvert  V_{22}\rvert & = \xi\cdot n,
  \end{align*}
  where $\xi\in(0,1)\in\QQ$ and $n\in\NN$ is large enough so that the quantities above are integers.

  Then for $U_1\df V_{11}\cup V_{12}$ and $U_2\df V_{21}\cup V_{22}$, we have
  \begin{gather*}
    d_G(U_1,U_2)
    =
    \frac{\lvert V_{11}\times V_{21}\rvert + \lvert V_{12}\times V_{22}\rvert}{\lvert U_1\rvert\cdot\lvert U_2\rvert}
    =
    (1 - \xi)^2 + \xi^2
    =
    1 - 2\cdot\xi + 2\cdot\xi^2,
  \end{gather*}
  so for $\epsilon'\in(0,1)$, the pair $(U_1,U_2)$ is $\epsilon'$-homogeneous in $G$ if and only if $\epsilon'\geq 2\cdot\xi -
  2\cdot\xi^2/2 = 2\cdot\xi - O_{\xi\to 0}(\xi^2)$.

  On the other hand, we note that for every $i\in[2]$ and every $v\in V_{i2}$, we have
  \begin{gather*}
    \frac{\lvert N_G(v)\cap U_{3-i}\rvert}{\lvert U_{3-i}\rvert}
    =
    \frac{\lvert V_{3-i,2}\rvert}{\lvert U_{3-i}\rvert}
    =
    \xi,
  \end{gather*}
  and for every $v\in V_{i1}$, we have
  \begin{gather*}
    \frac{\lvert N_G(v)\cap U_{3-i}\rvert}{\lvert U_{3-i}\rvert}
    =
    \frac{\lvert V_{3-i,1}\rvert}{\lvert U_{3-i}\rvert}
    =
    1-\xi.
  \end{gather*}
  Since there are no edges completely contained in $U_{3-i}$, it follows that for $\epsilon\in(0,1/2)$, the set $U_{3-i}$ is
  $\epsilon$-good in $G$ if and only if $\epsilon\geq\xi$.

  Thus, we conclude that for every $\epsilon\in(0,1/2)$, we have that $\epsilon$-goodness of both $U_1$ and $U_2$ in $G$
  \emph{cannot} imply $\epsilon'$-homogeneity for any $\epsilon' < 2\cdot\epsilon - 2\cdot\epsilon^2 = 2\epsilon -
  O_{\epsilon\to 0}(\epsilon^2)$.
\end{remark}

\begin{theorem}[(Weak version of) uniform convergence property]\label{thm:UC}
  There exists an absolute constant $C > 0$ such that for every $d\in\NN_+$, every countable set $X$, every family of sets
  $\cH\subseteq\cP(X)$ with $\VC(\cH)\leq d$, every $\delta,\epsilon > 0$, every probability measure $\mu$ over $X$ and every
  positive integer
  \begin{gather}\label{eq:UC:m}
    m \geq C\cdot\frac{d + \ln(1/\delta)}{\epsilon^2},
  \end{gather}
  we have
  \begin{gather*}
    \PP_{\rn{x}\sim\mu^m}\Bigl[
      \forall H\in\cH,
      \bigl\lvert
      \mu(H) - m^{-1}\lvert\{i\in[m]\mid \rn{x}_i\in H\}\rvert
      \bigr\rvert
      \leq
      \epsilon
      \Bigr]
    \geq
    1 - \delta.
  \end{gather*}
\end{theorem}

\begin{corollary}\label{cor:UC}
  There exists an absolute constant $C > 0$ such that for every $d\in\NN_+$, every countable set $X$, every family of sets
  $\cH\subseteq\cP(X)$ with $\VC(\cH)\leq d$, every $\epsilon > 0$, every probability measure $\mu$ over $X$ and every positive
  integer
  \begin{gather*}
    m\geq C\cdot\frac{d}{\epsilon^2},
  \end{gather*}
  there exists an $m$-tuple $(x_1,\ldots,x_m)\in X^m$ such that
  \begin{gather*}
    \bigl\lvert
    \mu(H) - m^{-1}\lvert\{i\in[m]\mid x_i\in H\}\rvert
    \bigr\rvert
    \leq
    \epsilon
  \end{gather*}
  for every $H\in\cH$.
\end{corollary}

\begin{proof}
  Apply Theorem~\ref{thm:UC} with any fixed $\delta < 1$ (say $\delta=1/2$) and let $(x_1,\ldots,x_m)$ be an outcome of the
  random $\rn{x}$ in the event of probability $1-\delta > 0$.
\end{proof}

\begin{lemma}\label{lem:atom->exc}
  Let $G$ be a graph, let $\delta,\epsilon > 0$ and let $W\subseteq V(G)$ be an $(\delta,\epsilon)$-atom of $G$. Then $W$ is
  $(\delta,\epsilon)$-excellent in $G$.
\end{lemma}

\begin{proof}
  For a $\delta$-good in $G$ set $U$, enumerate its elements as $U = \{x_1,\ldots,x_m\}$ ($m\df\lvert U\rvert$) and consider the
  tuple $x=(x_1,\ldots,x_m)$. Then for every $w\in V(G)$, we clearly have $t_G^\delta(w,U) = t_G^\delta(w,x)$.

  If there is no $(\delta,\epsilon)$-majority opinion of $W$ with respect to $U$ in $G$, then for
  every $b\in\{0,1\}$, we
  must have
  \begin{gather*}
    \frac{\lvert\{w\in W \mid t_G^\delta(w,x) = b\}\rvert}{\lvert W\rvert}
    =
    \frac{\lvert\{w\in W \mid t_G^\delta(w,U) = b\}\rvert}{\lvert W\rvert}
    >
    \epsilon,
  \end{gather*}
  so $W$ is $(\delta,\epsilon)$-split by $x$ in $G$, contradicting $(\delta,\epsilon)$-atomicity of $W$ in $G$.
\end{proof}

The next lemma informally says that in a graph of finite $\VC$-dimension, $(c\cdot\delta,\epsilon)$-atomicity follows from
$(\delta,\epsilon,m)$-atomicity for $m = O(\VC(G)/((1-c)\cdot\delta)^2)$, that is, we only need to check tuples of length that
is at most constant (i.e., the only dependence on $G$ is via its $\VC$-dimension, but not its size).

\begin{lemma}[Atomicity via constant-length tuples]\label{lem:atomconst}
  There exists an absolute constant $C > 0$ such that for every $d\in\NN_+$, every graph $G$ with $\VC(G)\leq d$, every
  $\delta,\epsilon > 0$, every $c_{\atom}\in(0,1)$ and every integer positive
  \begin{gather*}
    m \geq C\cdot\frac{d}{c_{\atom}^2\cdot\delta^2},
  \end{gather*}
  if $W$ is a $((1+c_{\atom})\cdot\delta,\epsilon,m)$-atom of $G$ set then $W$ is $(\delta,\epsilon)$-atom of $G$.
\end{lemma}

\begin{proof}
  Let $C$ be the absolute constant of Corollary~\ref{cor:UC} and let us prove the result by the contra-positive: suppose
  $W\subseteq V(G)$ is not a $(\delta,\epsilon)$-atom of $G$, that is, there exists a tuple $z=(z_1,\ldots,z_n)$ of vertices of
  $V(G)$ of some length $n\in\NN_+$ such that $W$ is $(\delta,\epsilon)$-split by $z$ in $G$.

  Consider the probability measure $\mu$ on $\im(z) = \{z_1,\ldots,z_n\}$ that assigns probability proportional to the number of
  times the vertex appears in the tuple $z$, that is, for each $v\in\im(z)$, we have $\mu(\{v\}) = \lvert\{i\in[n] \mid z_i =
  v\}\rvert/n$.

  We now apply Corollary~\ref{cor:UC} with $(\im(z),\cH_G\rest_{\im(z)},c_{\atom}\cdot\delta)$ as the parameters
  $(X,\cH,\epsilon)$ (hence our $m$ satisfies the conditions of the corollary) so that there exists an $m$-tuple
  $(x_1,\ldots,x_m)\in\im(z)^m$ such that for every $H\in\cH_G\rest_{\im(z)}$, we have
  \begin{gather*}
    \bigl\lvert
    \mu(H) - m^{-1}\lvert\{i\in[m]\mid x_i\in H\}\rvert
    \bigr\rvert
    \leq
    c_{\atom}\cdot\delta.
  \end{gather*}

  If we instantiate the above with $H=N_G(w)\cap\im(z)$ for some $w\in W$, we get $\lvert d_G(w,z) - d_G(w,x)\rvert \leq
  c_{\atom}\cdot\delta$, which implies $\lvert d_G^b(w,z) - d_G^b(w,x)\rvert \leq c_{\atom}\cdot\delta$ for every $b\in\{0,1\}$.
  In particular, if $w\in W$ is such that $z$ has a $\delta$-majority opinion $t_G^\delta(w,z)$ with respect to $w$ in $G$, then
  we have
  \begin{gather*}
    d_G^{t_G^{c_{\atom}\cdot\delta}(w,z)}(w,x)
    \geq
    d_G^{t_G^{c_{\atom}\cdot\delta}(w,z)}(w,z) - c_{\atom}\cdot\delta
    \geq
    1 - \delta - c_{\atom}\cdot\delta
    =
    1 - (1+c_{\atom})\cdot\delta,
  \end{gather*}
  that is, $x$ has a $(1+c_{\atom})\cdot\delta$-majority opinion $t_G^{(1+c_{\atom})\cdot\delta}(w,x)$ with respect to $w$ in
  $G$ and it is the same as the $\delta$-majority opinion of $z$ with respect to $w$ in $G$, i.e., we have
  $t_G^{(1+c_{\atom})\cdot\delta}(w,x) = t_G^\delta(w,z)$. Thus, we have shown the following containments
  \begin{align*}
    \{w\in W \mid t_G^\delta(w,z) = 0\} & \subseteq \{w\in W\mid t_G^{(1+c_{\atom})\cdot\delta}(w,x) = 0\},\\
    \{w\in W \mid t_G^\delta(w,z) = 1\} & \subseteq \{w\in W\mid t_G^{(1+c_{\atom})\cdot\delta}(w,x) = 1\}.
  \end{align*}

  Since we are assuming that $W$ is $(\delta,\epsilon)$-split by $z$, the sets on the left-hand sides of the above have size
  larger than $\epsilon\cdot\lvert W\rvert$ each, which means that the sets on the right-hand sides of the above have size
  larger than $\epsilon\cdot\lvert W\rvert$ each, hence $W$ is $(1+c_{\atom})\cdot\delta\cdot\lvert W\rvert$ each, that is, $W$
  is $((1+c_{\atom})\cdot\delta,\epsilon)$-split by $x$, hence $W$ is not a $((1+c_{\atom})\cdot\delta,\epsilon)$-atom of $G$.
\end{proof}

\begin{lemma}[Upward continuity]\label{lem:upcont}
  Let $\delta,\epsilon\in(0,1)$, let $m\in\NN_+$, let $G$ be a graph, let $U\subseteq\widetilde{U}\subseteq V(G)$ be non-empty
  sets and let
  \begin{gather*}
    \epsilon_+
    \df
    (\epsilon-1)\cdot\frac{\lvert U\rvert}{\lvert\widetilde{U}\rvert}
    +
    1
    \leq
    \epsilon + \frac{\lvert\widetilde{U}\setminus U\rvert}{\lvert\widetilde{U}\rvert}
  \end{gather*}
  Then the following hold:
  \begin{enumerate}
  \item\label{lem:upcont:good} If $U$ is $\epsilon$-good in $G$, then $\widetilde{U}$ is $\epsilon_+$-good in $G$.
  \item\label{lem:upcont:excellent} If $U$ is $(\delta,\epsilon)$-excellent in $G$, then $\widetilde{U}$ is
    $(\delta,\epsilon_+)$-excellent in $G$.
  \item\label{lem:upcont:atom} If $U$ is a $(\delta,\epsilon,m)$-atom of $G$, then $\widetilde{U}$ is a
    $(\delta,\epsilon_+,m)$-atom of $G$.
  \end{enumerate}
\end{lemma}

\begin{proof}
  We start by proving the bound on $\epsilon_+$:
  \begin{gather*}
    \epsilon_+
    \df
    (\epsilon-1)\cdot\frac{\lvert U\rvert}{\lvert\widetilde{U}\rvert}
    +
    1
    =
    \epsilon\cdot\frac{\lvert U\rvert}{\lvert\widetilde{U}\rvert}
    +
    \frac{\lvert\widetilde{U}\setminus U\rvert}{\lvert\widetilde{U}\rvert}
    \leq
    \epsilon
    +
    \frac{\lvert\widetilde{U}\setminus U\rvert}{\lvert\widetilde{U}\rvert}.
  \end{gather*}

  Note that item~\ref{lem:upcont:good} follows from item~\ref{lem:upcont:atom} with $m=1$ (and any $\delta\in(0,1)$). We prove
  the other two items essentially simultaneously.
  
  Let $X$ be either a $\delta$-good set in $G$ or an $m$-tuple of vertices of $G$ and for $b\in\{0,1\}$, let
  \begin{gather*}
    d_b^U(X) \df \frac{\lvert\{u\in U \;\bigm\vert\; t_G^\delta(u,X) = b\}\rvert}{\lvert U\rvert}
  \end{gather*}
  and define $d_b^{\widetilde{U}}(X)$ analogously.

  Note that
  \begin{gather*}
    d_b^{\widetilde{U}}(X)
    \leq
    \frac{d_b^U(X)\cdot\lvert U\rvert + \lvert\widetilde{U}\setminus U\rvert}{\lvert\widetilde{U}\rvert}
    =
    (d_b^U(X)-1)\cdot\frac{\lvert U\rvert}{\lvert\widetilde{U}\rvert} + 1.
  \end{gather*}
  Now, if $U$ is $(\delta,\epsilon)$-excellent in $G$ and $X$ is a $\delta$-good set in $G$, then for $b\df 1 -
  t_G^{\delta,\epsilon}(U,X)$ (i.e., the opposite of the $(\delta,\epsilon)$-majority opinion of $U$ with respect to $X$ in
  $G$), we have $d_b^U(X)\leq\epsilon$, which from the above implies $d_b^{\widetilde{U}}(X)\leq\epsilon_+$, so $1-b$ is also a
  $(\delta,\epsilon_+)$-majority opinion of $U$ with respect to $X$ in $G$. Hence $\widetilde{U}$ is
  $(\delta,\epsilon_+)$-excellent in $G$.

  On the other hand, if $U$ is a $(\delta,\epsilon,m)$-atom of $G$ and $X$ is an $m$-tuple of vertices of $G$, then there exists
  $b\in\{0,1\}$ with $d_b^U(X)\leq\epsilon$ (as $U$ is not $(\delta,\epsilon)$-split by $X$ in $G$), so the above implies
  $d_b^{\widetilde{U}}(X)\leq\epsilon_+$, hence $\widetilde{U}$ is not $(\delta,\epsilon_+)$-split by $X$ in $G$. Hence
  $\widetilde{U}$ is a $(\delta,\epsilon_+,m)$-atom of $G$.
\end{proof}

\begin{lemma}[Downward continuity]\label{lem:downcont}
  Let $\delta,\epsilon\in(0,1)$, let $m\in\NN_+$, let $G$ be a graph, let $\widetilde{U}\subseteq U\subseteq V(G)$ be non-empty
  sets and let
  \begin{gather*}
    \epsilon_-
    \df
    \epsilon\cdot\frac{\lvert U\rvert}{\lvert\widetilde{U}\rvert}
    \leq
    \epsilon + \frac{\lvert U\setminus\widetilde{U}\rvert}{\lvert\widetilde{U}\rvert}.
  \end{gather*}
  Then the following hold:
  \begin{enumerate}
  \item\label{lem:downcont:good} If $U$ is $\epsilon$-good in $G$, then $\widetilde{U}$ is $\epsilon_-$-good in $G$.
  \item\label{lem:downcont:excellent} If $U$ is $(\delta,\epsilon)$-excellent in $G$, then $\widetilde{U}$ is
    $(\delta,\epsilon_-)$-excellent in $G$.
  \item\label{lem:downcont:atom} If $U$ is a $(\delta,\epsilon,m)$-atom of $G$, then $\widetilde{U}$ is a
    $(\delta,\epsilon_-,m)$-atom of $G$.
  \end{enumerate}
\end{lemma}

\begin{proof}
  We start by proving the bound on $\epsilon_-$:
  \begin{gather*}
    \epsilon_-
    \df
    \epsilon\cdot\frac{\lvert U\rvert}{\lvert\widetilde{U}\rvert}
    =
    \epsilon + \epsilon\cdot\frac{\lvert U\setminus\widetilde{U}\rvert}{\lvert\widetilde{U}\rvert}
    \leq
    \epsilon + \frac{\lvert U\setminus\widetilde{U}\rvert}{\lvert\widetilde{U}\rvert}.
  \end{gather*}

  Similarly to Lemma~\ref{lem:upcont}, item~\ref{lem:downcont:good} follows from item~\ref{lem:downcont:atom} with $m=1$ (and any
  $\delta\in(0,1)$) and for $X$ being either a $\delta$-good set in $G$ or an $m$-tuple of vertices of $G$ and for
  $b\in\{0,1\}$, we let
  \begin{gather*}
    d_b^U(X) \df \frac{\lvert\{u\in U \;\bigm\vert\; t_G^\delta(u,X) = b\}\rvert}{\lvert U\rvert}
  \end{gather*}
  and define $d_b^{\widetilde{U}}(X)$ analogously.

  Then we have
  \begin{gather*}
    d_b^{\widetilde{U}}(X)
    \leq
    \frac{d_b^U(X)\cdot\lvert U\rvert}{\lvert\widetilde{U}\rvert}.
  \end{gather*}
  This means that if $U$ is $(\delta,\epsilon)$-excellent in $G$ and $X$ is a $\delta$-good set in $G$, then for $b\df 1 -
  t_G^{\delta,\epsilon}(U,X)$ we have $d_b^U(X)\leq\epsilon$ so that the above yields $d_b^{\widetilde{U}}\leq\epsilon_-$, hence
  $1-b$ is also a $(\delta,\epsilon_-)$-majority opinion of $\widetilde{U}$ with respect to $X$ in $G$, hence $\widetilde{U}$ is
  $(\delta,\epsilon_-)$-excellent in $G$.

  On the other hand, if $U$ is a $(\delta,\epsilon,m)$-atom of $G$, then there exists $b\in\{0,1\}$ with $d_b^U(X)\leq\epsilon$,
  so the above implies $d_b^{\widetilde{U}}(X)\leq\epsilon_-$ hence $\widetilde{U}$ is a $(\delta,\epsilon_-,m)$-atom of $G$.
\end{proof}

\begin{corollary}[Continuity]\label{cor:cont}
  Let $\delta,\epsilon\in(0,1)$, let $m\in\NN_+$, let $G$ be a graph, let $\widetilde{U},U\subseteq V(G)$ be non-empty sets that
  are comparable (i.e., either $U\subseteq\widetilde{U}$ or $U\supseteq\widetilde{U}$) and let
  \begin{gather*}
    \epsilon_{\symdiff}
    \df
    \epsilon + \frac{\lvert U\symdiff\widetilde{U}\rvert}{\lvert\widetilde{U}\rvert}.
  \end{gather*}
  Then the following hold:
  \begin{enumerate}
  \item\label{cor:cont:good} If $U$ is $\epsilon$-good in $G$, then $\widetilde{U}$ is $\epsilon_{\symdiff}$-good in $G$.
  \item\label{cor:cont:excellent} If $U$ is $(\delta,\epsilon)$-excellent in $G$, then $\widetilde{U}$ is
    $(\delta,\epsilon_{\symdiff})$-excellent in $G$.
  \item\label{cor:cont:atom} If $U$ is a $(\delta,\epsilon,m)$-atom of $G$, then $\widetilde{U}$ is a
    $(\delta,\epsilon_{\symdiff},m)$-atom of $G$.
  \end{enumerate}
\end{corollary}

\begin{proof}
  Follows by putting together Lemmas~\ref{lem:upcont} and~\ref{lem:downcont}.
\end{proof}


\section{The extraction process}

The next proposition says that every set in a bounded Littlestone dimension graph contains a linear-sized subset that is good,
excellent or atom; it also covers the algorithmic complexity of extracting one such subset in several different models of
computation.

\begin{proposition}[Extraction]\label{prop:ext}
  Let $n,m,d,\ell\in\NN_+$ with $d\leq\ell\leq\log_2(n)$, let
  $\widetilde{\delta},\widetilde{\epsilon},c_{\UC},\widetilde{\rho}_{\UC}\in(0,1)$, let $G$ be a graph with $\lvert G\rvert=n$,
  $\Lit(G)\leq\ell$ and $\VC(G)\leq d$ and let $U\subseteq V(G)$ be a non-empty set of vertices of $G$. Let further
  \begin{gather*}
    \gamma
    \df
    (1-c_{\UC})\cdot\widetilde{\epsilon}
  \end{gather*}
  and suppose $s_0,\ldots,s_\ell,s_{\ell+1}\in\NN_+$ are such that $\lvert U\rvert\geq s_0$
  \begin{gather}\label{eq:ext:si}
    s_i\leq\floor{\gamma\cdot s_{i-1}}+1
  \end{gather}
  for every $i\in[\ell+1]$. Finally, let $S\df\{s_0,\ldots,s_\ell\}$ (this does \emph{not} include $s_{\ell+1}$).

  Then the following hold:
  \begin{enumerate}
  \item\label{prop:ext:good} (Malliaris--Shelah~\cite{MS14}) There exists $W\subseteq U$ with $\lvert W\rvert\in S$ that is
    $\widetilde{\epsilon}$-good in $G$.
  \item\label{prop:ext:exc} (Malliaris--Shelah~\cite{MS14}) If $\widetilde{\delta} < 1/2^{\ell+1}$, then there exists $W\subseteq U$
    with $\lvert W\rvert\in S$ that is $(\widetilde{\delta},\widetilde{\epsilon})$-excellent in $G$.
  \item\label{prop:ext:atom} If $\widetilde{\delta} < 1/2^{\ell+1}$, then there exists $W\subseteq U$ with $\lvert W\rvert\in S$
    that is a $(\widetilde{\delta},\widetilde{\epsilon},m)$-atom of $G$.
  \item\label{prop:ext:goodalg} Suppose $c_{\UC} > 0$ and
    \begin{gather*}
      \widetilde{m}\df\Floor{\left(1 - \frac{c_{\UC}}{2}\right)\cdot\widetilde{\epsilon}\cdot m},
      \\
      m
      \geq
      \frac{C\cdot 4\cdot (d + \ln((2^{\ell+1}-1)\cdot n/\widetilde{\rho}_{\UC}))}{(c_{\UC}^2\cdot\widetilde{\epsilon}^2)},
    \end{gather*}
    where $C$ is the absolute constant of Theorem~\ref{thm:UC}. Then in the random-query model, Algorithm~\ref{alg:ext:good},
    with probability at least $1-\widetilde{\rho}_{\UC}$, computes $(W,s)$, where $W\subseteq U$ is a
    $\widetilde{\epsilon}$-good in $G$ with $\lvert W\rvert=s\in S$. The time-complexity of
    Algorithm~\ref{alg:ext:good} is:
    \begin{gather*}
      O(2^\ell\cdot n\cdot m)
    \end{gather*}
  \item\label{prop:ext:atomalg} If $\widetilde{\delta} < 1/2^{\ell+1}$ and $\widetilde{m}\df\floor{\widetilde{\delta}\cdot m}$,
    then in the deterministic model, Algorithm~\ref{alg:ext:atom} computes $(W,s)$, where $W\subseteq U$ is a
    $(\delta,\widetilde{\epsilon},m)$-atom of $G$ with $\lvert W\rvert=s\in S$. The time-complexity of
    Algorithm~\ref{alg:ext:atom} is:
    \begin{gather*}
      O\bigl((\ell+m)\cdot n^{m+1}\bigr).
    \end{gather*}
  \item\label{prop:ext:atomspace} If $\widetilde{\delta} < 1/2^{\ell+1}$, then in the query-into-oracle model,
    Algorithms~\ref{alg:ext:atomspace} and~\ref{alg:ext:atomspaceoracle} compute $(W,s)$, where $W\subseteq U$ is a
    $(\delta,\widetilde{\epsilon},m)$-atom of $G$ with $\lvert W\rvert=s\in S$. The space-complexities of
    Algorithms~\ref{alg:ext:atomspace} and~\ref{alg:ext:atomspaceoracle} are
    \begin{gather*}
      O\bigl(\ell\cdot m\cdot\log(n+1)\bigr),
      \\
      O\bigl(\log(\lvert\sigma\rvert+1) + \log(m+1) + \log(n+1)\bigr)
      \leq
      O\bigl(\log(\ell) + \log(m+1) + \log(n+1)\bigr),
    \end{gather*}
    respectively, the time-complexities are
    \begin{align*}
      O(2^\ell\cdot m\cdot n^m),
      & &
      O\bigl((\lvert\sigma\rvert+1)\cdot m\cdot n\bigr)
      \leq
      O(\ell\cdot m\cdot n),
    \end{align*}
    respectively, and the $U$-oracle-complexities are at most
    \begin{align*}
      (2^{\ell+1} - 1)\cdot n^{m+1},
      & &
      n,
    \end{align*}
    respectively.
  \end{enumerate}
\end{proposition}

Before we start the proof, let us mention that the assumption $d\leq\ell\leq\log_2(n)$ here is innocuous due to
$\VC(G)\leq\Lit(G)\leq\log_2(\lvert G\rvert)$ (see Lemma~\ref{lem:LitOPleqVC}). We also point out that although we do not make
this assumption, one could also assume that $\widetilde{\epsilon}\geq 1/n$ as the notions of $\widetilde{\epsilon}$-good,
$(\widetilde{\delta},\widetilde{\epsilon})$-good and $(\widetilde{\delta},\widetilde{\epsilon},m)$-atom remain unchanged if we
replace $\widetilde{\epsilon}$ with $\max\{\widetilde{\epsilon},1/\lvert G\rvert\}$.

Finally, let us note that none of the algorithms of Proposition~\ref{prop:ext} actually receive the values of
$\widetilde{\delta}$ or $\widetilde{\epsilon}$, they instead receive integer parameters that satisfy the correct properties with
respect to some choice of $\widetilde{\delta}$ and $\widetilde{\epsilon}$; this is important because sometimes we want to run
these algorithms with potentially irrational values of $\widetilde{\delta}$ and $\widetilde{\epsilon}$.

\begin{proof}
  Items~\ref{prop:ext:good} and~\ref{prop:ext:exc} are from Malliaris--Shelah~\cite{MS14} (in the hardest case when
  $c_{\UC}=0$), but they are consequences of item~\ref{prop:ext:atom}:

  To see how item~\ref{prop:ext:good} follows from item~\ref{prop:ext:atom}, we can take $m=1$ and any $\widetilde{\delta} <
  1/2^{\ell+1} < 1$ so that the notion of $(\widetilde{\delta},\widetilde{\epsilon},1)$-atom of $G$ coincides with the notion of
  $\widetilde{\epsilon}$-good in $G$.

  To see how item~\ref{prop:ext:exc} follows from item~\ref{prop:ext:atom}, we first note that since $\widetilde{\delta} <
  1/2^{\ell+1}$, we can pick $c_{\atom}\in(0,1)$ small enough so that
  $(1+c_{\atom})\cdot\widetilde{\delta} < 1/2^{\ell+1}$; we then let
  \begin{gather*}
    m \df \Ceil{C\cdot\frac{d}{c_{\atom}^2\cdot\widetilde{\delta}^2}},
  \end{gather*}
  where $C > 0$ is the absolute constant of Lemma~\ref{lem:atomconst} so that this lemma guarantees that every
  $((1+c_{\atom})\cdot\widetilde{\delta},\widetilde{\epsilon},m)$ atom of $G$ is a
  $(\widetilde{\delta},\widetilde{\epsilon})$-atom of $G$, and in turn is $(\widetilde{\delta},\widetilde{\epsilon})$-excellent
  in $G$ by Lemma~\ref{lem:atom->exc}. Now item~\ref{prop:ext:atom} gives the desired
  $((1+c_{\atom})\cdot\widetilde{\delta},\widetilde{\epsilon},m)$-atom of $G$ that is in particular
  $(\widetilde{\delta},\widetilde{\epsilon})$-excellent in $G$.

  Item~\ref{prop:ext:atom} itself follows from its algorithmic counterpart, item~\ref{prop:ext:atomalg}.

  \medskip

  Before we actually start the proof of the algorithmic items, let us observe that in all items except for
  item~\ref{prop:ext:goodalg} we can make the further assumption that $c_{\UC}=0$ (which also means that $\gamma =
  \widetilde{\epsilon}$). Indeed, the result for $c_{\UC} = 0$ implies the result for $c_{\UC} > 0$; this is because if
  $(s_i)_{i=1}^\ell$ satisfy~\eqref{eq:ext:si} for some $c_{\UC} > 0$, then they also satisfy~\eqref{eq:ext:si} for $c_{\UC}=0$
  (the $\gamma$ computed with $c_{\UC} > 0$ is less than $\widetilde{\epsilon}$, which is the $\gamma$ computed with
  $c_{\UC}=0$). As such, all of our algorithms (except for Algorithm~\ref{alg:ext:good} corresponding to
  item~\ref{prop:ext:goodalg}) will tacitly assume that $c_{\UC}=0$ and $\gamma\df\widetilde{\epsilon}$.

  One more comment is that our algorithms don't actually use the values of $c_{\UC}$, $\widetilde{\delta}$ or
  $\widetilde{\epsilon}$, so they do not receive them as input. In fact, except for Algorithm~\ref{alg:ext:good}, all other
  algorithms claim their atomicity or goodness for a particular value of $\widetilde{\epsilon}$. Namely, the value claimed is:
  \begin{gather*}
    \widehat{\epsilon} \df \max_{i\in[\ell+1]} \frac{s_i-1}{s_{i-1}}.
  \end{gather*}
  The reason for this is because we have
  \begin{align*}
    \forall i\in[\ell+1], s_i\geq\floor{\widetilde{\epsilon}\cdot s_{i-1}} + 1
    & \iff
    \forall i\in[\ell+1], s_i-1\geq \widetilde{\epsilon}\cdot s_{i-1}
    \\
    & \iff
    \widetilde{\epsilon} \geq \max_{i\in[\ell+1]}\frac{s_i-1}{s_{i-1}} = \widehat{\epsilon},
  \end{align*}
  so $\widehat{\epsilon}$-goodness ($(\widetilde{\delta},\widehat{\epsilon},m)$-atomicity, respectively) implies
  $\widetilde{\epsilon}$-goodness ($(\widetilde{\delta},\widetilde{\epsilon},m)$-atomicity, respectively). For the case of
  Algorithm~\ref{alg:ext:good}, we won't be able to tacitly an optimal $\widetilde{\epsilon}$, so we will instead assume that a
  correct $\widetilde{\epsilon}$ that is not given in the input exists. The fact that none of these algorithms actually receives
  $\widetilde{\epsilon}$ allows us to run them for implicit values of $\widetilde{\delta}$ or $\widetilde{\epsilon}$ that are
  not necessarily rational (as we have already done so in the proof of Item~\ref{prop:ext:atom}).

  We now prove the algorithmic items. We will do so slightly out of order for ease of presentation.
  \begin{description}[wide, itemsep={3ex}]
  \item[Item~\ref{prop:ext:atomalg}.] The idea of Algorithm~\ref{alg:ext:atom} is of a ``failed $(\ell+1)$-Littlestone tree''
    just as in the original Malliaris--Shelah argument~\cite{MS14}, except that since we are computing a
    $(\widetilde{\delta},\widetilde{\epsilon},m)$-atom of $G$, we need to consider all $m$-tuples of vertices of $G$.

    \begin{algorithm}[htbp]
      \caption{Deterministic model algorithm that returns $(W,s)$, where $W\subseteq U$ is a
        $(\widetilde{\delta},\widetilde{\epsilon},m)$-atom of $G$ with $\lvert W\rvert=s\in S\df\{s_0,\ldots,s_\ell\}$. The
        time-complexity is $O((\ell+m)\cdot n^{m+1})$.}
      \label{alg:ext:atom}
      \DontPrintSemicolon
      \KwIn{Numbers $n,m,\ell,s_0,\ldots,s_\ell,s_{\ell+1}\in\NN_+$, $\widetilde{m}\in\NN$ with $\ell\leq\log_2(n)$, a graph $G$
        with $\lvert G\rvert=n$ and $\Lit(G)\leq\ell$ and a non-empty set $U\subseteq V(G)$ with $\lvert U\rvert\geq s_0$. We
        further assume that there exists $\widetilde{\delta}\in(0,1/2^{\ell+1})$ such that $\widetilde{m} =
        \floor{\widetilde{\delta}\cdot m}$.}
      \KwOut{A pair $(W,s)$, where $W\subseteq U$ is a $(\widetilde{\delta},\widetilde{\epsilon},m)$-atom of $G$ with $\lvert
        W\rvert=s\in S\df\{s_0,\ldots,s_\ell\}$ and
        \begin{gather*}
          \widetilde{\epsilon}
          \df
          \max_{i\in[\ell+1]}\frac{s_i-1}{s_{i-1}}.
        \end{gather*}
      }
      Let $U_\varnothing$ be a subset of $U$ of size $s_0$.\;\label{alg:ext:atom:Uvarnothing}
      \For{$x\in V(G)^m$}{%
        \label{alg:ext:atom:prexFor}
        \For{$v\in V(G)$}{%
          \label{alg:ext:atom:prevFor}
          $d_G(v,x)\assign\lvert\{i\in[m]\mid x_i\in N_G(v)\}\rvert$\;
          \tcp*[l]{We use the names $t_G^{\widetilde{\delta}}$ for ease of proof; for the algorithm these are simply variables}
          \lIf{$m-d_G(v,x)\leq\widetilde{m}$}{$t_G^{\widetilde{\delta}}(v,x)\assign 1$}
          \lElseIf{$d_G(v,x)\leq\widetilde{m}$}{$t_G^{\widetilde{\delta}}(v,x)\assign 0$}
          \tcp*[l]{Else: leave $t_G^{\widetilde{\delta}}(v,x)$ undefined.}
        }
      }
      \For{$t\in\{0,\ldots,\ell\}$}{%
        \label{alg:ext:atom:tFor}
        \For{$\tau\in\{0,1\}^t$}{%
          \label{alg:ext:atom:tauFor}
          \For{$x\in V(G)^m$}{%
            \label{alg:ext:atom:xFor}
            $k^0\assign\lvert\{u\in U_\tau \mid t_G^{\widetilde{\delta}}(u,x) = 0\}\rvert$\;
            $k^1\assign\lvert\{u\in U_\tau \mid t_G^{\widetilde{\delta}}(u,x) = 1\}\rvert$\;
            \If(\tcp*[f]{$U_\tau$ seems $(\widetilde{\delta},\widetilde{\epsilon})$-split by $x$}){%
              $k^0\geq s_{t+1}$ and $k^1\geq s_{t+1}$
            }{%
              \label{alg:ext:atom:kIf}
              $x_\tau\assign x$\;
              \Break
            }
          }
          \uIf{$x_\tau$ is defined}{%
            \label{alg:ext:atom:xtauIf}
            Let $U_{\tau\conc 0}$ and $U_{\tau\conc 1}$ sets of size $s_{t+1}$ with
            \begin{gather*}
              U_{\tau\conc i}\subseteq\{u\in U_\tau \mid t_G^{\widetilde{\delta}}(u,x) = i\}
              \qquad (i\in\{0,1\}).
            \end{gather*}
          }
          \lElse{\label{alg:ext:atom:success}
            \Return{$(U_\tau,s_t)$.}
          }
        }
      }
      \Return{\Failure.}\tcp*[r]{This return statement never actually executes}
      \label{alg:ext:atom:totalfailure}
    \end{algorithm}

    We start by arguing correctness, that is, let us prove that Algorithm~\ref{alg:ext:atom} actually returns $(W,s)$ with
    $W\subseteq U$ a $(\widetilde{\delta},\widetilde{\epsilon},m)$-atom of $G$ and $\lvert W\rvert=s\in S$. First, we note that
    all sets of the form $U_\tau$ defined by the algorithm are defined either in line~\ref{alg:ext:atom:Uvarnothing} or in
    line~\ref{alg:ext:atom:xtauIf} and in both cases, the algorithm ensures that they have size $\lvert U_\tau\rvert =
    s_{\lvert\tau\rvert}$, so if the algorithm returns $(U_\tau,s_t)$, then we indeed have $\lvert U_\tau\rvert=s_t\in S$.

    Thus, we need to prove that the algorithm never reaches the failure return statement of line~\ref{alg:ext:atom:totalfailure}
    and that when it returns $(U_\tau,s_t)$, then $U_\tau$ is indeed a $(\widetilde{\delta},\widetilde{\epsilon})$-atom of $G$
    and $U_\tau\subseteq U$.

    We start by noting that all the definitions of the $U_\tau$ are indeed valid (i.e., we can indeed find such subsets of the
    required sizes): finding $U_\varnothing\subseteq U$ of size $s_0$ in line~\ref{alg:ext:atom:Uvarnothing} is possible due to
    our assumption that $\lvert U\rvert\geq s_0$ and finding the sets $U_{\tau\conc 0}$ and $U_{\tau\conc 1}$ of size $s_{t+1}$
    in line~\ref{alg:ext:atom:xtauIf} follows since for that line to execute, we must have $x_\tau$ defined, which only happens
    if the quantities
    \begin{align*}
      k^0
      & \df
      \lvert\{u\in U_\tau \mid t_G^{\widetilde{\delta}}(u,x) = 0\}\rvert
      \\
      k^1
      & \df
      \lvert\{u\in U_\tau \mid t_G^{\widetilde{\delta}}(u,x) = 1\}\rvert
    \end{align*}
    computed in the same iteration are both at least $s_{t+1}$, which means that the corresponding sets $U_{\tau\conc i}$ indeed
    exist as subsets of the sets above of size $s_{t+1}$.

    We now claim that the block of line~\ref{alg:ext:atom:prexFor} correctly computes all the $\widetilde{\delta}$-majority
    opinions $t_G^{\widetilde{\delta}}(v,x)$ of all $m$-tuples $x\in V(G)^m$ with respect to all vertices $v\in V(G)$ in $G$
    when they exist (and leaves $t_G^{\widetilde{\delta}}(v,x)$ undefined when there is no $\widetilde{\delta}$-majority
    opinion; recalling also that since $\widetilde{\delta} < 1/2^{\ell+1}\leq 1/2$ there can only be at most one such
    $\widetilde{\delta}$-majority opinion). Indeed, this is because we assume that
    \begin{gather*}
      \widetilde{m} = \floor{\widetilde{\delta}\cdot m}
    \end{gather*}
    and for $d_G(v,x)\df\lvert\{i\in[m]\mid x_i\in N_G(v)\}\rvert\in\NN$, we have
    \begin{gather*}
      d_G(v,x) \geq (1-\widetilde{\delta})\cdot m
      \iff
      m-d_G(v,x) \leq \widetilde{\delta}\cdot m
      \iff
      m-d_G(v,x) \leq \widetilde{m},
      \\
      d_G(v,x) \leq \widetilde{\delta}\cdot m
      \iff
      d_G(v,x) \leq \widetilde{m}.
    \end{gather*}

    We now argue that line~\ref{alg:ext:atom:totalfailure} never actually executes (which in particular means that the algorithm
    indeed returns one of the $(U_\tau,s_\tau)$ with $\lvert\tau\rvert\in\{0,\ldots,\ell\}$). Suppose not, then it must be that
    for every $\tau\in\{0,1\}^{<\ell+1}$, the $m$-tuple $x_\tau$ gets defined (otherwise line~\ref{alg:ext:atom:success} would
    execute and end the algorithm), which in turn means that the algorithm also defines a set $U_\tau$ for every
    $\tau\in\{0,1\}^{\ell+1}$. Since we have already shown that $\lvert U_\tau\rvert=s_{\ell+1}\geq 1$, it follows that $U_\tau$
    is non-empty, so we can let $u_\tau\in U_\tau$.

    By construction of the sets $U_\tau$, for every $t\in[\ell+1]$, we have $U_\tau\subseteq U_{\tau\rest_{[t]}}$, which in
    particular implies that
    \begin{gather*}
      t_G^{\widetilde{\delta}}(u_\tau, x_{\tau\rest_{[t-1]}}) = \tau_t
    \end{gather*}
    from the definition of the set $U_{\tau\rest_{[t]}}$. The $\widetilde{\delta}$-majority opinion above naturally means that
    \begin{gather}\label{eq:ext:atom:totalfailure}
      \Bigl\lvert\Bigl\{i\in[m] \;\Bigm\vert\;
      u_\tau\in N_G^{\tau_t}\bigl((x_{\tau\rest_{[t-1]}})_i\bigr)
      \Bigr\}
      \geq
      (1-\widetilde{\delta})\cdot m.
    \end{gather}

    We claim that for each $\sigma\in\{0,1\}^{<\ell+1}$, there exists an index $i_\sigma\in[m]$ such that for every
    $\tau\in\{0,1\}^{\ell+1}$ that extends $\sigma$, we have
    \begin{gather*}
      u_\tau\in N_G^{\tau_{\lvert\sigma\rvert+1}}\bigl((x_\sigma)_{i_\sigma}\bigr).
    \end{gather*}
    Indeed, we can simply consider the set of all $i\in[m]$ that satisfy the above conditions for each such $\tau$. Note that
    there are exactly $2^{\ell+1-\lvert\sigma\rvert}\leq 2^{\ell+1}$ many $\tau\in\{0,1\}^{\ell+1}$ extending $\sigma$, so by a
    union bound applied to~\eqref{eq:ext:atom:totalfailure}, the set of all $i\in[m]$ that satisfy the desired conditions has
    size at least
    \begin{gather*}
      \bigl(1 - 2^{\ell+1}\cdot\widetilde{\delta}\bigr)\cdot m,
    \end{gather*}
    which is positive since $\widetilde{\delta} < 1/2^{\ell+1}$ by assumption, guaranteeing the existence of our desired
    $i_\sigma\in[m]$.

    But then $(((x_\sigma)_{i_\sigma})_{\sigma\in\{0,1\}^{<\ell+1}},(N_G(u_\tau))_{\tau\in\{0,1\}^{\ell+1}})$ is an
    $(\ell+1)$-Littlestone tree in $\cH_G$, contradicting the fact that $\Lit(G) = \Lit(\cH_G)\leq\ell$.

    Thus, line~\ref{alg:ext:atom:totalfailure} never actually executes, which means that the algorithm
    line~\ref{alg:ext:atom:success} must get executed and the algorithm must return some $(U_\tau,s_\tau)$ with
    $\lvert\tau\rvert\in\{0,\ldots,\ell\}$. Let us now argue that this $U_\tau$ must be a
    $(\widetilde{\delta},\widetilde{\epsilon},m)$-atom of $G$. Suppose not, then it must be the case that there exists an
    $m$-tuple $x\in V(G)^m$ such that $U_\tau$ is $(\widetilde{\delta},\widetilde{\epsilon})$-split by $x$ in $G$, that is, for
    every $b\in\{0,1\}$, we must have
    \begin{gather*}
      \lvert\{u\in U_\tau \mid t_G^{\widetilde{\delta}}(u,x) = b\}\rvert
      >
      \widetilde{\epsilon}\cdot\lvert U_\tau\rvert
      =
      \widetilde{\epsilon}\cdot s_{\lvert\tau\rvert},
    \end{gather*}
    and since $s_{\lvert\tau\rvert+1}\leq\floor{\widetilde{\epsilon}\cdot s_{\lvert\tau\rvert}}+1$ and the left-hand side of the
    above is an integer, it follows that
    \begin{gather*}
      \lvert\{u\in U_\tau \mid t_G^{\widetilde{\delta}}(u,x) = b\}\rvert \geq s_{\lvert\tau\rvert+1}.
    \end{gather*}
    But this means that in the iteration of the loops of lines~\ref{alg:ext:atom:tauFor} and~\ref{alg:ext:atom:xFor}
    corresponding to $\tau$ and $x$, the test of line~\ref{alg:ext:atom:kIf} must have passed, which means that $x_\tau$ gets
    defined, but this contradicts the fact that $U_\tau$ gets returned by the algorithm on line~\ref{alg:ext:atom:success}
    executing due to $x_\tau$ not being defined.

    Therefore, Algorithm~\ref{alg:ext:atom} returns a $(\widetilde{\delta},\widetilde{\epsilon},m)$-atom of $G$.

    \smallskip

    We now argue the time-complexity of Algorithm~\ref{alg:ext:atom}. Then the loops of lines~\ref{alg:ext:atom:prexFor}
    and~\ref{alg:ext:atom:prevFor} run exactly $n^m$ and $n$ times, respectively. Each inner loop iteration takes time $O(m)$,
    so the total time-complexity of this initial block is at most $O(n^{m+1}\cdot m)$.

    The loops of lines~\ref{alg:ext:atom:tFor}, \ref{alg:ext:atom:tauFor} and~\ref{alg:ext:atom:xFor} have sizes $\ell+1$, $2^t$
    and $n^m$, respectively. Each iteration of the inner loop takes time $O(\lvert U_\tau\rvert)$.

    Now note that since $\widetilde{\delta} < 1/2^{\ell+1}\leq 1/2$, the definition of the sets $U_\tau$ in the block of
    line~\ref{alg:ext:atom:xtauIf} guarantees that $U_{\tau\conc 0}$ is disjoint from $U_{\tau\conc 1}$. With a simple
    induction, we conclude that for every $t\in\{0,\ldots,\ell+1\}$, the sets in $(U_\tau)_{\tau\in\{0,1\}^t}$ are pairwise
    disjoint, so we get $\sum_{\tau\in\{0,1\}^t}\lvert U_\tau\rvert\leq n$. This means that the total time-complexity of the
    loops of lines~\ref{alg:ext:atom:tFor}, \ref{alg:ext:atom:tauFor} and~\ref{alg:ext:atom:xFor} is at most
    \begin{gather*}
      O\left(\sum_{i=0}^\ell \sum_{\tau\in\{0,1\}^t} n^m\cdot\lvert U_\tau\rvert\right)
      \leq
      O\bigl((\ell+1)\cdot n^{m+1}\bigr)
    \end{gather*}
    so the final upper bound in the time-complexity is
    \begin{gather*}
      O(n^{m+1}\cdot m)
      +
      O\bigl((\ell+1)\cdot n^{m+1}\bigr)
      =
      O\bigl((\ell+m)\cdot n^{m+1}\bigr).
    \end{gather*}
  \item[Item~\ref{prop:ext:goodalg}.] The idea of Algorithm~\ref{alg:ext:good} is the same as the one in~\ref{alg:ext:atom}
    (with $m=1$ as we only want $\widetilde{\epsilon}$-good sets), except that since we cannot afford to spend time $O(n)$
    computing sizes of neighborhoods of vertices as that would yield a quadratic (in $n$) algorithm. Instead, we will estimate
    the size of neighborhoods with the help of uniform convergence (Theorem~\ref{thm:UC}) and the random-query oracle.
    \begin{algorithm}[htbp]
      \caption{Algorithm in the random-query model of computation that, with probability at least $1-\widetilde{\rho}_{\UC}$,
        returns $(W,s)$, where $W\subseteq U$ is $\widetilde{\epsilon}$-good in $G$ with $\lvert W\rvert=s\in
        S\df\{s_0,\ldots,s_\ell\}$. The time-complexity is:
        \begin{gather*}
          O(2^\ell\cdot n\cdot m).
        \end{gather*}
      }
      \label{alg:ext:good}
      \DontPrintSemicolon
      \KwIn{Numbers $n,m,\ell,s_0,\ldots,s_\ell,s_{\ell+1}\in\NN_+$, $\widetilde{m}\in\NN$ with $\ell\leq\log_2(n)$, a
        random-query oracle $\cO_G$ for a graph $G$ with $\lvert G\rvert=n$ and $\Lit(G)\leq\ell$ and a non-empty set
        $U\subseteq V(G)$ with $\lvert U\rvert\geq s_0$. We further assume that there exist
        $c_{\UC},\widetilde{\epsilon},\widetilde{\rho}_{\UC}\in(0,1)$ and $d\in[\ell]$ such that for
        $\gamma\df(1-c_{\UC})\cdot\widetilde{\epsilon}$, we have
        \begin{gather*}
          \VC(G)\leq d,
          \\
          \forall i\in[\ell+1], s_i\geq\floor{\gamma\cdot s_{i-1}}+1,
          \\
          \widetilde{m} = \Floor{\left(1 - \frac{c_{\UC}}{2}\right)\cdot\widetilde{\epsilon}\cdot m},
          \\
          m
          \geq
          \frac{C\cdot 4\cdot (d + \ln((2^{\ell+1}-1)\cdot n/\widetilde{\rho}_{\UC}))}{(c_{\UC}^2\cdot\widetilde{\epsilon}^2)},
        \end{gather*}
        where $C$ is the absolute constant of Theorem~\ref{thm:UC}.}
      \KwOut{A pair $(W,s)$, where $W\subseteq U$ is $\widetilde{\epsilon}$-good in $G$ with $\lvert W\rvert=s\in
        S\df\{s_0,\ldots,s_\ell\}$. With probability at most $\widetilde{\rho}_{\UC}$, may instead return ``\Failure'' or a pair
        $(W,s)$ in which $W$ is not $\widetilde{\epsilon}$-good (but we still have $\lvert W\rvert=s\in S$).}
      Let $U_\varnothing$ be a subset of $U$ of size $s_0$.\;\label{alg:ext:good:Uvarnothing}
      \For{$t\in\{0,\ldots,\ell\}$}{%
        \label{alg:ext:good:tFor}
        \For{$\tau\in\{0,1\}^t$}{%
          \label{alg:ext:good:tauFor}
          \For{$v\in V(G)$}{%
            \label{alg:ext:good:vFor}
            Consult $\cO_G(v,U_\tau)$ a total of $m$ times, let $k$ be the number of times it returns ``adjacent''.\;
            \label{alg:ext:good:oraclecall}
            \If(\tcc*[f]{$U_\tau$ seems to not have $\widetilde{\epsilon}$-majority opinion with respect to $v$}){%
              $\widetilde{m} < k < m - \widetilde{m}$
            }{%
              \label{alg:ext:good:kIf}
              $v_\tau\assign v$\;
              \Break
            }
          }
          \uIf{$v_\tau$ is defined}{%
            \label{alg:ext:good:vtauIf}
            Find sets $U_{\tau\conc 0}$ and $U_{\tau\conc 1}$ of size $s_{t+1}$ with $U_{\tau\conc b}\subseteq N_G^b(v_\tau)\cap U_\tau$\;
            \tcp*[l]{$N_G^b(v_\tau)$ can be computed consulting $\cO_G(v_\tau,\{u\})$}
            \If(\tcp*[f]{$N_G(v_\tau)\cap U_\tau$ incorrectly estimated}){%
              \label{alg:ext:good:incorrectnhd}
              At least one of $U_{\tau\conc 0}$ or $U_{\tau\conc 1}$ do not exist.
            }{%
              \Return{\Failure.}
            }
          }
          \lElse{\label{alg:ext:good:success}
            \Return{$(U_\tau,s_\tau)$.}
          }
        }
      }
      \Return{\Failure.}\tcp*[f]{Some neighborhood was incorrectly estimated along the way}
      \label{alg:ext:good:totalfailure}
    \end{algorithm}

    First, we argue the complexity of Algorithm~\ref{alg:ext:good}. The loops of lines~\ref{alg:ext:good:tFor},
    \ref{alg:ext:good:tauFor} and~\ref{alg:ext:good:vFor} have sizes $\ell+1$, $2^t$ and $n$, respectively. Each iteration of
    the inner loop takes time $O(m)$. Furthermore, the conditional block of line~\ref{alg:ext:good:vtauIf} takes time $O(n)$, so
    it follows that the time-complexity of the loop of line~\ref{alg:ext:good:tFor} (and hence of Algorithm~\ref{alg:ext:good})
    is
    \begin{gather*}
      O\left(
      \sum_{i=0}^\ell 2^t\cdot n\cdot m
      \right)
      =
      O(2^\ell\cdot n\cdot m).
    \end{gather*}

    \smallskip

    Now we argue correctness, that is, let us prove that with probability at least $1-\widetilde{\rho}_{\UC}$,
    Algorithm~\ref{alg:ext:good} returns $(W,s)$, where $W\subseteq U$ is $\widetilde{\epsilon}$-good in $G$ with $\lvert
    W\rvert=s\in S\df\{s_0,\ldots,s_\ell\}$. First we note that it is clear that every set $U_\tau$ defined by the algorithm
    must satisfy $\lvert U_\tau\rvert=s_{\lvert\tau\rvert}$. Furthermore, if the algorithm returns $(W,s)$, then it must be of
    the form $(U_\tau,s_{\lvert\tau\rvert})$ for some $\tau\in\{0,1\}^{<\ell+1}$.

    Now let us consider the event $E$ that each time line~\ref{alg:ext:good:oraclecall} is executed with some $(v,U_\tau)$, the
    number $k$ it computes satisfies
    \begin{gather}\label{eq:ext:goodoracle}
      \left\lvert\frac{k}{m} - \frac{\lvert N_G(v)\cap U_\tau\rvert}{\lvert U_\tau\rvert}\right\rvert
      \leq
      \frac{c_{\UC}}{2}\cdot\widetilde{\epsilon}.
    \end{gather}

    By Theorem~\ref{thm:UC} with parameters
    \begin{gather*}
      (d,\widetilde{\delta},\epsilon)
      \df
      \left(d,\frac{\widetilde{\rho}_{\UC}}{2^\ell\cdot n},\frac{c_{\UC}}{2}\cdot\widetilde{\epsilon}\right)
    \end{gather*}
    (and since our assumption on $m$ is precisely the hypothesis~\eqref{eq:UC:m} of Theorem~\ref{thm:UC} with these parameters
    and $\VC(\cH_G)=\VC(G)\leq d$), it follows that a single execution of line~\ref{alg:ext:good:oraclecall}
    satisfies~\eqref{eq:ext:goodoracle} with probability at least $1 - \widetilde{\rho}_{\UC}/((2^{\ell+1}-1)\cdot n)$. Since we
    have already argued that line~\ref{alg:ext:good:oraclecall} is executed at most $(2^{\ell+1}-1)\cdot n$ times, it follows by
    the union bound that the event $E$ has probability at least $1-\widetilde{\rho}_{\UC}$.

    \smallskip

    We now argue that within the event $E$ the algorithm never passes the test of line~\ref{alg:ext:good:incorrectnhd}. For
    this, note that for the line~\ref{alg:ext:good:incorrectnhd} to even execute, it must be that $v_\tau$ is defined. In turn,
    this means that it suffices to show that for every $b\in\{0,1\}$, the set $N_G^b(v_\tau)\cap U_\tau$ has size at least
    $s_{\lvert\tau\rvert+1}$, which itself is equivalent to
    \begin{gather}
      s_{\lvert\tau\rvert+1} \leq \lvert N_G(v_\tau)\cap U_\tau\rvert \leq s_{\lvert\tau\rvert} - s_{\lvert\tau\rvert+1}
    \end{gather}
    Since $s_{\lvert\tau\rvert+1}\leq\floor{\gamma\cdot s_{\lvert\tau\rvert}} + 1$, it suffices to show that
    \begin{gather}
      \gamma\cdot s_{\lvert\tau\rvert} < \lvert N_G(v_\tau)\cap U_\tau\rvert < (1-\gamma)\cdot s_{\lvert\tau\rvert}.
    \end{gather}
    Since $v_\tau$ is defined we know that the test of line~\ref{alg:ext:good:kIf} must have passed in the same iteration, that
    is, we must have had $\widetilde{m} < k < m - \widetilde{m}$. In turn, since we are in the event $E$, we know that we must
    have
    \begin{gather*}
      \left\lvert\frac{k}{m} - \frac{\lvert N_G(v)\cap U_\tau\rvert}{\lvert U_\tau\rvert}\right\rvert
      \leq
      \frac{c_{\UC}}{2}\cdot\widetilde{\epsilon},
    \end{gather*}
    which in turn means that
    \begin{gather*}
      \left(\frac{\widetilde{m}}{m} - \frac{c_{\UC}}{2}\cdot\widetilde{\epsilon}\right)\cdot\lvert U_\tau\rvert
      <
      \lvert N_G(v)\cap U_\tau\rvert
      <
      \left(1 - \frac{\widetilde{m}}{m} + \frac{c_{\UC}}{2}\cdot\widetilde{\epsilon}\right)\cdot\lvert U_\tau\rvert,
    \end{gather*}
    and since $\widetilde{m}=\floor{(1-c_{\UC}/2)\cdot\widetilde{\epsilon}\cdot m}$, we get
    \begin{gather*}
      \gamma\cdot s_{\lvert\tau\rvert}
      =
      (1-c_{\UC})\cdot\widetilde{\epsilon}\cdot\lvert U_\tau\rvert
      <
      \lvert N_G(v)\cap U_\tau\rvert
      <
      \bigl(1 - (1-c_{\UC})\cdot\widetilde{\epsilon}\bigr)\cdot\lvert U_{\lvert\tau\rvert}\rvert
      =
      (1-\gamma)\cdot s_{\lvert\tau\rvert},
    \end{gather*}
    as desired.

    \smallskip

    Now we claim that line~\ref{alg:ext:good:totalfailure} is never reached within the event $E$. Suppose not, then it must be
    the case that for every $\tau\in\{0,1\}^{<\ell+1}$, the vertex $v_\tau$ gets defined (otherwise
    line~\ref{alg:ext:good:success} would execute and end the algorithm), which in turn means that the algorithm also defines a
    set $U_\tau$ for every $\tau\in\{0,1\}^{\ell+1}$. Since $\lvert U_\tau\rvert=s_{\ell+1}\geq 1$, the set $U_\tau$ is
    non-empty, so we can let $u_\tau\in U_\tau$. Note that by construction of the sets $U_\tau$, for every $t\in[\ell+1]$, we
    have
    \begin{gather*}
      U_\tau\subseteq U_{\tau\rest_{[t]}}\subseteq N_G^{\tau_t}(v_{\tau\rest_{[t-1]}}),
    \end{gather*}
    which means that $((v_\sigma)_{\sigma\in\{0,1\}^{<\ell+1}}, (N_G(u_\tau))_{\tau\in\{0,1\}^{\ell+1}})$ is an
    $(\ell+1)$-Littlestone tree in $\cH_G$, contradicting the fact that $\Lit(G)=\Lit(\cH_G)\leq\ell$.

    Thus, line~\ref{alg:ext:good:totalfailure} is never reached within the event $E$, which means that within this event,
    line~\ref{alg:ext:good:success} must get executed and the algorithm must return some $(U_\tau,s_\tau)$. Let us now argue
    that this $U_\tau$ must be $\widetilde{\epsilon}$-good in $G$. Suppose not, then it must be the case that there exists $v\in
    V(G)$ such that there is no $\widetilde{\epsilon}$-majority opinion of $U$ with respect to $v$ in $G$, that is, we must have
    \begin{gather*}
      \widetilde{\epsilon} < \frac{\lvert N_G(v)\cap U_\tau\rvert}{\lvert U_\tau\rvert} < 1 - \widetilde{\epsilon}.
    \end{gather*}

    Let $k$ be the value computed calling the oracle $\cO_G$ in line~\ref{alg:ext:good:oraclecall} in the iteration
    corresponding to $(\tau,v)$ (i.e., the oracle got called $m$ times as $\cO_G(v,U_\tau)$). Since we are within the event $E$,
    we know that
    \begin{gather*}
      \left\lvert\frac{k}{m} - \frac{\lvert N_G(v)\cap U_\tau\rvert}{\lvert U_\tau\rvert}\right\rvert
      \leq
      \frac{c_{\UC}}{2}\cdot\widetilde{\epsilon},
    \end{gather*}
    from which we conclude that
    \begin{gather*}
      \left(1 - \frac{c_{\UC}}{2}\right)\cdot\widetilde{\epsilon}
      <
      \frac{k}{m}
      <
      1 - \left(1 - \frac{c_{\UC}}{2}\right)\cdot\widetilde{\epsilon},
    \end{gather*}
    hence (since $\widetilde{m}=\floor{(1-c_{\UC}/2)\cdot\widetilde{\epsilon}}$), we have
    \begin{gather*}
      \widetilde{m}
      <
      k
      <
      m - \widetilde{m}
    \end{gather*}
    and in turn implies that $v_\tau$ gets defined due to the test of line~\ref{alg:ext:good:kIf} passing, but this contradicts
    the fact that $(U_\tau,s_\tau)$ gets returned by the algorithm on line~\ref{alg:ext:good:success} executing due to $v_\tau$
    not being defined.

    Therefore, within the event $E$, Algorithm~\ref{alg:ext:good} returns $(W,s)$, where $W\subseteq U$ is
    $\widetilde{\epsilon}$-good in $G$ set with $\lvert W\rvert=s\in S=\{s_0,\ldots,s_\ell\}$.
  \item[Item~\ref{prop:ext:atomspace}.] Before we start the formal argument, let us mention that we will make the assumption
    that the vertex set of the graph is $[n]$; this has two important effects: names of vertices cost $O(\log(n+1))$ space and
    the vertex set comes equipped with the natural ordering of $[n]$, which will play a crucial role in our algorithm (as the
    algorithm will constantly think of the $i$th vertex in some set).

    As is typical with space-efficient algorithms, several subsets of $[n]$ are computed implicitly (as storing one such set in
    memory would take too much space, namely $\Theta(n)$) in the sense that the algorithm will only need to keep in memory a
    small ``code'' that describes it. Furthermore, we will need way more such sets (about $2^\ell$) than our space-complexity
    budget allows for (about $\ell+1$), so our algorithm will not even keep in memory at the same time the ``codes'' of all sets
    it uses.

    Algorithm~\ref{alg:ext:atomspace} uses a similar idea of ``failed $(\ell+1)$-Littlestone tree'' of Malliaris--Shelah~\cite{MS14}
    (which is the same idea as the one used in Algorithms~\ref{alg:ext:good} and~\ref{alg:ext:atom}), except that now we
    sacrifice a bit of time-complexity in favor of saving space. Just as in Algorithm~\ref{alg:ext:atom}, we would like to keep
    track of $m$-tuples $(x_\sigma)_{\sigma\in\{0,1\}^{<\ell+1}}$ and sets $(U_\tau)_{\tau\in\{0,1\}^{<\ell+2}}$, but since our
    space-complexity budget is $O((\ell+1)\cdot(m+1)\cdot\log(n+1))$, we don't have enough space even to keep all $x_\sigma$ in
    memory at the same time!

    The trick is to instead inspect one branch $\tau\in\{0,1\}^{\ell+1}$ of the failed Littlestone tree at a time\footnote{We
    don't actually need to construct the last vertex or sets of the tree, so the algorithm actually enumerates $\tau$ in
    $\{0,1\}^\ell$.}, which means that we will only need to keep at most $\ell+1$ many $x_\sigma$ in memory (namely,
    $(x_{\tau\rest_{[t]}})_{t=0}^\ell$), which costs $O((\ell+1)\cdot(m+1)\log(n+1))$ space as we can reuse the space; there is an
    extra space overhead of $O(\ell+1)$ that is used to keep track of which branch $\tau$ we are currently checking.

    As for the sets $(U_{\tau\rest_{[t]}})_{t=0}^{\ell+1}$ along this branch, they are completely determined by knowing $\tau$,
    $t$ and the sequence $(x_{\tau\rest_{[j]}})_{j=0}^{t-1}$, so we don't need to actually keep their explicit description in
    memory: the set $U_{\tau\rest_{[t]}}$ is simply the set of the first $s_t$ vertices $u$ such that
    $t_G^{\widetilde{\delta}}(u,x_{\tau\rest_{[j-1]}}) = \tau_j$ for every $j\in[t-1]$. We will also make use of constantly-many
    other variables to enumerate vertices, enumerate levels of the Littlestone tree, enumerate $m$-tuples of vertices or count
    vertices; this will incur an extra space overhead of $O(\log(\ell+1) + (m+1)\cdot\log(n+1))$.

    Finally, let us point out that the oracle algorithm, Algorithm~\ref{alg:ext:atomspaceoracle}, does not actually need several
    of the parameters of the computation algorithm, Algorithm~\ref{alg:ext:atomspace}; more specifically,
    Algorithm~\ref{alg:ext:atomspaceoracle} does not need the values of the $s_i$, nor the value of $\ell$. As such
    Algorithm~\ref{alg:ext:atomspaceoracle} does \emph{not} receive these parameters as input.

    \begin{algorithm}[htbp]
      \caption{Computation algorithm of query-into-oracle model with space-, time- and $U$-oracle-complexities
        \begin{align*}
          O\bigl(\ell\cdot m\cdot\log(n+1)\bigr),
          & &
          O(2^\ell\cdot m\cdot n^m),
          & &
          (2^{\ell+1} - 1)\cdot n^{m+1},
        \end{align*}
        respectively, that returns $(\sigma,(x_j)_{j=0}^{\lvert\sigma\rvert-1})$ such that when given to
        Algorithm~\ref{alg:ext:atomspaceoracle} produces an oracle for a $(\widetilde{\delta},\widetilde{\epsilon},m)$-atom of
        $G$ of size $s_{\lvert\sigma\rvert}$, where $\lvert\sigma\rvert\in\{0,\ldots,\ell\}$.
      }
      \label{alg:ext:atomspace}
      \DontPrintSemicolon
      \KwIn{Numbers $n,m,\ell,s_0,\ldots,s_\ell,s_{\ell+1}\in\NN_+$, $\widetilde{m}\in\NN$ with $\ell\leq\log_2(n)$, a
        query-oracle $\cO_G$ for a graph $G$ with $\lvert G\rvert=n$ and $\Lit(G)\leq\ell$, and a query-oracle $\cO_U$ for a set
        $U\subseteq V(G)$ with $\lvert U\rvert\geq s_0$. We further assume that there exists
        $\widetilde{\delta}\in(0,1/2^{\ell+1})$ such that $\widetilde{m} = \floor{\widetilde{\delta}\cdot m}$.}
      \KwOut{A tuple $(\sigma,(x_j)_{j=0}^{\lvert\sigma\rvert-1})$, where $\sigma\in\{0,1\}^{<\ell+1}$ and $x_j\in V(G)^m$, such
        that when given to Algorithm~\ref{alg:ext:atomspaceoracle} produces an oracle for a set $W\subseteq U$ that is a
        $(\widetilde{\delta},\widetilde{\epsilon},m)$-atom of $G$ with $\lvert W\rvert=s_{\lvert\sigma\rvert}$
        and
        \begin{gather*}
          \widetilde{\epsilon}
          \df
          \max_{i\in[\ell+1]}\frac{s_i-1}{s_{i-1}}.
        \end{gather*}
        }
      \tcc*[l]{Total memory required is for variables: $\tau\in\{0,1\}^\ell$, $t_0,t,j\in\{0,\ldots,\ell\}$, $w\in[n]$,
        $a,k\in\{0,\ldots,n\}$, $d\in\{0,\ldots,m\}$ and $z,x_0,\ldots,x_\ell\in[n]^m$}
      \For{$i\in[\ell+1]$}{%
        \If(\tcp*[f]{$\widetilde{\epsilon}\geq 1/2$}){$2\cdot (s_i-1)\geq s_{i-1}$}{%
          \label{alg:ext:atomspace:trivial}
          \Return{$(\varnothing,\varnothing)$}\tcp*[r]{Corresponds to first $s_0$ vertices of $U$}
        }
      }
      \uFor(\tcp*[f]{``Guess'' atom is in the $\tau$ branch}){$\tau\in\{0,1\}^\ell$ in lexicographic-order }{%
        \label{alg:ext:atomspace:tauFor}
        $t_0\assign\max\{i\in[\ell] \mid \tau_i = 1\}\cup\{0\}$.\;
        \label{alg:ext:atomspace:t0}
        \uFor(\tcp*[f]{``Guess'' atom is $U_{\tau\rest_{[t]}}$}){$t\in\{t_0,\ldots,\ell\}$}{%
          \label{alg:ext:atomspace:tFor}
          Undefine $x_t$.\;\label{alg:ext:atomspace:undefinext}
          \uFor(\tcp*[f]{$m$-tuple that potentially $(\widetilde{\delta},\widetilde{\epsilon})$-splits $U_{\tau_{[t]}}$}){%
            $z\in [n]^m$
          }{%
            \label{alg:ext:atomspace:zFor}
            $a\assign 0$\tcp*[r]{Vertices of $U_{\tau\rest_{[t]}}$ seen so far}
            $k^0\assign 0$\tcp*[r]{Vertices $w$ of $U_{\tau\rest_{[t]}}$ with $t_G^{\widetilde{\delta}}(w,z)=0$ seen so far}
            $k^1\assign 0$\tcp*[r]{Vertices $w$ of $U_{\tau\rest_{[t]}}$ with $t_G^{\widetilde{\delta}}(w,z)=1$ seen so far}
            \uFor(\tcp*[f]{Potential vertex of $U_{\tau\rest_{[t]}}$}){$w\in[n]$}{%
              \label{alg:ext:atomspace:wFor}
              \lIf{$\cO_U(w)=0$}{\label{alg:ext:atomspace:Utest}\Continue}
              \setcounter{algosplit}{\theAlgoLine}
            }
          }
        }
      }
    \end{algorithm}

    \begin{algorithm}[htbp]
      \caption*{(continued).}
      \DontPrintSemicolon
      \setcounter{AlgoLine}{\thealgosplit}
      \let\oldnl\nl
      \let\nl\relax
      \vspace{-\baselineskip}\InvisibleBegin{
        \vspace{-\baselineskip}\InvisibleBegin{
          \vspace{-\baselineskip}\InvisibleBegin{
            \vspace{-\baselineskip}\InvisibleBegin{
              \global\let\nl\oldnl 
              \uFor{$j\in\{0,\ldots,t-1\}$}{%
                \label{alg:ext:atomspace:jFor}
                $d\assign 0$\tcp*[r]{Indices $i\in[m]$ with $(x_j)_i\in N_G^{\tau_{j+1}}(w)$}  
                \For{$i\in[m]$}{%
                  \label{alg:ext:atomspace:iFor}
                  \lIf{$\cO_G(w,(x_j)_i)=\tau_{j+1}$}{$d\assign d+1$}
                }
                \If(\tcp*[f]{$w\notin U_{\tau\rest_{[t]}}$ since $t_G^{\widetilde{\delta}}(w,x_j)\neq\tau_{j+1}$}){%
                  $d < m - \widetilde{m}$
                }{%
                  \label{alg:ext:atomspace:BreakjFor}\Break
                }
              }
              \Else(\tcp*[f]{Executes only if ``for'' loop of line~\ref{alg:ext:atomspace:jFor} didn't break}){%
                \label{alg:ext:atomspace:jForElse}
                $a\assign a+1$\tcp*[r]{$w\in U_{\tau\rest_{[t]}}$, count it.}
                $d\assign 0$\tcp*[r]{Indices $i\in[m]$ with $z_i\in N_G(w)$ seen so far}
                \For{$i\in[m]$}{%
                  \label{alg:ext:atomspace:jForElseiFor}
                  \lIf{$\cO_G(w,z_i)=1$}{$d\assign d+1$}
                }
                \lIf(\tcp*[f]{$t_G^{\widetilde{\delta}}(w,z)=1$}){%
                  $d\geq m - \widetilde{m}$
                }{%
                  \label{alg:ext:atomspace:k1update}
                  $k^1\assign k^1+1$
                }
                \lElseIf(\tcp*[f]{$t_G^{\widetilde{\delta}}(w,z)=0$}){%
                  $d\leq\widetilde{m}$
                }{%
                  \label{alg:ext:atomspace:k0update}
                  $k^0\assign k^0+1$
                }
                \lIf(\tcp*[f]{All vertices of $U_{\tau\rest_{[t]}}$ already seen}){%
                  $a\geq s_t$}{\label{alg:ext:atomspace:mIf}\Break
                }
              }
            }
            \If(\tcp*[f]{$U_{\tau\rest_{[t]}}$ seems $(\widetilde{\delta},\widetilde{\epsilon})$-split by $z$}){%
              $k^0\geq s_{t+1}$ and $k^1\geq s_{t+1}$
            }{%
              \label{alg:ext:atomspace:kst+1If}
              $x_t\assign z$\;
              \label{alg:ext:atomspace:definext}
              \Break
            }
          }
          \If(\tcp*[f]{$U_{\tau\rest_{[t]}}$ is a $(\widetilde{\delta},\widetilde{\epsilon})$-atom}){$x_t$ is not defined}{%
            \label{alg:ext:atomspace:returnIf}
            \Return{$(\tau\rest_{[t]},(x_j)_{j=0}^{t-1})$}
          }
        }
      }
      \Return{\Failure.}\tcp*[r]{This return statement never actually executes}
      \label{alg:ext:atomspace:totalfailure}
    \end{algorithm}

    \begin{algorithm}[htbp]
      \caption{Oracle algorithm of query-into-oracle model with space-, time- and $U$-oracle-complexities
        \begin{align*}
          O\bigl(\log(\lvert\sigma\rvert+1) + \log(m+1) + \log(n+1)\bigr),
          & &
          O\bigl((\lvert\sigma\rvert+1)\cdot m\cdot n\bigr),
          & &
          n,
        \end{align*}
        respectively, that when given $(\sigma,(x_j)_{j=0}^{\lvert\sigma\rvert-1})$ computed via
        Algorithm~\ref{alg:ext:atomspace}, produces an oracle for a $(\widetilde{\delta},\widetilde{\epsilon},m)$-atom of $G$ of
        size $s_{\lvert\sigma\rvert}$. The algorithm does not actually use several of the parameters required by
        Algorithm~\ref{alg:ext:atomspace}, and as such, does not receive them as input.}
      \label{alg:ext:atomspaceoracle}
      \DontPrintSemicolon
      \KwIn{Numbers $n,m\in\NN_+$, $\widetilde{m}\in\NN$, a query-oracle $\cO_G$ for a graph $G$ with $V(G)=[n]$, a query-oracle
        $\cO_U$ for a set $U\subseteq V(G)$, a tuple $(\sigma,(x_j)_{j=0}^{\lvert\sigma\rvert-1})\in\{0,1\}^{\ell+1}\times
        (V(G)^m)^{\lvert\sigma\rvert-1}$ provided by Algorithm~\ref{alg:ext:atomspace} ran with the same parameters (and a
        choice of the other parameters required by it and satisfying the required hypotheses), and a vertex $u\in V(G)$. We
        further assume that there exists $\widetilde{\delta}\in(0,1/2^{\ell+1})$ such that $\widetilde{m} =
        \floor{\widetilde{\delta}\cdot m}$.}
      \KwOut{A bit $b\in\{0,1\}$ such that the set $W$ of $u\in V(G)$ for which the algorithm returns $1$ is a subset of $U$ of
        size $s_{\lvert\sigma\rvert}$ and that is a $(\widetilde{\delta},\widetilde{\epsilon},m)$-atom of $G$, where
        \begin{gather*}
          \widetilde{\epsilon}
          \df
          \max_{i\in[\ell+1]}\frac{s_i-1}{s_{i-1}},
        \end{gather*}
        where the $s_i$ are the ones used by Algorithm~\ref{alg:ext:atomspace}.}
      \tcc*[l]{Total memory required is for variables: $j\in\{0,\ldots,\lvert\sigma\rvert-1\}$, $w\in[n]$ and
        $d\in\{0,\ldots,m\}$}
      \lIf{$\cO_U(u)=0$}{\Return{$0$}\label{alg:ext:atomspaceoracle:uoracle}}
      \For{$j\in\{0,\ldots,\lvert\sigma\rvert-1\}$}{%
        \label{alg:ext:atomspaceoracle:ujFor}
        $d\assign 0$\tcp*[r]{Indices $i\in[m]$ with $(x_j)_i\in N_G^{\sigma_{j+1}}(u)$ seen so far}
        \For{$i\in[m]$}{%
          \label{alg:ext:atomspaceoracle:uiFor}
          \lIf{$\cO_G(u,(x_j)_i)=\sigma_{j+1}$}{$d\assign d+1$}
        }
        \lIf(\tcp*[f]{$u\notin W$ since $t_G^{\widetilde{\delta}}(u,x_j)\neq\sigma_{j+1}$}){%
          $d < m - \widetilde{m}$
        }{\label{alg:ext:atomspaceoracle:ureturn0}\Return{$0$}}
      }
      $a\assign 0$\tcp*[r]{Vertices of $W$ seen so far}
      \For(\tcp*[f]{Potential vertex of $W$}){$w\in[u-1]$}{%
        \label{alg:ext:atomspaceoracle:wFor}
        \lIf{$\cO_U(w)=0$}{\Continue}
        \uFor{$j\in\{0,\ldots,\lvert\sigma\rvert-1\}$}{%
          \label{alg:ext:atomspaceoracle:jFor}
          $d\assign 0$\tcp*[r]{Indices $i\in[m]$ with $(x_j)_i\in N_G^{\sigma_{j+1}}(w)$}  
          \For{$i\in[m]$}{%
            \label{alg:ext:atomspaceoracle:iFor}
            \lIf{$\cO_G(w,(x_j)_i)=\sigma_{j+1}$}{$d\assign d+1$}
          }
          \lIf(\tcp*[f]{$w\notin W$ since $t_G^{\widetilde{\delta}}(w,x_j)\neq\sigma_{j+1}$}){%
            $d < m - \widetilde{m}$
          }{\label{alg:ext:atomspaceoracle:BreakjFor}\Break}
        }
        \Else(\tcp*[f]{Executes only if ``for'' loop of line~\ref{alg:ext:atomspaceoracle:jFor} didn't break}){%
          \label{alg:ext:atomspaceoracle:jForElse}
          $a\assign a+1$\tcp*[r]{$w\in W$, count it.}
          \lIf(\tcp*[f]{All vertices of $W$ seen}){$a\geq s_{\lvert\sigma\rvert}$}{%
            \label{alg:ext:atomspaceoracle:return0}\Return{$0$}
          }
        }
      }
      \Return{$1$}\label{alg:ext:atomspaceoracle:return1}
    \end{algorithm}

    We start by arguing the space-complexity of Algorithm~\ref{alg:ext:atomspace}: the algorithm uses the variables
    $\tau\in\{0,1\}^\ell$, $t_0,t,j\in\{0,\ldots,\ell\}$, $w\in[n]$, $a,k\in\{0,\ldots,n\}$, $d\in\{0,\ldots,m\}$ and
    $z,x_0\ldots,x_\ell\in[n]^m$, so the number of bits needed to store these is\footnote{The algorithm also allows the $x_i$ to
    have an ``undefined state'' that is actually checked for, but this can be encoded by a special value encoding ``undefined
    state'' (which only incurs constant additive overhead). Also, it is easy to see that $t$ can actually reuse the space of
    $t_0$, as the latter is only used for initialization of the former, but this would only save $\ceil{\log_2(\ell+1)}$ space.}
    \begin{align*}
      \MoveEqLeft
      \ell + 3\cdot\ceil{\log_2(\ell+1)} + 3\cdot\ceil{\log_2(n+1)}
      + \ceil{\log_2(m+1)} + (\ell+2)\cdot\ceil{m\cdot\log_2(n+1)}
      \\
      & =
      O\bigl(\ell\cdot m\cdot\log(n+1)\bigr).
    \end{align*}

    The proof of the time- and $U$-oracle-complexities of Algorithm~\ref{alg:ext:atomspace} are more involved and require that we
    prove correctness first (as we will prove several structural properties of the algorithm along the way that will help in the
    analysis of the time- and $U$-oracle-complexities).

    \smallskip

    We now analyze the complexity of Algorithm~\ref{alg:ext:atomspaceoracle}. Space-complexity is easy: the algorithm uses the
    variables $j\in\{0,\ldots,\lvert\sigma\rvert-1\}$, $w\in[n]$ and $d\in\{0,\ldots,m\}$, so the number of bits needed to store
    these is
    \begin{gather*}
      \log(\lvert\sigma\rvert+1) + \log(m+1) + \log(n+1).
    \end{gather*}

    For time- and $U$-oracle-complexity, we note that the loops of lines~\ref{alg:ext:atomspaceoracle:ujFor},
    \ref{alg:ext:atomspaceoracle:uiFor}, \ref{alg:ext:atomspaceoracle:wFor}, \ref{alg:ext:atomspaceoracle:jFor}
    and~\ref{alg:ext:atomspaceoracle:iFor} execute at most $\lvert\sigma\rvert$, $m$, $n-1$, $\lvert\sigma\rvert$ and $m$ times,
    respectively, and there is exactly one $U$-oracle call in line~\ref{alg:ext:atomspaceoracle:uoracle} and one $U$-oracle call
    per iteration of the loop of line~\ref{alg:ext:atomspaceoracle:wFor}, so the $U$-oracle complexity is upper bounded by $n$.
    For the time-complexity, we note that each iteration of the loops of lines~\ref{alg:ext:atomspaceoracle:ujFor}
    and~\ref{alg:ext:atomspaceoracle:jFor} takes time $O(m)$ so the time-complexity is upper bounded by
    \begin{gather*}
      O\bigl((\lvert\sigma\rvert+1)\cdot m\cdot n\bigr).
    \end{gather*}

    \smallskip

    Let us now argue correctness of Algorithms~\ref{alg:ext:atomspace} and~\ref{alg:ext:atomspaceoracle}.

    We start by handling the trivial case: if the test of line~\ref{alg:ext:atomspace:trivial} passes, then we must have
    \begin{gather*}
      \widetilde{\epsilon}
      \df
      \max_{i\in[\ell+1]}\frac{s_i-1}{s_{i-1}}
      \geq
      \frac{1}{2},
    \end{gather*}
    which means that any set is a $(\widetilde{\delta},\widetilde{\epsilon},m)$-atom of $G$ (since
    $\widetilde{\delta}<1/2^{\ell+1}\leq 1/2$) and Algorithm~\ref{alg:ext:atomspace} returns $(\varnothing,\varnothing)$; once
    we feed $(\varnothing,\varnothing)$ to Algorithm~\ref{alg:ext:atomspaceoracle}, then the loops of
    lines~\ref{alg:ext:atomspaceoracle:ujFor} and~\ref{alg:ext:atomspace:jFor} have length $0$, so the ``Else'' block of
    line~\ref{alg:ext:atomspace:jForElse} always executes, making clear that the algorithm returns $1$ if and only if $u\in U$
    and it is among the first $s_0$ elements of $U$, this is a set of size exactly $s_0$ since $\lvert U\rvert\geq s_0$ (and as
    already argued is a $(\widetilde{\delta},\widetilde{\epsilon},m)$-atom of $G$).

    From now on, let us then suppose that $\widetilde{\epsilon} < 1/2$ and let us fix some notation for values of variables as
    the execution progresses and for sets that are being implicitly computed by the Algorithm~\ref{alg:ext:atomspace}.

    For each $\tau\in\{0,1\}^\ell$, let
    \begin{gather*}
      t_0[\tau]
      \df
      \max\{i\in[\ell] \mid \tau_i=1\}\cup\{0\}
    \end{gather*}
    be the value of $t_0$ computed in line~\ref{alg:ext:atomspace:t0} in the iteration $\tau$ of the loop of
    line~\ref{alg:ext:atomspace:tauFor} and note that $t_0[\tau]$ is exactly the unique position of $\tau$ that changed from $0$
    to $1$ from the previous iteration to the current one (and is $0$ for the first iteration, in which $\tau=(0,\ldots,0)$). We
    also note that in each iteration of the loop of line~\ref{alg:ext:atomspace:tauFor}, for each $k\in\{0,\ldots,\ell\}$, the
    value of $x_k$ either does not change (if $k < t_0[\tau]$ or if the test of line~\ref{alg:ext:atomspace:returnIf} passes for
    some $t < k$) or it gets undefined in line~\ref{alg:ext:atomspace:undefinext} and potentially redefined in
    line~\ref{alg:ext:atomspace:definext}, so we let $x_t[\tau]$ be the value of $x_t$ by the end of the iteration corresponding
    to $\tau$.

    Since $x_t$ only gets redefined when $t > t_0[\tau]$, it follows that for every $\tau,\tau'\in\{0,1\}^\ell$, if
    $k\in\{0,\ldots,\ell\}$ is such that $\tau\rest_{[k]} = \tau'\rest_{[k]}$, then $x_k[\tau] = x_k[\tau']$, this means that we
    can make the following definition: for every $\sigma\in\{0,1\}^{<\ell+1}$, let $x_\sigma\df x_{\lvert\sigma\rvert}[\tau]$
    for any $\tau\in\{0,1\}^\ell$ such that $x_{\lvert\sigma\rvert}[\tau]$ is defined and
    $\tau\rest_{[\lvert\sigma\rvert]}=\sigma$; if no such $\tau$ exists, then let us say that $x_\sigma$ is undefined.

    For each $\tau\in\{0,1\}^\ell$, each $k\in\{0,\ldots,\ell\}$, each $z\in V(G)^m$ and each $b\in\{0,1\}$, we let
    \begin{align*}
      U'[\tau,k]
      & \df
      \Bigl\{w\in U \;\Bigm\vert\;
      \forall j\in[k], t_G^{\widetilde{\delta}}\bigl(w,x_{j-1}[\tau]\bigr) = \tau_j
      \Bigr\}
      \\
      & =
      \bigl\{w\in U \;\bigm\vert\;
      \forall j\in[k], t_G^{\widetilde{\delta}}(w,x_{\tau\rest_{[j-1]}}) = \tau_j
      \bigr\},
      \\
      U'[\tau,k,z,b]
      & \df
      \bigl\{w\in U'[\tau,k] \;\bigm\vert\;
      t_G^{\widetilde{\delta}}(w,z) = b
      \bigr\}
      \\
      & =
      \bigl\{w\in U \;\bigm\vert\;
      \forall j\in[k], t_G^{\widetilde{\delta}}(w,x_{\tau\rest_{[j-1]}}) = \tau_j
      \land t_G^{\widetilde{\delta}}(w,z) = b
      \bigr\},
    \end{align*}
    provided all quantities above are defined.

    Again, it is clear from the definition above that if $\tau,\tau'\in\{0,1\}^\ell$ are such that
    $\tau\rest_{[k-1]}=\tau'\rest_{[k]}$, then $U'[\tau,k]=U'[\tau',k]$, so we can make the following definition: for every
    $\sigma\in\{0,1\}^{<\ell+1}$, let $U'_\sigma\df U'[\tau,\lvert\sigma\rvert]$ for any $\tau\in\{0,1\}^\ell$ such that
    $U'[\tau,\lvert\sigma\rvert]$ is defined and $\tau\rest_{[\lvert\sigma\rvert]} = \sigma$. For $\sigma\in\{0,1\}^{\ell+1}$,
    we also define $U'_\sigma\df U'[\sigma\rest_{[\ell]},\ell,\sigma_{\ell+1}]$.

    Note that this definition also guarantees that for every $\tau\in\{0,1\}^{\ell+1}$ and every $k\in\{0,\ldots,\ell\}$, we
    have $U'_{\tau\rest_{[k+1]}} = U'[\tau\rest_{[\ell]},k,x_{\tau\rest_{[k]}},\tau_{k+1}]$.

    Now, for each $\sigma\in\{0,1\}^{<\ell+2}$, we let $U_\sigma$ be the set of the first at most $s_{\lvert\sigma\rvert}$
    vertices of $U'_\sigma$ (in the natural order of $[n]=V(G)$).

    We claim that if $U'_\sigma$ is defined and $\lvert\sigma\rvert\leq\ell$, then $U_\sigma$ is exactly the set of vertices $w$
    for which the block of line~\ref{alg:ext:atomspace:jForElse} executes (i.e., the ones in which the loop of
    line~\ref{alg:ext:atomspace:jFor} does not break because the test of line~\ref{alg:ext:atomspace:BreakjFor} always fails) in
    any iteration of the loop of line~\ref{alg:ext:atomspace:tauFor} corresponding to $\tau\in\{0,1\}^\ell$ such that
    $\tau\rest_{[\lvert\sigma\rvert]} = \sigma$ and in the iteration of the loop of line~\ref{alg:ext:atomspace:tFor}
    corresponding to $t=\lvert\sigma\rvert$ (provided $\lvert\sigma\rvert\geq t_0[\tau]$).

    Indeed, in such an iteration we have $x_{j-1}[\tau] = x_{\sigma\rest_{j-1}}$ for every $j\in[\lvert\sigma\rvert]$. On the
    other hand, we note that in the iteration $j$ of the loop of line~\ref{alg:ext:atomspace:jFor}, after the loop of
    line~\ref{alg:ext:atomspace:iFor} is concluded, $d$ holds the value of
    \begin{gather*}
      \Bigl\lvert\Bigl\{i\in[m] \;\Bigm\vert\; \bigl(x_j[\tau]\bigr)_i\in N_G^{\tau_j}(w)\Bigr\}\Bigr\rvert
      =
      \bigl\lvert\bigl\{i\in[m] \;\bigm\vert\; (x_{\sigma\rest_{j-1}})_i\in N_G^{\sigma_j}(w)\bigr\}\bigr\rvert,
    \end{gather*}
    and the test of line~\ref{alg:ext:atomspace:BreakjFor} tests if the value above is less than $m - \widetilde{m} = m -
    \floor{\widetilde{\delta}\cdot m}$. Since these quantities are integers, this is equivalent to testing whether the set above
    has size strictly less than $(1-\widetilde{\delta})\cdot m$, which in turn is equivalent to
    \begin{gather*}
      t_G^{\widetilde{\delta}}(w,x_{\sigma\rest_{j-1}}) \neq \sigma_j,
    \end{gather*}
    so the set of vertices that always fail this test (hence do not break the loop of line~\ref{alg:ext:atomspace:jFor}) is
    exactly $U'_\sigma$. Finally, the test of line~\ref{alg:ext:atomspace:mIf} ensures that no vertices $w$ are checked after
    $s_t$ many vertices have executed the block of line~\ref{alg:ext:atomspace:jForElse}, making the final set of vertices that
    execute the block of line~\ref{alg:ext:atomspace:jForElse} exactly $U_\sigma$.

    We now claim that if $U'_\sigma$ is defined, then the size of $U_\sigma$ is exactly $s_{\lvert\sigma\rvert}$. By its
    definition, we know that the size of $U_\sigma$ is at most $s_{\lvert\sigma\rvert}$, so we only need to prove the other
    inequality, which amounts to showing that $U'_\sigma$ has size at least $s_{\lvert\sigma\rvert}$. The proof is by induction
    on $\lvert\sigma\rvert$. For $\lvert\sigma\rvert=0$, i.e., $\sigma=\varnothing$, this follows simply because $U'_\varnothing
    = U$, which has size at least $s_0$ by assumption. For $\lvert\sigma\rvert > 0$, since $U'_\sigma$ is defined, it must be
    the case that $U'_{\sigma\rest_{[\lvert\sigma\rvert-1]}}$ is also defined, which by inductive hypothesis implies that
    $U_{\sigma\rest_{[\lvert\sigma\rvert-1]}}$ has size $s_{\lvert\sigma\rvert-1}$.

    By our previous claim, we know that $U_{\sigma\rest_{[\lvert\sigma\rvert-1]}}$ is exactly the set of $w$ that execute the
    block of line~\ref{alg:ext:atomspace:jForElse} in the iteration of line~\ref{alg:ext:atomspace:tauFor} corresponding to
    \begin{gather*}
      \tau
      =
      \sigma\rest_{[\lvert\sigma\rvert-1]}
      \conc
      (\mathop{\underbrace{0,\ldots,0}}\limits_{\mathclap{\text{$\ell-\lvert\sigma\rvert+1$ times}}})
    \end{gather*}
    and in the iteration of the loop of line~\ref{alg:ext:atomspace:tFor} when $t=\lvert\sigma\rvert-1$ (note that this
    iteration exists as for the $\tau$ above we have $t_0[\tau]\leq\lvert\sigma\rvert-1$). Now, for one such $w\in
    U_{\sigma\rest_{[\lvert\sigma\rvert-1]}}$, it is clear that by the end of the loop of
    line~\ref{alg:ext:atomspace:jForElseiFor}, the variable $d$ holds the value
    \begin{gather*}
      \lvert\{i\in[m] \mid z_i\in N_G^{\tau_j}(w)\}\rvert.
    \end{gather*}

    Now, since $\widetilde{m} = \floor{\widetilde{\delta}\cdot m}$, we have
    \begin{align*}
      d\geq m - \widetilde{m}
      & \iff
      d\geq m - \floor{\widetilde{\delta}\cdot m}
      \iff
      d \geq (1-\delta)\cdot m
      \iff
      t_G^{\widetilde{\delta}}(w,z) = 1,
      \\
      d\leq\widetilde{m}
      & \iff
      d\leq\floor{\widetilde{\delta}\cdot m}
      \iff
      d\leq\widetilde{\delta}\cdot m
      \iff
      t_G^{\widetilde{\delta}}(w,z) = 0.
    \end{align*}
    Note also that since $\widetilde{\delta}<1/2^{\ell+1}\leq 1/2$, the conditions above are mutually exclusive. Thus, by the
    point line~\ref{alg:ext:atomspace:kst+1If} is executed, the variables $k^0$ and $k^1$ hold the values
    \begin{align*}
      \lvert\{w\in U_{\sigma\rest_{[\lvert\sigma\rvert-1]}} \mid t_G^{\widetilde{\delta}}(w,z) = 0\}\rvert,
      & &
      \lvert\{w\in U_{\sigma\rest_{[\lvert\sigma\rvert-1]}} \mid t_G^{\widetilde{\delta}}(w,z) = 1\}\rvert,
    \end{align*}
    respectively. Since $U'_\sigma$ being defined implies that $x_{\sigma\rest_{[\lvert\sigma\rvert-1]}}$ is
    defined, this must happen due to the line~\ref{alg:ext:atomspace:definext} in the iteration of the loop of
    line~\ref{alg:ext:atomspace:zFor} when $z = x_{\sigma\rest_{[\lvert\sigma\rvert-1]}}$, which in particular implies that for
    the value $k^{\sigma_{\lvert\sigma\rvert}}$ at that point, we have
    \begin{gather*}
      \lvert U'_\sigma\rvert
      =
      \lvert\{w\in U_{\sigma\rest_{[\lvert\sigma\rvert-1]}} \mid
      t_G^{\widetilde{\delta}}(w,z) = \sigma_{\lvert\sigma\rvert}\}\rvert
      =
      k^{\sigma_{\lvert\sigma\rvert}}
      \geq
      s_{t+1}
      =
      s_{\lvert\sigma\rvert},
    \end{gather*}
    so $\lvert U'_\sigma\rvert$ has size at least $s_{\lvert\sigma\rvert}$ as desired.

    \smallskip

    Let us now argue that Algorithm~\ref{alg:ext:atomspace} never executes line~\ref{alg:ext:atomspace:totalfailure}. Suppose
    not, then it must be the case that all $x_t[\tau]$ (hence all $x_\sigma$ and all $U_\sigma$) get defined. Since for every
    $\tau\in\{0,1\}^{\ell+1}$, we have $\lvert U_\tau\rvert=s_{\ell+1}\geq 1$, we can let $u_\tau\in U_\tau$.

    Note that the definition of the sets $U_\sigma$ imply for every $\tau\in\{0,1\}^\ell$ and every $t\in[\ell+1]$, we have
    $U_\tau\subseteq U_{\tau\rest_{[t]}}$, which in particular implies that
    \begin{gather}\label{eq:ext:atomspace:totalfailure}
      t_G^{\widetilde{\delta}}(u_\tau, x_{\tau\rest_{[t-1]}}) = \tau_t.
    \end{gather}

    We claim that for each $\sigma\in\{0,1\}^{<\ell+1}$, there exists an index $i_\sigma\in[m]$ such that for every
    $\tau\in\{0,1\}^{\ell+1}$ that extends $\sigma$, we have
    \begin{gather*}
      u_\tau\in N_G^{\tau_{\lvert\sigma\rvert+1}}\bigl((x_\sigma)_{i_\sigma}\bigr).
    \end{gather*}
    Indeed, we can simply consider the set of all $i\in[m]$ that satisfy the above conditions for each such $\tau$. Note that
    there are exactly $2^{\ell+1-\lvert\sigma\rvert}\leq 2^{\ell+1}$ many $\tau\in\{0,1\}^{\ell+1}$ extending $\sigma$, so by a
    union bound applied to~\eqref{eq:ext:atomspace:totalfailure}, the set of all $i\in[m]$ that satisfy the desired
    conditions has size at least
    \begin{gather*}
      \bigl(1 - 2^{\ell+1}\cdot\widetilde{\delta}\bigr)\cdot m,
    \end{gather*}
    which is positive since $\widetilde{\delta} < 1/2^{\ell+1}$ by assumption, guaranteing the existence of our desired
    $i_\sigma\in[m]$. But then $(((x_\sigma)_{i_\sigma})_{\sigma\in\{0,1\}^{<\ell+1}},(N_G(u_\tau))_{\tau\in\{0,1\}^{\ell+1}})$ is
    an $(\ell+1)$-Littlestone tree in $\cH_G$, contradicting the fact that $\Lit(G) = \Lit(\cH_G)\leq\ell$. Thus,
    Algorithm~\ref{alg:ext:atomspace} never executes line~\ref{alg:ext:atomspace:totalfailure}.

    \smallskip

    In particular, this means that the test of line~\ref{alg:ext:atomspace:returnIf} must pass in some iteration corresponding
    to $(\tau,t)$; we claim that for this $(\tau,t)$, if we set $\sigma\df\tau\rest_{[t]}$, then the set $U_\sigma$ is a
    $(\widetilde{\delta},\widetilde{\epsilon},m)$-atom of $G$. Suppose not, that is, suppose that there exists some $z\in
    V(G)^m$ that $(\widetilde{\delta},\widetilde{\epsilon})$-splits $U_\sigma$ in $G$, that is, we must have
    \begin{gather*}
      \frac{\lvert\{w\in U_\sigma \mid t_G^{\widetilde{\delta}}(w,z) = b\}\rvert}{\lvert U_\sigma\rvert} > \widetilde{\epsilon}
    \end{gather*}
    for every $b\in\{0,1\}$. But this means that in the final iteration corresponding to $(\tau,t)$, when the loop of
    line~\ref{alg:ext:atomspace:zFor} considered the tuple $z$, and the variables $k^0$ and $k^1$ ended up with the values
    \begin{align*}
      \lvert\{w\in U_\sigma \mid t_G^{\widetilde{\delta}}(w,z) = 0\}\rvert
      & >
      \widetilde{\epsilon}\cdot\lvert U_\sigma\rvert
      =
      \widetilde{\epsilon}\cdot s_{\lvert\sigma\rvert},
      \\      
      \lvert\{w\in U_\sigma \mid t_G^{\widetilde{\delta}}(w,z) = 1\}\rvert
      & >
      \widetilde{\epsilon}\cdot\lvert U_\sigma\rvert
      =
      \widetilde{\epsilon}\cdot s_{\lvert\sigma\rvert},
    \end{align*}
    respectively. Since $s_{\lvert\sigma\rvert+1}\leq\floor{\widetilde{\epsilon}\cdot s_{\lvert\sigma\rvert}}+1$, the sets above
    must have size at least $s_{\lvert\sigma\rvert+1}$, which would mean that the test of line~\ref{alg:ext:atomspace:kst+1If}
    would have passed, contradicting the fact that $x_t[\tau]$ is not defined. Thus $U_\sigma$ is a
    $(\widetilde{\delta},\widetilde{\epsilon},m)$-atom of $G$.

    \smallskip

    We now cover correctness of Algorithm~\ref{alg:ext:atomspaceoracle}. We already argued that
    Algorithm~\ref{alg:ext:atomspace} must pass the test of line~\ref{alg:ext:atomspace:returnIf} in some iteration
    corresponding to $(\tau,t)$, which means that it returns $(\sigma,(x_j)_{j=0}^{\lvert\sigma\rvert-1})$ given by
    $\sigma\df\tau\rest_{[t]}$ and $x_j\df x_{\tau\rest_{[j]}} = x_{\sigma\rest_{[j]}}$. We have also already argued that
    $U_\sigma$ is a $(\widetilde{\delta},\widetilde{\epsilon},m)$-atom of $G$, so now it suffices to prove that
    Algorithm~\ref{alg:ext:atomspaceoracle} is an oracle for the set $U_\sigma$.

    Indeed, we observe that the loops of lines~\ref{alg:ext:atomspaceoracle:ujFor} and~\ref{alg:ext:atomspaceoracle:jFor} of
    Algorithm~\ref{alg:ext:atomspaceoracle} have an analogous pattern to the loop of line~\ref{alg:ext:atomspace:jFor} of
    Algorithm~\ref{alg:ext:atomspace}, which means that the test of line~\ref{alg:ext:atomspaceoracle:ureturn0} passes if and
    only if $u\in U\setminus U'_\sigma$ and provided this test fails, the vertices $w$ that execute the ``Else'' block of
    line~\ref{alg:ext:atomspaceoracle:jForElse} are precisely the ones in $[u-1]\cap U_\sigma$.

    If $u\in U_\sigma$, then there are less than $s_{\lvert\sigma\rvert}$ vertices in $[u-1]\cap U_\sigma$ (as $\lvert
    U_\sigma\rvert=s_{\lvert\sigma\rvert}$), so the algorithm reaches line~\ref{alg:ext:atomspaceoracle:return0} and returns $1$.

    On the other hand, if $u\notin U_\sigma$, then exactly one of the following holds:
    \begin{itemize}
    \item $u\notin U$, in which case the algorithm returns $0$ in line~\ref{alg:ext:atomspaceoracle:uoracle}.
    \item $u\in U\setminus U'_\sigma$, in which case the algorithm returns $0$ in line~\ref{alg:ext:atomspaceoracle:ureturn0}.
    \item $u\in U'_\sigma\setminus U_\sigma$, in which case $[u-1]\cap U_\sigma = U_\sigma$, so the algorithm returns $0$ in
      line~\ref{alg:ext:atomspaceoracle:return0}.
    \end{itemize}

    Thus, Algorithm~\ref{alg:ext:atomspaceoracle} is indeed an oracle for the set $U_\sigma$, which we have already argued is a
    $(\widetilde{\delta},\widetilde{\epsilon},m)$-atom of $G$ and has size $s_{\lvert\sigma\rvert}\in S$, concluding the proof
    of correctness of Algorithms~\ref{alg:ext:atomspace} and~\ref{alg:ext:atomspaceoracle} in the query-into-oracle model.

    \medskip

    Finally, we analyze the time- and $U$-oracle-complexities of Algorithm~\ref{alg:ext:atomspace}. If $\widetilde{\epsilon}\geq
    1/2$, then the test of line~\ref{alg:ext:atomspace:trivial} passes and the algorithm takes time $O(\ell\cdot\log(n+1))$, so
    let us then suppose $\widetilde{\epsilon} < 1/2$.

    Now, we note that the loops of lines~\ref{alg:ext:atomspace:tauFor}, \ref{alg:ext:atomspace:tFor},
    \ref{alg:ext:atomspace:zFor}, \ref{alg:ext:atomspace:wFor}, \ref{alg:ext:atomspace:jFor}, \ref{alg:ext:atomspace:iFor}
    and~\ref{alg:ext:atomspace:jForElseiFor} execute at most $2^\ell$, $\ell-t_0[\tau]+1$, $n^m$, $n$, $t$, $m$ and $m$ times,
    respectively.

    These bounds suffice to analyze the $U$-oracle-complexity of Algorithm~\ref{alg:ext:atomspace}: the $U$-oracle is consulted
    only in line~\ref{alg:ext:atomspace:Utest}, so an upper bound on the $U$-oracle complexity is
    \begin{gather*}
      \sum_{\tau\in\{0,1\}^\ell} \sum_{t=t_0[\tau]}^\ell n^{m+1}
      =
      (2^{\ell+1} - \ell - 1)\cdot n^{m+1},
    \end{gather*}
    where the equality follows from
    \begin{gather}\label{eq:ext:tamortization}
      \sum_{\tau\in\{0,1\}^\ell} \bigl(\ell-t_0[\tau]+1\bigr)
      =
      \ell + 1 + \sum_{i=1}^\ell 2^{i-1}\cdot(\ell-i+1)
      =
      2^{\ell+1} - 1,
    \end{gather}
    since there for $i\in[\ell]$, there are exactly $2^{i-1}$ many $\tau\in\{0,1\}^\ell$ with $t_0[\tau]=i$, and there is
    exactly one $\tau\in\{0,1\}^\ell$ (namely $\tau=(0,\ldots,0)$) with $t_0[\tau]=0$.

    For time-complexity, if we make a naive analysis using the bounds above, we will obtain an upper bound that is larger than
    our goal by about a multiplicative factor of $\ell+1$. To get rid of this the extra factor of $\ell+1$, we make an analysis
    that amortizes both the $t$ of line~\ref{alg:ext:atomspace:tFor} and the $j$ of line~\ref{alg:ext:atomspace:jFor}. The
    amortization of $t$ will amount to a calculation similar to that of~\eqref{eq:ext:tamortization} and the amortization of
    $j$ will rely on the fact that the $j$ loop of line~\ref{alg:ext:atomspace:jFor} will run a total of $k$ times for
    approximately $\widetilde{\epsilon}^j\cdot n$ many vertices (except that the devil is in the details: the plus $1$ in the
    condition $s_i\leq\floor{\widetilde{\epsilon}\cdot s_{i-1}}+1$ makes the argument slightly more involved).

    For this amortized analysis, we consider again the tuples and sets $U'_\sigma$ and $U_\sigma$ and we let $(\tau_*,t_*)$ be
    the elements of the iterations of lines~\ref{alg:ext:atomspace:tauFor} and~\ref{alg:ext:atomspace:tFor} for which test of
    line~\ref{alg:ext:atomspace:returnIf} passes, causing the algorithm to terminate. This means that if $<_L$ is the strict
    lexicographic order on $\{0,1\}^\ell$, then $U'_\sigma$ and $U_\sigma$ is defined exactly for the $\sigma$ in the set
    \begin{multline*}
      \Sigma
      \df
      \Bigl\{\sigma\in\{0,1\}^{<\ell+2} \;\Bigm\vert\;
      \bigl(\lvert\sigma\rvert\leq t_*\land\tau_*\rest_{[\lvert\sigma\rvert]}=\sigma\bigr)
      \\
      \lor\exists\tau\in\{0,1\}^\ell,
      \bigl(\tau <_L \tau_*\land(\tau\rest_{[\lvert\sigma\rvert]} = \sigma\lor\sigma\rest_{[\ell]}=\tau)\bigr)
      \Bigr\},
    \end{multline*}
    that is, either $\sigma$ is a restriction or an extension of some $\tau\in\{0,1\}^\ell$ that is strictly before
    lexicographically than $\tau_*$, or $\sigma$ is a restriction of $\tau_*$ to a length of at most $t_*$.

    We now note that for every $w\in V(G)$ and every $t\in\{0,\ldots,\ell+1\}$, there exists at most one
    $\sigma^{t,w}\in\{0,1\}^t$ such that $w\in U_{\sigma^{t,w}}$ and we can let $\sigma^w$ be the $\sigma^{w,t}$ with maximum
    $t$ (if at least one such $\sigma^{w,t}$ exists, which happens precisely for $w\in U_\varnothing$). Given further
    $\tau\in\{0,1\}^\ell$ with $\tau\leq_L\tau_*$, we can then define
    \begin{gather*}
      j_w^{\tau,t} \df
      \begin{dcases*}
        0, & if $\sigma^w$ does not exist,
        \\
        t+1, & if $\sigma^w\rest_{[t]} = \tau\rest_{[t]}$,
        \\
        i, & if $\sigma^w$ exists, satisfies $\sigma^w\rest_{[t]}\neq\tau\rest_{[t]}$ and $i$ is the first position in which
        they differ.
      \end{dcases*}
    \end{gather*}

    Now we can make a better analysis of the loop of line~\ref{alg:ext:atomspace:jFor} together with its ``else'' block of
    line~\ref{alg:ext:atomspace:jForElse}: in an iteration of the outer loops corresponding to $(\tau,t,z,w)$, if we count the
    execution of the ``else'' block as one extra iteration, then the total number of iterations is exactly $j_w^{\tau,t}$, with
    each iteration costing $O(m)$ due to the loops of lines~\ref{alg:ext:atomspace:iFor}
    and~\ref{alg:ext:atomspace:jForElseiFor}.

    Note that even when $j_w^{\tau,t}=0$ and the loop of line~\ref{alg:ext:atomspace:jFor} does not execute at all due to the test
    of line~\ref{alg:ext:atomspace:Utest} failing, the iteration corresponding to $(\tau,t,z,w)$ still costs $O(1)$, so in total,
    the iteration corresponding to $(\tau,t,z,w)$ costs at most $O((j_w^\tau+1)\cdot m)$.

    Finally, due to the computation of line~\ref{alg:ext:atomspace:kst+1If}, each iteration corresponding to $(\tau,t,z)$ has an
    extra additive cost of $O(\log(n+1))$. This is completely dominated by the complexity of the loop of
    line~\ref{alg:ext:atomspace:wFor}, so we can ignore it.

    Thus, if we define for every $\tau\leq_L\tau_*$, the number
    \begin{gather*}
      t_1[\tau]
      \df
      \begin{dcases*}
        \ell, & if $\tau <_L \tau_*$,\\
        t_*, & if $\tau = \tau_*$,
      \end{dcases*}
    \end{gather*}
    then the loop of line~\ref{alg:ext:atomspace:tFor} iterates $t$ from $t_0[\tau]$ to $t_1[\tau]$, so we can upper bound the
    time-complexity of Algorithm~\ref{alg:ext:atomspace} by
    \begin{gather}\label{eq:ext:atomspace:timecompl}
      \begin{aligned}
        \MoveEqLeft
        \sum_{\substack{\tau\in\{0,1\}^\ell\\\tau\leq_L\tau_*}} \sum_{t=t_0[\tau]}^{t_1[\tau]}
        n^m\cdot
        \sum_{w\in[n]} O\bigl((j_w^{\tau,t} + 1)\cdot m\bigr).
        \\
        & =
        O(n^m\cdot m)\cdot
        \sum_{\substack{\tau\in\{0,1\}^\ell\\\tau\leq_L\tau_*}} \sum_{t=t_0[\tau]}^{t_1[\tau]}
        \left(
        \bigl(n - \lvert U_\varnothing\rvert\bigr)
        + \left(\sum_{j=1}^t j\cdot\lvert U_{\tau\rest_{[j-1]}}\setminus U_{\tau\rest_{[j]}}\rvert\right)
        + (t+1)\cdot\lvert U_{\tau\rest_{[t]}}\rvert
        \right)
        \\
        & =
        O(n^m\cdot m)\cdot
        \sum_{\substack{\tau\in\{0,1\}^\ell\\\tau\leq_L\tau_*}} \sum_{t=t_0[\tau]}^{t_1[\tau]}
        \left(
        \bigl(n - \lvert U_\varnothing\rvert\bigr)
        +
        \sum_{j=1}^{t+1} \lvert U_{\tau\rest_{[j-1]}}\rvert
        \right)
        \\
        & =
        O(n^m\cdot m)\cdot\sum_{\substack{\tau\in\{0,1\}^\ell\\\tau\leq_L\tau_*}} \sum_{t=t_0[\tau]}^{t_1[\tau]}
        \left(
        (n - s_0)
        +
        \sum_{j=1}^{t+1} s_{j-1}
        \right).
      \end{aligned}
    \end{gather}

    We now let
    \begin{gather*}
      i_* \df \min\{i\in\{0,\ldots,\ell\} \mid s_i = 1\}\cup\{\ell+1\}
    \end{gather*}
    and consider two cases.

    In the first case, we assume $i_* = \ell+1$. In this case, $s_\ell\geq 2$ and since $\widetilde{\epsilon} < 1/2$ and $2\leq
    s_\ell\leq\floor{\widetilde{\epsilon}\cdot s_{\ell-1}}+1$ (and $\ell > 0$), we conclude that $s_{\ell-1}\geq 3$. With a simple
    induction, we get that for every $i\in[\ell-1]$, we have
    \begin{gather*}
      3\leq s_i\leq\floor{\widetilde{\epsilon}\cdot s_{i-1}}+1\leq\widetilde{\epsilon}\cdot s_{i-1} + 1,
    \end{gather*}
    which in turn implies
    \begin{gather*}
      s_i \leq \frac{3}{2}\cdot\widetilde{\epsilon}\cdot s_{i-1}
    \end{gather*}
    for every $i\in[\ell-1]$, and now a simple induction yields
    \begin{gather*}
      s_i\leq \left(\frac{3}{2}\cdot\widetilde{\epsilon}\right)^i\cdot s_0
    \end{gather*}
    for every $i\in[\ell-1]$ (but this bound need not hold $s_\ell$ or $s_{\ell+1}$). Thus, since $\widetilde{\epsilon} < 1/2$
    (hence $(3/2)\cdot\widetilde{\epsilon} < 3/4$), we conclude that
    \begin{gather}\label{eq:ext:atomspace:ssum}
      \sum_{j=0}^\ell s_j
      \leq
      s_\ell + \sum_{j\in\NN} \left(\frac{3}{2}\cdot\widetilde{\epsilon}\right)^j\cdot s_0
      \leq
      s_\ell + \sum_{j\in\NN} \left(\frac{3}{4}\right)^j\cdot s_0
      \leq
      5\cdot n.
    \end{gather}

    Now we can upper bound the last expression in~\eqref{eq:ext:atomspace:timecompl} as follows:
    \begin{align*}
      O(n^m\cdot m)\cdot\sum_{\substack{\tau\in\{0,1\}^\ell\\\tau\leq_L\tau_*}} \sum_{t=t_0[\tau]}^{t_1[\tau]}
      \left(
      (n - s_0)
      +
      \sum_{j=1}^{t+1} s_{j-1}
      \right)
      & \leq
      O(n^m\cdot m)\cdot\sum_{\substack{\tau\in\{0,1\}^\ell\\\tau\leq_L\tau_*}} \sum_{t=t_0[\tau]}^{t_1[\tau]}
      \bigl(
      (n-s_0)
      +
      5\cdot n
      \bigr)
      \\
      & \leq
      O(n^m\cdot m)\cdot\sum_{\tau\in\{0,1\}^\ell} \sum_{t=t_0[\tau]}^\ell 6\cdot n
      \\
      & \leq
      O(n^m\cdot m)\cdot\sum_{\tau\in\{0,1\}^\ell} \bigl(\ell-t_0[\tau]+1\bigr) 6\cdot n
      \\
      & \leq
      O(2^\ell\cdot m\cdot n^{m+1}),
    \end{align*}
    where the first inequality follows from~\eqref{eq:ext:atomspace:ssum} and the last inequality follows
    from~\eqref{eq:ext:tamortization}.

    \smallskip

    We now consider the second case, in which $i_*\leq\ell$. In this case, since every set of size $1$ is a
    $(\widetilde{\delta},\widetilde{\epsilon},m)$-atom of $G$ (as $\widetilde{\delta} < 1/2^{\ell+1}\leq 1/2$), it follows that
    we must have $\tau_* = (0,\ldots,0)$ and $t_*\leq i_*$, so we can upper bound the last expression
    in~\eqref{eq:ext:atomspace:timecompl} as follows:
    \begin{align*}
      O(n^m\cdot m)\cdot\sum_{\substack{\tau\in\{0,1\}^\ell\\\tau\leq_L\tau_*}} \sum_{t=t_0[\tau]}^{t_1[\tau]}
      \left(
      (n - s_0)
      +
      \sum_{j=1}^{t+1} s_{j-1}
      \right)
      & =
      O(n^m\cdot m)\cdot\sum_{t=0}^{t_*}
      \left(
      (n - s_0)
      +
      \sum_{j=1}^{t+1} s_{j-1}
      \right)
      \\
      & \leq
      O(n^m\cdot m)\cdot\sum_{t=0}^\ell
      (t+2)\cdot\ell
      \leq
      O(\ell^2\cdot m\cdot n^m)
      \\
      & \leq
      O(2^\ell\cdot m\cdot n^m),
    \end{align*}
    concluding the time-complexity analysis of Algorithm~\ref{alg:ext:atomspace}.\qedhere
  \end{description}
\end{proof}

We conclude with a corollary about extracting large cliques/anticliques from graphs with bounded Littlestone dimension.

\begin{corollary}\label{cor:Ramsey}
  Let $n,\ell\in\NN_+$, let $\epsilon\in(0,1)$ and let $G$ be a graph with $\lvert G\rvert=n$ and $\Lit(G)\leq\ell$ and let
  $U\subseteq V(G)$ be non-empty. Then the following hold:
  \begin{enumerate}
  \item\label{cor:Ramsey:ext} There exists a subset $W\subseteq U$ such that either $\lvert
    W\rvert\geq\epsilon^{\ell-1}\cdot\lvert U\rvert$ and $W$ is $\epsilon$-good in $G$, or $W$ is $0$-good in $G$ and $\lvert
    W\rvert > \epsilon^\ell\cdot\lvert U\rvert$.
  \item\label{cor:Ramsey:0good} If $U$ is $0$-good, then $U$ is an anticlique of $G$.
  \item\label{cor:Ramsey:hom} If $(U,U)$ is $\epsilon$-homogeneous in $G$, then there exists $W\subseteq U$ with
    \begin{gather*}
      \lvert W\rvert \geq \Ceil{\left(\epsilon - \frac{1}{\lvert U\rvert}\right)^{-1}}
    \end{gather*}
    such that $W$ is either a clique or an anticlique of $G$.
  \item\label{cor:Ramsey:Ramsey} There exists $W\subseteq V(G)$ with
    \begin{gather*}
      \lvert W\rvert
      \geq
      \frac{n^{1/(\ell+1)}}{2^{1 - 1/(\ell+1)}} - O_{n\to\infty}(1)
    \end{gather*}
    such that $W$ is a clique or an anticlique of $G$.
  \end{enumerate}
\end{corollary}

\begin{proof}
  Item~\ref{cor:Ramsey:ext} is a slightly improved version of the argument of Proposition~\ref{prop:ext}: it essentially says
  that if the contruction of Proposition~\ref{prop:ext} has to look at the smallest sets (i.e., those whose size guarantee is
  only to be at least $\epsilon^\ell\cdot\lvert U\rvert$), then one of them must be $0$-good.

  Suppose for a contradiction that $U$ does not have any subset $W$ with the desired properties and let us construct sets
  $U_\tau$ for $\tau\in\{0,1\}^{<\ell+1}$ with
  \begin{gather*}
    \lvert U_\tau\rvert \geq \epsilon^{\lvert\tau\rvert}\cdot\lvert U\rvert.
  \end{gather*}
  and vertices $v_\sigma$ for $\sigma\in\{0,1\}^{<\ell}$ inductively as follows:
  \begin{enumerate}[label={\arabic*}]
  \item Let $U_\varnothing\df U$ (note that $\lvert U_\varnothing\rvert > \epsilon^0\cdot\lvert U\rvert -
    \One[\varnothing=\varnothing]$).
  \item Given $U_\tau$ for $\lvert\tau\rvert < \ell$, since
    \begin{gather*}
      \lvert U_\tau\rvert
      \geq
      \epsilon^{\lvert\tau\rvert}\cdot\lvert U\rvert
      \geq
      \epsilon^{\ell-1}\cdot\lvert U\rvert,
    \end{gather*}
    and $U_\tau\subseteq U$, our assumption toward a contradiction implies that $U_\tau$ is not $\epsilon$-good, so there exists
    $v_\tau\in V(G)$ such that the sets
    \begin{align*}
      U_{\tau\conc 0} & \df N_G^0(v_\tau)\cap U_\tau, &
      U_{\tau\conc 1} & \df N_G^1(v_\tau)\cap U_\tau,
    \end{align*}
    both have size at least $\epsilon\cdot\lvert U_\tau\rvert\geq\epsilon^{\lvert\tau\rvert+1}\cdot\lvert U\rvert$.
  \end{enumerate}

  We now note that for every $\tau\in\{0,1\}^\ell$, we have $\lvert U_\tau\rvert\geq\epsilon^\ell\cdot\lvert U\rvert$ and
  $U_\tau\subseteq U$, so our assumption toward a contradiction says that $U_\tau$ is not $0$-good, so there exists a vertex
  $v_\tau\in V(G)$ such that the sets
  \begin{align*}
    N_G^0(v_\tau)\cap U_\tau, & &
    N_G^1(v_\tau)\cap U_\tau,
  \end{align*}
  are both non-empty, so we can let $u_{\tau\conc 0}$ and $u_{\tau\conc 1}$ be any vertices in the sets above, respectively.

  It is now easy to see that $((v_\sigma)_{\sigma\in\{0,1\}^{<\ell+1}}, (N_G(u_\tau))_{\tau\in\{0,1\}^{\ell+1}})$ is an
  $(\ell+1)$-Littlestone tree in $\cH_G$, contradicting the fact that $\Lit(G)=\Lit(\cH_G)\leq\ell$.

  \medskip

  For item~\ref{cor:Ramsey:0good}, since $U$ is $0$-good, it follows that either every $v\in V(G)$ is adjacent to all vertices
  of $U$, or either every $v\in V(G)$ is non-adjacent to every vertex of $U$. The former case is impossible since for any $u\in
  U$ (which exists since $U\neq\varnothing$), the vertex $u$ is not adjacent to itself. Thus, every $v\in V(G)$ is non-adjacent
  to every vertex of $U$, which implies that $U$ is an anticlique of $G$.

  \medskip

  For item~\ref{cor:Ramsey:hom}, let $H\df G\rest_U$ and let $m\df\lvert U\rvert=\lvert H\rvert$. Since $(U,U)$ is
  $\epsilon$-homogeneous in $G$, we know that
  \begin{gather*}
    \frac{\lvert\vec{E}(H)\rvert}{m^2}  = d_H(U,U) = d_G(U,U) \in [0,\epsilon]\cup[1-\epsilon,1],
  \end{gather*}
  which in particular implies that either $\lvert\vec{E}(H)\rvert\geq (1-\epsilon)\cdot m^2$ or
  $\lvert\vec{E}(\overline{H})\rvert\geq (1-\epsilon)\cdot m^2 - m$ (the negative factor accounts for
  the fact that the diagonal consists of non-edges in both $H$ and $\overline{H}$). Thus, we get
  \begin{gather*}
    \max\{\lvert E(H)\rvert, \lvert E(\overline{H})\rvert\} \geq \left(1 - \epsilon - \frac{1}{m}\right)\cdot\frac{m^2}{2}.
  \end{gather*}
  Now let
  \begin{gather*}
    k \df \Ceil{\left(\epsilon - \frac{1}{m}\right)^{-1}} - 1 < \left(\epsilon - \frac{1}{m}\right)^{-1}
  \end{gather*}
  and consider the \Turan\ Graph $T(m,k)$ on $m$ vertices with $k$ parts, i.e., it is the complete $k$-partite graph on $m$
  vertices with parts as balanced as possible. Then
  \begin{gather*}
    \lvert E(T(m,k))\rvert
    \leq
    \binom{k}{2}\cdot\left(\frac{m}{k}\right)^2
    =
    \frac{m^2}{2}\cdot\left(1 - \frac{1}{k}\right)
    <
    \frac{m^2}{2}\cdot\left(1 - \epsilon - \frac{1}{m}\right),
  \end{gather*}
  so by \Turan's Theorem, either $H$ or $\overline{H}$ has a clique of size
  \begin{gather*}
    k + 1 = \Ceil{\left(\epsilon - \frac{1}{m}\right)^{-1}};
  \end{gather*}
  in the former case, $U$ has a clique of size $k+1$ and in the latter case, $U$ has an anticlique of size $k+1$.

  \medskip

  For item~\ref{cor:Ramsey:Ramsey}, let
  \begin{gather*}
    \epsilon \df \left(\frac{1}{2\cdot n}\right)^{1/(\ell+1)}.
  \end{gather*}
  By item~\ref{cor:Ramsey:ext}, there exists $U\subseteq V(G)$ such that either $\lvert U\rvert\geq\epsilon^{\ell-1}\cdot n$ and
  $U$ is $\epsilon$-good in $G$, or $\lvert U\rvert\geq\epsilon^\ell\cdot n$ and $U$ is $0$-good in $G$.

  In the latter case, item~\ref{cor:Ramsey:0good} guarantees $U$ is an anticlique and its size is
  \begin{gather*}
    \lvert U\rvert\geq\epsilon^\ell\cdot n = \frac{n^{1/(\ell+1)}}{2^{1 - 1/(\ell+1)}}.
  \end{gather*}

  In the former case, we have
  \begin{gather*}
    \lvert U\rvert\geq\epsilon^{\ell-1}\cdot n = \frac{n^{2/(\ell+1)}}{2^{1 - 2/(\ell+1)}}.
  \end{gather*}
  Let
  \begin{align*}
    \epsilon'
    & \df
    \begin{dcases*}
      \frac{1 - \sqrt{1 - 8\cdot\epsilon\cdot(1-\epsilon)}}{2}, & if $\epsilon < (2-\sqrt{2})/4$,\\
      1/2, & otherwise,
    \end{dcases*}
    \\
    & =
    2\cdot\epsilon + O_{\epsilon\to 0}(\epsilon^2)
    =
    \frac{2^{1 - 1/(\ell+1)}}{n^{1/(\ell+1)}} + O_{n\to\infty}(n^{-2/(\ell+1)})
  \end{align*}
  be as in Lemma~\ref{lem:good->hom} so that its item~\ref{lem:good->hom:hom} implies that $(U,U)$ is $\epsilon'$-homogeneous in
  $G$, which in turn by item~\ref{cor:Ramsey:hom} of this corollary implies that there exists $W\subseteq U$ that is either a
  clique or an anticlique of $G$ of size
  \begin{align*}
    \lvert W\rvert
    & \geq
    \Ceil{\left(\epsilon' - \frac{1}{\lvert U\rvert}\right)^{-1}}
    \\
    & \geq
    \Ceil{
      \left(
      \frac{2^{1 - 1/(\ell+1)}}{n^{1/(\ell+1)}} + O_{n\to\infty}(n^{-2/(\ell+1)}) - \frac{2^{1 - 2/(\ell+1)}}{n^{2/(\ell+1)}}
      \right)^{-1}
    }
    =
    \frac{n^{1/(\ell+1)}}{2^{1 - 1/(\ell+1)}} - O_{n\to\infty}(1).
    \qedhere
  \end{align*}
\end{proof}


\section{Non-equitable partitions}

\begin{lemma}[Calculations for non-equitable partitions]\label{lem:nonequicalc}
  Let $\ell\in\NN_+$, let $c_{\UC}\in[0,1)$ and let $c_{\ZZ},c_{\abs}\in(0,1)$. Define
  \begin{align*}
    \widetilde{c}_{\UC} & \df \frac{1}{(1-c_{\UC})^\ell} - 1,
    &
    \widetilde{c}_{\ZZ} & \df \frac{1}{(1-c_{\ZZ})^\ell} - 1,
    &
    \widetilde{c}_{\abs} & \df 2\cdot c_{\abs},
  \end{align*}
  and for $\epsilon\in(0,1)$, let
  \begin{gather*}
    \begin{aligned}
      \zeta
      & \df
      (1-c_{\ZZ})\cdot\epsilon,
      \\
      \widetilde{\epsilon}
      & \df
      (1 - \zeta^{c_{\abs}})\cdot\zeta,
      \\
      \gamma
      & \df
      (1 - c_{\UC})\cdot\widetilde{\epsilon},
    \end{aligned}
    \\
    K
    \df
    K_{\nonequi}(\ell,c_{\UC},c_{\ZZ},c_{\abs},\epsilon)
    \df
    1
    +
    \frac{
      \ln(\zeta-\widetilde{\epsilon}) - \ln(1-\widetilde{\epsilon})
    }{
      \ln(1 - \gamma^\ell)
    }.
  \end{gather*}

  Define further
  \begin{align*}
    \MoveEqLeft
    \epsilon_{\nonequi}(\ell,c_{\ZZ},c_{\abs})
    \\
    & \df
    \sup\left\{\epsilon'\in(0,1)
    \;\middle\vert\;
    \forall\epsilon\in(0,\epsilon'),
    \frac{\ln(1/\zeta)}{\zeta^\ell}\cdot\frac{1+c_{\abs}}{(1-\zeta^{c_{\abs}})^\ell}
    \leq
    (1+\widetilde{c}_{\ZZ})\cdot(1+2\cdot c_{\abs})\cdot\frac{\ln(1/\epsilon)}{\epsilon^\ell}
    \right\}
  \end{align*}
  (note that the above does not depend on $c_{\UC}$) and for $\epsilon\in(0,\epsilon_{\nonequi}(\ell,c_{\ZZ},c_{\abs}))$, define
  \begin{gather*}
    n_{\nonequi}(\ell,c_{\UC},c_{\ZZ},c_{\abs},\epsilon)
    \df
    \Ceil{
      \frac{1 - \zeta}{(\epsilon - \zeta)\cdot(1+(1 - \gamma^\ell)^{K-1}\cdot\gamma^\ell)}
    }.
  \end{gather*}

  Then the following hold:
  \begin{enumerate}
  \item\label{lem:nonequicalc:epsilonnonequi} We have
    \begin{gather*}
      \epsilon_{\nonequi}(\ell,c_{\ZZ},c_{\abs})
      \geq
      \bigl(1+o_{c_{\abs}\to 0,\ell}(1)\bigr)\cdot
      \left(\frac{c_{\abs}}{1 + 2\cdot c_{\abs}}\right)^{1/(\ell\cdot c_{\abs})}
      \cdot\frac{1}{1-c_{\ZZ}}.
    \end{gather*}
  \item\label{lem:nonequicalc:Knonequi} For $\epsilon\in(0,\epsilon_{\nonequi}(\ell,c_{\ZZ},c_{\abs}))$, we have
    \begin{gather*}
      K_{\nonequi}(\ell,c_{\UC},c_{\ZZ},c_{\abs},\epsilon)
      \leq
      (1+\widetilde{c}_{\UC})\cdot(1+\widetilde{c}_{\ZZ})\cdot(1+\widetilde{c}_{\abs})
      \cdot\frac{\ln(1/\epsilon)}{\epsilon^\ell}.
    \end{gather*}
  \item\label{lem:nonequicalc:nnonequiimplies} For every $\epsilon\in(0,\epsilon_{\nonequi}(\ell,c_{\ZZ},c_{\abs}))$ and every
    $n\in\NN_+$ with $n\geq n_{\nonequi}(\ell,c_{\UC},c_{\ZZ},c_{\abs},\epsilon)$, we have
    \begin{gather}\label{eq:nonequicalc:nnonequiimplies}
      (\zeta-1)\cdot\left(1 - \frac{1}{1 + (1-\gamma^\ell)^{K-1}\cdot\gamma^\ell\cdot n}\right)
      +
      1
      \leq
      \epsilon.
    \end{gather}
  \item\label{lem:nonequicalc:nnonequi} For $\epsilon\in(0,\epsilon_{\nonequi}(\ell,c_{\ZZ},c_{\abs}))$, we have
    \begin{gather*}
      n_{\nonequi}(\ell,c_{\UC},c_{\ZZ},c_{\abs},\epsilon)
      \leq
      \frac{1}{c_{\ZZ}\cdot\epsilon} + 1.
    \end{gather*}
  \end{enumerate}
\end{lemma}

\begin{proof}
  For item~\ref{lem:nonequicalc:epsilonnonequi}, we note that for $\epsilon\in(0,1)$, we have
  \begin{align*}
    \frac{\ln(1/\zeta)}{\zeta^\ell}\cdot\frac{1+c_{\abs}}{(1-\zeta^{c_{\abs}})^\ell}
    & =
    \frac{
      \ln(1/\epsilon) + \ln(1/(1-c_{\ZZ}))
    }{
      (1-c_{\ZZ})^\ell\cdot\epsilon^\ell
    }\cdot\frac{1+c_{\abs}}{(1-(1-c_{\ZZ})^{c_{\abs}}\epsilon^{c_{\abs}})^\ell}
    \\
    & =
    \frac{1}{(1-(1-c_{\ZZ})^{c_{\abs}}\epsilon^{c_{\abs}})^\ell}
    \cdot
    \frac{1+\widetilde{c}_{\ZZ}}{\epsilon^\ell}
    \cdot
    \left(
    \ln\left(\frac{1}{\epsilon}\right)
    +
    \ln\left(\frac{1}{1-c_{\ZZ}}\right)
    \right)
    \\
    &
    \leq
    \frac{1}{(1-(1-c_{\ZZ})^{c_{\abs}}\epsilon^{c_{\abs}})^\ell}
    \cdot
    \frac{1+\widetilde{c}_{\ZZ}}{\epsilon^\ell}
    \cdot
    \left(
    \ln\left(\frac{1}{\epsilon}\right)
    +
    \ln\left(\frac{1}{1-c_{\ZZ}}\right)
    \right).
  \end{align*}
  Thus, to show that the above is at most
  $(1+\widetilde{c}_{\ZZ})\cdot(1+2\cdot c_{\abs})\cdot\epsilon^{-\ell}\cdot\ln(1/\epsilon)$, it suffices to show that
  \begin{align*}
    \ln\left(\frac{1}{1-c_{\ZZ}}\right)
    & \leq
    c_{\abs}\cdot\ln\left(\frac{1}{\epsilon}\right),
    \\
    \frac{1+c_{\abs}}{(1-(1-c_{\ZZ})^{c_{\abs}}\epsilon^{c_{\abs}})^\ell}
    & \leq
    1 + 2\cdot c_{\abs}.
  \end{align*}
  In turn, the above are equivalent respectively to:
  \begin{align*}
    \epsilon
    & \leq
    (1-c_{\ZZ})^{1/c_{\abs}}
    =
    1 + o_{c_{\abs}\to 0,\ell}(1)
    =
    1 + o_{c_{\abs}\to 0,\ell}(1)
    \\
    \epsilon
    & \leq
    \left(\frac{c_{\abs}}{1 + 2\cdot c_{\abs}}\right)^{1/(\ell\cdot c_{\abs})}
    \cdot\frac{1}{1-c_{\ZZ}}
  \end{align*}
  where the first equality follows since $1 - c_{\ZZ} = (1+\widetilde{c}_{\ZZ})^{-1/\ell} > 2^{-1/\ell} > 0$, so the second
  condition dominates the asymptotic behavior of $\epsilon_{\nonequi}(\ell,c_{\ZZ},c_{\abs})$ as $c_{\abs}\to 0$ with $\ell$
  fixed, yielding the result.

  \medskip

  We now prove item~\ref{lem:nonequicalc:Knonequi}. Note that
  \begin{align*}
    K_{\nonequi}(\ell,c_{\UC},c_{\ZZ},c_{\abs},\epsilon)
    & \df
    1
    +
    \frac{
      \ln(\zeta-\widetilde{\epsilon}) - \ln(1-\widetilde{\epsilon})
    }{
      \ln(1 - \gamma^\ell)
    }
    =
    1
    +
    \frac{
      \ln(\zeta-\widetilde{\epsilon}) - \ln(1-\widetilde{\epsilon})
    }{
      \ln(1 - (1-c_{\UC})^\ell\cdot\widetilde{\epsilon}^\ell)
    }
    \\
    & \leq
    1
    +
    \frac{
      \ln(\zeta-\widetilde{\epsilon}) - \ln(1-\widetilde{\epsilon})
    }{
       -(1-c_{\UC})^\ell\cdot\widetilde{\epsilon}^\ell
    }
    \leq
    1
    +
    \frac{
      \ln(\zeta-\widetilde{\epsilon}) + \widetilde{\epsilon}
    }{
       -(1-c_{\UC})^\ell\cdot\widetilde{\epsilon}^\ell
    }
    \\
    & =
    1
    +
    (1+\widetilde{c}_{\UC})\cdot
    \frac{
      (1+c_{\abs})\cdot\ln(\zeta) + (1-\zeta^{c_{\abs}})\cdot\zeta
    }{
       -(1-\zeta^{c_{\abs}})^\ell\cdot\zeta^\ell
    }
    \\
    & =
    (1+\widetilde{c}_{\UC})\cdot
    \frac{\ln(1/\zeta)}{\zeta^\ell}\cdot\frac{1+c_{\abs}}{(1-\zeta^{c_{\abs}})^\ell}
    +
    1
    -
    (1+\widetilde{c}_{\UC})\cdot(1-\zeta^{c_{\abs}})^{1-\ell}\cdot\zeta^{1-\ell}
    \\
    & \leq
    (1+\widetilde{c}_{\UC})\cdot
    \frac{\ln(1/\zeta)}{\zeta^\ell}\cdot\frac{1+c_{\abs}}{(1-\zeta^{c_{\abs}})^\ell}
    \\
    & \leq
    (1+\widetilde{c}_{\UC})\cdot(1+\widetilde{c}_{\ZZ})\cdot(1+2\cdot c_{\abs})
    \cdot\frac{\ln(1/\epsilon)}{\epsilon^\ell}
    \\
    & =
    (1+\widetilde{c}_{\UC})\cdot(1+\widetilde{c}_{\ZZ})\cdot(1+\widetilde{c}_{\abs})
    \cdot\frac{\ln(1/\epsilon)}{\epsilon^\ell},
  \end{align*}
  where the first two inequalities follow from $\ln(1+x)\leq x$ (note that both the numerator and denominator are negative) and
  the last inequality follows since $\epsilon < \epsilon_{\nonequi}(\ell,c_{\ZZ},c_{\abs})$.

  \medskip

  For item~\ref{lem:nonequicalc:nnonequi}, we note that the left-hand side of~\eqref{eq:nonequicalc:nnonequiimplies} is a
  decreasing function of $n$ (when all the other parameters are fixed). Thus, it suffices to show only the case
  $n=n_{\nonequi}(\ell,c_{\UC},c_{\ZZ},c_{\abs},\epsilon)$. In turn, this case amounts to a straightforward computation that
  shows that~\eqref{eq:nonequicalc:nnonequiimplies} is true exactly when $n$ is at least the expression under the ceiling in the
  definition of $n_{\nonequi}(\ell,c_{\UC},c_{\ZZ},c_{\abs},\epsilon)$.

  \medskip

  For item~\ref{lem:nonequicalc:nnonequi}, we simply note that
  \begin{gather*}
    \Ceil{
      \frac{1 - \zeta}{(\epsilon - \zeta)\cdot(1+(1 - \gamma^\ell)^{K-1}\cdot\gamma^\ell)}
    }
    \leq
    \Ceil{\frac{1 - \zeta}{\epsilon - \zeta}}
    \leq
    \Ceil{\frac{1}{\epsilon-\zeta}}
    =
    \Ceil{\frac{1}{c_{\ZZ}\cdot\epsilon}}
    \leq
    \frac{1}{c_{\ZZ}\cdot\epsilon}
    \qedhere
  \end{gather*}
\end{proof}

\begin{remark}\label{rmk:nonequicalc}
  It is straightforward to check that the constants $\widetilde{c}_{\UC},\widetilde{c}_{\ZZ},\widetilde{c}_{\abs}$ of
  Lemma~\ref{lem:nonequicalc} are bijective increasing functions of the respective constants $c_{\UC},c_{\ZZ},c_{\abs}$ and
  their inverses are given by
  \begin{align*}
    c_{\UC} & \df 1 - \frac{1}{(1+\widetilde{c}_{\UC})^{1/\ell}},
    &
    c_{\ZZ} & \df 1 - \frac{1}{(1+\widetilde{c}_{\ZZ})^{1/\ell}},
    &
    c_{\abs} & \df \frac{\widetilde{c}_{\abs}}{2},
  \end{align*}
\end{remark}

\begin{theorem}[Non-equitable partitions]\label{thm:nonequi}
  Let $\ell,d\in\NN_+$ with $d\leq\ell$, let $\delta,\epsilon\in(0,1)$, let $c_{\atom},\rho_{\UC}\in(0,1)\cap\QQ$ be rationals, let
  $c_{\UC}\in[0,1)\cap\QQ$ be rational and let $G$ be a graph with $\lvert G\rvert=n$, $\Lit(G)\leq\ell$ and $\VC(G)\leq d$.

  Given further $c_{\ZZ},c_{\abs}\in(0,1)$, let
  \begin{gather*}
    \epsilon_{\nonequi}(\ell,c_{\ZZ},c_{\abs})
    \geq
    \bigl(1+o_{c_{\abs}\to 0,\ell}(1)\bigr)\cdot
    \left(\frac{c_{\abs}}{2 + 2\cdot c_{\abs}}\right)^{2/(\ell\cdot c_{\abs})}
    \cdot(1+c_{\ZZ})^{1/\ell}
  \end{gather*}
  be as in Lemma~\ref{lem:nonequicalc} and for $\epsilon\in(0,\epsilon_{\nonequi}(\ell,c_{\ZZ},c_{\abs}))$, let
  \begin{align*}
    K_{\nonequi}(\ell,c_{\UC},c_{\ZZ},c_{\abs},\epsilon)
    & \leq
    (1+\widetilde{c}_{\UC})\cdot(1+\widetilde{c}_{\ZZ})\cdot(1+\widetilde{c}_{\abs})
    \cdot\frac{\ln(1/\epsilon)}{\epsilon^\ell},
    \\
    n_{\nonequi}(\ell,c_{\UC},c_{\ZZ},c_{\abs},\epsilon)
    & \leq
    \frac{1}{c_{\ZZ}\cdot\epsilon} + 1
  \end{align*}
  also be as in Lemma~\ref{lem:nonequicalc}.

  Then for $\epsilon\in(0,\epsilon_{\nonequi}(\ell,c_{\ZZ},c_{\abs}))$ and $n\in\NN_+$ with $n\geq
  n_{\nonequi}(\ell,c_{\UC},c_{\ZZ},c_{\abs},\epsilon)$, the following hold:
  \begin{enumerate}
  \item\label{thm:nonequi:exists} If $c_{\UC}=0$ and $\delta < 1/2^{\ell+1}$, there there exists a $(\delta,\epsilon)$-atomic
    partition of $G$ into at most $K_{\nonequi}(\ell,0,c_{\ZZ},c_{\abs},\epsilon)$ parts.
  \end{enumerate}
  For the items below, we further assume that $\delta$, $\epsilon$, $c_{\ZZ}$ and $c_{\abs}$ are rational (the other constants
  are already assumed to be rational):
  \begin{enumerate}[resume*]
  \item\label{thm:nonequi:good} If $c_{\UC} > 0$, then in the random-query model, Algorithm~\ref{alg:nonequi:good}, with
    probability at least $1-\rho_{\UC}$, computes an $\epsilon$-good partition of $G$ into at most
    $K_{\nonequi}(\ell,c_{\UC},c_{\ZZ},c_{\abs},\epsilon)$ parts in time
    \begin{multline*}
      O\bigggl(
      \frac{2^\ell\cdot\ell^2\cdot(\log(1/\epsilon))^2}{c_{\UC}^2\cdot(1-c_{\ZZ})^2\cdot\epsilon^{\ell+2}}
      \cdot n\cdot\log\left(\frac{n}{\rho_{\UC}}\right)
      \\
      +
      \frac{\log(1/\epsilon)}{\epsilon^\ell}
      \cdot\ell^2\cdot\log(\ell+1)
      \cdot 2^{3\cdot\len(c_{\abs})}\cdot\bigl(\len(c_{\ZZ})+\len(\epsilon)\bigr)^3\cdot\len(c_{\UC})^6\cdot\len(\rho_{\UC})\cdot
      \bigl(\log(n+1)\bigr)^3
      \bigggr).
    \end{multline*}
  \item\label{thm:nonequi:gooddet} If $c_{\UC}=0$, then in the deterministic model, Algorithm~\ref{alg:nonequi:gooddet} computes
    an $\epsilon$-good partition of $G$ into at most $K_{\nonequi}(\ell,0,c_{\ZZ},c_{\abs},\epsilon)$ parts in time
    \begin{gather*}
      O\left(
      \ell\cdot\frac{\log(1/\epsilon)}{\epsilon^\ell}\cdot n^2
      +
      \frac{\log(1/\epsilon)}{\epsilon^\ell}
      \cdot\ell
      \cdot
      2^{3\cdot\len(c_{\abs})}\cdot\bigl(\len(c_{\ZZ})+\len(\epsilon)\bigr)^3\cdot\bigl(\log(n+1)\bigr)^3
      \right).
    \end{gather*}
  \item\label{thm:nonequi:atom} If $c_{\UC}=0$ and $(1+c_{\atom})\cdot\delta < 1/2^{\ell+1}$, then in the deterministic model,
    Algorithm~\ref{alg:nonequi:atom} computes a $(\delta,\epsilon)$-atomic partition of $G$ into at most
    $K_{\nonequi}(\ell,0,c_{\ZZ},c_{\abs},\epsilon)$ parts in time
    \begin{multline*}
      \frac{\log(1/\epsilon)}{\epsilon^\ell}\cdot n^{O(d/(c_{\atom}^2\cdot\delta^2))}
      +
      O\bigggl(
      \len(c_{\atom})^2\cdot\len(\delta)^2\cdot\bigl(\log(d+1)\bigr)^2
      \\
      +
      \frac{\log(1/\epsilon)}{\epsilon^\ell}
      \cdot\ell
      \cdot
      2^{3\cdot\len(c_{\abs})}\cdot\bigl(\len(c_{\ZZ})+\len(\epsilon)\bigr)^3\cdot\bigl(\log(n+1)\bigr)^3
      \bigggr).
    \end{multline*}
  \item\label{thm:nonequi:goodspace} If $c_{\UC}=0$, then in the query-into-oracle model, Algorithms~\ref{alg:nonequi:goodspace}
    and~\ref{alg:nonequi:atomspaceoracle} compute an $\epsilon$-good partition of $G$ into at most
    $K_{\nonequi}(\ell,0,c_{\ZZ},c_{\abs},\epsilon)$ parts. The space-complexities of Algorithms~\ref{alg:nonequi:goodspace}
    and~\ref{alg:nonequi:atomspaceoracle} are
    \begin{gather*}
      O\left(
      \frac{\ell\cdot(\log(1/\epsilon))^2\cdot\log(n+1)}{\epsilon^\ell}
      +
      2^{2\cdot\len(c_{\abs})}\cdot\bigl(\len(c_{\ZZ}) + \len(\epsilon) + \log(n+1)\bigr)
      \right),
      \\
      O\Bigl(K\cdot\bigl(\log(K+1) + \log(m+1) + \log(n+1)\bigr)\Bigr)
      \leq
      O\left(
      \frac{\log(1/\epsilon)}{\epsilon^\ell}\cdot\left(\log\left(\frac{1}{\epsilon}\right) + \log(n+1)\right)
      \right),
    \end{gather*}
    respectively, and the time-complexities are
    \begin{gather*}
      n^{O(\epsilon^{-\ell}\cdot\log(1/\epsilon))}
      +
      O\left(
      \frac{\log(1/\epsilon)}{\epsilon^\ell}
      \cdot\ell
      \cdot
      2^{3\cdot\len(c_{\abs})}\cdot\bigl(\len(c_{\ZZ})+\len(\epsilon)\bigr)^3\cdot\bigl(\log(n+1)\bigr)^3
      \right),
      \\
      O\bigl(K^2\cdot n^{K+1}\bigr)
      \leq
      n^{O(\epsilon^{-\ell}\cdot\log(1/\epsilon))},
    \end{gather*}
    respectively.
  \item\label{thm:nonequi:atomspace} If $c_{\UC}=0$ and $(1+c_{\atom})\cdot\delta < 1/2^{\ell+1}$, then in the query-into-oracle
    model, Algorithms~\ref{alg:nonequi:atomspace} and~\ref{alg:nonequi:atomspaceoracle} compute an $(\delta,\epsilon)$-atomic
    partition of $G$ into at most $K_{\nonequi}(\ell,0,c_{\ZZ},c_{\abs},\epsilon)$ parts. The space-complexities of
    Algorithms~\ref{alg:nonequi:atomspace} and~\ref{alg:nonequi:atomspaceoracle} are
    \begin{gather*}
      \begin{multlined}[t]
        O\bigggl(
        \frac{\ell\cdot d\cdot(\log(1/\epsilon))^2\cdot\log(n+1)}{c_{\atom}^2\cdot\delta^2\cdot\epsilon^\ell}
        +
        2^{2\cdot\len(c_{\abs})}\cdot\bigl(\len(c_{\ZZ}) + \len(\epsilon) + \log(n+1)\bigr)
        \\
        +
        \len(c_{\atom}) + \len(\delta)
        \bigggr),
      \end{multlined}
      \\
      O\Bigl(K\cdot\bigl(\log(K+1) + \log(m+1) + \log(n+1)\bigr)\Bigr)
      \leq
      O\left(
      \frac{\log(1/\epsilon)}{\epsilon^\ell}\cdot
      \log\left(\frac{n}{c_{\atom}\cdot\delta\cdot\epsilon^\ell}\right)
      \right),
    \end{gather*}
    respectively, and the time-complexities are
    \begin{gather*}
      \begin{multlined}[t]
        n^{O(d/(c_{\atom}^2\cdot\delta^2) + \epsilon^{-\ell}\cdot\log(1/\epsilon))}
        +
        O\bigggl(
        \len(c_{\atom})^2\cdot\len(\delta)^2\cdot\bigl(\log(d+1)\bigr)^2
        \\
        +
        \frac{\log(1/\epsilon)}{\epsilon^\ell}
        \cdot\ell
        \cdot
        2^{3\cdot\len(c_{\abs})}\cdot\bigl(\len(c_{\ZZ})+\len(\epsilon)\bigr)^3\cdot\bigl(\log(n+1)\bigr)^3
        \bigggr),
      \end{multlined}
      \\
      O\bigl(K^2\cdot m\cdot n^{K+1}\bigr)
      \leq
      \frac{n^{O(\epsilon^{-\ell}\cdot\log(1/\epsilon))}}{c_{\atom}^2\cdot\delta^2}.
    \end{gather*}
    respectively.
  \end{enumerate}

  Furthermore, for
  \begin{align*}
    \epsilon'
    & \df
    \begin{dcases*}
      \frac{1 - \sqrt{1 - 8\cdot\epsilon\cdot(1-\epsilon)}}{2}, & if $\epsilon < (2-\sqrt{2})/4$,\\
      1/2, & otherwise,
    \end{dcases*}
    \\
    & =
    2\cdot\epsilon + O_{\epsilon\to 0}(\epsilon^2),
    \\
    \epsilon''
    & \df
    \frac{\epsilon'}{1-\epsilon}
    =
    \frac{1 - \sqrt{1 - 8\cdot\epsilon\cdot(1-\epsilon)}}{2\cdot(1-\epsilon)}
    =
    2\cdot\epsilon + O_{\epsilon\to 0}(\epsilon^2),
  \end{align*}
  all $(\delta,\epsilon)$-atomic partitions are $(\delta,\epsilon)$-excellent and $\epsilon$-good; and all $\epsilon$-good
  partitions of $G$ are in particular totally $\epsilon'$-homogeneous and $(\epsilon,\epsilon'')$-excellent.
\end{theorem}

\begin{proof}
  Before we start the proof, let us point out that the final that the corresponding $(\delta,\epsilon)$-atomic partitions are
  also $(\delta,\epsilon)$-excellent partitions and $\epsilon$-good partitions follow from Lemma~\ref{lem:atom->exc} (and the
  fact that since $\delta < 1$, $(\delta,\epsilon)$-excellence implies $\epsilon$-goodness).

  Similarly, the final assertions of all items that the corresponding $\epsilon$-good partitions are also totally
  $\epsilon'$-homogeneous partitions follow from Lemma~\ref{lem:good->hom}\ref{lem:good->hom:partition}.

  Item~\ref{thm:nonequi:exists} follows from its algorithmic counterpart, item~\ref{thm:nonequi:atom}, in the following way:
  since $\delta < 1/2^{\ell+1}$, this allows us to pick $c_{\atom}\in(0,1)\cap\QQ$ small enough so that
  $(1+c_{\atom})\cdot\delta < 1/2^{\ell+1}$ (note that $c_{\atom}$ does not affect any size bounds; it only affects algorithmic
  complexities). We can then apply item~\ref{thm:nonequi:atom} to obtain our desired $(\delta,\epsilon)$-atomic partition,
  except that we need to do a small ``rationalization of parameters trick'': if any of the parameters $\delta$, $\epsilon$,
  $c_{\ZZ}$ and $c_{\abs}$ are irrational, then we can approximate them to a rational to enough precision so that the notion of
  atomicity does not change on a fixed graph $G$ and the size guarantee does not change; alternatively, one can think of running
  the algorithm on a computational model that ``accepts irrational values'' (i.e., simply follow the algorithm instructions to
  make a mathematical construction, then follow its proof of correctness for a proof that the construction has the desired
  properties).

  Before we prove the algorithmic items, let us point out that they all follow the same strategy:
  \begin{enumerate}[label={\arabic*.}]
  \item For items about atomicity, we will start by letting $m$ be given by Lemma~\ref{lem:atomconst}.
  \item For all items except for item~\ref{thm:nonequi:good}, we will set $c_{\UC}\df 0$ (the excluded item already receives
    $c_{\UC}\in(0,1)\cap\QQ$ as input).
  \item Using the parameters of Lemma~\ref{lem:nonequicalc}, we will repeatedly extract $\widetilde{\epsilon}$-good set or
    $((1+c_{\atom})\cdot\delta,\widetilde{\epsilon})$-atoms of relative size larger than $\gamma^\ell$ until we are left with a
    set $R$ of size at most
    \begin{gather}\label{eq:nonequi:R}
      \lvert R\rvert
      \leq
      \frac{\zeta-\widetilde{\epsilon}}{1-\widetilde{\epsilon}}\cdot n.
    \end{gather}
    Our choice of parameters using Lemma~\ref{lem:atomconst} ensures that the number of sets extracted is at most
    $K_{\nonequi}(\ell,c_{\UC},c_{\ZZ},c_{\abs},\epsilon)$ (but we still have a remainder set $R$).
  \item We will then redistribute the vertices of the remainder $R$ among the already generated parts proportionally to their
    size; this will deteriorate our goodness/atomicity from $\widetilde{\epsilon}$ to a value that we want to be at most
    $\zeta$. Using Lemma~\ref{lem:upcont}, we will argue that a part of size $s$ can absorb at most
    \begin{gather}\label{eq:nonequi:absorption}
      \frac{\zeta-\widetilde{\epsilon}}{1-\zeta}\cdot s
    \end{gather}
    vertices while having the deterioration be from $\widetilde{\epsilon}$ to at most $\zeta$. Comparing with our stopping
    condition~\eqref{eq:nonequi:R}, it might seem that we can absorb the whole remainder in the parts
    because~\eqref{eq:nonequi:R} is equivalent to
    \begin{gather*}
      \lvert R\rvert
      \leq
      \frac{\zeta-\widetilde{\epsilon}}{1-\widetilde{\epsilon}}\cdot(n-\lvert R\rvert).
    \end{gather*}
    However, the issue is that the value in~\eqref{eq:nonequi:absorption} need not be an integer, so a part $s$ can only
    realistically absorb at most
    \begin{gather*}
      \Floor{\frac{\zeta-\widetilde{\epsilon}}{1-\zeta}\cdot s}
    \end{gather*}
    many vertices while having its goodness/atomicity deteriorate to at most $\zeta$. This means that when we use these integer
    values instead, we get a partition with goodness/atomicity parameter $\zeta$, at most
    $K_{\nonequi}(\ell,c_{\UC},c_{\ZZ},c_{\abs},\epsilon)$ parts, and a new remainder set $R'$ of size at most
    $K_{\nonequi}(\ell,c_{\UC},c_{\ZZ},c_{\abs},\epsilon)$ (as we can only have absorbed at most $1$ less than we wanted to absorb for
    each part).
  \item Since $\lvert R'\rvert\leq K_{\nonequi}(\ell,c_{\UC},c_{\ZZ},c_{\abs},\epsilon)$. We can now add at most $1$ of these vertices to
    each part and use Lemma~\ref{lem:upcont} again to claim that the goodness/atomicity will deteriorate from $\zeta$ to at most
    $\epsilon$; this will be a consequence of the fact that the parts are large since $n\geq
    n_{\nonequi}(\ell,c_{\UC},c_{\ZZ},c_{\abs},\epsilon)$. Note that this two-step absorption is equivalent to making each part
    of size $s$ absorb at most
    \begin{gather*}
      \Ceil{\frac{\zeta-\widetilde{\epsilon}}{1-\zeta}\cdot s}
    \end{gather*}
    vertices.
  \item Finally, for the atomicity items, we actually obtained a partition into $((1+c_{\atom})\cdot\delta,\epsilon,m)$-atoms of
    $G$; our choice of $m$ as in Lemma~\ref{lem:atomconst} then ensures that these are $(\delta,\epsilon)$-atoms of $G$.
  \end{enumerate}

  Let us then prove the algorithmic items (and we will prove them slightly out of order for simplicity sake):
  \begin{description}[wide, itemsep={3ex}]
  \item[Item~\ref{thm:nonequi:atom}.] We start by arguing correctness of Algorithm~\ref{alg:nonequi:atom}.

    \begin{algorithm}[htbp]
      \caption{Deterministic model algorithm that returns a $(\delta,\epsilon)$-atomic partition of $G$ into at most
        $K_{\nonequi}(\ell,0,c_{\ZZ},c_{\abs},\epsilon)$ parts. The time-complexity is:
        \begin{multline*}
          \frac{\log(1/\epsilon)}{\epsilon^\ell}\cdot n^{O(d/(c_{\atom}^2\cdot\delta^2))}
          +
          O\bigggl(
          \len(c_{\atom})^2\cdot\len(\delta)^2\cdot\bigl(\log(d+1)\bigr)^2
          \\
          +
          \frac{\log(1/\epsilon)}{\epsilon^\ell}
          \cdot\ell
          \cdot
          2^{3\cdot\len(c_{\abs})}\cdot\bigl(\len(c_{\ZZ})+\len(\epsilon)\bigr)^3\cdot\bigl(\log(n+1)\bigr)^3
          \bigggr).
        \end{multline*}
      }
      \label{alg:nonequi:atom}
      \DontPrintSemicolon
      \KwIn{Numbers $n,d,\ell\in\NN_+$ with $d\leq\ell\leq\log_2(n)$, $c_{\atom},c_{\ZZ},c_{\abs},\delta\in(0,1)\cap\QQ$ with
        $(1+c_{\atom})\cdot\delta < 1/2^{\ell+1}$, $\epsilon\in(0,\epsilon_{\nonequi}(\ell,c_{\ZZ},c_{\abs}))\cap\QQ$, and a
        graph $G$ with $\lvert G\rvert=n$, $\Lit(G)\leq\ell$ and $\VC(G)\leq d$. We further assume $n\geq
        n_{\nonequi}(\ell,0,c_{\ZZ},c_{\abs},\epsilon)$.}
      \KwOut{A $(\delta,\epsilon)$-atomic partition of $G$ into at most $K_{\nonequi}(\ell,0,c_{\ZZ},c_{\abs},\epsilon)$ parts.}
      $\widetilde{\delta}\assign(1+c_{\atom})\cdot\delta$\;
      Let $m\assign\ceil{C\cdot d/(c_{\atom}^2\cdot\delta^2)}$, where $C$ is the absolute constant of
      Lemma~\ref{lem:atomconst}.\;\label{alg:nonequi:atom:m}
      $\widetilde{m}\assign\floor{\widetilde{\delta}\cdot m}$\;
      $c_{\UC}\assign 0$\;
      $\zeta\assign (1-c_{\ZZ})\cdot\epsilon$\;
      Compute a sufficiently precise rational approximation of $\widetilde{\epsilon}\df (1 - \zeta^{c_{\abs}})\cdot\zeta$ to
      compute $S\df\floor{(\zeta-\widetilde{\epsilon})\cdot n/(1-\widetilde{\epsilon})}$.\;
      \label{alg:nonequi:atom:rational:S}
      $\cQ\assign\varnothing$\;
      $R\assign V(G)$\;
      $r\assign n$\;
      \While{$r > S$}{%
        \label{alg:nonequi:atom:While}
        $s_0\assign r$\;
        Compute sufficiently precise rational approximations of $\widetilde{\epsilon}\df (1 - \zeta^{c_{\abs}})\cdot\zeta$ and
        $\gamma\df (1 - c_{\UC})\cdot\widetilde{\epsilon}$ to compute $s_i\df\floor{\gamma\cdot s_{i-1}} + 1$ for every
        $i\in[\ell+1]$.\;
        \label{alg:nonequi:atom:rational:si}
        Run Algorithm~\ref{alg:ext:atom} with parameters:
        \begin{gather*}
          (n,m,\ell,s_0,\ldots,s_\ell,s_{\ell+1},\widetilde{m},G,U)
          \df
          (n,m,\ell,s_0,\ldots,s_\ell,s_{\ell+1},\widetilde{m},G,R).
        \end{gather*}
        \;\label{alg:nonequi:atom:call}
        Let $(W,s)$ be the pair returned by the algorithm.\;
        $\cQ\assign\cP\cup\{(W,s)\}$\;
        $R\assign R\setminus W$\;
        $r\assign r - s$\;
      }
      \setcounter{algosplit}{\theAlgoLine}
    \end{algorithm}

    \begin{algorithm}[htbp]
      \caption*{(continued).}
      \DontPrintSemicolon
      \setcounter{AlgoLine}{\thealgosplit}
      \let\oldnl\nl
      \let\nl\relax
      \global\let\nl\oldnl 
      $\cP\assign\varnothing$\;
      \For{$(W,s)\in\cQ$}{%
        \label{alg:nonequi:atom:For}
        Compute a sufficently precise approximation of $\widetilde{\epsilon}\df (1-\zeta^{c_{\abs}})\cdot\zeta$ to compute
        $\widetilde{s}\df\ceil{(\zeta-\widetilde{\epsilon})\cdot s/(1-\zeta)}$.\;
        \label{alg:nonequi:atom:rational:ws}
        Let $A\subseteq R$ be a set of size as large as possible with $\lvert A\rvert\leq\widetilde{s}$.\;
        $\cP\assign\cP\cup\{W\cup A\}$\;
        $R\assign R\setminus A$\;
      }
      \Return{$\cP$}
    \end{algorithm}

    Note first that due to the definition of the parameters $s_i$ of line~\ref{alg:nonequi:atom:rational:si}, the parameters $s_i$
    always satisfy\footnote{Even though $\gamma=\widetilde{\epsilon}$, we opt to use the notation $\gamma$ in some places of
    this proof and $\widetilde{\epsilon}$ as it makes the application of Lemma~\ref{lem:nonequicalc} more evident and makes the
    transition to the proof of item~\ref{thm:nonequi:good} simpler.}
    \begin{gather*}
      s_i \geq \floor{\gamma\cdot s_{i-1}} + 1
    \end{gather*}
    for every $i\in[\ell+1]$, which in particular implies that
    \begin{gather*}
      s_i \geq s_\ell \geq \gamma^\ell\cdot s_0 = \gamma^\ell\cdot r,
    \end{gather*}
    for every $i\in\{0,\ldots,\ell\}$.
    
    By Proposition~\ref{prop:ext:atomalg}, we know that each execution of Algorithm~\ref{alg:ext:atom} in
    line~\ref{alg:nonequi:atom:call} returns a pair $(W,s)$, where $W$ is a $(\widetilde{\delta},\widetilde{\epsilon},m)$-atom
    of $G$ of size $s=\lvert W\rvert\in S\df\{s_0,\ldots,s_\ell\}$, which in particular means that
    \begin{gather*}
      s\geq s_\ell\geq\gamma^\ell\cdot r = \gamma^\ell\cdot\lvert R\rvert
    \end{gather*}
    (it is straightforward to check that the variable $r$ keeps track of $\lvert R\rvert$). This means that after $i$ executions
    of the loop of line~\ref{alg:nonequi:atom:While}, the remainder set $R$ has size
    \begin{gather*}
      r = \lvert R\rvert\leq(1-\gamma^\ell)^i\cdot n.
    \end{gather*}
    In particular, the loop of line~\ref{alg:nonequi:atom:While} must execute a finite number of times, say $K$. Furthermore, in
    the beginning of the last iteration of the loop, we must have
    \begin{gather*}
      S
      =
      \Floor{\frac{\zeta-\widetilde{\epsilon}}{1-\widetilde{\epsilon}}\cdot n}
      <
      r
      =
      \lvert R\rvert
      \leq
      (1-\gamma^\ell)^{K-1}\cdot n,
    \end{gather*}
    from which we conclude that
    \begin{gather*}
      K
      \leq
      1 + \frac{\ln(\zeta-\widetilde{\epsilon}) - \ln(1-\widetilde{\epsilon})}{\ln(1-\gamma^\ell)}
      =
      K_{\nonequi}(\ell,0,c_{\ZZ},c_{\abs},\epsilon)
    \end{gather*}
    (recall that $c_{\UC}=\widetilde{c}_{\UC}=0$).

    In particular, by the end of the loop of line~\ref{alg:nonequi:atom:While}, we know that $\lvert\cQ\rvert = K\leq
    K_{\nonequi}(\ell,0,c_{\ZZ},c_{\abs},\epsilon)$, $\lvert R\rvert < (\zeta-\widetilde{\epsilon})\cdot
    n/(1-\widetilde{\epsilon})$ and
    \begin{gather*}
      \cQ' \df \{W \mid (W,s)\in\cQ\}
    \end{gather*}
    is a partition of $V(G)\setminus R$ into $(\widetilde{\delta},\widetilde{\epsilon},m)$-atoms of $G$. Now it is clear that
    the block of line~\ref{alg:nonequi:atom:For} generates a partition $\cP$ from $\cQ'$ by adding at most
    \begin{gather*}
      \Ceil{\frac{\zeta-\widetilde{\epsilon}}{1-\zeta}\cdot\lvert Q\rvert}
    \end{gather*}
    vertices to each part $Q\in\cQ'$. The size of the final partition $\cP$ is the same: $\lvert\cP\rvert=\lvert\cQ'\rvert =
    K\leq K_{\nonequi}(\ell,0,c_{\ZZ},c_{\abs},\epsilon)$. Let us then prove that each set of $P$ is a
    $(\widetilde{\delta},\epsilon,m)$-atom of $G$. To do this, for each $Q\in\cQ'$, we let $P_Q$ be the part of $\cP$ generated
    from $Q$ and we let $P'_Q$ be a set with $Q\subseteq P'_Q\subseteq P_Q$ and
    \begin{gather*}
      \lvert P'_Q\setminus Q\rvert
      =
      \Floor{\frac{\zeta-\widetilde{\epsilon}}{1-\zeta}\cdot\lvert Q\rvert}
      \leq
      \frac{\zeta-\widetilde{\epsilon}}{1-\zeta}\cdot\lvert Q\rvert
    \end{gather*}
    so that $\lvert P_Q\setminus P'_Q\rvert\leq 1$.

    Since $Q$ is a $(\widetilde{\delta},\widetilde{\epsilon},m)$-atom of $G$, by Lemma~\ref{lem:upcont}\ref{lem:upcont:atom}, we
    know that $P'_Q$ is a $(\widetilde{\delta},\epsilon_+,m)$-atom of $G$, where
    \begin{gather*}
      \epsilon_+
      \df
      (\widetilde{\epsilon}-1)\cdot\frac{\lvert Q\rvert}{\lvert P'_Q\rvert} + 1
      \leq
      (\widetilde{\epsilon}-1)\cdot\left(1 + \frac{\zeta-\widetilde{\epsilon}}{1-\zeta}\right)^{-1}
      + 1
      =
      \zeta,
    \end{gather*}
    hence $P'_Q$ is a $(\widetilde{\delta},\zeta,m)$-atom of $G$.

    On the other hand, note that by construction, we know that
    \begin{gather*}
      \lvert P'_Q\rvert
      \geq
      \lvert Q\rvert
      \geq
      (1-\gamma^\ell)^{K-1}\cdot\gamma^\ell\cdot n
      \geq
      (1-\gamma^\ell)^{K_{\nonequi}(\ell,0,c_{\ZZ},c_{\abs},\epsilon)-1}\cdot\gamma^\ell\cdot n.
    \end{gather*}

    Since $\lvert P_Q\setminus P'_Q\rvert\leq 1$, another application of Lemma~\ref{lem:upcont}\ref{lem:upcont:atom}, this time
    to $P'_Q$, says that $P_Q$ is a $(\widetilde{\delta},\epsilon_{++},m)$-atom of $G$, where
    \begin{align*}
      \epsilon_{++}
      & \df
      (\zeta-1)\cdot\frac{\lvert P'_Q\rvert}{\lvert P_Q\rvert} + 1
      \leq
      (\zeta-1)\cdot\left(1 - \frac{1}{1+\lvert P'_Q\rvert}\right) + 1
      \\
      & \leq
      (\zeta-1)\cdot\left(1 - \frac{1}{1 + (1-\gamma^\ell)^{K_{\nonequi}(\ell,c_{\UC},c_{\ZZ},c_{\abs},\epsilon)-1}\cdot\gamma^\ell\cdot n}\right)
      + 1
      \leq
      \epsilon,
    \end{align*}
    where the last inequality follows from Lemma~\ref{lem:nonequicalc}\ref{lem:nonequicalc:nnonequiimplies} (as $n\geq
    n_{\nonequi}(\ell,0,c_{\ZZ},c_{\abs},\epsilon)$), hence $P_Q$ is a $(\widetilde{\delta},\epsilon,m)$ atom of $G$.

    Finally, since $\widetilde{\delta}=(1+c_{\atom})\cdot\delta$, by Lemma~\ref{lem:atomconst} (and our definition of $m$ in
    line~\ref{alg:nonequi:atom:m} as exactly the ceiling of the right-hand side of the condition of the lemma), it follows that
    $P_Q$ is a $(\delta,\epsilon)$-atom of $G$, hence the returned partition $\cP$ is a $(\delta,\epsilon)$-atomic partition of
    $G$ into at most
    \begin{gather*}
      K_{\nonequi}(\ell,0,c_{\ZZ},c_{\abs},\epsilon)
      \leq
      (1+\widetilde{c}_{\ZZ})\cdot(1+\widetilde{c}_{\abs})\cdot\frac{\ln(1/\epsilon)}{\epsilon^\ell}
    \end{gather*}
    parts.

    \medskip

    We now analyze the time-complexity of Algorithm~\ref{alg:nonequi:atom}. We have already argued that the loop of
    line~\eqref{alg:nonequi:atom:While} executes at most
    \begin{gather*}
      K_{\nonequi}(\ell,0,c_{\ZZ},c_{\abs},\epsilon)
      \leq
      O\left(\frac{\log(1/\epsilon)}{\epsilon^\ell}\right)
    \end{gather*}
    times. The same is true for the loop of line~\ref{alg:nonequi:atom:For}.

    For the former loop of line~\ref{alg:nonequi:atom:While}, we claim that if we exclude the time it takes to make the rational
    estimations of line~\ref{alg:nonequi:atom:rational:si} (we will account for these later), then the time-complexity of each
    iteration is at most
    \begin{gather*}
      O\bigl(
      (\ell+m)\cdot n^{m+1}
      \bigr)
      \leq
      O\bigl((\ell+m)\cdot n^{m+1}\bigr)
    \end{gather*}
    due to the call to Algorithm~\ref{alg:ext:atom} from Proposition~\ref{prop:ext}\ref{prop:ext:atomalg}. So the total
    time-complexity of the loop of line~\ref{alg:nonequi:atom:While} (excluding the rational approximations) is at most
    \begin{gather}\label{eq:nonequi:atom:While:complexity}
      O\left((\ell+m)\cdot\frac{\log(1/\epsilon)}{\epsilon^\ell}\cdot n^{m+1}\right).
    \end{gather}

    Similarly, for the loop of line~\ref{alg:nonequi:atom:For}, if we exclude the time it takes to make the rational estimation
    of line~\ref{alg:nonequi:atom:rational:ws}, then it is clear that each iteration costs $O(n)$, amounting to a total
    time-complexity of
    \begin{gather}\label{eq:nonequi:atom:For:complexity}
      O\left(\frac{\log(1/\epsilon)}{\epsilon^\ell}\cdot n\right).
    \end{gather}

    \smallskip

    Let us now analyze the time-complexity of the calculations of parameters. We can give the following crude upper bounds for the
    time-complexities of the calculations of $\widetilde{\delta}$, $m$, $\widetilde{m}$ and $\zeta$, respectively:
    \begin{gather*}
      \begin{aligned}
        O\bigl(\len(c_{\atom})\cdot\len(\delta)\bigr),
        & &
        O\Bigl(
        \log(d+1)\cdot\len(c_{\atom})
        + \bigl(\log(d+1) + \len(c_{\atom})\bigr)\cdot\len(\delta)
        \Bigr),
      \end{aligned}
      \\
      \begin{aligned}
        O\Bigl(\bigl(\log(d+1) + \len(c_{\atom}) + \len(\delta)\bigr)\cdot\bigl(\len(c_{\atom})+\len(\delta)\bigr)\Bigr),
        & &
        O\bigl(\len(c_{\ZZ})\cdot\len(\epsilon)\bigr),
      \end{aligned}
    \end{gather*}
    which together are at most
    \begin{gather}\label{eq:nonequi:atom:param:complexity}
      O\Bigl(
      \len(c_{\atom})^2\cdot\len(\delta)^2\cdot\bigl(\log(d+1)\bigr)^2
      +
      \len(c_{\ZZ})\cdot\len(\epsilon)
      \Bigr).
    \end{gather}

    \smallskip

    We now analyze the time-complexity of the calculations of line~\ref{alg:nonequi:atom:rational:S}. For these rational
    approximations (and the ones that follow), we assume $c_{\abs}$ and $\zeta$ are presented as $c_{\abs}=a_{\abs}/b_{\abs}$
    and $\zeta = a_\zeta/b_\zeta$ for $a_{\abs},b_{\abs},a_\zeta,b_\zeta\in\NN_+$. Thus, we have
    \begin{gather*}
      \gamma \df \widetilde{\epsilon} \df (1-\zeta^{c_{\abs}})\cdot\zeta = \zeta - \zeta^{1+a_{\abs}/b_{\abs}},
    \end{gather*}
    to compute
    \begin{align*}
      S
      & \df
      \Floor{\frac{\zeta-\widetilde{\epsilon}}{1-\widetilde{\epsilon}}\cdot n}
      =
      \Floor{\frac{\zeta^{1 + c_{\abs}}\cdot n}{1-\zeta+\zeta^{1+c_{\abs}}}}
      \\
      & =
      \Floor{\frac{n}{\zeta^{-1}\cdot(\zeta^{-c_{\abs}}-1) + 1}}
      =
      \Floor{
        n\cdot
        \left(
        \frac{b_\zeta}{a_\zeta}\cdot
        \left(\left(\frac{b_\zeta}{a_\zeta}\right)^{a_{\abs}/b_{\abs}} - 1\right)
        + 1
        \right)^{-1}
      }.
    \end{align*}
    Let $c_{\den}$ be the expression under the last parentheses so that we are interested in $\floor{n/c_{\den}}$. Note that it
    suffices to compute the ceiling of $c_{\den}$ then make an integer division of $n$ by it, which will take time at most
    \begin{gather*}
      O\bigl(\log(n+1)\cdot\log(\ceil{c_{\den}}+1)\bigr)
      \leq
      O\bigl(\log(n+1)\cdot\len(c_{\den})\bigr)
      \leq
      O\bigl(\log(n+1)\cdot 2^{\len(c_{\abs})}\cdot\len(\zeta)\bigr)
    \end{gather*}
    In turn, to compute such a ceiling $\ceil{c_{\den}}$, it suffices to compute the rational $\zeta^{-1}=b_\zeta/a_\zeta$
    then determine the ceiling of the following:
    \begin{align*}
      \frac{b_\zeta^{1+a_{\abs}/b_{\abs}}}{a_\zeta^{a_{\abs}/b_{\abs}}} = \frac{b_\zeta}{\zeta^{c_{\abs}}},
    \end{align*}
    that is, we need to determine the unique $f\in\NN_+$ such that
    \begin{gather*}
      a_\zeta^{a_{\abs}/b_{\abs}}\cdot f
      <
      b_\zeta^{1+a_{\abs}/b_{\abs}}
      \leq
      a_\zeta^{a_{\abs}/b_{\abs}}\cdot (f+1),
    \end{gather*}
    which is equivalent to
    \begin{gather}\label{eq:nonequi:atom:rational:S:ceil}
      a_\zeta^{a_{\abs}}\cdot f^{b_{\abs}}
      <
      b_\zeta^{b_{\abs}+a_{\abs}}
      \leq
      a_\zeta^{a_{\abs}}\cdot (f+1)^{b_{\abs}}.
    \end{gather}

    Since we know that $f\in\{0,\ldots,b_\zeta\cdot\zeta^{-1}+1\}$, we can do a simple binary search for it (which takes at most
    $\log(b_\zeta\cdot\zeta^{-1}+2)\leq O(\len(\zeta))$ iterations). Using again that $f\leq b_\zeta\cdot\zeta^{-1}$ so
    $\len(f)\leq O(\len(\zeta))$, we can give the following crude upper bound for the time-complexity of computing the quantities
    in~\eqref{eq:nonequi:atom:rational:S:ceil}:
    \begin{align*}
      \MoveEqLeft
      \begin{multlined}[t]
        O\Bigl(
        \log(a_{\abs}+1)\cdot a_{\abs}^2\cdot\bigl(\log(a_\zeta+1)\bigr)^2
        +
        \log(b_{\abs}+1)\cdot b_{\abs}^2\cdot\len(\zeta)^2
        \\
        +
        a_{\abs}\cdot\log(a_\zeta+1)\cdot b_{\abs}\cdot\len(\zeta)
        +
        \log(b_{\abs}+a_{\abs})\cdot\bigl(b_{\abs}+a_{\abs}\bigr)^2\cdot\bigl(\log(b_\zeta+1)\bigr)^2
        \Bigr)
      \end{multlined}
      \\
      & \leq
      O\bigl(2^{2\cdot\len(c_{\abs})}\cdot\len(c_{\abs})\cdot\len(\zeta)^2\bigr)
      \\
      & \leq
      O\bigl(2^{3\cdot\len(c_{\abs})}\cdot\len(\zeta)^2\bigr)
    \end{align*}
    (Note that this completely dominates the time it takes to do the other computations after this to determine
    $\ceil{c_{\den}}$.) Since we need to make the binary search and then the integer division of $n$ by $\ceil{c_{\den}}$, the
    total time-complexity of line~\ref{alg:nonequi:atom:rational:S} is upper bounded by
    \begin{gather}\label{eq:nonequi:atom:rational:S:complexity}
      \begin{aligned}
        \MoveEqLeft
        O\bigl(
        2^{3\cdot\len(c_{\abs})}\cdot\len(\zeta)^3
        +
        \log(n+1)\cdot 2^{\len(c_{\abs})}\cdot\len(\zeta)
        \bigr)
        \\
        & \leq
        O\bigl(
        2^{3\cdot\len(c_{\abs})}\cdot\len(\zeta)^3\cdot\log(n+1)
        \bigr)
        \\
        & \leq
        O\bigl(
        2^{3\cdot\len(c_{\abs})}\cdot\bigl(\len(c_{\ZZ})+\len(\epsilon)\bigr)^3\cdot\log(n+1)
        \bigr).
      \end{aligned}
    \end{gather}

    \smallskip

    We now analyze the time-complexity of making the appropriate rational approximations of
    line~\ref{alg:nonequi:atom:rational:si}. We note that all the $s_i$ are always numbers in $[n]$, so we need to determine the
    time-complexity of computing the following for some $s\in[n]$:
    \begin{gather*}
      \floor{\zeta\cdot s - \zeta^{1+a_{\abs}/b_{\abs}}\cdot s}
      =
      \Floor{\frac{a_\zeta\cdot s}{b_\zeta} - \left(\frac{a_\zeta}{b_\zeta}\right)^{1+a_{\abs}/b_{\abs}}\cdot s}.
    \end{gather*}
    From the last expression, note that, it suffices to compute the rational $a_\zeta\cdot s/b_\zeta$ (which takes time at most
    $O(\len(\zeta)\cdot\log(n+1))$), then determine the ceiling of the following:
    \begin{gather*}
      \frac{a_\zeta^{1+a_{\abs}/b_{\abs}}\cdot s}{b_\zeta^{a_{\abs}/b_{\abs}}},
    \end{gather*}
    that is, we need to determine the unique $f\in\NN$ such that
    \begin{gather*}
      b_\zeta^{a_{\abs}/b_{\abs}}\cdot f
      <
      a_\zeta^{1+a_{\abs}/b_{\abs}}\cdot s
      \leq
      b_\zeta^{a_{\abs}/b_{\abs}}\cdot (f+1),
    \end{gather*}
    which in turn is equivalent to
    \begin{gather}\label{eq:nonequi:atom:rational:si:ceil}
      b_\zeta^{a_{\abs}}\cdot f^{b_{\abs}}
      <
      a_\zeta^{a_{\abs} + b_{\abs}}\cdot s^{b_{\abs}}
      \leq
      b_\zeta^{a_{\abs}}\cdot (f+1)^{b_{\abs}}.
    \end{gather}

    We know that such $f$ will be in $\{0,\ldots,b_\zeta\cdot n+1\}$, so we can do a simple binary search for it (which takes at
    most $\log_2(b_\zeta\cdot n+2)\leq\len(\zeta)+\log(n)$ iterations). Using $\len(s)\leq\log(n+1)$ and that $f\leq
    b_\zeta\cdot n+1$ implies $\len(f)\leq O(\len(\zeta) + \log(n))$, we can give the following crude upper bound for the
    time-complexity of computing the quantities in~\eqref{eq:nonequi:atom:rational:si:ceil}:
    \begin{align*}
      \MoveEqLeft
      \begin{multlined}[t]
        O\Bigl(
        \log(a_{\abs}+1)\cdot a_{\abs}^2\cdot\bigl(\log(b_\zeta+1)\bigr)^2
        +
        \log(b_{\abs}+1)\cdot b_{\abs}^2\cdot\bigl(\len(\zeta)+\log(n)\bigr)^2
        \\
        +
        a_{\abs}\cdot\log(b_{\abs}+1)\cdot b_{\abs}\cdot\bigl(\len(\zeta)+\log(n)\bigr)
        +
        \log(a_{\abs} + b_{\abs})\cdot(a_{\abs}+b_{\abs})^2\cdot\bigl(\log(a_\zeta+1)\bigr)^2
        \\
        +
        \log(b_{\abs}+1)\cdot b_{\abs}^2\cdot\bigl(\log(n+1)\bigr)^2
        +
        (a_{\abs}+b_{\abs})\cdot\log(a_\zeta+1)\cdot b_{\abs}\cdot\log(n+1)
        \Bigr)
      \end{multlined}
      \\
      & \leq
      O\bigl(
      2^{3\cdot\len(c_{\abs})}\cdot\len(\zeta)^2\cdot\bigl(\log(n+1)\bigr)^2
      \bigr).
    \end{align*}
    (Note that this completely dominates the time it takes after this to make the final calculations and determine
    $\floor{\zeta\cdot s - \zeta^{1+a_{\abs}/b_{\abs}}\cdot s}$.) Since we need to make the binary search, the total time to
    make such rational estimation for a single $s\in[n]$ is upper bounded by
    \begin{gather*}
      O\Bigl(
      \bigl(\len(\zeta) + \log(n)\bigr)
      \cdot
      2^{3\cdot\len(c_{\abs})}\cdot\len(\zeta)^2\cdot\bigl(\log(n+1)\bigr)^2
      \Bigr)
      \leq
      O\Bigl(
      2^{3\cdot\len(c_{\abs})}\cdot\len(\zeta)^3\cdot\bigl(\log(n+1)\bigr)^3
      \Bigr).
    \end{gather*}

    Finally, we make $\ell+1$ such estimations per iteration of the loop of line~\ref{alg:nonequi:atom:While} and since we have
    already established that the loop executes at most
    \begin{gather*}
      K_{\nonequi}(\ell,0,c_{\ZZ},c_{\abs},\epsilon)
      \leq
      O\left(\frac{\log(1/\epsilon)}{\epsilon^\ell}\right)
    \end{gather*}
    times, the total time-complexity of line~\ref{alg:nonequi:atom:rational:si}, accounting for every time it executes, is at
    most
    \begin{gather}\label{eq:nonequi:atom:rational:si:complexity}
      \begin{aligned}
        \MoveEqLeft
        O\left(
        \frac{\log(1/\epsilon)}{\epsilon^\ell}
        \cdot\ell
        \cdot
        2^{3\cdot\len(c_{\abs})}\cdot\len(\zeta)^3\cdot\bigl(\log(n+1)\bigr)^3
        \right)
        \\
        & \leq
        O\left(
        \frac{\log(1/\epsilon)}{\epsilon^\ell}
        \cdot\ell
        \cdot
        2^{3\cdot\len(c_{\abs})}\cdot\bigl(\len(c_{\ZZ})+\len(\epsilon)\bigr)^3\cdot\bigl(\log(n+1)\bigr)^3
        \right).
      \end{aligned}
    \end{gather}

    \smallskip

    We now analyze the time-complexity of making the appropriate rational approximations of
    line~\ref{alg:nonequi:atom:rational:ws}. We note that $s$ is always a number in $[n]$ and we need to determine the
    time-complexity of computing
    \begin{gather*}
      \widetilde{s}
      \df
      \Ceil{\frac{\zeta-\widetilde{\epsilon}}{1-\zeta}\cdot s}
      =
      \Ceil{\frac{\zeta^{1+c_{\abs}}\cdot s}{1-\zeta}}
      =
      \Ceil{\frac{a_\zeta^{a_{\abs}/b_{\abs}}\cdot s}{b_\zeta^{a_{\abs}/b_{\abs}}\cdot(b_\zeta - a_\zeta)}},
    \end{gather*}
    that is, we need to determine the unique $f\in\NN_+$ such that
    \begin{gather*}
      b_\zeta^{a_{\abs}/b_{\abs}}\cdot(b_\zeta-a_\zeta)\cdot f
      <
      a_\zeta^{a_{\abs}/b_{\abs}}\cdot s
      \leq
      b_\zeta^{a_{\abs}/b_{\abs}}\cdot(b_\zeta-a_\zeta)\cdot (f+1),
    \end{gather*}
    which is equivalent to
    \begin{gather}\label{eq:nonequi:atom:rational:ws:ceil}
      b_\zeta^{a_{\abs}}\cdot(b_\zeta-a_\zeta)^{b_{\abs}}\cdot f^{b_{\abs}}
      <
      a_\zeta^{a_{\abs}}\cdot s^{b_{\abs}}
      \leq
      b_\zeta^{a_{\abs}}\cdot(b_\zeta-a_\zeta)^{b_{\abs}}\cdot (f+1)^{b_{\abs}}.
    \end{gather}

    We now such $f$ will be in $\{0,\ldots,n\}$, so we can do a simple binary search for it (which takes at most $\log(n+1)$
    iterations). Using again that $f\leq n$ so $\len(f)\leq\log(n+1)$, we can give the following crude upper bound for the
    time-complexity of computing the quantities in~\eqref{eq:nonequi:atom:rational:ws:ceil}:
    \begin{align*}
      \MoveEqLeft
      \begin{multlined}[t]
        O\Bigl(
        \log(a_{\abs}+1)\cdot a_{\abs}^2\cdot\bigl(\log(b_\zeta+1)\bigr)^2
        +
        \log(b_{\abs}+1)\cdot b_{\abs}^2\cdot\len(\zeta)^2
        \\
        +
        \log(b_{\abs}+1)\cdot b_{\abs}^2\cdot\bigl(\log(n+1)\bigr)^2
        +
        a_{\abs}\cdot\log(b_\zeta+1)\cdot b_{\abs}\cdot\len(\zeta)\cdot b_{\abs}\cdot\log(n+1)
        \\
        +
        \log(a_{\abs}+1)\cdot a_{\abs}^2\cdot\bigl(\log(a_\zeta+1)\bigr)^2
        +
        \log(b_{\abs}+1)\cdot b_{\abs}^2\cdot\bigl(\log(n+1)\bigr)^2
        \\
        +
        a_{\abs}\cdot\log(a_\zeta+1)\cdot b_{\abs}\cdot\log(n+1)
        \Bigr)
      \end{multlined}
      \\
      & \leq
      O\Bigl(
      2^{3\cdot\len(c_{\abs})}\cdot\len(\zeta)^2\cdot\bigl(\log(n+1)\bigr)^2
      \Bigr).
    \end{align*}
    Since we need to make the binary search and the loop of line~\ref{alg:nonequi:atom:For} executes at most
    \begin{gather*}
      K_{\nonequi}(\ell,0,c_{\ZZ},c_{\abs},\epsilon)
      \leq
      O\left(\frac{\log(1/\epsilon)}{\epsilon^\ell}\right),
    \end{gather*}
    the total time-complexity of line~\ref{alg:nonequi:atom:rational:ws}, accounting for every time it executes, is at most
    \begin{gather}\label{eq:nonequi:atom:rational:ws:complexity}
      \begin{aligned}
        \MoveEqLeft
        O\left(
        \frac{\log(1/\epsilon)}{\epsilon^\ell}
        \cdot
        2^{3\cdot\len(c_{\abs})}\cdot\len(\zeta)^2\cdot\bigl(\log(n+1)\bigr)^3
        \right)
        \\
        & \leq
        O\left(
        \frac{\log(1/\epsilon)}{\epsilon^\ell}
        \cdot
        2^{3\cdot\len(c_{\abs})}\cdot\bigl(\len(c_{\ZZ})+\len(\epsilon)\bigr)^2\cdot\bigl(\log(n+1)\bigr)^3
        \right)
      \end{aligned}
    \end{gather}
    
    \smallskip

    Putting everything together (i.e., equations~\eqref{eq:nonequi:atom:While:complexity},
    \eqref{eq:nonequi:atom:For:complexity}, \eqref{eq:nonequi:atom:param:complexity},
    \eqref{eq:nonequi:atom:rational:S:complexity} and~\eqref{eq:nonequi:atom:rational:si:complexity}), we bound the
    time-complexity of Algorithm~\ref{alg:nonequi:atom} by:
    \begin{align*}
      \MoveEqLeft
      \begin{multlined}[t]
        O\bigggl(
        (\ell+m)\cdot\frac{\log(1/\epsilon)}{\epsilon^\ell}\cdot n^{m+1}
        +
        \frac{\log(1/\epsilon)}{\epsilon^\ell}\cdot n
        +
        \len(c_{\atom})^2\cdot\len(\delta)^2\cdot\bigl(\log(d+1)\bigr)^2
        \\
        +
        \len(c_{\ZZ})\cdot\len(\epsilon)
        +
        2^{3\cdot\len(c_{\abs})}\cdot\bigl(\len(c_{\ZZ})+\len(\epsilon)\bigr)^3\cdot\log(n+1)
        \\
        +
        \frac{\log(1/\epsilon)}{\epsilon^\ell}
        \cdot\ell
        \cdot
        2^{3\cdot\len(c_{\abs})}\cdot\bigl(\len(c_{\ZZ})+\len(\epsilon)\bigr)^3\cdot\bigl(\log(n+1)\bigr)^3
        \\
        +
        \frac{\log(1/\epsilon)}{\epsilon^\ell}
        \cdot
        2^{3\cdot\len(c_{\abs})}\cdot\bigl(\len(c_{\ZZ})+\len(\epsilon)\bigr)^2\cdot\bigl(\log(n+1)\bigr)^3
        \bigggr)
      \end{multlined}
      \\
      & \leq
      \begin{multlined}[t]
        O\bigggl(
        (\ell+m)\cdot\frac{\log(1/\epsilon)}{\epsilon^\ell}\cdot n^{m+1}
        +
        \len(c_{\atom})^2\cdot\len(\delta)^2\cdot\bigl(\log(d+1)\bigr)^2
        \\
        +
        \frac{\log(1/\epsilon)}{\epsilon^\ell}
        \cdot\ell
        \cdot
        2^{3\cdot\len(c_{\abs})}\cdot\bigl(\len(c_{\ZZ})+\len(\epsilon)\bigr)^3\cdot\bigl(\log(n+1)\bigr)^3
        \bigggr)
      \end{multlined}
      \\
      & \leq
      \begin{multlined}[t]
        \frac{\log(1/\epsilon)}{\epsilon^\ell}\cdot n^{O(d/(c_{\atom}^2\cdot\delta^2))}
        +
        O\bigggl(
        \len(c_{\atom})^2\cdot\len(\delta)^2\cdot\bigl(\log(d+1)\bigr)^2
        \\
        +
        \frac{\log(1/\epsilon)}{\epsilon^\ell}
        \cdot\ell
        \cdot
        2^{3\cdot\len(c_{\abs})}\cdot\bigl(\len(c_{\ZZ})+\len(\epsilon)\bigr)^3\cdot\bigl(\log(n+1)\bigr)^3
        \bigggr),
      \end{multlined}
    \end{align*}
    where the last inequality follows since
    \begin{gather*}
      m
      =
      \Ceil{C\cdot\frac{d}{c_{\atom}^2\cdot\delta^2}}
      \leq
      O\left(\frac{d}{c_{\atom}^2\cdot\delta^2}\right).
    \end{gather*}
  \item[Item~\ref{thm:nonequi:gooddet}.] The deterministic algorithm to produce $\epsilon$-good sets,
    Algorithm~\ref{alg:nonequi:gooddet}, is the same as Algorithm~\ref{alg:nonequi:atom}, except that we use $m\df 1$ (and an
    arbitrary $\delta\df 1/2^{\ell+2}$) along with the fact that $\epsilon$-goodness is equivalent to
    $(\delta,\epsilon,1)$-atomicity (provided $\delta<1$).

    \begin{algorithm}[htbp]
      \caption{Deterministic model algorithm that returns an $\epsilon$-good partition of $G$ into at most
        $K_{\nonequi}(\ell,0,c_{\ZZ},c_{\abs},\epsilon)$ parts. The time-complexity is:
        \begin{gather*}
          O\left(
          \ell\cdot\frac{\log(1/\epsilon)}{\epsilon^\ell}\cdot n^2
          +
          \frac{\log(1/\epsilon)}{\epsilon^\ell}
          \cdot\ell
          \cdot
          2^{3\cdot\len(c_{\abs})}\cdot\bigl(\len(c_{\ZZ})+\len(\epsilon)\bigr)^3\cdot\bigl(\log(n+1)\bigr)^3
          \right).
        \end{gather*}
      }
      \label{alg:nonequi:gooddet}
      \DontPrintSemicolon
      \KwIn{Numbers $n,\ell\in\NN_+$ with $\ell\leq\log_2(n)$, $c_{\ZZ},c_{\abs}\in(0,1)\cap\QQ$,
        $\epsilon\in(0,\epsilon_{\nonequi}(\ell,c_{\ZZ},c_{\abs}))\cap\QQ$, and a graph $G$ with $\lvert G\rvert=n$ and
        $\Lit(G)\leq\ell$. We further assume $n\geq n_{\nonequi}(\ell,0,c_{\ZZ},c_{\abs},\epsilon)$.}
      \KwOut{An $\epsilon$-good partition of $G$ into at most $K_{\nonequi}(\ell,0,c_{\ZZ},c_{\abs},\epsilon)$ parts.}
      $m\assign 1$\;
      $\widetilde{m}\assign 0$\tcp*[r]{$\widetilde{m} = \floor{\widetilde{\delta}\cdot m}$ for $\widetilde{\delta} < 1$}
      $c_{\UC}\assign 0$\;
      $\zeta\assign (1-c_{\ZZ})\cdot\epsilon$\;
      Compute a sufficiently precise rational approximation of $\widetilde{\epsilon}\df (1 - \zeta^{c_{\abs}})\cdot\zeta$ to
      compute $S\df\floor{(\zeta-\widetilde{\epsilon})\cdot n/(1-\widetilde{\epsilon})}$\;
      \label{alg:nonequi:gooddet:rational:S}
      $\cQ\assign\varnothing$\;
      $R\assign V(G)$\;
      $r\assign n$\;
      \While{$r > S$}{%
        \label{alg:nonequi:gooddet:While}
        $s_0\assign r$\;
        Compute sufficiently precise rational approximations of $\widetilde{\epsilon}\df (1 - \zeta^{c_{\abs}})\cdot\zeta$ and
        $\gamma\df (1 - c_{\UC})\cdot\widetilde{\epsilon}$ to compute $s_i\df\floor{\gamma\cdot s_{i-1}} + 1$ for every
        $i\in[\ell+1]$.\;
        \label{alg:nonequi:gooddet:rational:si}
        Run Algorithm~\ref{alg:ext:atom} with parameters:
        \begin{gather*}
          (n,m,\ell,s_0,\ldots,s_\ell,s_{\ell+1},\widetilde{m},G,U)
          \df
          (n,m,\ell,s_0,\ldots,s_\ell,s_{\ell+1},\widetilde{m},G,R).
        \end{gather*}
        \;\label{alg:nonequi:gooddet:call}
        Let $(W,s)$ be the pair returned by the algorithm.\;
        $\cQ\assign\cP\cup\{(W,s)\}$\;
        $R\assign R\setminus W$\;
        $r\assign r - s$\;
      }
      $\cP\assign\varnothing$\;
      \For{$(W,s)\in\cQ$}{%
        \label{alg:nonequi:gooddet:For}
        Compute a sufficently precise approximation of $\widetilde{\epsilon}\df (1-\zeta^{c_{\abs}})\cdot\zeta$ to compute
        $\widetilde{s}\df\ceil{(\zeta-\widetilde{\epsilon})\cdot s/(1-\zeta)}$.\;
        \label{alg:nonequi:gooddet:rational:ws}
        Let $A\subseteq R$ be a set of size as large as possible with $\lvert A\rvert\leq\widetilde{s}$.\;
        $\cP\assign\cP\cup\{W\cup A\}$\;
        $R\assign R\setminus A$\;
      }
      \Return{$\cP$}
    \end{algorithm}

    The complexity analysis of Algorithm~\ref{alg:nonequi:gooddet} is also analogous to that of Algorithm~\ref{alg:nonequi:atom}
    and yields the time-complexity bound:
    \begin{align*}
      \MoveEqLeft
      O\left(
      (\ell+m)\cdot\frac{\log(1/\epsilon)}{\epsilon^\ell}\cdot n^{m+1}
      +
      \frac{\log(1/\epsilon)}{\epsilon^\ell}
      \cdot\ell
      \cdot
      2^{3\cdot\len(c_{\abs})}\cdot\bigl(\len(c_{\ZZ})+\len(\epsilon)\bigr)^3\cdot\bigl(\log(n+1)\bigr)^3
      \right)
      \\
      & =
      O\left(
      \ell\cdot\frac{\log(1/\epsilon)}{\epsilon^\ell}\cdot n^2
      +
      \frac{\log(1/\epsilon)}{\epsilon^\ell}
      \cdot\ell
      \cdot
      2^{3\cdot\len(c_{\abs})}\cdot\bigl(\len(c_{\ZZ})+\len(\epsilon)\bigr)^3\cdot\bigl(\log(n+1)\bigr)^3
      \right).
    \end{align*}
  \item[Item~\ref{thm:nonequi:good}.] As mentioned before, the proof of this item is very similar to that of
    item~\ref{thm:nonequi:atom}, except that we will replace Algorithm~\ref{alg:ext:atom} with Algorithm~\ref{alg:ext:good} and
    we will have to control the probability of failure.

    \begin{algorithm}[htbp]
      \caption{Random-query model algorithm that, with probability at least $1-\rho_{\UC}$, returns an $\epsilon$-good partition
        of $G$ into at most $K_{\nonequi}(\ell,c_{\UC},c_{\ZZ},c_{\abs},\epsilon)$ parts. The time-complexity is:
        \begin{multline*}
          O\bigggl(
          \frac{2^\ell\cdot\ell^2\cdot(\log(1/\epsilon))^2}{c_{\UC}^2\cdot(1-c_{\ZZ})^2\cdot\epsilon^{\ell+2}}
          \cdot n\cdot\log\left(\frac{n}{\rho_{\UC}}\right)
          \\
          +
          \ell^2\cdot\log(\ell+1)
          \cdot 2^{3\cdot\len(c_{\abs})}\cdot\bigl(\len(c_{\ZZ})+\len(\epsilon)\bigr)^3\cdot\len(c_{\UC})^6\cdot\len(\rho_{\UC})
          \cdot\frac{\log(1/\epsilon)}{\epsilon^\ell}\cdot\bigl(\log(n+1)\bigr)^3
          \bigggr).
        \end{multline*}
      }
      \label{alg:nonequi:good}
      \DontPrintSemicolon
      \KwIn{Numbers $n,d,\ell\in\NN_+$ with $d\leq\ell\leq\log_2(n)$, $c_{\UC},c_{\ZZ},c_{\abs}\in(0,1)\cap\QQ$,
        $\epsilon\in(0,\epsilon_{\nonequi}(\ell,c_{\ZZ},c_{\abs}))\cap\QQ$, and a random-query oracle $\cO_G$ for a graph $G$
        with $\lvert G\rvert=n$, $\Lit(G)\leq\ell$ and $\VC(G)\leq d$. We further assume $n\geq
        n_{\nonequi}(\ell,c_{\UC},c_{\ZZ},c_{\abs},\epsilon)$.}
      \KwOut{An $\epsilon$-good partition of $G$ into at most $K_{\nonequi}(\ell,c_{\UC},c_{\ZZ},c_{\abs},\epsilon)$ parts. With
        probability at most $\rho_{\UC}$, may instead return an incorrect partition or ``\Failure''.}
      $\widetilde{c}_{\UC}\assign(1+c_{\UC})^{-\ell}-1$\;
      $\widetilde{c}_{\ZZ}\assign(1+c_{\ZZ})^{-\ell}-1$\;
      $\widetilde{c}_{\abs}\assign 2\cdot c_{\abs}$\;
      $K\assign\floor{(1+\widetilde{c}_{\UC})\cdot(1+\widetilde{c}_{\ZZ})\cdot(1+\widetilde{c}_{\abs})\cdot\epsilon^{-\ell-1}}$\;
      $\zeta\assign(1-c_{\ZZ})\cdot\epsilon$\;
      $\widetilde{\rho}_{\UC}\assign\rho_{\UC}/K$\;
      Compute a sufficiently precise rational approximation of $\widetilde{\epsilon}\df (1 - \zeta^{c_{\abs}})\cdot\zeta$ to
      compute:
      \begin{gather*}
        S \df \Floor{\frac{\zeta-\widetilde{\epsilon}}{1-\widetilde{\epsilon}}\cdot n},
        \\
        m
        \df
        \Ceil{\frac{C\cdot 4\cdot (d + \ell + 1 + \ceil{\log_2(n/\widetilde{\rho}_{\UC})})}{c_{\UC}^2\cdot\widetilde{\epsilon}^2}},
        \\
        \widetilde{m}\df\Floor{\left(1 - \frac{c_{\UC}}{2}\right)\cdot\widetilde{\epsilon}\cdot m},
      \end{gather*}
      where $C$ is the absolute constant of Theorem~\ref{thm:UC}.\;\label{alg:nonequi:good:rational:initial}
      $\cQ\assign\varnothing$\;
      $R\assign V(G)$\;
      $r\assign n$\;
      \uWhile{$r > S$}{%
        \label{alg:nonequi:good:While}
        $s_0\assign r$\;
        Compute sufficiently precise rational approximations of $\widetilde{\epsilon}\df (1 - \zeta^{c_{\abs}})\cdot\zeta$ and
        $\gamma\df (1 - c_{\UC})\cdot\widetilde{\epsilon}$ to compute $s_i\df\floor{\gamma\cdot s_{i-1}} + 1$ for every
        $i\in[\ell+1]$.\;
        \label{alg:nonequi:good:rational:si}
        \setcounter{algosplit}{\theAlgoLine}
      }
    \end{algorithm}

    \begin{algorithm}[htbp]
      \caption*{(continued).}
      \DontPrintSemicolon
      \setcounter{AlgoLine}{\thealgosplit}
      \let\oldnl\nl
      \let\nl\relax
      \vspace{-\baselineskip}\InvisibleBegin{
        \global\let\nl\oldnl 
        %
        Run Algorithm~\ref{alg:ext:good} with parameters:
        \begin{gather*}
          (n,m,\ell,s_0,\ldots,s_\ell,s_{\ell+1},\widetilde{m},\cO_G,U)
          \df
          (n,m,\ell,s_0,\ldots,s_\ell,s_{\ell+1},\widetilde{m},\cO_G,R).
        \end{gather*}
        \;\label{alg:nonequi:good:call}
        \lIf{The algorithm returned ``\Failure''}{\Return{\Failure}}
        Let $(W,s)$ be the pair returned by the algorithm.\;
        $\cQ\assign\cP\cup\{(W,s)\}$\;
        $R\assign R\setminus W$\;
        $r\assign r - s$\;
      }
      $\cP\assign\varnothing$\;
      \For{$(W,s)\in\cQ$}{%
        \label{alg:nonequi:good:For}
        Compute a sufficiently precise approximation of $\widetilde{\epsilon}\df (1-\zeta^{c_{\abs}})\cdot\zeta$ to compute
        $\widetilde{s}\df\ceil{(\zeta-\widetilde{\epsilon})\cdot s/(1-\zeta)}$.\;
        \label{alg:nonequi:good:rational:ws}
        Let $A\subseteq R$ be a set of size as large as possible with $\lvert A\rvert\leq\widetilde{s}$.\;
        $\cP\assign\cP\cup\{W\cup A\}$\;
        $R\assign R\setminus A$\;
      }
      \Return{$\cP$}
    \end{algorithm}

    The proof of correctness of Algorithm~\ref{alg:nonequi:good} is completely analogous to that of
    Algorithm~\ref{alg:nonequi:atom}, provided all calls to Algorithm~\ref{alg:ext:good} in line~\ref{alg:nonequi:good:call} are
    successful. We also note that our particular calculation of the parameter $m$ of Algorithm~\ref{alg:nonequi:good} is
    \begin{gather*}
      m
      \df
      \Ceil{\frac{C\cdot 4\cdot (d + \ell + 1 + \ceil{\log_2(n/\widetilde{\rho}_{\UC})})}{c_{\UC}^2\cdot\widetilde{\epsilon}^2}}
      \geq
      \frac{C\cdot 4\cdot (d + \ln((2^{\ell+1}-1)\cdot n/\widetilde{\rho}_{\UC}))}{(c_{\UC}^2\cdot\widetilde{\epsilon}^2)},
    \end{gather*}
    so we can indeed apply Proposition~\ref{prop:ext}\ref{prop:ext:goodalg} and use Algorithm~\ref{alg:ext:good} (correctness
    would still work if we had defined $m$ as the last expression, but we opt to define a slightly larger $m$ to simplify our
    analysis of the time-complexity of computing it).

    Thus, we need only to analyze the probability of such an event. We know that there will be at most
    \begin{gather*}
      K_{\nonequi}(\ell,c_{\UC},c_{\ZZ},c_{\abs},\epsilon)
      \leq
      (1+c_{\UC})\cdot(1+c_{\ZZ})\cdot(1+c_{\abs})\cdot\frac{\ln(1/\epsilon)}{\epsilon^\ell}
    \end{gather*}
    calls to Algorithm~\ref{alg:ext:good}, each with success probability at least $1-\widetilde{\rho}_{\UC}$ (by
    Proposition~\ref{prop:ext}\ref{prop:ext:goodalg}). Note that Algorithm~\ref{alg:nonequi:good} computes
    \begin{gather*}
      K
      \df
      \floor{(1+c_{\UC})\cdot(1+c_{\ZZ})\cdot(1+c_{\abs})\cdot\epsilon^{-\ell-1}}
      \geq
      \floor{K_{\nonequi}(\ell,c_{\UC},c_{\ZZ},c_{\abs},\epsilon)},
    \end{gather*}
    which itself is responsible for the computation of $\widetilde{\rho}_{\UC}$ and $m$. Thus, the final probability of success
    is at least\footnote{There are two sources of non-optimality in our computation of $K$ in Algorithm~\ref{alg:ext:good}: The
    first is that we replace a $\ln(1/\epsilon)$ by $\epsilon^{-1}$, this is to avoid having to deal with the time-complexity of
    approximating $\ln$ to enough precision. The second is that our calculation of $\widetilde{\rho}_{\UC}$ in terms of
    $\rho_{\UC}$ is toward an application of a union bound of the probabilities of failures. However, since each run of
    Algorithm~\ref{alg:ext:good} has the same conditional probability of failure given that it is actually reached, the actual
    probability of success is at least $(1-\widetilde{\rho}_{\UC})^K$ rather than just at least $1 - \widetilde{\rho}_{\UC}/K$;
    the reason why we do not use the former bound for the computation of $\widetilde{\rho}_{\UC}$ in terms of $\rho_{\UC}$ is
    because it would yield the formula $\widetilde{\rho}_{\UC} = 1 - (1-\rho_{\UC})^{1/K}$ and then we would have to make
    another sufficiently good rational approximation to a potentially irrational value.}
    \begin{gather*}
      (1-\widetilde{\rho}_{\UC})^K \geq 1 - \frac{\widetilde{\rho}_{\UC}}{K} = 1 - \rho_{\UC}.
    \end{gather*}

    \medskip

    The complexity analysis of Algorithm~\ref{alg:nonequi:good} is mostly analogous to that of Algorithm~\ref{alg:nonequi:atom},
    with the major differences being in the fact that it calls Algorithm~\ref{alg:ext:good} (instead of
    Algorithm~\ref{alg:ext:atom}) and in the time-complexity of the parameter computations.

    The loops of lines~\ref{alg:nonequi:good:While} and~\ref{alg:nonequi:good:For} run at most
    \begin{gather*}
      K_{\nonequi}(\ell,c_{\UC},c_{\ZZ},c_{\abs},\epsilon)
      \leq
      O\left(\frac{\log(1/\epsilon)}{\epsilon^\ell}\right)
    \end{gather*}
    times. The time-complexity of each iteration of the loop of line~\ref{alg:nonequi:good:While}, excluding the rational
    approximations of line~\ref{alg:nonequi:good:rational:si} (which we will account for later), is at most
    \begin{align*}
      O\bigl(
      2^\ell\cdot n\cdot m
      \bigr)
      & \leq
      O\left(
      2^\ell\cdot n\cdot
      \frac{\ell\cdot\ln(n/\widetilde{\rho}_{\UC})}{c_{\UC}^2\cdot\widetilde{\epsilon}^2}
      \right)
      \\
      & \leq
      O\left(
      \frac{2^\ell\cdot\ell^2\cdot\log(1/\epsilon)}{c_{\UC}^2\cdot(1-c_{\ZZ})^2\cdot\epsilon^2}
      \cdot n\cdot\log\left(\frac{n}{\rho_{\UC}}\right)
      \right),
    \end{align*}
    due to the call to Algorithm~\ref{alg:ext:good} from Proposition~\ref{prop:ext}\ref{prop:ext:goodalg}. So the total
    time-complexity of the loop of line~\ref{alg:nonequi:good:While} (excluding the rational approximations) is at most
    \begin{gather}\label{eq:nonequi:good:While:complexity}
      O\left(
      \frac{2^\ell\cdot\ell^2\cdot(\log(1/\epsilon))^2}{c_{\UC}^2\cdot(1-c_{\ZZ})^2\cdot\epsilon^{\ell+2}}
      \cdot n\cdot\log\left(\frac{n}{\rho_{\UC}}\right)
      \right).
    \end{gather}

    Similarly, for the loop of line~\ref{alg:nonequi:good:For}, excluding the time it takes to make the rational approximations
    of line~\ref{alg:nonequi:good:rational:ws}, each iteration costs $O(n)$, so we get a total time-complexity of
    \begin{gather}\label{eq:nonequi:good:For:complexity}
      O\left(\frac{\log(1/\epsilon)}{\epsilon^\ell}\cdot n\right).
    \end{gather}

    \smallskip

    It remains to analyze the time-complexity of the calculation of parameters. The following are crude bounds for the
    time-complexities of the calculations of $\widetilde{c}_{\UC}$, $\widetilde{c}_{\ZZ}$, $\widetilde{c}_{\abs}$, $K$, $\zeta$ and
    $\widetilde{\rho}_{\UC}$, respectively:
    \begin{gather*}
      \begin{aligned}
        O\bigl(\ell\cdot\log(\ell+1)\cdot\len(c_{\UC})^2\bigr),
        & &
        O\bigl(\ell\cdot\log(\ell+1)\cdot\len(c_{\ZZ})^2\bigr),
        & &
        O\bigl(\len(c_{\abs})\bigr),
      \end{aligned}
      \\
      \begin{aligned}
        O\bigl(\ell^2\cdot\log(\ell+1)\len(c_{\UC})\cdot\len(c_{\ZZ})\cdot\len(c_{\abs})\cdot\len(\epsilon)^2\bigr),
        & &
        O\bigl(\len(c_{\ZZ})\cdot\len(\epsilon)\bigr),
      \end{aligned}
      \\
      O\Bigl(
      \len(\rho_{\UC})\cdot\ell\cdot\bigl(\len(c_{\UC}) + \len(c_{\ZZ}) + \len(c_{\abs})\bigr)\cdot\len(\epsilon)
      \Bigr),
    \end{gather*}
    which together are at most
    \begin{gather}\label{eq:nonequi:good:param:complexity}
      O\Bigl(
      \ell^2\cdot\log(\ell+1)\cdot\len(c_{\UC})\cdot\len(c_{\abs})\cdot\len(c_{\ZZ})\cdot\len(\epsilon)^2\cdot\len(\rho_{\UC})
      \Bigr).
    \end{gather}

    \smallskip

    We now analyze the time-complexity of the rational approximations of line~\ref{alg:nonequi:good:rational:initial}. For these
    rational approximations (and the ones that follow), we assume $c_{\abs}$, $c_{\UC}$ and $\zeta$ are presented as
    \begin{align*}
      c_{\abs} & = \frac{a_{\abs}}{b_{\abs}}, &
      c_{\UC} & = \frac{a_{\UC}}{b_{\UC}}, &
      \zeta & = \frac{a_\zeta}{b_\zeta},
    \end{align*}
    for $a_{\abs},b_{\abs},a_{\UC},b_{\UC},a_\zeta,b_\zeta\in\NN_+$. Since
    \begin{gather*}
      \widetilde{\epsilon}
      =
      (1 - \zeta^{c_{\abs}})\cdot\zeta
      =
      (\zeta - \zeta^{1+a_{\abs}/b_{\abs}}),
    \end{gather*}
    we want to compute the above to enough precision so that the following are determined:
    \begin{gather*}
      S \df \Floor{\frac{\zeta-\widetilde{\epsilon}}{1-\widetilde{\epsilon}}\cdot n},
      \\
      m
      \df
      \Ceil{\frac{C\cdot 4\cdot (d + \ell + 1 + \ceil{\log_2(n/\widetilde{\rho}_{\UC}))}}{c_{\UC}^2\cdot\widetilde{\epsilon}^2}},
      \\
      \widetilde{m}\df\Floor{\left(1 - \frac{c_{\UC}}{2}\right)\cdot\widetilde{\epsilon}\cdot m},
    \end{gather*}
    where $C$ is the absolute constant of Theorem~\ref{thm:UC}.

    The analysis of the computation of $S$ is the same as that of Algorithm~\ref{alg:nonequi:atom}, hence yields the
    time-complexity as in~\eqref{eq:nonequi:atom:rational:S:complexity}:
    \begin{gather}\label{eq:nonequi:good:rational:S:complexity}
      \begin{aligned}
        \MoveEqLeft
        O\bigl(
        2^{3\cdot\len(c_{\abs})}\cdot\len(\zeta)^3
        +
        \log(n+1)\cdot 2^{\len(c_{\abs})}\cdot\len(\zeta)
        \bigr)
        \\
        & \leq
        O\bigl(
        2^{3\cdot\len(c_{\abs})}\cdot\len(\zeta)^3\cdot\log(n+1)
        \bigr)
        \\
        & \leq
        O\bigl(
        2^{3\cdot\len(c_{\abs})}\cdot\bigl(\len(c_{\ZZ})+\len(\epsilon)\bigr)^3\cdot\log(n+1)
        \bigr).
      \end{aligned}
    \end{gather}

    \smallskip

    Let us now analyze the computation of $m$. For this, we compute the numerator $c_{\num}\in\NN_+$ of the expression under the
    ceiling for $m$. Recalling that $d\leq\ell$, this takes time at most
    \begin{gather}\label{eq:nonequi:good:rational:num}
      O\bigl(\log(\ell) + \log(n) + \len(\widetilde{\rho}_{\UC})\bigr)
      \leq
      O\bigl(
      \log(n) + \len(\rho_{\UC}) + \ell\cdot\len(\epsilon)
      \bigr).
    \end{gather}
    Then we need to determine
    \begin{align*}
      m
      & =
      \ceil{c_{\num}\cdot c_{\UC}^2\cdot\widetilde{\epsilon}^2}
      =
      \Ceil{
        \frac{
          c_{\num}
        }{
          (a_{\UC}/b_{\UC})^2\cdot(1-(a_{\UC}/b_{\UC}))^2
          \cdot(a_\zeta/b_\zeta - (a_\zeta/b_\zeta)^{1+a_{\abs}/b_{\abs}})
        }
      }
      \\
      & =
      \Ceil{
        \frac{
          c_{\num}\cdot b_{\UC}^4
        }{
          a_{\UC}^2\cdot(b_{\UC}-a_{\UC})^2
          \cdot(a_\zeta/b_\zeta - (a_\zeta/b_\zeta)^{1+a_{\abs}/b_{\abs}})
        }
      }
    \end{align*}
    It then suffices to determine the floor of the denominator above and in turn, it suffices to determine the floor of the
    following:
    \begin{gather*}
      \frac{a_{\UC}^2\cdot(b_{\UC}-a_{\UC})^2\cdot a_\zeta^{1+a_{\abs}/b_{\abs}}}{b_\zeta^{a_{\abs}/b_{\abs}}},
    \end{gather*}
    that is, we need to find the unique $f\in\NN$ such that
    \begin{gather*}
      \frac{f}{a_{\UC}^2\cdot(b_{\UC}-a_{\UC})^2}
      \leq
      \frac{a_\zeta^{1+a_{\abs}/b_{\abs}}}{b_\zeta^{a_{\abs}/b_{\abs}}}
      <
      \frac{f+1}{a_{\UC}^2\cdot(b_{\UC}-a_{\UC})^2},
    \end{gather*}
    which is equivalent to
    \begin{gather}\label{eq:nonequi:good:rational:m:floor}
      \left(\frac{f}{a_{\UC}^2\cdot(b_{\UC}-a_{\UC})^2}\right)^{b_{\abs}}
      \leq
      \frac{a_\zeta^{b_{\abs}+a_{\abs}}}{b_\zeta^{a_{\abs}}}
      <
      \left(\frac{f+1}{a_{\UC}^2\cdot(b_{\UC}-a_{\UC})^2}\right)^{b_{\abs}}.
    \end{gather}
    Note that such $f$ is guaranteed to be a non-negative integer with
    \begin{gather*}
      f
      \leq
      a_{\UC}^2\cdot b_{\UC}^2\cdot a_\zeta
    \end{gather*}
    so we can do a simple binary search for it. Note also that $\log(f+1)\leq O(\len(c_{\UC})+\len(\zeta))$.

    Now we can give the following crude bound for the time-complexity of computing the quantities
    in~\eqref{eq:nonequi:good:rational:m:floor}:
    \begin{align*}
      \MoveEqLeft
      \begin{multlined}[t]
        O\bigl(
        b_{\abs}^2\cdot\log(b_{\abs}+1)\cdot
        \bigl(\len(c_{\UC})^4 + \len(\zeta)\bigr)^2
        \cdot\bigl(\log(a_{\UC}+1)\bigr)^4\cdot\bigl(\log(b_{\UC}+1)\bigr)^4
        \\
        +
        (b_{\abs}+a_{\abs})^2\cdot\log(b_{\abs}+a_{\abs})\cdot\bigl(\log(a_\zeta+1)\bigr)^2
        +
        a_{\abs}^2\cdot\log(a_{\abs}+1)\cdot\bigl(\log(b_\zeta+1)\bigr)^2
        \bigr)
      \end{multlined}
      \\
      & \leq
      O\bigl(
      2^{3\cdot\len(c_{\abs})}
      \cdot\len(c_{\UC})^5\cdot\len(\zeta)^2
      \bigr).
    \end{align*}
    Since we need to make the binary search, compute $c_{\num}\cdot b_{\UC}^4$ and make an integer division of it by the floor
    $f$, the total time to make such rational estimation is upper bounded by
    \begin{gather}\label{eq:nonequi:good:rational:m:complexity}
      \begin{aligned}
        \MoveEqLeft
        \begin{multlined}[t]
          O\Bigl(
          2^{3\cdot\len(c_{\abs})}
          \cdot\len(c_{\UC})^6\cdot\len(\zeta)^3
          \\
          +
          \bigl(\log(n) + \len(\rho_{\UC}) + \ell\cdot\len(\epsilon)\bigr)\cdot\len(c_{\UC})
          \cdot\len(c_{\UC})\cdot\len(\zeta)
          \Bigr)
        \end{multlined}
        \\
        & \leq
        O\Bigl(
        2^{3\cdot\len(c_{\abs})}
        \cdot\len(c_{\UC})^6\cdot\bigl(\len(c_{\ZZ})+\len(\epsilon)\bigr)^3
        \cdot\len(\rho_{\UC})\cdot\ell\cdot\log(n+1)
        \Bigr).
      \end{aligned}
    \end{gather}

    \smallskip

    We now analyze the rational estimation of
    \begin{align*}
      \widetilde{m}
      & \df
      \Floor{\left(1 - \frac{c_{\UC}}{2}\right)\cdot\widetilde{\epsilon}\cdot m}
      \\
      & =
      \Floor{
        \left(1 - \frac{a_{\UC}}{2\cdot b_{\UC}}\right)
        \cdot\left(\frac{a_\zeta}{b_\zeta} - \left(\frac{a_\zeta}{b_\zeta}\right)^{1+a_{\abs}/b_{\abs}}\right)
      }.
      \\
      & =
      \Floor{
        \frac{2\cdot b_{\UC} + a_{\UC}}{2\cdot b_{\UC}}
        \cdot\left(\frac{a_\zeta}{b_\zeta} - \left(\frac{a_\zeta}{b_\zeta}\right)^{1+a_{\abs}/b_{\abs}}\right)
      }.
    \end{align*}
    From the last expression, it suffices to compute the rational $(2\cdot b_{\UC} + a_{\UC})\cdot a_\zeta/(2\cdot b_{\UC}\cdot
    b_\zeta)$ (which takes time at most $O(\len(c_{\UC})^2\cdot\len(\zeta))$) and determine the ceiling of
    \begin{gather*}
      \frac{(2\cdot b_{\UC} + a_{\UC})\cdot a_\zeta^{1+a_{\abs}/b_{\abs}}}{2\cdot b_{\UC}\cdot b_\zeta^{a_{\abs}/b_{\abs}}},
    \end{gather*}
    that is, we need to determine the unique $f\in\NN$ such that
    \begin{gather*}
      2\cdot b_{\UC}\cdot b_\zeta^{a_{\abs}/b_{\abs}}\cdot f
      <
      (2\cdot b_{\UC} + a_{\UC})\cdot a_\zeta^{1+a_{\abs}/b_{\abs}}
      \leq
      2\cdot b_{\UC}\cdot b_\zeta^{a_{\abs}/b_{\abs}}\cdot (f+1),
    \end{gather*}
    which in turn is equivalent to
    \begin{gather}\label{eq:nonequi:good:rational:wm:ceil}
      2^{b_{\abs}}\cdot b_{\UC}^{b_{\abs}}\cdot b_\zeta^{a_{\abs}}\cdot f^{b_{\abs}}
      <
      (2\cdot b_{\UC} + a_{\UC})\cdot a_\zeta^{b_{\abs}+a_{\abs}}
      \leq
      2^{b_{\abs}}\cdot b_{\UC}^{b_{\abs}}\cdot b_\zeta^{a_{\abs}}\cdot (f+1)^{b_{\abs}}.
    \end{gather}

    We know that such $f$ will be in $\{0,\ldots,2\cdot a_\zeta+1\}$, so we can do a simple binary search for it (which takes at
    most $\log_2(2\cdot a_\zeta+2)\leq O(\len(\zeta))$ iterations). Using $f\leq 2\cdot a_\zeta+1$ again, we can give the
    following crude bound for the time-complexity of computing the quantities in~\eqref{eq:nonequi:good:rational:wm:ceil}:
    \begin{align*}
      \MoveEqLeft
      \begin{multlined}[t]
        O\Bigl(
        \log(b_{\abs}+1)\cdot b_{\abs}^2\cdot\bigl(\log(b_{\UC}+1)\bigr)^2
        +
        \log(a_{\abs}+1)\cdot a_{\abs}^2\cdot\bigl(\log(b_\zeta+1)\bigr)^2
        \\
        +
        \log(b_{\abs}+1)\cdot b_{\abs}^2\cdot\len(\zeta)^2
        +
        b_{\abs}\cdot\log(b_{\UC}+1)\cdot a_{\abs}\cdot\log(b_\zeta+1)\cdot b_{\abs}\cdot\len(\zeta)
        \\
        +
        \log(b_{\abs}+a_{\abs})\cdot(b_{\abs}+a_{\abs})^2\cdot\log(a_\zeta)
        +
        \log(2\cdot b_{\UC} + a_{\UC})\cdot(b_{\abs}+a_{\abs})\cdot\log(a_\zeta)
        \Bigr)
      \end{multlined}
      \\
      & \leq
      O\Bigl(
      2^{3\cdot\len(c_{\abs})}\cdot\len(c_{\UC})^2\cdot\len(\zeta)^2
      \Bigr).
    \end{align*}
    (Note that this completely dominates the time it takes after this to make the final calculations and determine
    $\widetilde{m}$.) Since we need to make the binary search, the total time-complexity it takes to compute $\widetilde{m}$ is
    at most
    \begin{gather}\label{eq:nonequi:good:rational:wm:complexity}
      O\Bigl(
      2^{3\cdot\len(c_{\abs})}\cdot\len(c_{\UC})^2\cdot\len(\zeta)^3
      \Bigr)
      \leq
      O\Bigl(
      2^{3\cdot\len(c_{\abs})}\cdot\len(c_{\UC})^2\cdot\bigl(\len(c_{\ZZ})+\len(\epsilon)\bigr)^3
      \Bigr).
    \end{gather}

    \smallskip

    We now analyze the rational estimations of line~\ref{alg:nonequi:good:rational:si}. This is very similar to the analysis
    made for Algorithm~\ref{alg:nonequi:atom}, except that now our $c_{\UC}$ is not $0$, so our computation of $\gamma$ is
    slightly more expensive. Since all $s_i$ are always numbers in $[n]$, we need to determine the time-complexity of computing
    the following for some $s\in[n]$:
    \begin{gather*}
      \floor{\gamma\cdot s}
      =
      \Floor{
        \left(1-\frac{a_{\UC}}{b_{\UC}}\right)
        \cdot\left(\frac{a_\zeta\cdot s}{b_\zeta} - \left(\frac{a_\zeta}{b_\zeta}\right)^{1+a_{\abs}/b_{\abs}}\cdot s
        \right)
      }.
    \end{gather*}
    From the last expression, it suffices to compute the rational $(1-a_{|UC}/b_{\UC})\cdot(a_\zeta\cdot s/b_\zeta)$ (which
    takes time at most $O(\len(c_{\UC})\cdot\len(\zeta)\cdot\log(n+1))$) and determine the ceiling of
    \begin{gather*}
      \frac{a_\zeta^{1+a_{\abs}/b_{\abs}}\cdot s\cdot b_{\UC}}{b_\zeta^{a_{\abs}/b_{\abs}}},
    \end{gather*}
    that is, we need to determine the unique $f\in\NN$ such that
    \begin{gather*}
      b_\zeta^{a_{\abs}/b_{\abs}}\cdot f
      <
      a_\zeta^{1+a_{\abs}/b_{\abs}}\cdot s\cdot b_{\UC}
      \leq
      b_\zeta^{a_{\abs}/b_{\abs}}\cdot (f+1),
    \end{gather*}
    which in turn is equivalent to
    \begin{gather}\label{eq:nonequi:good:rational:si:ceil}
      b_\zeta^{a_{\abs}}\cdot f^{b_{\abs}}
      <
      a_\zeta^{a_{\abs} + b_{\abs}}\cdot s^{b_{\abs}}\cdot b_{\UC}^{b_{\abs}}
      \leq
      b_\zeta^{a_{\abs}}\cdot (f+1)^{b_{\abs}}.
    \end{gather}

    Such $f$ will be in $\{0,\ldots,b_{\UC}\cdot b_\zeta\cdot n\}$, so we can do a simple binary search for it (which takes at
    most $\log_2(b_{\UC}\cdot b_\zeta\cdot n+1)\leq\len(c_{\UC})+\len(\zeta)+\log(n)$ iterations). We can give the following
    crude bound for the time-complexity of computing the quantities in~\eqref{eq:nonequi:good:rational:si:ceil}:
    \begin{align*}
      \MoveEqLeft
      \begin{multlined}[t]
        O\bigl(
        \log(a_{\abs}+1)\cdot a_{\abs}^2\cdot\bigl(\log(b_\zeta+1)\bigr)^2
        +
        \log(b_{\abs}+1)\cdot b_{\abs}^2\cdot\bigl(\log(b_\zeta\cdot n+1)\bigr)^2
        \\
        +
        a_{\abs}\cdot\log(b_{\abs}+1)\cdot b_{\abs}\cdot\log(b_\zeta\cdot n+1)
        +
        \log(a_{\abs} + b_{\abs})\cdot(a_{\abs}+b_{\abs})^2\cdot\bigl(\log(a_\zeta+1)\bigr)^2
        \\
        +
        \log(b_{\abs}+1)\cdot b_{\abs}^2\cdot\bigl(\log(n+1)\bigr)^2
        +
        \log(b_{\abs}+1)\cdot b_{\abs}^2\cdot\bigl(\log(b_{\UC}+1)\bigr)^2
        \\
        +
        (a_{\abs}+b_{\abs})\cdot\log(a_\zeta+1)\cdot b_{\abs}\cdot\log(n+1)\cdot b_{\abs}\cdot\log(b_{\UC}+1)
        \bigr)
      \end{multlined}
      \\
      & \leq
      O\bigl(
      2^{3\cdot\len(c_{\abs})}\cdot\len(\zeta)^2\cdot\len(c_{\UC})^2\cdot\bigl(\log(n+1)\bigr)^2
      \bigr).
    \end{align*}
    (Note that this completely dominates the time it took to compute the rational $(1-a_{|UC}/b_{\UC})\cdot(a_\zeta\cdot
    s/b_\zeta)$.) Since we need to make the binary search, we make $\ell+1$ such estimations per iteration of the loop of
    line~\ref{alg:nonequi:good:While} and this loop executes at most
    \begin{gather*}
      K_{\nonequi}(\ell,c_{\UC},c_{\ZZ},c_{\abs},\epsilon)
      \leq
      O\left(\frac{\log(1/\epsilon)}{\epsilon^\ell}\right)
    \end{gather*}
    times, the total time-complexity of line~\ref{alg:nonequi:good:rational:si}, accounting for every time it executes, is at
    most
    \begin{gather}\label{eq:nonequi:good:rational:si:complexity}
      \begin{aligned}
        \MoveEqLeft
        O\left(
        \frac{\log(1/\epsilon)}{\epsilon^\ell}
        \cdot\ell
        \cdot
        2^{3\cdot\len(c_{\abs})}\cdot\len(\zeta)^2\cdot\len(c_{\UC})^2\cdot\bigl(\log(n+1)\bigr)^2
        \right)
        \\
        & \leq
        O\left(
        \frac{\log(1/\epsilon)}{\epsilon^\ell}
        \cdot\ell
        \cdot
        2^{3\cdot\len(c_{\abs})}\cdot\bigl(\len(c_{\ZZ})+\len(\epsilon)\bigr)^2\cdot\len(c_{\UC})^2\cdot\bigl(\log(n+1)\bigr)^2
        \right).
      \end{aligned}
    \end{gather}

    \smallskip

    The time-complexity analysis of the rational approximations of line~\ref{alg:nonequi:good:rational:ws} is the same as that
    of Algorithm~\ref{alg:nonequi:atom}, hence yields the time-complexity as in~\eqref{eq:nonequi:atom:rational:ws:complexity}:
    \begin{gather}\label{eq:nonequi:good:rational:ws:complexity}
      \begin{aligned}
        \MoveEqLeft
        O\left(
        \frac{\log(1/\epsilon)}{\epsilon^\ell}
        \cdot
        2^{3\cdot\len(c_{\abs})}\cdot\len(\zeta)^2\cdot\bigl(\log(n+1)\bigr)^3
        \right)
        \\
        & \leq
        O\left(
        \frac{\log(1/\epsilon)}{\epsilon^\ell}
        \cdot
        2^{3\cdot\len(c_{\abs})}\cdot\bigl(\len(c_{\ZZ})+\len(\epsilon)\bigr)^2\cdot\bigl(\log(n+1)\bigr)^3
        \right).
      \end{aligned}
    \end{gather}

    Putting everything together (i.e., equations~\eqref{eq:nonequi:good:While:complexity},
    \eqref{eq:nonequi:good:For:complexity}, \eqref{eq:nonequi:good:param:complexity},
    \eqref{eq:nonequi:good:rational:S:complexity}, \eqref{eq:nonequi:good:rational:m:complexity},
    \eqref{eq:nonequi:good:rational:wm:complexity}, \eqref{eq:nonequi:good:rational:si:complexity}
    and~\eqref{eq:nonequi:good:rational:ws:complexity}), we bound the time-complexity of Algorithm~\ref{alg:nonequi:good} by:
    \begin{align*}
      \MoveEqLeft
      \begin{multlined}[t]
        O\bigggl(
        \frac{2^\ell\cdot\ell^2\cdot(\log(1/\epsilon))^2}{c_{\UC}^2\cdot(1-c_{\ZZ})^2\cdot\epsilon^{\ell+2}}
        \cdot n\cdot\log\left(\frac{n}{\rho_{\UC}}\right)
        +
        \frac{\log(1/\epsilon)}{\epsilon^\ell}\cdot n
        \\
        +
        \ell^2\cdot\log(\ell+1)\cdot\len(c_{\UC})\cdot\len(c_{\abs})\cdot\len(c_{\ZZ})\cdot\len(\epsilon)^2\cdot\len(\rho_{\UC})
        \\
        +
        2^{3\cdot\len(c_{\abs})}\cdot\bigl(\len(c_{\ZZ})+\len(\epsilon)\bigr)^3\cdot\log(n+1)
        \\
        +
        2^{3\cdot\len(c_{\abs})}
        \cdot\len(c_{\UC})^6\cdot\bigl(\len(c_{\ZZ})+\len(\epsilon)\bigr)^3
        \cdot\len(\rho_{\UC})\cdot\ell\cdot\log(n+1)
        \\
        +
        2^{3\cdot\len(c_{\abs})}\cdot\len(c_{\UC})^2\cdot\bigl(\len(c_{\ZZ})+\len(\epsilon)\bigr)^3
        \\
        +
        \frac{\log(1/\epsilon)}{\epsilon^\ell}
        \cdot\ell
        \cdot
        2^{3\cdot\len(c_{\abs})}\cdot\bigl(\len(c_{\ZZ})+\len(\epsilon)\bigr)^2\cdot\len(c_{\UC})^2\cdot\bigl(\log(n+1)\bigr)^2
        \\
        +
        \frac{\log(1/\epsilon)}{\epsilon^\ell}
        \cdot
        2^{3\cdot\len(c_{\abs})}\cdot\bigl(\len(c_{\ZZ})+\len(\epsilon)\bigr)^2\cdot\bigl(\log(n+1)\bigr)^3
        \bigggr)
      \end{multlined}
      \\
      & \leq
      \begin{multlined}[t]
        O\bigggl(
        \frac{2^\ell\cdot\ell^2\cdot(\log(1/\epsilon))^2}{c_{\UC}^2\cdot(1-c_{\ZZ})^2\cdot\epsilon^{\ell+2}}
        \cdot n\cdot\log\left(\frac{n}{\rho_{\UC}}\right)
        \\
        +
        \ell^2\cdot\log(\ell+1)
        \cdot 2^{3\cdot\len(c_{\abs})}\cdot\bigl(\len(c_{\ZZ})+\len(\epsilon)\bigr)^3\cdot\len(c_{\UC})^6\cdot\len(\rho_{\UC})
        \\
        \cdot
        \frac{\log(1/\epsilon)}{\epsilon^\ell}
        \cdot
        \bigl(\log(n+1)\bigr)^3
        \bigggr).
      \end{multlined}
    \end{align*}
  \item[Item~\ref{thm:nonequi:atomspace}.] The idea behind the space-efficient version of Algorithm~\ref{alg:nonequi:atom} is
    the same, except that we call Algorithm~\ref{alg:ext:atomspace} and we will need an auxiliary oracle algorithm responsible
    for keeping track of the ``remainder'' set. Furthermore, let us already point out that both the main oracle algorithm and
    the auxiliary algorithm will also be used as the oracle algorithms of item~\ref{thm:nonequi:goodspace}.

    \begin{algorithm}[htbp]
      \caption{Computation algorithm of query-into-oracle model with space- and time-complexities
        \begin{gather*}
          \begin{multlined}[t]
            O\bigggl(
            \frac{\ell\cdot d\cdot(\log(1/\epsilon))^2\cdot\log(n+1)}{c_{\atom}^2\cdot\delta^2\cdot\epsilon^\ell}
            +
            2^{2\cdot\len(c_{\abs})}\cdot\bigl(\len(c_{\ZZ}) + \len(\epsilon) + \log(n+1)\bigr)
            \\
            +
            \len(c_{\atom}) + \len(\delta)
            \bigggr),
          \end{multlined}
          \\
          \begin{multlined}[t]
            n^{O(d/(c_{\atom}^2\cdot\delta^2) + \epsilon^{-\ell}\cdot\log(1/\epsilon))}
            +
            O\bigggl(
            \len(c_{\atom})^2\cdot\len(\delta)^2\cdot\bigl(\log(d+1)\bigr)^2
            \\
            +
            \frac{\log(1/\epsilon)}{\epsilon^\ell}
            \cdot\ell
            \cdot
            2^{3\cdot\len(c_{\abs})}\cdot\bigl(\len(c_{\ZZ})+\len(\epsilon)\bigr)^3\cdot\bigl(\log(n+1)\bigr)^3
            \bigggr).
          \end{multlined}
        \end{gather*}
        respectively, that returns $(m,\widetilde{m},(\sigma^k,x^k,s^k,\widetilde{s}^k)_{k=1}^K)$ such that when given to
        Algorithm~\ref{alg:nonequi:atomspaceoracle} produces an oracle for a $(\delta,\epsilon)$-atomic partition of $G$ into
        $K\leq K_{\nonequi}(\ell,0,c_{\ZZ},c_{\abs},\epsilon)$ parts.}
      \label{alg:nonequi:atomspace}
      \DontPrintSemicolon
      \KwIn{Numbers $n,d,\ell\in\NN_+$ with $d\leq\ell\leq\log_2(n)$, $c_{\atom},c_{\ZZ},c_{\abs},\delta\in(0,1)\cap\QQ$ with
        $(1+c_{\atom})\cdot\delta < 1/2^{\ell+1}$, $\epsilon\in(0,\epsilon_{\nonequi}(\ell,c_{\ZZ},c_{\abs}))\cap\QQ$, and a
        query-oracle $\cO_G$ for a graph $G$ with $\lvert G\rvert=n$, $\Lit(G)\leq\ell$ and $\VC(G)\leq d$. We further assume
        $n\geq n_{\nonequi}(\ell,0,c_{\ZZ},c_{\abs},\epsilon)$.}
      \KwOut{A tuple $(m,\widetilde{m},(\sigma^k,x^k,s^k,\widetilde{s}^k)_{k=1}^K)$, with $K\leq
        K_{\nonequi}(\ell,0,c_{\ZZ},c_{\abs},\epsilon)$, where $m\in\NN_+$, $\widetilde{m}\in\NN$ and for each $k\in[K]$, we
        have $\sigma^k\in\{0,1\}^{<\ell+1}$, $x^k\in (V(G)^m)^{\lvert\sigma\rvert}$ and $s^k,\widetilde{s}^k\in[n]$, such that
        when given to Algorithm~\ref{alg:nonequi:atomspaceoracle} produces an oracle for a $(\delta,\epsilon)$-atomic partition
        of $G$ into $K$ parts.}
      $\widetilde{\delta}\assign(1+c_{\atom})\cdot\delta$\;
      Let $m\assign\ceil{C\cdot d/(c_{\atom}^2\cdot\delta^2)}$, where $C$ is the absolute constant of
      Lemma~\ref{lem:atomconst}.\;\label{alg:nonequi:atomspace:m}
      $\widetilde{m}\assign\floor{\widetilde{\delta}\cdot m}$\;
      $c_{\UC}\assign 0$\;
      $\zeta\assign (1-c_{\ZZ})\cdot\epsilon$\;
      Compute a sufficiently precise rational approximation of $\widetilde{\epsilon}\df (1 - \zeta^{c_{\abs}})\cdot\zeta$ to
      compute $S\df\floor{(\zeta-\widetilde{\epsilon})\cdot n/(1-\widetilde{\epsilon})}$.\;
      \label{alg:nonequi:atomspace:rational:S}
      $K\assign 0$\;
      $r\assign n$\;
      \uWhile{$r > S$}{%
        \label{alg:nonequi:atomspace:While}
        $s_0\assign r$\;
        Compute sufficiently precise rational approximations of $\widetilde{\epsilon}\df (1 - \zeta^{c_{\abs}})\cdot\zeta$ and
        $\gamma\df (1 - c_{\UC})\cdot\widetilde{\epsilon}$ to compute $s_i\df\floor{\gamma\cdot s_{i-1}} + 1$ for every
        $i\in[\ell+1]$.\;
        \label{alg:nonequi:atomspace:rational:si}
        %
        \setcounter{algosplit}{\theAlgoLine}
      }
    \end{algorithm}

    \begin{algorithm}[htbp]
      \caption*{(continued).}
      \DontPrintSemicolon
      \setcounter{AlgoLine}{\thealgosplit}
      \let\oldnl\nl
      \let\nl\relax
      \vspace{-\baselineskip}\InvisibleBegin{
        \global\let\nl\oldnl 
        Let $\cO$ be Algorithm~\ref{alg:nonequi:atomspaceauxiliary} with parameters:
        \begin{gather*}
          \Bigl(n,\cO_G,\bigl(m,\widetilde{m},(\sigma^k,x^k,s^k)_{k=1}^K\bigr),u\Bigr)
          \df
          \Bigl(n,\cO_G,\bigl(m,\widetilde{m},(\sigma^k,x^k,s^k)_{k=1}^K\bigr),\place\Bigr).
        \end{gather*}
        \;
        Run Algorithm~\ref{alg:ext:atomspace} with parameters:
        \begin{gather*}
          (n,m,\ell,s_0,\ldots,s_\ell,s_{\ell+1},\widetilde{m},\cO_G,\cO_U)
          \df
          (n,m,\ell,s_0,\ldots,s_\ell,s_{\ell+1},\widetilde{m},\cO_G,\cO).
        \end{gather*}
        \;\label{alg:nonequi:atomspace:call}
        Let $(\sigma,(x_j)_{j=0}^{\lvert\sigma\rvert-1})$ be the tuple returned by the algorithm.\;
        $K\assign K+1$\;
        $\sigma^K\assign\sigma$\;
        $x^K\assign (x_j)_{j=0}^{\lvert\sigma\rvert-1}$\;
        $s^K\assign s_{\lvert\sigma\rvert}$\;
        $r\assign r - s^K$\;
        Compute a sufficiently precise rational approximation of $\widetilde{\epsilon}\df (1 - \zeta^{c_{\abs}})\cdot\zeta$ to
        compute $\widetilde{s}^K\df\ceil{(\zeta-\widetilde{\epsilon})\cdot s^K/(1-\zeta)}$.\;
        \label{alg:nonequi:atomspace:rational:ws}
      }
      \Return{$(m,\widetilde{m},(\sigma^k,x^k,s^k,\widetilde{s}^k)_{k=1}^K)$}
    \end{algorithm}

    \begin{algorithm}[htbp]
      \caption{Oracle algorithm of query-into-oracle model with space- and time-complexities
        \begin{align*}
          O\bigl(\log(K+1) + \log(m+1) + \log(n+1)\bigr),
          & &
          O\bigl(K^2\cdot m\cdot n^{K+1}\bigr).
        \end{align*}
        respectively, that when given $(m,\widetilde{m},(\sigma^k,x^k,s^k,\widetilde{s}^k)_{k=1}^K)$ computed via
        Algorithm~\ref{alg:nonequi:atomspace}, produces an oracle for a $(\delta,\epsilon)$-atomic partition of $G$ with $K\leq
        K_{\nonequi}(\ell,0,c_{\ZZ},c_{\abs},\epsilon)$ parts. If instead
        $(m,\widetilde{m},(\sigma^k,x^k,s^k,\widetilde{s}^k)_{k=1}^K)$ is computed via Algorithm~\ref{alg:nonequi:goodspace},
        then the oracle produced is for an $\epsilon$-good partition of $G$ with $K\leq
        K_{\nonequi}(\ell,0,c_{\ZZ},c_{\abs},\epsilon)$ parts.}
      \label{alg:nonequi:atomspaceoracle}
      \DontPrintSemicolon
      \KwIn{A number $n\in\NN_+$, a query-oracle $\cO_G$ for a graph $G$ with $V(G)=[n]$, a tuple
        $(m,\widetilde{m},(\sigma^k,x^k,s^k,\widetilde{s}^k)_{k=1}^K)$ provided by Algorithm~\ref{alg:nonequi:atomspace} or
        Algorithm~\ref{alg:nonequi:goodspace} ran with the same parameters (and a choice of the other parameters required by it
        and satisfying the required hypotheses), a vertex $u\in V(G)$ and an index $k\in[K]$.}
      \KwOut{A bit $b\in\{0,1\}$ such that if $P_k$ is the set of $u\in V(G)$ for which the algorithm returns $1$ for $k\in[K]$,
        then $\cP = \{P_k \mid k\in[K]\}$ is a partition of $G$. If the input came from Algorithm~\ref{alg:nonequi:atomspace},
        then $\cP$ is guaranteed to be $(\delta,\epsilon)$-atomic, where $\delta$ and $\epsilon$ are the parameters used by
        Algorithm~\ref{alg:nonequi:atomspace}. If the input came from Algorithm~\ref{alg:nonequi:goodspace}, then $\cP$ is
        guaranteed to be $\epsilon$-good, where $\epsilon$ is the parameter used by Algorithm~\ref{alg:nonequi:goodspace}.}
      Run Algorithm~\ref{alg:nonequi:atomspaceauxiliary} with parameters:
      \begin{gather*}
        \Bigl(n,\cO_G,\bigl(m,\widetilde{m},(\sigma^k,x^k,s^k)_{k=1}^K\bigr),u\Bigr)
        \df
        \Bigl(n,\cO_G,
        \bigl(m,\widetilde{m},(\sigma^{\widetilde{k}},x^{\widetilde{k}},s^{\widetilde{k}})_{\widetilde{k}=1}^{k-1}\bigr),u\Bigr).
      \end{gather*}
      \;
      \lIf(\tcp*[f]{$u$ is in the main portion of the first $k-1$ parts}){The algorithm returned $0$}{\Return{$0$}}
      Run Algorithm~\ref{alg:nonequi:atomspaceauxiliary} with parameters:
      \begin{gather*}
        \Bigl(n,\cO_G,\bigl(m,\widetilde{m},(\sigma^k,x^k,s^k)_{k=1}^K\bigr),u\Bigr)
        \df
        \Bigl(n,\cO_G,
        \bigl(m,\widetilde{m},(\sigma^{\widetilde{k}},x^{\widetilde{k}},s^{\widetilde{k}})_{\widetilde{k}=1}^k\bigr),u\Bigr).
      \end{gather*}
      \;
      \lIf(\tcp*[f]{$u$ is in the main portion of the $k$th part}){The algorithm returned $0$}{\Return{$1$}}
      \tcc*[l]{At this point, we know $u$ is not in the main portion of the $k$th part, but could be an absorbed vertex of
        remainder.}
      Run Algorithm~\ref{alg:nonequi:atomspaceauxiliary} with parameters:
      \begin{gather*}
        \Bigl(n,\cO_G,\bigl(m,\widetilde{m},(\sigma^k,x^k,s^k)_{k=1}^K\bigr),u\Bigr)
        \df
        \Bigl(n,\cO_G,
        \bigl(m,\widetilde{m},(\sigma^{\widetilde{k}},x^{\widetilde{k}},s^{\widetilde{k}})_{\widetilde{k}=1}^K\bigr),u\Bigr).
      \end{gather*}
      \;
      \lIf(\tcp*[f]{$u$ is not in the remainder}){The algorithm returned $0$}{\Return{$0$}}
      \tcp*[l]{Now $u$ is in the remainder; need to find which part absorbed it}
      $i\assign 1$\tcp*[r]{Index of part absorbing current vertex}
      $a\assign 0$\tcp*[r]{Number of vertices of remainder absorbed by current part so far}
      \setcounter{algosplit}{\theAlgoLine}
    \end{algorithm}

    \begin{algorithm}[htbp]
      \caption*{(continued).}
      \DontPrintSemicolon
      \setcounter{AlgoLine}{\thealgosplit}
      \let\oldnl\nl
      \let\nl\relax
      \global\let\nl\oldnl 
      \For(\tcp*[f]{Vertex that potentially is in the remainder}){$w\in[u-1]$}{%
        Run Algorithm~\ref{alg:nonequi:atomspaceauxiliary} with parameters:
        \begin{gather*}
          \Bigl(n,\cO_G,\bigl(m,\widetilde{m},(\sigma^k,x^k,s^k)_{k=1}^K\bigr),u\Bigr)
          \df
          \Bigl(n,\cO_G,
          \bigl(m,\widetilde{m},(\sigma^{\widetilde{k}},x^{\widetilde{k}},s^{\widetilde{k}})_{\widetilde{k}=1}^K\bigr),w\Bigr).
        \end{gather*}
        \;
        \If(\tcp*[f]{$w$ is in remainder and is absorbed by $i$th part}){The algorithm returned $1$}{%
          $a\assign a+1$\;
          \If(\tcp*[f]{$i$th part already absorbed all its remainder vertices}){$a\geq\widetilde{s}^i$}{%
            $i\assign i+1$\;
            \lIf(\tcp*[f]{$u$ is absorbed by a part after the $k$th}){$i > k$}{\Return{$0$}}
            $a\assign 0$\;
          }
        }
      }
      \Return{$1$}\tcp*[f]{$u$ is absorbed by $k$th part}
    \end{algorithm}

    \begin{algorithm}[htbp]
      \caption{Auxiliary oracle algorithm that when given $(m,\widetilde{m},(\sigma^k,x^k,s^k)_{k=1}^K)$ provided within the
        execution of one of Algorithms~\ref{alg:nonequi:atomspace} or~\ref{alg:nonequi:goodspace}, produces an oracle for the
        set of vertices that have not yet been covered by the partition. The space- and time-complexities are
        \begin{align*}
          O\bigl(\log(K+1) + \log(m+1) + \log(n+1)\bigr),
          & &
          O\bigl(1 + K^2\cdot m\cdot n^K\bigr).
        \end{align*}
        respectively.}
      \label{alg:nonequi:atomspaceauxiliary}
      \DontPrintSemicolon
      \KwIn{A number $n\in\NN_+$, a query-oracle $\cO_G$ for a graph $G$ with $V(G)=[n]$, a tuple
        $(m,\widetilde{m},(\sigma^k,x^k,s^k)_{k=1}^K)$ provided by one of Algorithms~\ref{alg:nonequi:atomspace}
        or~\ref{alg:nonequi:goodspace} ran with the same parameters (and a choice of the other parameters required by it and
        satisfying the required hypotheses), and a vertex $u\in V(G)$.}
      \KwOut{A bit $b\in\{0,1\}$ such that if $W$ is the set of $u\in V(G)$ for which the algorithm returns $1$, then the sets
        defined by $(\sigma^k,x^k,s^k)_{k=1}^K$ cover exactly $V(G)\setminus W$ (so $W$ is the remainder set).}
      \lIf{$K=0$}{\Return{$1$}}
      \For{$\widetilde{k}\in[K]$}{%
        Let $\cO$ be Algorithm~\ref{alg:nonequi:atomspaceauxiliary} with parameters:
        \begin{gather*}
          \Bigl(n,\cO_G,\bigl(m,\widetilde{m},(\sigma^k,x^k,s^k)_{k=1}^K\bigr),u\Bigr)
          \df
          \Bigl(n,\cO_G,\bigl(m,\widetilde{m},(\sigma^k,x^k,s^k)_{k=1}^{\widetilde{k}-1}\bigr),\place\Bigr).
        \end{gather*}
        \;
        Run Algorithm~\ref{alg:ext:atomspaceoracle} with parameters:
        \begin{gather*}
          \bigl(n,m,\widetilde{m},\cO_G,\cO_U,(\sigma,(x_j)_{j=0}^{\lvert\sigma\rvert-1}),u\bigr)
          \df
          \bigl(n,m,\widetilde{m},\cO_G,\cO,(\sigma^{\widetilde{k}},x^{\widetilde{k}}),u\bigr).
        \end{gather*}
        \;
        \lIf(\tcp*[f]{$u$ is in $\widetilde{k}$th part}){The algorithm returned $1$}{\Return{$0$}}
      }
      \Return{$1$}
    \end{algorithm}

    Correctness of Algorithms~\ref{alg:nonequi:atomspace} and~\ref{alg:nonequi:atomspaceoracle} amount to observing that they
    are essentially mimicking Algorithm~\ref{alg:nonequi:atom} but using the space-efficient versions of the algorithms from
    Proposition~\ref{prop:ext}. Namely, they are repeatedly extracting $\widetilde{\epsilon}$-good parts until the remainder is
    of size at most $S\df\floor{(\zeta-\widetilde{\epsilon})\cdot n/(1-\widetilde{\epsilon})}$, then adding at most
    $\widetilde{s}^i\df\ceil{(\zeta-\widetilde{\epsilon})\cdot s^i/(1-\zeta)}$ vertices to the $i$th part, which was guaranteed
    to have size $s^i$. One only needs to check that the auxiliary oracle algorithm,
    Algorithm~\ref{alg:nonequi:atomspaceauxiliary}, indeed encodes an oracle of the set of vertices that have not yet been
    covered by the growing partition (i.e., not covered by the first $K$ parts); this itself is also straightforward to prove
    inductively on $K$.

    \medskip

    We now analyze the space- and time-complexities of the algorithms. We start with
    Algorithm~\ref{alg:nonequi:atomspaceauxiliary}. Since the algorithm passes itself as an oracle parameter to
    Algorithm~\ref{alg:ext:atomspaceoracle}, our analysis is recursive:

    \begin{claim}\label{clm:nonequi:atomspaceauxiliary}
      For an input $(n,\cO_G,(m,\widetilde{m},(\sigma^k,x^k,s^k)_{k=1}^K),u)$ and for $i\in\{0,\ldots,K\}$, let $\cO_i$ be
      Algorithm~\ref{alg:nonequi:atomspaceauxiliary} when ran with $(n,\cO_G,(m,\widetilde{m},(\sigma^k,x^k,s^k)_{k=1}^i),u)$
      and let $S(i)$ and $T(i)$ its space- and time-complexities. Then we have
      \begin{align*}
        S(K)
        & \leq
        O\Bigl(K\cdot\bigl(\log(K+1) + \log(m+1) + \log(n+1)\bigr)\Bigr),
        \\
        T(i)
        & \leq
        O\bigl(1 + i^2\cdot m\cdot n^i\bigr),
      \end{align*}
      for every $i\in\{0,\ldots,K\}$.
    \end{claim}

    \begin{proofof}{Claim~\ref{clm:nonequi:atomspaceauxiliary}}
      Let $\widetilde{S}(i)$, $\widetilde{T}(i)$ and $\widetilde{U}(i)$ be the space-, time- and $U$-oracle-complexities of
      Algorithm~\ref{alg:ext:atomspaceoracle} when ran with $(n,m,\widetilde{m},\cO_G,\cO_i,(\sigma^i,x^i),u)$. By
      Proposition~\ref{prop:ext}\ref{prop:ext:atomspace}, we have
      \begin{align*}
        \widetilde{S}(i) & = O\bigl(\log(\lvert\sigma\rvert+1) + \log(m+1) + \log(n+1)\bigr)
        &
        \widetilde{T}(i) & = O\bigl((i+1)\cdot m\cdot n\bigr),
        &
        \widetilde{U}(i) & \leq n.
      \end{align*}
      Thus, we get the following recursive formulas for $i\in[K]$:
      \begin{align*}
        S(i)
        & =
        \log_2(i+1) + \max_{j\in\{0,\ldots,i-1\}} \bigl(\widetilde{S}(j) + S(j)\bigr)
        \\
        & \leq
        \log_2(i+1) + O\bigl(\log(i) + \log(m+1) + \log(n+1)\bigr) + \max_{j\in\{0,\ldots,i-1\}} S(j),
        \\
        T(i)
        & \leq
        \sum_{j=0}^{i-1} \bigl(\widetilde{T}(j) + \widetilde{U}(j)\cdot T(j) + O(1)\bigr)
        \leq
        \sum_{j=0}^{i-1} \Bigl(O\bigl((j+1)\cdot m\cdot n\bigr) + n\cdot T(j)\Bigr)
        \\
        & \leq
        O\left(\binom{i+1}{2}\cdot m\cdot n\right) + n\cdot\sum_{j=0}^{i-1} T(j).
      \end{align*}

      Therefore, we conclude that
      \begin{gather*}
        S(K)
        \leq
        O\Bigl(K\cdot\bigl(\log(K+1) + \log(m+1) + \log(n+1)\bigr)\Bigr).
      \end{gather*}
      (Note that this also upper bounds $S(i)$.)

      For $T$, it is straightforward to prove by induction in $i$ that
      \begin{gather*}
        T(i)
        \leq
        O\left(
        \binom{i+1}{2}\cdot m\cdot n
        +
        \sum_{j=1}^{i-1}\binom{j+1}{2}\cdot m\cdot n^2\cdot(n+1)^{i-j-1}
        \right),
      \end{gather*}
      from which we get
      \begin{gather*}
        T(i)
        \leq
        O\bigl(1 + i^2\cdot m\cdot n^i\bigr)
      \end{gather*}
      for every $i\in\{0,\ldots,K\}$.
    \end{proofof}

    We now proceed to the analysis of oracle Algorithm~\ref{alg:nonequi:atomspaceoracle}. Note that it calls the auxiliary
    algorithm, Algorithm~\ref{alg:nonequi:atomspaceauxiliary}, at most $u+2\leq n+2$ times. It is also clear that it needs
    internal memory to store the variables $i\in[K]$, $a\in\{0,\ldots,n\}$ and $w\in[n]$, which takes space
    \begin{gather*}
      O\bigl(\log(K+1) + \log(n+1)\bigr).
    \end{gather*}

    Taking into account our earlier analysis of the auxiliary Algorithm~\ref{alg:nonequi:atomspaceauxiliary} in
    Claim~\ref{clm:nonequi:atomspaceauxiliary}, we conclude that the space- and time-complexities of
    Algorithm~\ref{alg:nonequi:atomspaceoracle} are at most
    \begin{align*}
      O\Bigl(K\cdot\bigl(\log(K+1) + \log(m+1) + \log(n+1)\bigr)\Bigr),
      & &
      O\bigl(K^2\cdot m\cdot n^{K+1}\bigr),
    \end{align*}
    respectively.

    Taking into account the fact that this will be executed with parameters from the computation
    Algorithm~\ref{alg:nonequi:atomspace}, we know that
    \begin{align*}
      K
      & \leq
      K_{\nonequi}(\ell,0,c_{\ZZ},c_{\abs},\epsilon)
      \leq
      O\left(\frac{\log(1/\epsilon)}{\epsilon^\ell}\right),
      \\
      m
      & \leq
      O\left(\frac{d}{c_{\atom}^2\cdot\delta^2}\right)
      \leq
      O\left(\frac{\ell}{c_{\atom}^2\cdot\delta^2}\right),
    \end{align*}
    so the space- and time-complexities become
    \begin{align*}
      O\left(
      \frac{\log(1/\epsilon)}{\epsilon^\ell}\cdot
      \log\left(\frac{n}{c_{\atom}\cdot\delta\cdot\epsilon^\ell}\right)
      \right),
      & &
      \frac{n^{O(\epsilon^{-\ell}\cdot\log(1/\epsilon))}}{c_{\atom}^2\cdot\delta^2}.
    \end{align*}

    \smallskip

    Finally, we analyze the computation Algorithm~\ref{alg:nonequi:atomspace}. We know that $K$ is always at most
    \begin{gather*}
      K_{\nonequi}(\ell,0,c_{\ZZ},c_{\abs},\epsilon)
      \leq
      O\left(\frac{\log(1/\epsilon)}{\epsilon^\ell}\right)
    \end{gather*}
    and that this also bounds the number of times that the loop of line~\ref{alg:nonequi:atomspace:While} executes.

    We also know that Algorithm~\ref{alg:nonequi:atomspace} needs space to store the following variables:
    \begin{gather*}
      \begin{aligned}
        \widetilde{\delta}, & &
        m, & &
        \widetilde{m},
        \\
        \zeta, & &
        K, & &
        r,
      \end{aligned}
      \\
      \begin{aligned}
        s_0,\ldots,s_{\ell+1}, & &
        (\sigma^k,x^k,s^k,\widetilde{s}^k)_{k=1}^K.
      \end{aligned}
    \end{gather*}
    These require at most the following space, respectively:
    \begin{gather*}
      \begin{aligned}
        O\bigl(\len(c_{\atom}) + \len(\delta)\bigr),
        & &
        O\bigl(\log(d) + \len(c_{\atom}) + \len(\delta)\bigr),
        & &
        O\bigl(\log(d) + \len(c_{\atom}) + \len(\delta)\bigr),
        \\
        O\bigl(\len(c_{\ZZ}) + \len(\epsilon)\bigr),
        & &
        O\left(\ell\cdot\log\left(\frac{1}{\epsilon}\right)\right),
        & &
        \log(n+1),
      \end{aligned}
      \\
      \begin{aligned}
        \ell\cdot\log(n+1),
        & &
        O\left(
        \frac{\log(1/\epsilon)}{\epsilon^\ell}
        \cdot
        \ell\cdot m\cdot\log(n+1)
        \right).
      \end{aligned}
    \end{gather*}
    These can be upper bounded together by
    \begin{gather}\label{eq:nonequi:atomspace:varspace}
      \begin{aligned}
        \MoveEqLeft
        O\left(
        \len(c_{\atom}) + \len(\delta) + \len(c_{\ZZ}) + \len(\epsilon)
        + \ell\cdot\frac{\log(1/\epsilon)}{\epsilon^\ell}\cdot m\cdot\log(n+1)
        \right)
        \\
        & \leq
        O\left(
        \len(c_{\atom}) + \len(\delta) + \len(c_{\ZZ}) + \len(\epsilon)
        +
        \frac{\ell\cdot d\cdot\log(1/\epsilon)\cdot\log(n+1)}{c_{\atom}^2\cdot\delta^2\cdot\epsilon^\ell}
        \right).
      \end{aligned}
    \end{gather}

    Furthermore, Algorithm~\ref{alg:nonequi:atomspace} also needs space to compute the rational approximations of
    lines~\ref{alg:nonequi:atomspace:rational:S}, \ref{alg:nonequi:atomspace:rational:si}
    and~\ref{alg:nonequi:atomspace:rational:ws}. These approximations are computed just as in Algorithm~\ref{alg:nonequi:atom}
    (see~\eqref{eq:nonequi:atom:rational:S:complexity}, \eqref{eq:nonequi:atom:rational:si:complexity}
    and~\eqref{eq:nonequi:atom:rational:ws:complexity}). This means that we make in place calculations and binary searches,
    meaning that we only need to keep in memory a constant amount of variables with which we are performing calculations. This
    amounts to the following respective space-complexities:
    \begin{gather*}
      O\Bigl(
      2^{2\cdot\len(c_{\abs})}\cdot\bigl(\log(n+1) + \len(c_{\ZZ}) + \len(\epsilon)\bigr)
      \Bigr),
      \\
      O\Bigl(
      2^{2\cdot\len(c_{\abs})}\cdot\bigl(\len(c_{\ZZ}) + \len(\epsilon) + \log(n+1)\bigr)
      \Bigr),
      \\
      O\Bigl(
      2^{\len(c_{\abs})}\cdot\bigl(\len(c_{\ZZ}) + \len(\epsilon) + \log(n+1)\bigr)
      \Bigr),
    \end{gather*}
    which can be upper bounded together by
    \begin{gather}\label{eq:nonequi:atomspace:rationalspace}
      O\Bigl(
      2^{2\cdot\len(c_{\abs})}\cdot\bigl(\len(c_{\ZZ}) + \len(\epsilon) + \log(n+1)\bigr)
      \Bigr).
    \end{gather}

    Finally, since the algorithm calls Algorithm~\ref{alg:ext:atomspace} with Algorithm~\ref{alg:nonequi:atomspaceauxiliary} as
    an oracle, we need to account for the space-complexity of these algorithms, which are:
    \begin{gather*}
      O\bigl(\ell\cdot m\cdot\log(n+1)\bigr)
      \leq
      O\left(\frac{\ell\cdot d\cdot\log(n+1)}{c_{\atom}^2\cdot\delta^2}\right),
      \\
      O\Bigl(K\cdot\bigl(\log(K+1) + \log(m+1) + \log(n+1)\bigr)\Bigr)
      \leq
      O\left(
      \frac{\log(1/\epsilon)}{\epsilon^\ell}\cdot
      \log\left(\frac{n}{c_{\atom}\cdot\delta\cdot\epsilon^\ell}\right)
      \right),
    \end{gather*}
    respectively.

    Putting these together with~\eqref{eq:nonequi:atomspace:varspace} and~\eqref{eq:nonequi:atomspace:rationalspace}, the total
    space-complexity of Algorithm~\ref{alg:nonequi:atomspace} is at most
    \begin{align*}
      \MoveEqLeft
      \begin{multlined}[t]
        O\bigggl(
        \len(c_{\atom}) + \len(\delta) + \len(c_{\ZZ}) + \len(\epsilon)
        +
        \frac{\ell\cdot d\cdot\log(1/\epsilon)\cdot\log(n+1)}{c_{\atom}^2\cdot\delta^2\cdot\epsilon^\ell}
        \\
        +
        2^{2\cdot\len(c_{\abs})}\cdot\bigl(\len(c_{\ZZ}) + \len(\epsilon) + \log(n+1)\bigr)
        \\
        +
        \frac{\ell\cdot d\cdot\log(n+1)}{c_{\atom}^2\cdot\delta^2}
        +
        \frac{\log(1/\epsilon)}{\epsilon^\ell}\cdot
        \log\left(\frac{n}{c_{\atom}\cdot\delta\cdot\epsilon^\ell}\right)
        \bigggr)
      \end{multlined}
      \\
      & \leq
      \begin{multlined}[t]
        O\bigggl(
        \frac{\ell\cdot d\cdot(\log(1/\epsilon))^2\cdot\log(n+1)}{c_{\atom}^2\cdot\delta^2\cdot\epsilon^\ell}
        +
        2^{2\cdot\len(c_{\abs})}\cdot\bigl(\len(c_{\ZZ}) + \len(\epsilon) + \log(n+1)\bigr)
        \\
        +
        \len(c_{\atom}) + \len(\delta)
        \bigggr).
      \end{multlined}
    \end{align*}

    \smallskip

    Finally, we analyze the time-complexity of Algorithm~\ref{alg:nonequi:atomspace}. First, the analyses of the parameter
    computations and rational approximations are exactly as in Algorithm~\ref{alg:nonequi:atom}; these account for a
    time-complexity of at most
    \begin{gather}\label{eq:nonequi:atomspace:param:complexity}
      \begin{multlined}
        O\bigggl(
        \len(c_{\atom})^2\cdot\len(\delta)^2\cdot\bigl(\log(d+1)\bigr)^2
        \\
        +
        \frac{\log(1/\epsilon)}{\epsilon^\ell}
        \cdot\ell
        \cdot
        2^{3\cdot\len(c_{\abs})}\cdot\bigl(\len(c_{\ZZ})+\len(\epsilon)\bigr)^3\cdot\bigl(\log(n+1)\bigr)^3
        \bigggr).
      \end{multlined}
    \end{gather}

    It remains to analyze the contribution of the call to Algorithm~\ref{alg:ext:atomspace} to the cost of the loop of
    line~\ref{alg:nonequi:atomspace:While} (it completely dominates all other costs of the loop, except for rational
    approximations, which we already accounted for). Letting $\widehat{T}(i)$ and $\widehat{U}(i)$ be the time- and $U$-oracle
    complexities of the call to Algorithm~\ref{alg:ext:atomspace} in the $i$th iteration,
    Proposition~\ref{prop:ext}\ref{prop:ext:atomspace} gives
    \begin{align*}
      \widehat{T}(i) & =  O(2^\ell\cdot m\cdot n^m),
      &
      \widehat{U}(i) & \leq (2^{\ell+1} - 1)\cdot n^{m+1}.
    \end{align*}

    We also recall that the time-complexity of Algorithm~\ref{alg:nonequi:atomspaceauxiliary} in the $i$th iteration is
    \begin{gather*}
      T(i-1)
      \leq
      O\bigl(1 + (i-1)^2\cdot m\cdot n^{i-1}\bigr).
    \end{gather*}

    Then the total time-complexity of these calls is upper bounded by
    \begin{align*}
      \sum_{i=1}^K \bigl(\widehat{T}(i) + \widehat{U}(i)\cdot T(i-1)\bigr)
      & \leq
      \sum_{i=1}^K O\Bigl(
      2^\ell\cdot m\cdot n^m + (2^{\ell+1}-1)\cdot n^{m+1}\cdot\bigl(1 + (i-1)^2\cdot m\cdot n^{i-1}\bigr)
      \Bigr)
      \\
      & \leq
      O\bigl(K^2\cdot 2^\ell\cdot m\cdot n^{m+K}\bigr)
      \\
      & \leq
      n^{O(d/(c_{\atom}^2\cdot\delta^2) + \epsilon^{-\ell}\cdot\log(1/\epsilon))}.
    \end{align*}

    Putting this together with~\eqref{eq:nonequi:atomspace:param:complexity}, the final time-complexity of
    Algorithm~\ref{alg:nonequi:atomspace} is upper bounded by
    \begin{multline*}
      n^{O(d/(c_{\atom}^2\cdot\delta^2) + \epsilon^{-\ell}\cdot\log(1/\epsilon))}
      +
      O\bigggl(
      \len(c_{\atom})^2\cdot\len(\delta)^2\cdot\bigl(\log(d+1)\bigr)^2
      \\
      +
      \frac{\log(1/\epsilon)}{\epsilon^\ell}
      \cdot\ell
      \cdot
      2^{3\cdot\len(c_{\abs})}\cdot\bigl(\len(c_{\ZZ})+\len(\epsilon)\bigr)^3\cdot\bigl(\log(n+1)\bigr)^3
      \bigggr).
    \end{multline*}
  \item[Item~\ref{thm:nonequi:goodspace}.] Just as Algorithm~\ref{alg:nonequi:gooddet} is essentially
    Algorithm~\ref{alg:nonequi:atom} but with $m\df 1$ (and an arbitrary $\delta\df 1/2^{\ell+2}$),
    Algorithm~\ref{alg:nonequi:goodspace} is essentially Algorithm~\ref{alg:nonequi:atomspace} with the same choice $m\df 1$ and
    $\delta\df 1/2^{\ell+2}$, and we use the fact that $\epsilon$-goodness is equivalent to $(\delta,\epsilon,1)$-atomicity
    provided $\delta < 1$.

    \begin{algorithm}[htbp]
      \caption{Computation algorithm of query-into-oracle model with space- and time-complexities
        \begin{gather*}
          O\left(
          \frac{\ell\cdot(\log(1/\epsilon))^2\cdot\log(n+1)}{\epsilon^\ell}
          +
          2^{2\cdot\len(c_{\abs})}\cdot\bigl(\len(c_{\ZZ}) + \len(\epsilon) + \log(n+1)\bigr)
          \right),
          \\
          n^{O(\epsilon^{-\ell}\cdot\log(1/\epsilon))}
          +
          O\left(
          \frac{\log(1/\epsilon)}{\epsilon^\ell}
          \cdot\ell
          \cdot
          2^{3\cdot\len(c_{\abs})}\cdot\bigl(\len(c_{\ZZ})+\len(\epsilon)\bigr)^3\cdot\bigl(\log(n+1)\bigr)^3
          \right).
        \end{gather*}
        respectively, that returns $(\widetilde{m},(\sigma^k,x^k,s^k,\widetilde{s}^k)_{k=1}^K)$ such that when given to
        Algorithm~\ref{alg:nonequi:atomspaceoracle} produces an oracle for an $\epsilon$-good partition of $G$ into
        $K\leq K_{\nonequi}(\ell,0,c_{\ZZ},c_{\abs},\epsilon)$ parts.}
      \label{alg:nonequi:goodspace}
      \DontPrintSemicolon
      \KwIn{Numbers $n,\ell\in\NN_+$ with $\ell\leq\log_2(n)$, $c_{\ZZ},c_{\abs}\in(0,1)\cap\QQ$,
        $\epsilon\in(0,\epsilon_{\nonequi}(\ell,c_{\ZZ},c_{\abs}))\cap\QQ$, and a query-oracle $\cO_G$ for a graph $G$ with
        $\lvert G\rvert=n$ and $\Lit(G)\leq\ell$. We further assume $n\geq n_{\nonequi}(\ell,0,c_{\ZZ},c_{\abs},\epsilon)$.}
      \KwOut{A tuple $(\widetilde{m},(\sigma^k,x^k,s^k,\widetilde{s}^k)_{k=1}^K)$, with $K\leq
        K_{\nonequi}(\ell,0,c_{\ZZ},c_{\abs},\epsilon)$, where $\widetilde{m} = 1$ and for each $k\in[K]$, we have
        $\sigma^k\in\{0,1\}^{<\ell+1}$, $x^k\in (V(G)^1)^{\lvert\sigma\rvert}$ and $s^k,\widetilde{s}^k\in[n]$, such that when
        given to Algorithm~\ref{alg:nonequi:atomspaceoracle} produces an oracle for an $\epsilon$-good partition of
        $G$ into $K$ parts.}
      $m\assign 1$\;
      $\widetilde{m}\assign 0$\tcp*[r]{$\widetilde{m} = \floor{\widetilde{\delta}\cdot m}$ for $\widetilde{\delta} < 1$}
      $c_{\UC}\assign 0$\;
      $\zeta\assign (1-c_{\ZZ})\cdot\epsilon$\;
      Compute a sufficiently precise rational approximation of $\widetilde{\epsilon}\df (1 - \zeta^{c_{\abs}})\cdot\zeta$ to
      compute $S\df\floor{(\zeta-\widetilde{\epsilon})\cdot n/(1-\widetilde{\epsilon})}$.\;
      \label{alg:nonequi:goodspace:rational:S}
      $K\assign 0$\;
      $r\assign n$\;
      \uWhile{$r > S$}{%
        \label{alg:nonequi:goodspace:While}
        $s_0\assign r$\;
        Compute sufficiently precise rational approximations of $\widetilde{\epsilon}\df (1 - \zeta^{c_{\abs}})\cdot\zeta$ and
        $\gamma\df (1 - c_{\UC})\cdot\widetilde{\epsilon}$ to compute $s_i\df\floor{\gamma\cdot s_{i-1}} + 1$ for every
        $i\in[\ell+1]$.\;
        \label{alg:nonequi:goodspace:rational:si}
        Let $\cO$ be Algorithm~\ref{alg:nonequi:atomspaceauxiliary} with parameters:
        \begin{gather*}
          \Bigl(n,\cO_G,\bigl(m,\widetilde{m},(\sigma^k,x^k,s^k)_{k=1}^K\bigr),u\Bigr)
          \df
          \Bigl(n,\cO_G,\bigl(m,\widetilde{m},(\sigma^k,x^k,s^k)_{k=1}^K\bigr),\place\Bigr).
        \end{gather*}
        \;
        Run Algorithm~\ref{alg:ext:atomspace} with parameters:
        \begin{gather*}
          (n,m,\ell,s_0,\ldots,s_\ell,s_{\ell+1},\widetilde{m},\cO_G,\cO_U)
          \df
          (n,m,\ell,s_0,\ldots,s_\ell,s_{\ell+1},\widetilde{m},\cO_G,\cO).
        \end{gather*}
        \;\label{alg:nonequi:goodspace:call}
        Let $(\sigma,(x_j)_{j=0}^{\lvert\sigma\rvert-1})$ be the tuple returned by the algorithm.\;
        \setcounter{algosplit}{\theAlgoLine}
      }
    \end{algorithm}

    \begin{algorithm}[htbp]
      \caption*{(continued).}
      \DontPrintSemicolon
      \setcounter{AlgoLine}{\thealgosplit}
      \let\oldnl\nl
      \let\nl\relax
      \vspace{-\baselineskip}\InvisibleBegin{
        \global\let\nl\oldnl 
        $K\assign K+1$\;
        $\sigma^K\assign\sigma$\;
        $x^K\assign (x_j)_{j=0}^{\lvert\sigma\rvert-1}$\;
        $s^K\assign s_{\lvert\sigma\rvert}$\;
        $r\assign r - s^K$\;
        Compute a sufficiently precise rational approximation of $\widetilde{\epsilon}\df (1 - \zeta^{c_{\abs}})\cdot\zeta$ to
        compute $\widetilde{s}^K\df\ceil{(\zeta-\widetilde{\epsilon})\cdot s^K/(1-\zeta)}$.\;
      }
      \Return{$(\sigma^k,x^k,s^k,\widetilde{s}^k)_{k=1}^K$}
    \end{algorithm}

    The complexity analysis of Algorithm~\ref{alg:nonequi:goodspace} is also analogous to that of
    Algorithm~\ref{alg:nonequi:atomspace} and yields the space-complexity bound
    \begin{align*}
      \MoveEqLeft
      \begin{multlined}[t]
        O\bigggl(
        \len(c_{\ZZ}) + \len(\epsilon)
        + \ell\cdot\frac{\log(1/\epsilon)}{\epsilon^\ell}\cdot m\cdot\log(n+1)
        \\
        +
        2^{2\cdot\len(c_{\abs})}\cdot\bigl(\len(c_{\ZZ}) + \len(\epsilon) + \log(n+1)\bigr)
        \\
        +
        \ell\cdot m\cdot\log(n+1)
        +
        K\cdot\bigl(\log(K+1) + \log(m+1) + \log(n+1)\bigr)
        \bigggr)
      \end{multlined}
      \\
      & \leq
      O\left(
      \frac{\ell\cdot(\log(1/\epsilon))^2\cdot\log(n+1)}{\epsilon^\ell}
      +
      2^{2\cdot\len(c_{\abs})}\cdot\bigl(\len(c_{\ZZ}) + \len(\epsilon) + \log(n+1)\bigr)
      \right),
    \end{align*}
    and the time-complexity bound
    \begin{align*}
      \MoveEqLeft
      \begin{multlined}[t]
        O\bigggl(
        \frac{\log(1/\epsilon)}{\epsilon^\ell}
        \cdot\ell
        \cdot
        2^{3\cdot\len(c_{\abs})}\cdot\bigl(\len(c_{\ZZ})+\len(\epsilon)\bigr)^3\cdot\bigl(\log(n+1)\bigr)^3
        \\
        +
        K^2\cdot 2^\ell\cdot m\cdot n^{m+K}\bigr)
      \end{multlined}
      \\
      & \leq
      n^{O(\epsilon^{-\ell}\cdot\log(1/\epsilon))}
      +
      O\left(
      \frac{\log(1/\epsilon)}{\epsilon^\ell}
      \cdot\ell
      \cdot
      2^{3\cdot\len(c_{\abs})}\cdot\bigl(\len(c_{\ZZ})+\len(\epsilon)\bigr)^3\cdot\bigl(\log(n+1)\bigr)^3
      \right).
      \qedhere
    \end{align*}
  \end{description}
\end{proof}

\begin{remark}\label{rmk:MCLV}
  Let us point out that Algorithm~\ref{alg:nonequi:good} is a randomized algorithm of Monte Carlo type (i.e., it has a fixed
  running time, but a positive probability of failure), which in the RAM model of computation (in which arithmetic operations
  cost $O(1)$) has running time $O_{\epsilon,\ell}(n\cdot\log(n))$. We also point out that it is not trivial to produce a Las
  Vegas randomized algorithm (i.e., one that always succeeds, but only has a bound on the expected running time) from
  Algorithm~\ref{alg:nonequi:good} with expected running time $O_{\epsilon,\ell}(n\cdot\log(n))$. This is because checking
  whether a set is $\epsilon$-good naively takes time $O_{\epsilon,\ell}(n^2)$, so the standard Monte Carlo to Las Vegas
  conversion trick of repeatedly running and checking the output would yield an expected running time of at least
  $O_{\epsilon,\ell}(n^2)$; at this complexity regime, we could instead use Algorithm~\ref{alg:nonequi:gooddet}, which in the
  RAM deterministic model is time $O_{\epsilon,\ell}(n^2)$.
\end{remark}

\begin{remark}[Choosing $c_{\atom}$ and $\delta$]\label{rmk:delta}
  Due to Lemma~\ref{lem:mono}, we know that $(\delta,\epsilon)$-excellence/atomicity implies
  $(\delta',\epsilon)$-excellence/atomicity when $\delta\geq\delta'$, which means that we are free to run any of the atomicity
  algorithms of Theorem~\ref{thm:nonequi} with a larger $\delta$ than the one desired. On the other hand, we also see that
  making $\delta$ larger also has another advantage: all bounds in algorithmic complexity (both space and time) get smaller.
  However, the algorithms do require that $(1+c_{\atom})\cdot\delta < 1/2^{\ell+1}$, so the question becomes how to optimize the
  choice of $(c_{\atom},\delta)$?

  All complexity bounds actually factor as decreasing functions of the product $c_{\atom}\cdot\delta$, so this amounts to
  the calculation
  \begin{gather*}
    \sup\left\{c_{\atom}\cdot\delta \;\middle\vert\;
    c_{\atom},\delta\in(0,1)\land
    (1+c_{\atom})\cdot\delta < \frac{1}{2^{\ell+1}}
    \right\}
    =
    \frac{1}{2^{\ell+2}}.
  \end{gather*}
  However, the supremum above is not attained at any valid $(c_{\atom},\delta)$, but we can get arbitrarily close by setting
  $c_{\atom}\df 1 - a$ for $a > 0$ small and $\delta\df 2^{-\ell-2}$ so that $c_{\atom}\cdot\delta = (1-a)/2^{\ell+2}$. This
  effectively means that in all complexity bounds we can replace $c_{\atom}\cdot\delta$ by $(1-o(1))\cdot 2^{-\ell-2}$.
\end{remark}


\section{Equitable partitions, using randomness}

\begin{proposition}\label{prop:randsubset}
  Let $r,s,d,n,m\in\NN_+$ with $r\leq s$, let $G$ be a graph with $\lvert G\rvert=n$, $\VC(G)\leq d$, let $W\subseteq V(G)$ be a
  set with $\lvert W\rvert=s$, let $\widetilde{\delta},\widetilde{\epsilon},\xi\in(0,1)$ with $\widetilde{\epsilon} < \xi$, let
  $\cH\subseteq\cP(W)$ be a family of subsets of $W$ such that $\lvert H\rvert\leq\widetilde{\epsilon}\cdot s$ and let
  $\rn{\widetilde{W}}$ be a uniformly at random subset of $W$ of size $r$. Then the following hold:
  \begin{enumerate}
  \item\label{prop:randsubset:conc} We have
    \begin{gather*}
      \PP\left[
        \forall H\in\cH,
        \lvert H\cap\rn{\widetilde{W}}\rvert
        \leq
        \xi\cdot r
        \right]
      \geq
      1 - \lvert\cH\rvert\cdot\exp\bigl(-2\cdot(\xi-\widetilde{\epsilon})^2\cdot r\bigr).
    \end{gather*}
  \item\label{prop:randsubset:good} (essentially from Cameron--Freer--Patel~\cite{AFP18}) If $W$ is $\widetilde{\epsilon}$-good
    in $G$, then
    \begin{gather*}
      \PP[\rn{\widetilde{W}}\text{ is $\xi$-good in $G$}]
      \geq
      1 - n\cdot\exp\bigl(-2\cdot(\xi-\widetilde{\epsilon})^2\cdot r\bigr).
    \end{gather*}
  \item\label{prop:randsubset:exc} (essentially from Cameron--Freer--Patel~\cite{AFP18}) If $\widetilde{\delta}<1/(d+1)$ and $W$
    is $(\widetilde{\delta},\widetilde{\epsilon})$-excellent in $G$, then
    \begin{gather*}
      \PP[\rn{\widetilde{W}}\text{ is $(\widetilde{\delta},\xi)$-excellent in $G$}]
      \geq
      1 - 2\cdot(1+s^d)\cdot\exp\bigl(-2\cdot(\xi-\widetilde{\epsilon})^2\cdot r\bigr).
    \end{gather*}
  \item\label{prop:randsubset:atom} If $W$ is a $(\widetilde{\delta},\widetilde{\epsilon},m)$-atom of $G$, then
    \begin{gather*}
      \PP[\rn{\widetilde{W}}\text{ is a $(\widetilde{\delta},\xi,m)$-atom of $G$}]
      \geq
      1 - n^m\cdot\exp\bigl(-2\cdot(\xi-\widetilde{\epsilon})^2\cdot r\bigr).
    \end{gather*}
  \end{enumerate}
\end{proposition}

\begin{proof}
  For item~\ref{prop:randsubset:conc}, for a fixed $H\in\cH$, we let
  \begin{gather*}
    \mu_H
    \df
    \EE[\lvert H\cap\rn{\widetilde{W}}\rvert]
    =
    \lvert H\rvert\cdot\frac{r}{s}
    \leq
    \widetilde{\epsilon}\cdot r,
  \end{gather*}
  so by a standard tail bound of the hypergeometric distribution (see e.g.~\cite{Ska13}), we get
  \begin{align*}
    \PP\left[
      \lvert H\cap\rn{\widetilde{W}}\rvert
      \geq
      \xi\cdot r
      \right]
    & \leq
    \PP\left[
      \lvert H\cap\rn{\widetilde{W}}\rvert
      \geq
      \mu_H + (\xi-\widetilde{\epsilon})\cdot r
      \right]
    \\
    & \leq
    \exp\bigl(-2\cdot(\xi-\widetilde{\epsilon})^2\cdot r\bigr).
  \end{align*}

  Applying a union bound to all elements of $\cH$ and taking the complementary event, we get
  \begin{gather*}
    \PP\left[
      \lvert H\cap\rn{\widetilde{W}}\rvert
      \leq
      \xi\cdot r
      \right]
    \geq
    1 - \lvert\cH\rvert\cdot\exp\bigl(-2\cdot(\xi-\widetilde{\epsilon})^2\cdot r\bigr).
  \end{gather*}

  \medskip

  For item~\ref{prop:randsubset:good}, since $W$ is $\widetilde{\epsilon}$-good in $G$, we know that for every $v\in V(G)$,
  there exists an $\widetilde{\epsilon}$-majority opinion $t_G^{\widetilde{\epsilon}}(v,W)$ of $W$ with respect to $v$ in $G$,
  which means that taking
  \begin{gather*}
    \cH \df \{W\cap N_G^{1-t_G^{\widetilde{\epsilon}}(v,W)}(v) \mid v\in V(G)\},
  \end{gather*}
  we have $\lvert H\rvert\leq\widetilde{\epsilon}\cdot s$ for every $H\in\cH$. On the other hand, it is clear that
  $\lvert\cH\rvert\leq\lvert G\rvert = n$, so applying item~\ref{prop:randsubset:conc}, it then suffices to show that if
  \begin{gather*}
    \lvert H\cap\rn{\widetilde{W}}\rvert
    \leq
    \xi\cdot r,
  \end{gather*}
  for every $H\in\cH$, then $\rn{\widetilde{W}}$ is $\xi$-good in $G$. But indeed, for $v\in V(G)$, we note
  that the set $H\df W\cap N_G^{1-t_G^{\widetilde{\epsilon}}(v,W)}(v)$ is in $\cH$, so the above implies
  \begin{gather*}
    \lvert N_G^{1-t_G^{\widetilde{\epsilon}}(v,W)}(v)\cap\rn{\widetilde{W}}\rvert
    \leq
    \xi\cdot r
    =
    \xi\cdot\lvert\rn{\widetilde{W}}\rvert,
  \end{gather*}
  i.e., $t_G^{\widetilde{\epsilon}}(v,W)$ is also a $\xi$-majority opinion of $\rn{\widetilde{W}}$ with respect to $v$ in $G$,
  so $\rn{\widetilde{W}}$ is $\xi$-good in $G$.

  \medskip

  Item~\ref{prop:randsubset:exc} has an analogous proof to that of item~\ref{prop:randsubset:good}, except that we need a few
  more ingredients (which were also present in Malliaris--Shelah~\cite{MS14} and Cameron--Freer--Patel~\cite{AFP18}): since $W$
  is $(\widetilde{\delta},\widetilde{\epsilon})$-excellent in $G$, for every $\widetilde{\delta}$-good set $U$, there exists a
  $(\widetilde{\delta},\widetilde{\epsilon})$-majority opinion $t_G^{\widetilde{\delta},\widetilde{\epsilon}}(W,U)$ of $W$ with
  respect to $U$ in $G$, which means that taking
  \begin{gather*}
    \cH
    \df
    \{W\cap N_{G,\widetilde{\delta}}^{1-t_G^{\widetilde{\delta},\widetilde{\epsilon}}(W,U)}(U)
    \mid U\text{ is $\widetilde{\delta}$-good in $G$}\},
  \end{gather*}
  where
  \begin{gather*}
    N_{G,\widetilde{\delta}}^b(U)
    \df
    \{v\in V(G) \mid t_G^{\widetilde{\delta}}(v,U) = b\}
  \end{gather*}
  we have $\lvert H\rvert\leq\widetilde{\epsilon}\cdot s$ for every $H\in\cH$.

  On the other hand, for
  \begin{gather*}
    \cH_{G,\widetilde{\delta}}
    \df
    N_{G,\widetilde{\delta}}^1(U)
    \mid U\text{ is $\widetilde{\delta}$-good in $G$}\},
  \end{gather*}
  Lemma~\ref{lem:VCgood} guarantees that $\VC(\cH_{G,\widetilde{\delta}})\leq d$ (as $\widetilde{\delta}<1/(d+1)$ and
  $\VC(G)\leq d$). Thus, we get
  \begin{gather*}
    \lvert\cH\rvert
    \leq
    2\cdot\lvert\cH_{G,\widetilde{\delta}}\rest_W\rvert
    \leq
    2\cdot\sum_{i=0}^d \binom{s}{i}
    \leq
    2\cdot(1+s^d),
  \end{gather*}
  where the first inequality follows from the Sauer--Shelah--Perles Lemma (Lemma~\ref{lem:SSP}).

  Applying item~\ref{prop:randsubset:conc}, it then suffices to show that if
  \begin{gather*}
    \lvert H\cap\rn{\widetilde{W}}\rvert
    \leq
    \xi\cdot r,
  \end{gather*}
  for every $H\in\cH$, then $\rn{\widetilde{W}}$ is $(\widetilde{\delta},\xi)$-excellent in $G$. But indeed,
  for a $\widetilde{\delta}$-good set $U$ in $G$, we note that the set $H\df W\cap
  N_{G,\widetilde{\delta}}^{1-t_G^{\widetilde{\delta},\widetilde{\epsilon}}(W,U)}(U)$ is in $\cH$, so the above implies
  \begin{gather*}
    \lvert N_{G,\widetilde{\delta}}^{1-t_G^{\widetilde{\delta},\widetilde{\epsilon}}(W,U)}(U)\cap\rn{\widetilde{W}}\rvert
    \leq
    \xi\cdot r
    =
    \xi\cdot\lvert\rn{\widetilde{W}}\rvert,
  \end{gather*}
  i.e., $t_G^{\widetilde{\delta},\widetilde{\epsilon}}(W,U)$ is also a $(\widetilde{\delta},\xi)$-majority opinion of
  $\rn{\widetilde{W}}$ with respect to $U$ in $G$, so $\rn{\widetilde{W}}$ is $(\widetilde{\delta},\xi)$-excellent in $G$.

  \medskip

  Finally, we prove item~\ref{prop:randsubset:atom}. Since $W$ is a $(\widetilde{\delta},\widetilde{\epsilon},m)$-atom of $G$,
  we know that it is not $(\widetilde{\delta},\widetilde{\epsilon})$-split by any $m$-tuple $x\in V(G)^m$ in $G$, which means
  that there exists $b_x\in\{0,1\}$ such that
  \begin{gather*}
    \lvert N_{G,\widetilde{\delta}}^{b_x}(x)\cap W\rvert\leq\widetilde{\epsilon}\cdot\lvert W\rvert,
  \end{gather*}
  where
  \begin{gather*}
    N_{G,\widetilde{\delta}}^b(x)
    \df
    \{v\in V(G) \mid t_G^{\widetilde{\delta}}(v,x) = b\}.
  \end{gather*}
  This means that if we set
  \begin{gather*}
    \cH \df \{W\cap N_{G,\widetilde{\delta}}^{b_x}(x) \mid x\in V(G)^m\},
  \end{gather*}
  we have $\lvert H\rvert\leq\widetilde{\epsilon}\cdot s$ for every $H\in\cH$.

  It is clear that $\lvert\cH\rvert\leq n^m$, so applying item~\ref{prop:randsubset:conc}, it then suffices to show that if
  \begin{gather*}
    \lvert H\cap\rn{\widetilde{W}}\rvert \leq \xi\cdot r,
  \end{gather*}
  for every $H\in\cH$, then $\rn{\widetilde{W}}$ is a $(\widetilde{\delta},\xi,m)$-atom of $G$. But indeed, for $x\in V(G)^m$,
  we note that the set $H\df W\cap N_{G,\widetilde{\delta}}^{b_x}(x)$ is in $\cH$, so the above implies
  \begin{gather*}
    \lvert N_{G,\widetilde{\delta}}^{b_x}(x)\cap\rn{\widetilde{W}}\rvert
    \leq
    \xi\cdot r
    =
    \xi\cdot\lvert\rn{\widetilde{W}}\rvert,
  \end{gather*}
  i.e., $x$ does not $(\widetilde{\delta},\xi)$-split $\rn{\widetilde{W}}$ in $G$, hence $\rn{\widetilde{W}}$ is a
  $(\widetilde{\delta},\xi,m)$-atom of $G$.
\end{proof}

\begin{lemma}[Base calculations for equitable partitions]\label{lem:equicalc}
  Let $\ell\in\NN_+$, let $c_{\UC},c_{\adj},c_{\Root}\in[0,1)$ and let $c_{\ZZ},c_{\samp},c_{\Sum}\in(0,1)$. Define
  \begin{align*}
    \widetilde{c}_{\UC} & \df \frac{1}{(1-c_{\UC})^\ell} - 1,
    &
    \widetilde{c}_{\adj} & \df \frac{1}{(1-c_{\adj})^\ell} - 1,
    &
    \widetilde{c}_{\samp} & \df \frac{1}{(1-c_{\samp})^\ell} - 1,
    \\
    \widetilde{c}_{\ZZ} & \df \frac{1}{(1-c_{\ZZ})^{\ell+1}} - 1,
    &
    \widetilde{c}_{\Root} & \df \frac{1}{(1-c_{\Root})^{\ell+1}} - 1,
  \end{align*}
  and for $\epsilon\in(0,1)$, let
  \begin{align*}
    \zeta
    & \df
    (1-c_{\ZZ})\cdot\epsilon.
    \\
    K
    & \df
    K_{\equi}(\ell,c_{\UC},c_{\ZZ},c_{\samp},c_{\adj},c_{\Root},c_{\Sum},\epsilon)
    \\
    & \df
    \ceil{
      (1+\widetilde{c}_{\UC})\cdot(1+\widetilde{c}_{\samp})\cdot(1+\widetilde{c}_{\adj})\cdot(1+\widetilde{c}_{\Root})
      \cdot(1+c_{\Sum})\cdot\zeta^{-\ell-1}
    }
    \\
    & =
    \ceil{
      (1+\widetilde{c}_{\UC})\cdot(1+\widetilde{c}_{\ZZ})\cdot(1+\widetilde{c}_{\samp})\cdot(1+\widetilde{c}_{\adj})
      \cdot(1+\widetilde{c}_{\Root})\cdot(1+c_{\Sum})\cdot\epsilon^{-\ell-1}
    },
  \end{align*}
  and let
  \begin{multline*}
    \epsilon_{\equi}(\ell,c_{\UC},c_{\ZZ},c_{\samp},c_{\adj},c_{\Root},c_{\Sum})
    \df
    \sup\bigggl\{\epsilon'\in(0,1)
    \;\bigggm\vert\;
    \forall\epsilon\in(0,\epsilon'),
    \\
    \Biggl(
    K\cdot(1-c_{\Root})\cdot\zeta
    -
    \frac{1}{(1-c_{\UC})\cdot(1-c_{\samp})\cdot(1-c_{\adj})}\cdot
    \sum_{t=3}^K t^{-1/\ell}
    \\
    >
    \bigl(
    (1-c_{\UC})\cdot(1-c_{\samp})\cdot(1-c_{\adj})\cdot(1-c_{\Root})\cdot\zeta
    \bigr)^{-\ell}
    \Biggr)
    \bigggr\}.
  \end{multline*}

  Suppose now that $\epsilon < \epsilon_{\equi}(\ell,c_{\UC},c_{\ZZ},c_{\samp},c_{\adj},c_{\Root},c_{\Sum})$ and define
  \begin{align*}
    \MoveEqLeft
    n_{\equi}(\ell,c_{\UC},c_{\ZZ},c_{\samp},c_{\adj},c_{\Root},c_{\Sum},\epsilon)
    \\
    & \df
    \begin{multlined}[t]
      \max\bigggl\{
      K\cdot\frac{1-\epsilon}{\epsilon-\zeta},
      K\cdot
      \bigggl(
      1
      +
      (K-2)\cdot
      \Biggl(K\cdot(1-c_{\Root})\cdot\zeta
      \\
      -
      \frac{1}{(1-c_{\UC})\cdot(1-c_{\samp})\cdot(1-c_{\adj})}\cdot
      \sum_{t=3}^K t^{-1/\ell}
      \\
      -
      \bigl(
      (1-c_{\UC})\cdot(1-c_{\samp})\cdot(1-c_{\adj})\cdot(1-c_{\Root})\cdot\zeta
      \bigr)^{-\ell}
      \Biggr)^{-1}
      \biggr)
      \biggr\}.
    \end{multlined}
  \end{align*}

  Finally, $n\in\NN_+$ with $n\geq n_{\equi}(\ell,c_{\UC},c_{\ZZ},c_{\samp},c_{\adj},c_{\Root},c_{\Sum},\epsilon)$, define
  \begin{gather*}
    r \df \Floor{\frac{n}{K}}
    \\
    n_0 \df n,
  \end{gather*}
  and make the following definitions inductively in $i\in[K]$:
  \begin{enumerate}[label={\arabic*.}]
  \item Given $n_{i-1}$, let $\xi_i$ satisfy
    \begin{gather}
      0 < \xi_i\leq (1-c_{\Root})\cdot\zeta,
      \label{eq:equicalc:xii1}
      \\
      0\leq -1 + \zeta - \xi_i + \gamma_i^\ell\cdot\frac{n_{i-1}}{r} \leq c_{\Root}\cdot\zeta,
      \label{eq:equicalc:xii2}
    \end{gather}
    where
    \begin{align*}
      \gamma_i & \df (1-c_{\UC})\cdot\widetilde{\epsilon}_i,
      &
      \widetilde{\epsilon}_i & \df (1-c_{\samp})\cdot\widetilde{\xi}_i,
      &
      \widetilde{\xi}_i & \df (1-c_{\adj})\cdot\xi_i.
    \end{align*}
  \item Given $n_{i-1}$ and $\xi_i$ as above (which also defines $\gamma_i$), let
    \begin{align*}
      r_i & \df \ceil{\gamma_i^\ell\cdot n_{i-1}},
      &
      n_i & \df n_{i-1} - r_i.
    \end{align*}
  \end{enumerate}

  Then the following hold for every
  $\epsilon\in(0,\epsilon_{\equi}(\ell,c_{\UC},c_{\ZZ},c_{\samp},c_{\adj},c_{\Root},c_{\Sum}))$, every $n\in\NN_+$ with $n\geq
  n_{\equi}(\ell,c_{\UC},c_{\ZZ},c_{\samp},c_{\adj},c_{\Root},c_{\Sum},\epsilon)$, and every $i\in[K]$:
  \begin{enumerate}
  \item\label{lem:equicalc:epsilonequi} We have
    \begin{align*}
      \MoveEqLeft
      \epsilon_{\equi}(\ell,c_{\UC},c_{\ZZ},c_{\samp},c_{\adj},c_{\Root},c_{\Sum})
      \\
      & \geq
      \begin{dcases*}
        \bigl(1+o_{c_{\Sum}\to 0,c_{\UC}\to 0,c_{\samp}\to 0,c_{\adj}\to 0,c_{\ZZ},c_{\Root}}(1)\bigr)\cdot
        \frac{2\cdot(1-c_{\ZZ})\cdot\ln(1/c_{\Sum})}{c_{\Sum}},
        & if $\ell=1$,
        \\
        \bigl(1 + o_{c_{\UC}\to 0, c_{\samp}\to 0, c_{\adj}\to 0}(1)\bigr)\cdot
        \frac{
          (1+c_{\Sum})^{\ell-1}\cdot(1-c_{\ZZ})
        }{
          (1-c_{\Root})^\ell\cdot(1-1/\ell)^\ell\cdot c_{\Sum}^\ell
        },
        & if $\ell\geq 2$.
      \end{dcases*}
    \end{align*}
  \item\label{lem:equicalc:nequi} We have
    \begin{align*}
      \MoveEqLeft
      n_{\equi}(\ell,c_{\UC},c_{\ZZ},c_{\samp},c_{\adj},c_{\Root},c_{\Sum},\epsilon)
      \\
      & \leq
      \bigl(1+o_{\epsilon\to 0,c_{\UC},c_{\ZZ},c_{\samp},c_{\adj},c_{\Root},c_{\Sum}}(1)\bigr)
      \max\left\{\frac{1-c_{\ZZ}}{c_{\ZZ}}, \frac{1+c_{\Sum}}{(1-c_{\Root})\cdot c_{\Sum}}\right\}
      \cdot
      \frac{K}{\zeta}
      \\
      & =
      \bigl(1+o_{\epsilon\to 0,c_{\UC},c_{\ZZ},c_{\samp},c_{\adj},c_{\Root},c_{\Sum}}(1)\bigr)
      \cdot
      \max\left\{\frac{1-c_{\ZZ}}{c_{\ZZ}}, \frac{1+c_{\Sum}}{(1-c_{\Root})\cdot c_{\Sum}}\right\}
      \cdot\frac{K}{(1-c_{\ZZ})\cdot\epsilon}.
    \end{align*}
  \item\label{lem:equicalc:xii} For every $i\in[K]$, a number $\xi_i$ satisfying~\eqref{eq:equicalc:xii1}
    and~\eqref{eq:equicalc:xii2} indeed exists. Furthermore, if $c_{\Root} > 0$, then one can find a rational $\xi_i$
    satisfying~\eqref{eq:equicalc:xii1} and~\eqref{eq:equicalc:xii2} and such that $\gamma_i$ is dyadic (i.e., an element of
    $\{a/2^b \mid a\in\ZZ, b\in\NN\}$) provided $c_{\UC}$, $c_{\samp}$ and $c_{\adj}$ are rationals. In fact, the middle term
    of~\eqref{eq:equicalc:xii2} is negative when evaluated at $\xi_i=0$ and is at least $c_{\Root}\cdot\zeta$ when evaluated at
    $\xi_i=(1-c_{\Root})\cdot\zeta$.
  \item\label{lem:equicalc:eval} For every $i\in[K]$, we have
    \begin{gather*}
      n_{i-1}
      \geq
      \bigl(
      (1-c_{\UC})\cdot(1-c_{\samp})\cdot(1-c_{\adj})\cdot(1-c_{\Root})\cdot\zeta
      \bigr)^{-\ell}\cdot r.
    \end{gather*}
  \item\label{lem:equicalc:ri} For every $i\in[K]$, we have
    \begin{gather*}
      \gamma_i^\ell\cdot n_{i-1} \leq r_i \leq r,
    \end{gather*}
  \item\label{lem:equicalc:xiibound} For every $i\in[K]$, we have
    \begin{gather*}
      \xi_i \leq \frac{(K-i+1)^{-1/\ell}}{(1+c_{\UC})\cdot(1-c_{\samp})\cdot(1-c_{\adj})}
    \end{gather*}
  \item\label{lem:equicalc:cont} For every $i\in[K]$, we have
    \begin{gather*}
      (\xi_i-1)\cdot\frac{r_i}{r} + 1\leq\zeta.
    \end{gather*}
  \item\label{lem:equicalc:ZZ} For every $n\in\NN_+$ with $n\geq
    n_{\equi}(\ell,c_{\UC},c_{\ZZ},c_{\samp},c_{\adj},c_{\Root},c_{\Sum},\epsilon)$, we have
    \begin{gather*}
      (\zeta - 1)\cdot\frac{r}{r+1} + 1\leq\epsilon.
    \end{gather*}
  \item\label{lem:equicalc:xiilowerbound} For every $i\in[K]$, we have
    \begin{gather*}
      \xi_i \geq \widetilde{\xi}_i \geq \widetilde{\epsilon}_i
      \geq
      \gamma_i \geq \left(\frac{(1 - \zeta)\cdot r}{n}\right)^{1/\ell}
    \end{gather*}
  \item\label{lem:equicalc:rilowerbound} For every $i\in[K]$, we have
    \begin{gather*}
      r_i \geq (1-\zeta)\cdot r.
    \end{gather*}
  \end{enumerate}
\end{lemma}

\begin{proof}
  For item~\ref{lem:equicalc:epsilonequi}, we consider two cases. First suppose $\ell=1$, so that the condition in the
  definition of $\epsilon_{\equi}(\ell,c_{\UC},c_{\ZZ},c_{\samp},c_{\adj},c_{\Root},c_{\Sum})$ becomes
  \begin{gather}\label{eq:equicalc:epsilonequiell=1}
    \begin{aligned}
      \MoveEqLeft
      K\cdot(1-c_{\Root})\cdot\zeta
      -
      \frac{1}{(1-c_{\UC})\cdot(1-c_{\samp})\cdot(1-c_{\adj})}\cdot
      \sum_{t=3}^K t^{-1}
      \\
      & >
      \bigl(
      (1-c_{\UC})\cdot(1-c_{\samp})\cdot(1-c_{\adj})\cdot(1-c_{\Root})\cdot\zeta
      \bigr)^{-1}.
    \end{aligned}
  \end{gather}
  Since
  \begin{gather*}
    \sum_{t=3}^K t^{-1}
    \leq
    \int_2^K t^{-1}\ dt
    =
    \ln(K) - \ln(2)
    \leq
    \ln(K),
  \end{gather*}
  the condition in~\eqref{eq:equicalc:epsilonequiell=1} is implied by
  \begin{align*}
    \MoveEqLeft
    K\cdot(1-c_{\Root})\cdot\zeta
    -
    \frac{\ln(K)}{(1-c_{\UC})\cdot(1-c_{\samp})\cdot(1-c_{\adj})}
    \\
    & >
    \bigl(
    (1-c_{\UC})\cdot(1-c_{\samp})\cdot(1-c_{\adj})\cdot(1-c_{\Root})\cdot\zeta
    \bigr)^{-1}.
  \end{align*}
  Noting that the left-hand side of the above is increasing in $K$ provided $c_{\UC}$, $c_{\samp}$ and $c_{\adj}$ are
  sufficiently close to $0$ and recalling that
  $K\df\ceil{(1+\widetilde{c}_{\UC})\cdot(1+\widetilde{c}_{\samp})\cdot(1+\widetilde{c}_{\adj})\cdot(1+\widetilde{c}_{\Root})
    \cdot(1+c_{\Sum})\cdot\zeta^{-2}}$, the above is implied by
  \begin{align*}
    \MoveEqLeft
    \begin{multlined}[t]
      (1+\widetilde{c}_{\UC})\cdot(1+\widetilde{c}_{\samp})\cdot(1+\widetilde{c}_{\adj})\cdot(1+\widetilde{c}_{\Root})
      \cdot(1+c_{\Sum})\cdot(1-c_{\Root})\cdot\zeta^{-1}
      \\
      -
      \frac{
        \ln((1+\widetilde{c}_{\UC})\cdot(1+\widetilde{c}_{\samp})\cdot(1+\widetilde{c}_{\adj})\cdot(1+\widetilde{c}_{\Root})
        \cdot(1+c_{\Sum})\cdot\zeta^{-2})
      }{
        (1-c_{\UC})\cdot(1-c_{\samp})\cdot(1-c_{\adj})
      }
    \end{multlined}
    \\
    & >
    \bigl(
    (1-c_{\UC})\cdot(1-c_{\samp})\cdot(1-c_{\adj})\cdot(1-c_{\Root})\cdot\zeta
    \bigr)^{-1},
  \end{align*}
  which itself is equivalent to
  \begin{align*}
    \MoveEqLeft
    \begin{multlined}[t]
      \frac{
        (1+\widetilde{c}_{\UC})\cdot(1+\widetilde{c}_{\samp})\cdot(1+\widetilde{c}_{\adj})\cdot(1+\widetilde{c}_{\Root})^{1/2}
        \cdot(1+c_{\Sum})
      }{
        \zeta
      }
      \\
      -
      (1+\widetilde{c}_{\UC})\cdot(1+\widetilde{c}_{\samp})\cdot(1+\widetilde{c}_{\adj})\cdot
      \\
      \ln((1+\widetilde{c}_{\UC})\cdot(1+\widetilde{c}_{\samp})\cdot(1+\widetilde{c}_{\adj})\cdot(1+\widetilde{c}_{\Root})^{1/2}
      \cdot(1+c_{\Sum})\cdot\zeta^{-2})
    \end{multlined}
    \\
    & >
    \frac{(1+\widetilde{c}_{\UC})\cdot(1+\widetilde{c}_{\samp})\cdot(1+\widetilde{c}_{\adj})\cdot(1+\widetilde{c}_{\Root})^{1/2}}{\zeta},
  \end{align*}
  hence equivalent to
  \begin{gather*}
    \frac{c_{\Sum}}{\zeta}
    >
    \ln((1+\widetilde{c}_{\UC})\cdot(1+\widetilde{c}_{\samp})\cdot(1+\widetilde{c}_{\adj})\cdot(1+\widetilde{c}_{\Root})
    \cdot(1+c_{\Sum})\cdot\zeta^{-2})
  \end{gather*}

  It is clear that we can satisfy the above with $\zeta$ satisfying:
  \begin{gather*}
    \zeta
    =
    \bigl(1+o_{c_{\Sum}\to 0,c_{\UC},c_{\ZZ},c_{\samp},c_{\adj,c_{\Root}}}(1)\bigr)\cdot
    \frac{2\cdot\ln(1/c_{\Sum})}{c_{\Sum}}.
  \end{gather*}
  and recalling that $\zeta\df(1-c_{\ZZ})\cdot\epsilon = \epsilon/(1+\widetilde{c}_{\ZZ})^{1/2}$ (and that we already made the
  assumption that $c_{\UC}$, $c_{\samp}$ and $c_{\adj}$ were small enough), we get
  \begin{align*}
    \MoveEqLeft
    \epsilon_{\equi}(\ell,c_{\UC},c_{\ZZ},c_{\samp},c_{\adj},c_{\Root},c_{\Sum})
    \\
    & \geq
    \bigl(1+o_{c_{\Sum}\to 0,c_{\UC}\to 0,c_{\samp}\to 0,c_{\adj}\to 0,c_{\ZZ},c_{\Root}}(1)\bigr)\cdot
    \frac{2\cdot\ln(1/c_{\Sum})}{(1+\widetilde{c}_{\ZZ})^{1/2}\cdot c_{\Sum}}
    \\
    & =
    \bigl(1+o_{c_{\Sum}\to 0,c_{\UC}\to 0,c_{\samp}\to 0,c_{\adj}\to 0,c_{\ZZ},c_{\Root}}(1)\bigr)\cdot
    \frac{2\cdot(1-c_{\ZZ})\cdot\ln(1/c_{\Sum})}{c_{\Sum}}.
  \end{align*}

  \smallskip

  We now consider the case $\ell\geq 2$. Since the condition of
  $\epsilon_{\equi}(\ell,c_{\UC},c_{\ZZ},c_{\samp},c_{\adj},c_{\Root},c_{\Sum})$ is
  \begin{gather}\label{eq:equicalc:epsilonequiellgeq2}
    \begin{aligned}
      \MoveEqLeft
      K\cdot(1-c_{\Root})\cdot\zeta
      -
      \frac{1}{(1-c_{\UC})\cdot(1-c_{\samp})\cdot(1-c_{\adj})}\cdot
      \sum_{t=3}^K t^{-1/\ell}
      \\
      & >
      \bigl(
      (1-c_{\UC})\cdot(1-c_{\samp})\cdot(1-c_{\adj})\cdot(1-c_{\Root})\cdot\zeta
      \bigr)^{-\ell}
    \end{aligned}
  \end{gather}
  and since
  \begin{gather*}
    \sum_{t=3}^K t^{-1/\ell}
    \leq
    \int_2^K t^{-1/\ell}\ dt
    =
    \frac{K^{1-1/\ell} - 2^{1-1/\ell}}{1-1/\ell}
    \leq
    \frac{K^{(\ell-1)/\ell}}{1-1/\ell},
  \end{gather*}
  the condition in~\eqref{eq:equicalc:epsilonequiellgeq2} is implied by
  \begin{align*}
    \MoveEqLeft
    K\cdot(1-c_{\Root})\cdot\zeta
    -
    \frac{1}{(1-c_{\UC})\cdot(1-c_{\samp})\cdot(1-c_{\adj})}\cdot
    \frac{K^{(\ell-1)/\ell}}{1-1/\ell},
    \\
    & >
    \bigl(
    (1-c_{\UC})\cdot(1-c_{\samp})\cdot(1-c_{\adj})\cdot(1-c_{\Root})\cdot\zeta
    \bigr)^{-\ell}.
  \end{align*}
  Noting that the left-hand side of the above is increasing in $K$ provided $c_{\UC}$, $c_{\samp}$ and $c_{\adj}$ are
  sufficiently close to $0$ and recalling that
  $K\df\ceil{(1+\widetilde{c}_{\UC})\cdot(1+\widetilde{c}_{\samp})\cdot(1+\widetilde{c}_{\adj})\cdot(1+\widetilde{c}_{\Root})
  \cdot(1+c_{\Sum})\cdot\zeta^{-\ell-1}}$, the above is implied by
  \begin{align*}
    \MoveEqLeft
    \begin{multlined}[t]
      (1+\widetilde{c}_{\UC})\cdot(1+\widetilde{c}_{\samp})\cdot(1+\widetilde{c}_{\adj})\cdot(1+\widetilde{c}_{\Root})
      \cdot(1+c_{\Sum})\cdot(1-c_{\Root})\cdot\zeta^{-\ell}
      \\
      -
      \frac{
        ((1+\widetilde{c}_{\UC})\cdot(1+\widetilde{c}_{\samp})\cdot(1+\widetilde{c}_{\adj})\cdot(1+\widetilde{c}_{\Root})
        \cdot(1+c_{\Sum})\cdot(1-c_{\Root})
        \cdot\zeta^{-\ell-1})^{(\ell-1)/\ell}
      }{
        (1-c_{\UC})\cdot(1-c_{\samp})\cdot(1-c_{\adj})\cdot(1-1/\ell)
      }
    \end{multlined}
    \\
    & >
    \bigl(
    (1-c_{\UC})\cdot(1-c_{\samp})\cdot(1-c_{\adj})\cdot(1-c_{\Root})\cdot\zeta
    \bigr)^{-\ell},
  \end{align*}
  which, using the definitions of $c_{\UC}$, $c_{\samp}$, $c_{\adj}$ and
  $c_{\Root}$, is equivalent to
  \begin{align*}
    \MoveEqLeft
    \begin{multlined}[t]
      \frac{
        (1+\widetilde{c}_{\UC})\cdot(1+\widetilde{c}_{\samp})\cdot(1+\widetilde{c}_{\adj})
        \cdot(1+\widetilde{c}_{\Root})^{1-1/(\ell+1)}\cdot(1+c_{\Sum})
      }{
        \zeta^\ell
      }
      \\
      -
      \frac{
        (1+\widetilde{c}_{\UC})\cdot(1+\widetilde{c}_{\samp})\cdot(1+\widetilde{c}_{\adj})\cdot(1+\widetilde{c}_{\Root})\cdot
        (1+c_{\Sum})^{1-1/\ell}\cdot\zeta^{-\ell+1/\ell}
      }{
        (1-1/\ell)
      }
    \end{multlined}
    \\
    & >
    \frac{
      (1+\widetilde{c}_{\UC})\cdot(1+\widetilde{c}_{\samp})\cdot(1+\widetilde{c}_{\adj})\cdot(1+\widetilde{c}_{\Root})^{1-1/(\ell+1)}
    }{
      \zeta^\ell
    }
  \end{align*}
  hence equivalent to
  \begin{gather*}
    \frac{c_{\Sum}}{\zeta^\ell}
    >
    \frac{(1+\widetilde{c}_{\Root})^{1/(\ell+1)}\cdot(1+c_{\Sum})^{1-1/\ell}}{(1-1/\ell)\cdot\zeta^{\ell-1/\ell}},
  \end{gather*}
  which is satisfied precisely when
  \begin{gather*}
    \zeta
    <
    \frac{(1+\widetilde{c}_{\Root})^{1 - 1/(\ell+1)}\cdot(1+c_{\Sum})^{\ell-1}}{(1-1/\ell)^\ell\cdot c_{\Sum}^\ell}.
  \end{gather*}
  Recalling that $\zeta\df(1-c_{\ZZ})\cdot\epsilon = \epsilon/(1+\widetilde{c}_{\ZZ})^{1/(\ell+1)}$ (and that we already made
  the assumption that $c_{\UC}$, $c_{\samp}$ and $c_{\adj}$ were small enough), we get
  \begin{align*}
    \MoveEqLeft
    \epsilon_{\equi}(\ell,c_{\UC},c_{\ZZ},c_{\samp},c_{\adj},c_{\Root},c_{\Sum})
    \\
    & \geq
    \bigl(1 + o_{c_{\UC}\to 0, c_{\samp}\to 0, c_{\adj}\to 0}(1)\bigr)\cdot
    \frac{
      (1+\widetilde{c}_{\Root})^{1 - 1/(\ell+1)}\cdot(1+c_{\Sum})^{\ell-1}
    }{
      (1+\widetilde{c}_{\ZZ})^{1/(\ell+1)}\cdot(1-1/\ell)^\ell\cdot c_{\Sum}^\ell
    }
    \\
    & =
    \bigl(1 + o_{c_{\UC}\to 0, c_{\samp}\to 0, c_{\adj}\to 0}(1)\bigr)\cdot
    \frac{
      (1+c_{\Sum})^{\ell-1}\cdot(1-c_{\ZZ})
    }{
      (1-c_{\Root})^\ell\cdot(1-1/\ell)^\ell\cdot c_{\Sum}^\ell
    }.
  \end{align*}

  \medskip

  Let us now prove item~\ref{lem:equicalc:nequi}. Since $\zeta=(1-c_{\ZZ})\cdot\epsilon$ and
  \begin{align*}
    K
    & =
    \ceil{(1+\widetilde{c}_{\UC})\cdot(1+\widetilde{c}_{\samp})\cdot(1+\widetilde{c}_{\adj})\cdot(1+\widetilde{c}_{\Root})
    \cdot(1+c_{\Sum})\cdot\zeta^{-\ell-1}}
    \\
    & =
    \ceil{(1+\widetilde{c}_{\UC})\cdot(1+\widetilde{c}_{\ZZ})\cdot(1+\widetilde{c}_{\samp})\cdot(1+\widetilde{c}_{\adj})
    \cdot(1+\widetilde{c}_{\Root})\cdot(1+c_{\Sum})\cdot\epsilon^{-\ell-1}},
  \end{align*}
  we get
  \begin{gather*}
    K\cdot\frac{1-\epsilon}{\epsilon-\zeta}
    =
    \bigl(1+o_{\epsilon\to 0,c_{\UC},c_{\ZZ},c_{\samp},c_{\adj},c_{\Root},c_{\Sum}}(1)\bigr)\cdot
    \frac{1-c_{\ZZ}}{c_{\ZZ}}\cdot\frac{K}{\zeta}
  \end{gather*}
  and
  \begin{align*}
    \MoveEqLeft
    \begin{multlined}[t]
      K\cdot
      \bigggl(
      1
      +
      (K-2)\cdot
      \Biggl(K\cdot(1-c_{\Root})\cdot\zeta
      -
      \frac{1}{(1-c_{\UC})\cdot(1-c_{\samp})\cdot(1-c_{\adj})}\cdot
      \sum_{t=3}^K t^{-1/\ell}
      \\
      -
      \bigl(
      (1-c_{\UC})\cdot(1-c_{\samp})\cdot(1-c_{\adj})\cdot(1-c_{\Root})\cdot\zeta
      \bigr)^{-\ell}
      \Biggr)^{-1}
      \bigggr)
    \end{multlined}
    \\
    & =
    \bigl(1+o_{\epsilon\to 0,c_{\UC},c_{\ZZ},c_{\samp},c_{\adj},c_{\Root},c_{\Sum}}(1)\bigr)\cdot
    \frac{K^2}{K\cdot(1-c_{\Root})\cdot\zeta - K\cdot(1-c_{\Root})\cdot\zeta/(1+c_{\Sum})}
    \\
    & =
    \bigl(1+o_{\epsilon\to 0,c_{\UC},c_{\ZZ},c_{\samp},c_{\adj},c_{\Root},c_{\Sum}}(1)\bigr)\cdot
    \frac{1+c_{\Sum}}{(1-c_{\Root})\cdot c_{\Sum}}\cdot\frac{K}{\zeta}.
  \end{align*}
  Now the result follows straightforwardly.

  \medskip
  
  We prove items~\ref{lem:equicalc:eval}, \ref{lem:equicalc:xii}, \ref{lem:equicalc:ri} and~\ref{lem:equicalc:xiibound} by
  induction in $i\in[K]$. The logic of the
  induction is as follows:
  \begin{enumerate}[label={\Roman*.}, ref={(\Roman*)}]
  \item\label{lem:equicalc:step} Items~\ref{lem:equicalc:xii} and~\ref{lem:equicalc:xiibound} for all $j\in[i-1]$ imply
    item~\ref{lem:equicalc:eval} for $i$.
  \item\label{lem:equicalc:def} Item~\ref{lem:equicalc:eval} for $i$ and item~\ref{lem:equicalc:xii} for all $j\in[i-1]$ imply
    item~\ref{lem:equicalc:xii} for $i$.
  \item\label{lem:equicalc:size} Item~\ref{lem:equicalc:xii} for all $j\in[i]$ implies item~\ref{lem:equicalc:ri} for $i$.
  \item\label{lem:equicalc:bound} Items~\ref{lem:equicalc:xii} and~\ref{lem:equicalc:ri} for all $j\in[i]$ imply
    item~\ref{lem:equicalc:xiibound} for $i$.
  \end{enumerate}
  In all cases above, the assumption that item~\ref{lem:equicalc:xii} holds for (at least) all $j\in[i-1]$ is simply because we
  need all previous $\xi_j$, etc.\ to be defined.

  The hardest item of the above is~\ref{lem:equicalc:step}, so we leave it for last.

  For item~\ref{lem:equicalc:def}, we want to show that a number $\xi_i$ satisfying conditions~\eqref{eq:equicalc:xii1}
  and~\eqref{eq:equicalc:xii2} exists (and if $c_{\Root} > 0$, then a rational $\xi_i$ exists that makes $\gamma_i$ dyadic). For
  this, let us expand the definition of $\gamma_i$ in the condition~\eqref{eq:equicalc:xii2}:
  \begin{gather*}
    0
    \leq
    -1 + \zeta - \xi_i + 
    \bigl(
      (1-c_{\UC})\cdot(1-c_{\samp})\cdot(1-c_{\adj})
    \cdot\xi_i\bigr)^\ell
    \cdot\frac{n_{i-1}}{r}
    \leq
    c_{\Root}\cdot\zeta,
  \end{gather*}
  and we note that the middle term is a continuous function of $\xi_i$, call it $f(\xi_i)$. We now evaluate $f$ at the two
  extremities of the interval of condition~\eqref{eq:equicalc:xii1}:
  \begin{align*}
    f(0)
    & =
    -1 + \zeta
    <
    0,
    \\
    f\bigl((1-c_{\Root})\cdot\zeta\bigr)
    & =
    c_{\Root}\cdot\zeta - 1
    +
    \bigl(
      (1-c_{\UC})\cdot(1-c_{\samp})\cdot(1-c_{\adj})\cdot(1-c_{\Root})
    \cdot\zeta\bigr)^\ell
    \cdot\frac{n_{i-1}}{r}.
  \end{align*}
  Since $f$ is continuous, the existence of $\xi_i$ satisfying~\eqref{eq:equicalc:xii1} and~\eqref{eq:equicalc:xii2} and making
  $\gamma_i$ dyadic provided $c_{\UC}$, $c_{\samp}$ and $c_{\adj}$ are rationals (which makes $\gamma_i$ be a rational multiple
  of $\xi_i$) will be a consequence of the Intermediate Value Theorem provided $f((1-c_{\Root})\cdot\zeta)\geq
  c_{\Root}\cdot\zeta$ (furthermore if $c_{\Root} > 0$, then we will be guaranteed that a rational $\xi_i$ exists), that is, it
  suffices to prove
  \begin{gather*}
    n_{i-1}
    \geq
    \bigl(
    (1-c_{\UC})\cdot(1-c_{\samp})\cdot(1-c_{\adj})\cdot(1-c_{\Root})\cdot\zeta
    \bigr)^{-\ell}\cdot r,
  \end{gather*}
  which is precisely item~\ref{lem:equicalc:eval} for $i$, which is our inductive hypothesis.

  \smallskip

  For item~\ref{lem:equicalc:size}, we need to show that $\gamma_i^\ell\cdot n_{i-1}\leq r_i\leq r$. Note that since $\xi_i$ is
  defined (by inductive hypothesis that item~\ref{lem:equicalc:xii} holds for $i$), from~\eqref{eq:equicalc:xii2}, we have
  \begin{gather*}
    \gamma_i^\ell\cdot n_{i-1}
    \leq
    r\cdot(1 + c_{\Root}\cdot\zeta - \zeta + \xi_i)
    \leq
    r,
  \end{gather*}
  where the second inequality follows from~\eqref{eq:equicalc:xii1}. Since $r$ is an integer, we get
  \begin{gather*}
    \gamma_i^\ell\cdot n_{i-1}
    \leq
    r_i
    =
    \ceil{\gamma_i^\ell\cdot n_{i-1}}
    \leq
    r,
  \end{gather*}
  as desired.

  \smallskip

  For item~\ref{lem:equicalc:bound}, we need to show that
  \begin{gather*}
    \xi_i \leq \frac{(K-i+1)^{-1/\ell}}{(1+c_{\UC})\cdot(1-c_{\samp})\cdot(1-c_{\adj})}.
  \end{gather*}
  Since our inductive hypothesis says that $r_j\leq r$ for every $j\in[i]$ and we have $n_j = n_{j-1} - r_j$, a simple induction
  (along with $n_0=n$) gives
  \begin{gather*}
    n_{i-1} \geq n - (i-1)\cdot r \geq K\cdot r - (i-1)\cdot r = (K-i+1)\cdot r.
  \end{gather*}
  In turn, using our inductive hypothesis that $\gamma_i^\ell\cdot n_{i-1}\leq r$, we get
  \begin{align*}
    \xi_i
    & =
    \frac{\gamma_i}{(1-c_{\UC})\cdot(1-c_{\samp})\cdot(1-c_{\adj})}
    \leq
    \frac{(r/n_{i-1})^{1/\ell}}{(1-c_{\UC})\cdot(1-c_{\samp})\cdot(1-c_{\adj})}
    \\
    & \leq
    \frac{(K-i+1)^{-1/\ell}}{(1-c_{\UC})\cdot(1-c_{\samp})\cdot(1-c_{\adj})},
  \end{align*}
  as desired.

  \smallskip

  We now prove the hardest item~\ref{lem:equicalc:step}: note that
  \begin{align*}
    n_{i-1}
    & =
    n - \sum_{j=1}^{i-2} r_j
    =
    n - \sum_{j=1}^{i-2}\ceil{\gamma_j^\ell\cdot n_{j-1}}
    \geq
    n - \sum_{j=1}^{i-2}(\gamma_j^\ell\cdot n_{j-1} + 1)
    \\
    & \geq
    K\cdot r - \sum_{j=1}^{i-2} \Bigl(r\cdot\bigl(1 + \xi_j - (1-c_{\Root})\cdot\zeta\bigr) + 1\Bigr)
    \\
    & =
    K\cdot r - (i-2) - (i-2)\cdot r\bigl(1 - (1-c_{\Root})\cdot\zeta\bigr) - r\cdot\sum_{j=1}^{i-2} \xi_j
    \\
    & \geq
    K\cdot r - (i-2) - (i-2)\cdot r\bigl(1 - (1-c_{\Root})\cdot\zeta\bigr)
    - r\cdot\sum_{j=1}^{i-2}
    \frac{(K-i+2)^{-1/\ell}}{(1-c_{\UC})\cdot(1-c_{\samp})\cdot(1-c_{\adj})}
    \\
    & \geq
    K\cdot r\cdot(1-c_{\Root})\cdot\zeta - (K-2)
    - r\cdot\sum_{j=1}^{K-2}
    \frac{(K-j+1)^{-1/\ell}}{(1-c_{\UC})\cdot(1-c_{\samp})\cdot(1-c_{\adj})}
    \\
    & =
    r\cdot\left(
    K\cdot(1-c_{\Root})\cdot\zeta
    -
    \frac{K-2}{r}
    - 
    \frac{1}{(1-c_{\UC})\cdot(1-c_{\samp})\cdot(1-c_{\adj})}\cdot
    \sum_{t=3}^K t^{-1/\ell}
    \right)
    \\
    & \geq
    r\cdot\left(
    K\cdot(1-c_{\Root})\cdot\zeta
    -
    \frac{K-2}{n/K-1}
    - 
    \frac{1}{(1-c_{\UC})\cdot(1-c_{\samp})\cdot(1-c_{\adj})}\cdot
    \sum_{t=3}^K t^{-1/\ell}
    \right),
  \end{align*}
  where the second inequality follows since $r = \floor{n/K}\leq n/K$ and from the condition~\eqref{eq:equicalc:xii2} for all
  $\xi_j$ with $j\in[i-1]$ (from our inductive hypothesis of item~\ref{lem:equicalc:xii}), the third inequality follows from
  item~\ref{lem:equicalc:xiibound} (from our inductive hypothesis), the fourth inequality follows since $i\leq K$ and the last
  inequality follows since $r = \floor{n/K}\geq n/K-1$.

  Since our goal is to prove that
  \begin{gather*}
    n_{i-1}
    \geq
    \bigl(
    (1-c_{\UC})\cdot(1-c_{\samp})\cdot(1-c_{\adj})\cdot(1-c_{\Root})\cdot\zeta
    \bigr)^{-\ell}\cdot r,
  \end{gather*}
  we see that it suffices to show
  \begin{align*}
    \MoveEqLeft
    K\cdot(1-c_{\Root})\cdot\zeta
    -
    \frac{K-2}{n/K-1}
    - 
    \frac{1}{(1-c_{\UC})\cdot(1-c_{\samp})\cdot(1-c_{\adj})}\cdot
    \sum_{t=3}^K t^{-1/\ell}
    \\
    & \geq
    \bigl(
    (1-c_{\UC})\cdot(1-c_{\samp})\cdot(1-c_{\adj})\cdot(1-c_{\Root})\cdot\zeta
    \bigr)^{-\ell},
  \end{align*}
  which itself is equivalent to
  \begin{multline*}
    n
    \geq
    K\cdot
    \bigggl(
    1
    +
    (K-2)\cdot
    \Biggl(K\cdot(1-c_{\Root})\cdot\zeta
    -
    \frac{1}{(1-c_{\UC})\cdot(1-c_{\samp})\cdot(1-c_{\adj})}\cdot
    \sum_{t=3}^K t^{-1/\ell}
    \\
    -
    \bigl(
    (1-c_{\UC})\cdot(1-c_{\samp})\cdot(1-c_{\adj})\cdot(1-c_{\Root})\cdot\zeta
    \bigr)^{-\ell}
    \Biggr)^{-1}
    \bigggr),
  \end{multline*}
  which in turn follows from $n\geq n_{\equi}(\ell,c_{\UC},c_{\ZZ},c_{\samp},c_{\adj},c_{\Root},c_{\Sum},\epsilon)$.
  
  \medskip

  We now prove item~\ref{lem:equicalc:cont}. Note that
  \begin{align*}
    (\xi_i-1)\cdot\frac{r_i}{r} + 1
    & \leq
    \xi_i + 1 - \frac{r_i}{r}
    =
    \xi_i + 1 - \frac{\ceil{\gamma_i^\ell\cdot n_{i-1}}}{r}
    \\
    & \leq
    \xi_i + 1 - \frac{\gamma_i^\ell\cdot n_{i-1}}{r}
    \leq
    \zeta
  \end{align*}
  where the first inequality follows since $r_i\leq r$ from item~\ref{lem:equicalc:ri} and the last inequality follows from the
  first inequality of~\eqref{eq:equicalc:xii2}.

  \medskip

  Let us now prove item~\ref{lem:equicalc:ZZ}. We note that
  \begin{gather*}
    (\zeta - 1)\cdot\frac{r}{r+1} + 1 \leq \epsilon.
  \end{gather*}
  is equivalent to
  \begin{gather*}
    1 - \frac{1}{r+1} \geq \frac{1-\epsilon}{1-\zeta},
  \end{gather*}
  which in turn is equivalent to
  \begin{gather*}
    r \geq \frac{1-\epsilon}{\epsilon-\zeta}.
  \end{gather*}
  Since $r = \floor{n/K}$, the above is a consequence of
  \begin{gather*}
    n
    \geq
    K\cdot\frac{1-\epsilon}{\epsilon-\zeta},
  \end{gather*}
  which in turn follows from $n\geq n_{\equi}(\ell,c_{\UC},c_{\ZZ},c_{\samp},c_{\adj},c_{\Root},c_{\Sum},\epsilon)$.

  \medskip

  We now prove item~\ref{lem:equicalc:xiilowerbound}. From the definitions of $\gamma_i$, $\widetilde{\epsilon}_i$ and
  $\widetilde{\xi}_i$, we have
  \begin{align*}
    \xi_i & \geq \widetilde{\xi}_i \geq \widetilde{\epsilon}_i
    \geq
    \gamma_i
    \\
    & \geq
    \frac{(1 - \zeta + \xi_i)\cdot r}{n_{i-1}}
    \geq
    \frac{(1 - \zeta)\cdot r}{n},
  \end{align*}
  where the penultimate inequality follows from~\eqref{eq:equicalc:xii2} and the last inequality follows since $n_{i-1}\leq n$
  and $\xi_i\geq 0$.

  \medskip

  Finally, we prove item~\ref{lem:equicalc:rilowerbound}. From item~\ref{lem:equicalc:cont} (and the fact that $\xi_i < 1$), we have
  \begin{gather*}
    r_i
    \geq
    \frac{1 - \zeta}{1-\xi_i}\cdot r
    \geq
    (1 - \zeta)\cdot r.
    \qedhere
  \end{gather*}
\end{proof}

\begin{remark}\label{rmk:equicalc}
  It is straightforward to check that the constants
  $\widetilde{c}_{\UC},\widetilde{c}_{\adj},\widetilde{c}_{\samp},\widetilde{c}_{\ZZ},\widetilde{c}_{\Root}$ of
  Lemma~\ref{lem:nonequicalc} are bijective increasing functions of the respective constants
  $c_{\UC},c_{\adj},c_{\samp},c_{\ZZ},c_{\Root}$ and their inverses are given by
  \begin{align*}
    c_{\UC} & \df 1 - \frac{1}{(1+\widetilde{c}_{\UC})^{1/\ell}},
    &
    c_{\adj} & \df 1 - \frac{1}{(1+\widetilde{c}_{\adj})^{1/\ell}},
    &
    c_{\samp} & \df 1 - \frac{1}{(1+\widetilde{c}_{\samp})^{1/\ell}},
    \\
    c_{\ZZ} & \df 1 - \frac{1}{(1+\widetilde{c}_{\ZZ})^{1/(\ell+1)}},
    &
    c_{\Root} & \df 1 - \frac{1}{(1+\widetilde{c}_{\Root})^{1/(\ell+1)}}.
  \end{align*}
\end{remark}

\begin{lemma}[Calculations for equitable partitions using randomness]\label{lem:equirandcalc}
  Let $d,\ell\in\NN_+$, let $c_{\UC},c_{\Root}\in[0,1)$ and let $c_{\atom},c_{\ZZ},c_{\samp},c_{\Sum},\rho_{\samp}\in(0,1)$.
  Set $c_{\adj}\df 0$ and as in Lemma~\ref{lem:equicalc}, define the parameters:
  \begin{align*}
    \widetilde{c}_{\UC} & \df \frac{1}{(1-c_{\UC})^\ell} - 1,
    &
    \widetilde{c}_{\adj} & \df \frac{1}{(1-c_{\adj})^\ell} - 1 = 0,
    &
    \widetilde{c}_{\samp} & \df \frac{1}{(1-c_{\samp})^\ell} - 1,
    \\
    \widetilde{c}_{\ZZ} & \df \frac{1}{(1-c_{\ZZ})^{\ell+1}} - 1,
    &
    \widetilde{c}_{\Root} & \df \frac{1}{(1-c_{\Root})^{\ell+1}} - 1,
  \end{align*}
  and for $\epsilon\in(0,1)$, let
  \begin{align*}
    \zeta
    & \df
    (1-c_{\ZZ})\cdot\epsilon.
    \\
    K
    \df
    K_{\equirand}(\ell,c_{\UC},c_{\ZZ},c_{\samp},c_{\Root},c_{\Sum},\epsilon)
    & \df
    K_{\equi}(\ell,c_{\UC},c_{\ZZ},c_{\samp},0,c_{\Root},c_{\Sum},\epsilon),
    \\
    \epsilon_{\equirand}(\ell,c_{\UC},c_{\ZZ},c_{\samp},c_{\Root},c_{\Sum})
    & \df
    \epsilon_{\equi}(\ell,c_{\UC},c_{\ZZ},c_{\samp},0,c_{\Root},c_{\Sum}).
  \end{align*}
  For $n\in\NN_+$ with $n\geq n_{\equi}(\ell,c_{\UC},c_{\ZZ},c_{\samp},0,c_{\Root},c_{\Sum},\epsilon)$ as in
  Lemma~\ref{lem:equicalc}, define the quantities $r$, $\xi_i$, $\widetilde{\xi}_i$, $\widetilde{\epsilon}_i$, $\gamma_i$, $r_i$
  and $n_i$ as in Lemma~\ref{lem:equicalc}.

  Given $\delta\in(0,1)$, let further
  \begin{gather*}
    m \df \Ceil{C\cdot\frac{d}{c_{\atom}^2\cdot\delta^2}},
  \end{gather*}
  where $C$ is the absolute constant of Lemma~\ref{lem:atomconst}.

  Suppose now that $\epsilon < \epsilon_{\equirand}(\ell,c_{\UC},c_{\ZZ},c_{\samp},c_{\Root},c_{\Sum})$ and let
  \begin{align*}
    \widetilde{\rho}_{\samp}
    & \df
    1 -  (1-\rho_{\samp})^{1/K},
    \\
    \MoveEqLeft
    n_{\equirand}^{\good}(\ell,c_{\UC},c_{\ZZ},c_{\samp},c_{\Root},c_{\Sum},\epsilon,\rho_{\samp})
    \\
    & \df
    \begin{multlined}[t]
      \min\bigggl\{n\in\NN_+ \;\bigggm\vert\;
      n\geq n_{\equi}(\ell,c_{\UC},c_{\ZZ},c_{\samp},0,c_{\Root},c_{\Sum},\epsilon)
      \\
      \land
      2\cdot c_{\samp}^2\cdot(1-\zeta)^{1+2/\ell}\cdot\frac{r^{1+2/\ell}}{n^{2/\ell}}
      \geq
      \ln(n) + \ln\left(\frac{1}{\widetilde{\rho}_{\samp}}\right)
      \bigggr\}.
    \end{multlined}
    \\
    \MoveEqLeft
    n_{\equirand}^{\atom}(\ell,d,c_{\atom},c_{\ZZ},c_{\samp},c_{\Root},c_{\Sum},\epsilon,\delta,\rho_{\samp})
    \\
    & \df
    \begin{multlined}[t]
      \min\bigggl\{n\in\NN_+ \;\bigggm\vert\;
      n\geq n_{\equi}(\ell,c_{\UC},c_{\ZZ},c_{\samp},0,c_{\Root},c_{\Sum},\epsilon)
      \\
      \land
      2\cdot c_{\samp}^2\cdot(1-\zeta)^{1+2/\ell}\cdot\frac{r^{1+2/\ell}}{n^{2/\ell}}
      \geq
      m\cdot\ln(n) + \ln\left(\frac{1}{\widetilde{\rho}_{\samp}}\right)
      \bigggr\}.
    \end{multlined}
  \end{align*}

  Then the following hold for every
  $\epsilon\in(0,\epsilon_{\equirand}(\ell,c_{\UC},c_{\ZZ},c_{\samp},c_{\Root},c_{\Sum}))$ and $n\in\NN_+$:
  \begin{enumerate}
  \item\label{lem:equirandcalc:good} If $n\geq
    n_{\equirand}^{\good}(\ell,c_{\UC},c_{\ZZ},c_{\samp},c_{\Root},c_{\Sum},\epsilon,\rho_{\samp})$, then for every
    $i\in[K]$, we have
    \begin{gather*}
      n\cdot\exp\bigl(-2\cdot(\xi_i-\widetilde{\epsilon}_i)^2\cdot r_i\bigr) \leq \widetilde{\rho}_{\samp}.
    \end{gather*}
  \item\label{lem:equirandcalc:ngood} We have
    \begin{align*}
      \MoveEqLeft
      n_{\equirand}^{\good}(\ell,c_{\UC},c_{\ZZ},c_{\samp},c_{\Root},c_{\Sum},\epsilon,\rho_{\samp})
      \\
      & \leq
      \bigl(1 + o_{\rho_{\samp}\to 0, \epsilon\to 0, c_{\UC},c_{\ZZ},c_{\samp},c_{\Root},c_{\Sum}, \ell}(1)\bigr)\cdot
      \frac{K^{1+2/\ell}\cdot\ln(K/\rho_{\samp})}{2\cdot c_{\samp}^2\cdot c_{\ZZ}\cdot c_{\Sum}}.
    \end{align*}
  \item\label{lem:equirandcalc:atom} If $n\geq
    n_{\equirand}^{\atom}(\ell,d,c_{\atom},c_{\ZZ},c_{\samp},c_{\Root},c_{\Sum},\epsilon,\delta,\rho_{\samp})$, then for every
    $i\in[K]$, we have
    \begin{gather*}
      n^m\cdot\exp\bigl(-2\cdot(\xi_i-\widetilde{\epsilon}_i)^2\cdot r_i\bigr) \leq \widetilde{\rho}_{\samp}.
    \end{gather*}
  \item\label{lem:equirandcalc:natom} We have
    \begin{align*}
      \MoveEqLeft
      n_{\equirand}^{\atom}(\ell,d,c_{\atom},c_{\ZZ},c_{\samp},c_{\Root},c_{\Sum},\epsilon,\delta,\rho_{\samp})
      \\
      & =
      \bigl(1 + o_{\rho_{\samp}\to 0, \epsilon\to 0, \delta, c_{\atom}, c_{\UC},c_{\ZZ},c_{\samp},c_{\Root},c_{\Sum}, d, \ell}(1)\bigr)\cdot
      \frac{K^{1+2/\ell}\cdot\ln(K/\rho_{\samp})}{2\cdot c_{\samp}^2\cdot c_{\ZZ}\cdot c_{\Sum}}.
    \end{align*}
  \end{enumerate}
\end{lemma}

\begin{proof}
  Items~\ref{lem:equirandcalc:good} and~\ref{lem:equirandcalc:atom} follow from the definitions of $n_{\equirand}^{\good}$ and
  $n_{\equirand}^{\atom}$ and fact that
  \begin{align*}
    \xi_i - \widetilde{\epsilon}_i
    & =
    c_{\samp}\cdot\xi_i
    \geq
    c_{\samp}\cdot\left(\frac{(1-\zeta)\cdot r}{n}\right)^{1/\ell},
    \\
    r_i & \geq (1-\zeta)\cdot r,
  \end{align*}
  where the inequalities follow from Lemma~\ref{lem:equicalc}, items~\ref{lem:equicalc:xiilowerbound}
  and~\ref{lem:equicalc:rilowerbound}, respectively.

  \medskip

  For item~\ref{lem:equirandcalc:ngood}, we first note the following asymptotic behaviors:
  \begin{gather}\label{eq:equirandcalc:asymp}
    \begin{aligned}
      (1-\zeta)^{1+2/\ell}
      & =
      \bigl(1 + o_{\epsilon\to 0, c_{\ZZ}, \ell}(1)\bigr),
      \\
      r
      & \df
      \Floor{\frac{n}{K}}
      =
      \bigl(1 + o_{n\to\infty, \epsilon, c_{\UC},c_{\ZZ},c_{\samp},c_{\Root},c_{\Sum}, \ell}(1)\bigr)\cdot
      \frac{n}{K},
      \\
      \widetilde{\rho}_{\samp}
      & \df
      1 -  (1-\rho_{\samp})^{1/K}
      =
      \bigl(1 + o_{\rho_{\samp}\to 0, \epsilon, c_{\UC},c_{\ZZ},c_{\samp},c_{\Root},c_{\Sum}, \ell}(1)\bigr)\cdot
      \frac{\rho_{\samp}}{K}.
    \end{aligned}
  \end{gather}

  Thus, the second condition in the definition of $n_{\equirand}^{\good}$ is asymptotically equivalent to
  \begin{gather*}
    \bigl(1 + o_{n\to\infty, \rho_{\samp}\to 0, \epsilon\to 0, c_{\UC},c_{\ZZ},c_{\samp},c_{\Root},c_{\Sum}, \ell}(1)\bigr)\cdot
    \frac{2\cdot c_{\samp}^2\cdot n}{K^{1+2/\ell}}
    \geq
    \ln(n) + \ln\left(\frac{K}{\rho_{\samp}}\right)
  \end{gather*}
  and is therefore satisfied with
  \begin{align*}
    n
    & \leq
    \begin{multlined}[t]
      \bigl(1 + o_{\rho_{\samp}\to 0, \epsilon\to 0, c_{\UC},c_{\ZZ},c_{\samp},c_{\Root},c_{\Sum}, \ell}(1)\bigr)\cdot
      \frac{K^{1+2/\ell}}{2\cdot c_{\samp}^2}
      \\
      \cdot\left(
      2\cdot\ln\left(\frac{1}{c_{\samp}}\right)
      +
      \left(1 + \frac{2}{\ell} + \ln\left(\frac{1}{\rho_{\samp}}\right)\right)\cdot\ln(K)
      \right)
    \end{multlined}
    \\
    & \leq
    \bigl(1 + o_{\rho_{\samp}\to 0, \epsilon\to 0, c_{\UC},c_{\ZZ},c_{\samp},c_{\Root},c_{\Sum}, \ell}(1)\bigr)\cdot
    \frac{K^{1+2/\ell}\cdot\ln(K/\rho_{\samp})}{2\cdot c_{\samp}^2}.
  \end{align*}

  Recalling from Lemma~\ref{lem:equicalc}\ref{lem:equicalc:nequi} that
  \begin{align*}
    \MoveEqLeft
    n_{\equi}(\ell,c_{\UC},c_{\ZZ},c_{\samp},c_{\adj},c_{\Root},c_{\Sum},\epsilon)
    \\
    & \leq
    \bigl(1+o_{\epsilon\to 0,c_{\UC},c_{\ZZ},c_{\samp},c_{\adj},c_{\Root},c_{\Sum}}(1)\bigr)
    \cdot
    \max\left\{\frac{1-c_{\ZZ}}{c_{\ZZ}}, \frac{1+c_{\Sum}}{(1-c_{\Root})\cdot c_{\Sum}}\right\}
    \cdot\frac{K}{(1-c_{\ZZ})\cdot\epsilon}
    \\
    & \leq
    \bigl(1+o_{\epsilon\to 0,c_{\UC},c_{\ZZ},c_{\samp},c_{\adj},c_{\Root},c_{\Sum}}(1)\bigr)
    \cdot
    \max\left\{\frac{1-c_{\ZZ}}{c_{\ZZ}}, \frac{1+c_{\Sum}}{c_{\Sum}}\right\}
    \cdot K^{1+1/(\ell+1)},
  \end{align*}
  we get
  \begin{align*}
    \MoveEqLeft
    n_{\equirand}^{\good}(\ell,c_{\UC},c_{\ZZ},c_{\samp},c_{\Root},c_{\Sum},\epsilon,\rho_{\samp})
    \\
    & \leq
    \bigl(1 + o_{\rho_{\samp}\to 0, \epsilon\to 0, c_{\UC},c_{\ZZ},c_{\samp},c_{\Root},c_{\Sum}, \ell}(1)\bigr)\cdot
    \frac{K^{1+2/\ell}\cdot\ln(K/\rho_{\samp})}{2\cdot c_{\samp}^2\cdot c_{\ZZ}\cdot c_{\Sum}}.
  \end{align*}

  \medskip

  Finally, for item~\ref{lem:equirandcalc:natom}, using again the asymptotic expressions of~\eqref{eq:equirandcalc:asymp} and
  the fact that $m = \ceil{C\cdot d/(c_{\atom}^2\cdot\delta^2)}$, we note that the second condition in the definition of
  $n_{\equirand}^{\atom}$ is asymptotically equivalent to
  \begin{gather*}
    \bigl(1 + o_{n\to\infty, \rho_{\samp}\to 0, \epsilon\to 0, c_{\UC},c_{\ZZ},c_{\samp},c_{\Root},c_{\Sum}, \ell}(1)\bigr)\cdot
    \frac{2\cdot c_{\samp}^2\cdot n}{K^{1+2/\ell}}
    \geq
    \frac{C\cdot d\cdot\ln(n)}{c_{\atom}^2\cdot\delta^2}
    + \ln\left(\frac{K}{\rho_{\samp}}\right)
  \end{gather*}
  and is therefore satisfied with
  \begin{align*}
    n
    & \leq
    \begin{multlined}[t]
      \bigl(1 + o_{\rho_{\samp}\to 0, \epsilon\to 0, c_{\UC},c_{\ZZ},c_{\samp},c_{\Root},c_{\Sum}, \ell}(1)\bigr)
      \cdot
      \frac{K^{1+2/\ell}}{2\cdot c_{\samp}^2}
      \\
      \cdot\left(
      \frac{2\cdot C\cdot d}{c_{\atom}^2\cdot\delta^2}\cdot\ln\left(\frac{1}{c_{\samp}}\right)
      +
      \left(\frac{C\cdot d}{c_{\atom}^2\cdot\delta^2}\cdot\ln\left(\frac{1}{c_{\samp}}\right)
      \cdot\left(1 + \frac{2}{\ell}\right)
      +
      \ln\left(\frac{1}{\rho_{\samp}}\right)
      \right)\cdot\ln(K)
      \right)
    \end{multlined}
    \\
    & \leq
    \bigl(1 + o_{\rho_{\samp}\to 0, \epsilon\to 0, \delta, c_{\atom}, c_{\UC},c_{\ZZ},c_{\samp},c_{\Root},c_{\Sum}, d, \ell}(1)\bigr)\cdot
    \frac{K^{1+2/\ell}\cdot\ln(K/\rho_{\samp})}{2\cdot c_{\samp}^2}.
  \end{align*}

  Thus, accounting for the asymptotic behavior of $n_{\equi}$ again, we conclude that
  \begin{align*}
    \MoveEqLeft
    n_{\equirand}^{\good}(\ell,d,c_{\UC},c_{\ZZ},c_{\samp},c_{\Root},c_{\Sum},\epsilon,\delta,\rho_{\samp})
    \\
    & \leq
    \bigl(1 + o_{\rho_{\samp}\to 0, \epsilon\to 0, \delta, c_{\atom}, c_{\UC},c_{\ZZ},c_{\samp},c_{\Root},c_{\Sum}, d, \ell}(1)\bigr)\cdot
    \frac{K^{1+2/\ell}\cdot\ln(K/\rho_{\samp})}{2\cdot c_{\samp}^2\cdot c_{\ZZ}\cdot c_{\Sum}}.
    \qedhere
  \end{align*}
\end{proof}

\begin{theorem}[Equitable partitions with randomness]\label{thm:equirand}
  Let $\ell,d\in\NN_+$ with $d\leq\ell$, let $\delta,\epsilon,\rho_{\samp}\in(0,1)$, let $c_{\atom},\rho_{\UC}\in(0,1)\cap\QQ$
  be rationals, let $c_{\Root},c_{\UC}\in[0,1)\cap\QQ$ be rational and let $G$ be a graph with $\lvert G\rvert=n$,
  $\Lit(G)\leq\ell$ and $\VC(G)\leq d$.

  Given further $c_{\ZZ},c_{\samp},c_{\Sum}\in(0,1)$, let
  \begin{align*}
    \MoveEqLeft
    \epsilon_{\equirand}(\ell,c_{\UC},c_{\ZZ},c_{\samp},c_{\Root},c_{\Sum})
    \\
    & =
    \epsilon_{\equi}(\ell,c_{\UC},c_{\ZZ},c_{\samp},0,c_{\Root},c_{\Sum})
    \\
    & \geq
    \begin{dcases*}
      \bigl(1+o_{c_{\Sum}\to 0,c_{\UC}\to 0,c_{\samp}\to 0,c_{\adj}\to 0,c_{\ZZ},c_{\Root}}(1)\bigr)\cdot
      \frac{2\cdot(1-c_{\ZZ})\cdot\ln(1/c_{\Sum})}{c_{\Sum}},
      & if $\ell=1$,
      \\
      \bigl(1 + o_{c_{\UC}\to 0, c_{\samp}\to 0, c_{\adj}\to 0}(1)\bigr)\cdot
      \frac{
        (1+c_{\Sum})^{\ell-1}\cdot(1-c_{\ZZ})
      }{
        (1-c_{\Root})^\ell\cdot(1-1/\ell)^\ell\cdot c_{\Sum}^\ell
      },
      & if $\ell\geq 2$.
    \end{dcases*}
  \end{align*}
  be as in Lemmas~\ref{lem:equicalc} and~\ref{lem:equirandcalc} and for
  $\epsilon\in(0,\epsilon_{\equirand}(\ell,c_{\UC},c_{\ZZ},c_{\samp},c_{\Root},c_{\Sum}))$, let
  \begin{align*}
    K
    & \df
    K_{\equirand}(\ell,c_{\UC},c_{\ZZ},c_{\samp},c_{\Root},c_{\Sum},\epsilon)
    \df
    K_{\equi}(\ell,c_{\UC},c_{\ZZ},c_{\samp},0,c_{\Root},c_{\Sum},\epsilon)
    \\
    & \df
    \ceil{
      (1+\widetilde{c}_{\UC})\cdot(1+\widetilde{c}_{\ZZ})\cdot(1+\widetilde{c}_{\samp})\cdot(1+\widetilde{c}_{\adj})
      \cdot(1+\widetilde{c}_{\Root})\cdot(1+c_{\Sum})\cdot\epsilon^{-\ell-1}
    },
    \\
    \MoveEqLeft
    n_{\equirand}^{\good}(\ell,c_{\UC},c_{\ZZ},c_{\samp},c_{\Root},c_{\Sum},\epsilon,\rho_{\samp})
    \\
    & \leq
    \bigl(1 + o_{\rho_{\samp}\to 0, \epsilon\to 0, c_{\UC},c_{\ZZ},c_{\samp},c_{\Root},c_{\Sum}, \ell}(1)\bigr)\cdot
    \frac{K^{1+2/\ell}\cdot\ln(K/\rho_{\samp})}{2\cdot c_{\samp}^2\cdot c_{\ZZ}\cdot c_{\Sum}}
    \\
    \MoveEqLeft
    n_{\equirand}^{\atom}(\ell,d,c_{\atom},c_{\ZZ},c_{\samp},c_{\Root},c_{\Sum},\epsilon,\delta,\rho_{\samp})
    \\
    & \leq
    \bigl(1 + o_{\rho_{\samp}\to 0, \epsilon\to 0, \delta, c_{\atom}, c_{\UC},c_{\ZZ},c_{\samp},c_{\Root},c_{\Sum}, d, \ell}(1)\bigr)\cdot
    \frac{K^{1+2/\ell}\cdot\ln(K/\rho_{\samp})}{2\cdot c_{\samp}^2\cdot c_{\ZZ}\cdot c_{\Sum}}
  \end{align*}
  also be as in Lemmas~\ref{lem:equicalc} and~\ref{lem:equirandcalc}. Let $n\in\NN_+$. For the items below about good
  equipartitions, we assume $n\geq
  n_{\equirand}^{\good}(\ell,c_{\UC},c_{\ZZ},c_{\samp},c_{\Root},c_{\Sum},\epsilon,\rho_{\samp})$ and for the items below
  about atomic equipartitions, we assume $n\geq
  n_{\equirand}^{\atom}(\ell,d,c_{\atom},c_{\ZZ},c_{\samp},c_{\Root},c_{\Sum},\epsilon,\delta,\rho_{\samp})$. Then the
  following hold:
  \begin{enumerate}
  \item\label{thm:equirand:good} If $c_{\UC}=c_{\Root}=0$, then there exists an $\epsilon$-good equipartition of $G$ into $K$
    parts.
  \item\label{thm:equirand:atom} If $c_{\UC}=c_{\Root}=0$ and $\delta < 1/2^{\ell+1}$, then there exists a
    $(\delta,\epsilon)$-atomic equipartition of $G$ into $K$ parts.
  \item\label{thm:equirand:goodalg} If $c_{\UC},c_{\Root} > 0$, then in the random-query model,
    Algorithm~\ref{alg:equirand:good}, with probability at least $(1-\rho_{\UC})\cdot(1-\rho_{\samp})$, computes an
    $\epsilon$-good equipartition of $G$ into $K$ parts in time
    \begin{multline*}
      O\bigggl(
      \frac{2^\ell\cdot\log(n/\rho_{\UC})\cdot n}{c_{\UC}^2\cdot\epsilon^{\ell+5}}
      \\
      +
      \frac{
        \len(\epsilon)^3\cdot\len(c_{\ZZ})^3\cdot\len(c_{\samp})^4\cdot\len(c_{\UC})^2
        \cdot\len(c_{\Root})^3\cdot\len(c_{\Sum})\cdot\len(\rho_{\UC})
      }{
        \epsilon^{\ell+1}
      }
      \\
      \cdot\ell^5\cdot\log(\ell+1)\cdot(\log(n+1))^2
      \bigggr).
    \end{multline*}
  \item\label{thm:equirand:atomalg} If $c_{\UC} = 0$, $c_{\Root} > 0$ and $(1+c_{\atom})\cdot\delta < 1/2^{\ell+1}$, then in the
    randomized model, Algorithm~\ref{alg:equirand:atom}, with probability at least $1-\rho_{\samp}$, computes a
    $(\delta,\epsilon)$-atomic equipartition of $G$ into $K$ parts in time
    \begin{multline*}
      \frac{\ell}{\epsilon^{\ell+1}}\cdot n^{O(d/(c_{\atom}^2\cdot\delta^2))}
      +
      \log(d+1)\cdot\len(c_{\atom})\cdot\len(\delta)
      \\
      +
      \frac{
        \ell^5\cdot\log(\ell+1)\cdot\len(\epsilon)^3\cdot\len(c_{\ZZ})^3\cdot\len(c_{\samp})^4
        \cdot\len(c_{\Root})^3\cdot\len(c_{\Sum})\cdot(\log(n+1))^2
      }{
        \epsilon^{\ell+1}
      }.
    \end{multline*}
  \end{enumerate}

  Furthermore, for
  \begin{align*}
    \epsilon'
    & \df
    \begin{dcases*}
      \frac{1 - \sqrt{1 - 8\cdot\epsilon\cdot(1-\epsilon)}}{2}, & if $\epsilon < (2-\sqrt{2})/4$,\\
      1/2, & otherwise,
    \end{dcases*}
    \\
    & =
    2\cdot\epsilon + O_{\epsilon\to 0}(\epsilon^2),
    \\
    \epsilon''
    & \df
    \frac{\epsilon'}{1-\epsilon}
    =
    \frac{1 - \sqrt{1 - 8\cdot\epsilon\cdot(1-\epsilon)}}{2\cdot(1-\epsilon)}
    =
    2\cdot\epsilon + O_{\epsilon\to 0}(\epsilon^2),
  \end{align*}
  all $(\delta,\epsilon)$-atomic partitions are $(\delta,\epsilon)$-excellent and $\epsilon$-good; and all $\epsilon$-good
  partitions of $G$ are in particular totally $\epsilon'$-homogeneous and $(\epsilon,\epsilon'')$-excellent.
\end{theorem}

\begin{proof}
  Just as in Theorem~\ref{thm:nonequi}, the final assertions on how to translate atomic into excellent, into good, into totally
  homogeneous and back into excellent follow from Lemmas~\ref{lem:atom->exc}
  and~\ref{lem:good->hom}\ref{lem:good->hom:partition}.

  However, this time, we cannot exactly derive the existence items~\ref{thm:equirand:good} and~\ref{thm:equirand:atom} 
  from their algorithmic counterparts, items~\ref{thm:equirand:goodalg} and~\ref{thm:equirand:atomalg}, as the existence items
  take $c_{\Root}=0$; this yields better bounds than taking the limit $c_{\Root}\to 0$. Nevertheless, we point out that all
  items will follow the same strategy:
  \begin{enumerate}[label={\arabic*.}]
  \item For items about atomicity, we will start by letting $m$ be given by Lemma~\ref{lem:atomconst}.
  \item For all items except for item~\ref{thm:equirand:good}, we will set $c_{\UC}\df 0$ (the excluded item already receives
    $c_{\UC}\in(0,1)\cap\QQ$ as input).
  \item We compute the parameters $K$, $r$, $\xi_i$, $\widetilde{\xi}_i$, $\widetilde{\epsilon}_i$, $\gamma_i$, $r_i$, $n_i$ and
    $\widetilde{\rho}_{\samp}$ via Lemmas~\ref{lem:equicalc} and~\ref{lem:equirandcalc}, then we will use
    Proposition~\ref{prop:ext} a total of $K$ times so that in the $i$th stage, we extract a $\widetilde{\epsilon}_i$-good set
    or a $((1+c_{\atom})\cdot\delta,\widetilde{\epsilon}_i,m)$-atom. The extraction process is guaranteed to return a set of
    relative size at least $\gamma_i^\ell$. By sampling uniformly at random a subset, we will normalize the size to be of
    relative size exactly $\gamma_i^\ell$, rounded up to the nearest integer. This will ensure that after the $i$th extraction,
    we have exactly $n_i$ vertices left and the $i$th extracted set is of size exactly $r_i\df\ceil{\gamma_i^\ell\cdot
      n_{i-1}}$. Proposition~\ref{prop:randsubset} and our choice of parameters will ensure that our goodness/atomicity quality
    deteriorates at most from $\widetilde{\epsilon}_i$ to $\xi_i=\widetilde{\xi}_i$.
  \item After all $K$ parts are extracted, we will take the remainder vertices and add them to the parts to make them all have
    size $r\df\floor{n/K}$. Our choice of parameters will also ensure that the goodness/atomicity quality of the $i$th part
    deteriorates at most from $\xi_i$ to $\zeta$.
  \item Since now all parts have size $r\df\floor{n/K}$, we can add at most one vertex of the remainder (of size $n\bmod K$) to
    each of the parts to get an equipartition. Our choice of parameters will ensure that the goodness/atomicity quality
    deteriorates at most from $\zeta$ to $\epsilon$.
  \item Finally, for the atomicity items, we actually obtained a partition into $((1+c_{\atom})\cdot\delta,\epsilon,m)$-atoms of
    $G$; our choice of $m$ as in Lemma~\ref{lem:atomconst} then ensures that these are $(\delta,\epsilon)$-atoms of $G$.
  \end{enumerate}

  Let us then prove the items (and we will prove them out of order for simplicity sake):
  \begin{description}[wide, itemsep={3ex}]
  \item[Item~\ref{thm:equirand:atom}.] Let $m$ be given by Lemma~\ref{lem:atomconst} and let
    $\widetilde{\delta}\df(1+c_{\atom})\cdot\delta$. We will follow the inductive construction of the parameters of
    Lemma~\ref{lem:equicalc} as we simultaneously construct $(\widetilde{\delta},\widetilde{\epsilon}_i,m)$-atoms $Q_i$ of $G$.
    Namely, we set $r\df\floor{n/K}$, $n_0\df n = \lvert G\rvert$, $R_0\df V(G)$ and for $i\in\NN_+$, given a subset
    $R_{i-1}\subseteq V(G)$ of size $n_{i-1}$, we do the following:
    \begin{enumerate}[label={\Roman*.}, ref={(\Roman*)}]
    \item Let $\xi_i$ be such that
      \begin{align}\label{eq:equirand:atom:xii}
        0 & < \xi_i\leq\zeta,
        &
        -1 + \zeta - \xi_i + \gamma_i^\ell\cdot\frac{n_{i-1}}{r} & = 0,
      \end{align}
      where
      \begin{align*}
        \gamma_i & \df \widetilde{\epsilon}_i \df (1-c_{\samp})\cdot\widetilde{\xi}_i,
        &
        \widetilde{\xi}_i & \df \xi_i.
      \end{align*}
      Lemma~\ref{lem:equicalc}\ref{lem:equicalc:xii} (and the fact that $c_{\UC}=c_{\Root}=c_{\adj}=0$) guarantees that such $\xi_i$
      exists.
    \item\label{thm:equirand:atom:extcall} By taking $s_0\df n_{i-1}$ and $s_j\df\floor{\gamma_i\cdot s_{j-1}} + 1$, we
      can apply Proposition~\ref{prop:ext}\ref{prop:ext:atom} and obtain a $(\widetilde{\delta},\widetilde{\epsilon}_i,m)$-atom
      $\widetilde{Q}_i\subseteq R_{i-1}$ of $G$ with $\lvert\widetilde{Q}_i\rvert\in\{s_0,\ldots,s_\ell\}$, which in particular
      implies that for $r_i\df\ceil{\gamma_i^\ell\cdot n_{i-1}}$, we have $\lvert\widetilde{Q}_i\rvert\geq r_i$.
    \item\label{thm:equirand:atom:randsubsetcall} We claim that there exists $Q_i\subseteq\widetilde{Q}_i$ of size exactly $r_i$
      that is a $(\widetilde{\delta},\xi_i,m)$-atom of $G$. To prove this, it suffices to prove that if $\rn{W}$ is a uniformly
      at random subset of $\widetilde{Q}_i$ of size $r_i$, then with positive probability $\rn{W}$ is a
      $(\widetilde{\delta},\xi_i,m)$-atom of $G$. For this, we apply
      Proposition~\ref{prop:randsubset}\ref{prop:randsubset:atom}:
      \begin{gather*}
        \PP[\rn{W}\text{ is a $(\widetilde{\delta},\xi_i,m)$-atom of $G$}]
        \geq
        1 - n^m\cdot\exp\bigl(-2\cdot(\xi_i-\widetilde{\epsilon}_i)^2\cdot r_i\bigr)
        \geq
        1 - \widetilde{\rho}_{\samp}
        >
        0,
      \end{gather*}
      where $\widetilde{\rho}_{\samp}\df 1 - (1-\rho_{\samp})^{-1/K}$ and the second inequality follows from
      Lemma~\ref{lem:equirandcalc}\ref{lem:equirandcalc:atom}. Thus such $Q_i$ exists.
    \item We repeat such extraction process a total of $K$ times (as per guaranteed to be possible by Lemma~\ref{lem:equicalc})
      to get $\{Q_i\mid i\in[K]\}$ pairwise disjoint with $\lvert Q_i\rvert=r_i=\ceil{\gamma_i^\ell\cdot n_{i-1}}$ and $Q_i$ a
      $(\widetilde{\delta},\xi_i,m)$-atom of $G$ for every $i\in[K]$. Since Lemma~\ref{lem:equicalc}\ref{lem:equicalc:ri}
      guarantees that $r_i\leq r = \floor{n/K}$, we know that from the remaining $n_K$ vertices, we can add enough vertices to
      each of the $Q_i$ to make them all have size $r$ (and remain pairwise disjoint). Let $Q'_i$ be the part obtained from
      $Q_i$ and observe since $Q_i$ is a $(\widetilde{\delta},\xi_i,m)$-atom of $G$, then
      Lemma~\ref{lem:upcont}\ref{lem:upcont:atom} guarantees that $Q'_i$ is a $(\widetilde{\delta},\epsilon_+,m)$-atom of $G$,
      where
      \begin{gather*}
        \epsilon_+
        \df
        (\xi_i-1)\cdot\frac{\lvert Q_i\rvert}{\lvert Q'_i\rvert} + 1
        =
        (\xi_i-1)\cdot\frac{r_i}{r} + 1
        \leq
        \zeta,
      \end{gather*}
      where the last inequality follows from Lemma~\ref{lem:equicalc}\ref{lem:equicalc:cont}.
    \item\label{thm:equirand:atom:final} We now know that $\{Q'_i \mid i\in[K]\}$ is a pairwise disjoint collection of
      $(\widetilde{\delta},\zeta,m)$-atoms of $G$, all of size $r\df\floor{n}{K}$, so for the remaining $n\bmod K$ vertices, we
      can at most one of these to each part to get an equipartition $\cP$ of $V(G)$. By
      Lemma~\ref{lem:upcont}\ref{lem:upcont:atom}, we know that the resulting partition is of
      $(\widetilde{\delta},\epsilon_{++},m)$-atoms of $G$, where
      \begin{gather*}
        \epsilon_{++}
        \df
        (\zeta-1)\cdot\frac{r}{r+1} + 1
        \leq
        \epsilon,
      \end{gather*}
      where the last inequality follows from Lemma~\ref{lem:equicalc:ZZ}, that is, these are
      $(\widetilde{\delta},\epsilon,m)$-atoms of $G$, which by Lemma~\ref{lem:atomconst} and our choice of
      $\widetilde{\delta}\df(1+c_{\atom})\cdot\delta$ and $m\df\ceil{C\cdot d/(c_{\atom}^2\cdot\delta^2)}$ implies that these
      are $(\delta,\epsilon)$-atoms of $G$, that is, $\cP$ is a $(\delta,\epsilon)$-atomic equipartition of $G$ into $K$ parts.
    \end{enumerate}
  \item[Item~\ref{thm:equirand:good}.] The construction of an $\epsilon$-good equipartition of $G$ is similar to the
    construction of a $(\delta,\epsilon)$-atomic equipartition of $G$ of item~\ref{thm:equirand:atom}, except for the following
    changes (which yield the slightly better bounds on the starting $n$, i.e., we can use $n_{\equirand}^{\good}$ instead of
    $n_{\equirand}^{\atom}$):
    \begin{itemize}
    \item In step~\ref{thm:equirand:atom:extcall}, we naturally invoke item~\ref{prop:ext:good} of Proposition~\ref{prop:ext}
      instead of item~\ref{prop:ext:atom}.
    \item In step~\ref{thm:equirand:atom:randsubsetcall}, we naturally invoke item~\ref{prop:randsubset:good} of
      Proposition~\ref{prop:randsubset} instead of item~\ref{prop:randsubset:atom} and item~\ref{lem:equirandcalc:good} of
      Lemma~\ref{lem:equirandcalc} instead of item~\ref{lem:equirandcalc:atom}.
    \item In the final step~\ref{thm:equirand:atom:final}, there is obviously no need to invoke Lemma~\ref{lem:atomconst}, as we
      already get an $\epsilon$-good equipartition into $K$ parts.
    \end{itemize}
  \item[Item~\ref{thm:equirand:atomalg}.] Algorithm~\ref{alg:equirand:atom} follows the same construction as
    item~\ref{thm:equirand:atom}, except that it cannot pick the number $\xi_i$ satisfying~\eqref{eq:equirand:atom:xii} as it
    might be an irrational number. This is why in the algorithmic items we have $c_{\Root} > 0$ so that now
    Lemma~\ref{lem:equicalc}\ref{lem:equicalc:xii} guarantees that we can find rational $\xi_i$ satisfying instead
    \begin{align*}
      0 & < \xi_i\leq (1-c_{\Root})\cdot\zeta,
      &
      0 & \leq -1 + \zeta - \xi_i + \gamma_i^\ell\cdot\frac{n_{i-1}}{r} \leq c_{\Root}\cdot\zeta,
    \end{align*}
    where
    \begin{align*}
      \gamma_i & \df \widetilde{\epsilon}_i \df (1-c_{\samp})\cdot\widetilde{\xi}_i,
      &
      \widetilde{\xi}_i & \df \xi_i.
    \end{align*}
    All other aspects of the proof of correctness of Algorithm~\ref{alg:equirand:atom} are exactly as in
    item~\ref{thm:equirand:atom}.

    \begin{algorithm}[htbp]
      \caption{Randomized model algorithm that returns a $(\delta,\epsilon)$-atomic equipartition of $G$ into
        $K_{\equirand}(\ell,0,c_{\ZZ},c_{\samp},c_{\Root},c_{\Sum},\epsilon)$ parts. The success probability is at least
        $1-\rho_{\samp}$ and the time-complexity is:
        \begin{multline*}
          \frac{\ell}{\epsilon^{\ell+1}}\cdot n^{O(d/(c_{\atom}^2\cdot\delta^2))}
          +
          \log(d+1)\cdot\len(c_{\atom})\cdot\len(\delta)
          \\
          +
          \frac{
            \ell^5\cdot\log(\ell+1)\cdot\len(\epsilon)^3\cdot\len(c_{\ZZ})^3\cdot\len(c_{\samp})^4
            \cdot\len(c_{\Root})^3\cdot\len(c_{\Sum})\cdot(\log(n+1))^2
          }{
            \epsilon^{\ell+1}
          }.
        \end{multline*}
      }
      \label{alg:equirand:atom}
      \DontPrintSemicolon
      \KwIn{Numbers $n,d,\ell\in\NN_+$ with $d\leq\ell\leq\log_2(n)$,
        $c_{\atom},c_{\ZZ},c_{\samp},c_{\Root},c_{\Sum},\delta\in(0,1)\cap\QQ$ with $(1+c_{\atom})\cdot\delta < 1/2^{\ell+1}$,
        and $\epsilon\in(0,\epsilon_{\equirand}(\ell,0,c_{\ZZ},c_{\samp},c_{\Root},c_{\Sum}))$, and a graph $G$ with
        $\lvert G\rvert=n$, $\Lit(G)\leq\ell$ and $\VC(G)\leq d$. We further assume that $n\geq
        n_{\equirand}^{\atom}(\ell,d,c_{\atom},c_{\ZZ},c_{\samp},c_{\Root},c_{\Sum},\epsilon,\delta,\rho_{\samp})$ for some
        choice of $\rho_{\samp}\in(0,1)$.}
      \KwOut{A $(\delta,\epsilon)$-atomic equipartition of $G$ into
        $K_{\equirand}(\ell,0,c_{\ZZ},c_{\samp},c_{\Root},c_{\Sum},\epsilon)$ parts. With probability at most $\rho_{\samp}$,
        the algorithm may return an incorrect partition.}
      $\widetilde{\delta}\assign(1+c_{\atom})\cdot\delta$\;
      Let $m\assign\ceil{C\cdot d/(c_{\atom}^2\cdot\delta^2)}$, where $C$ is the absolute constant of
      Lemma~\ref{lem:atomconst}.\;\label{alg:equirand:atom:m}
      $\widetilde{m}\assign\floor{\widetilde{\delta}\cdot m}$\;
      $\zeta\assign (1-c_{\ZZ})\cdot\epsilon$\;
      $\widetilde{c}_{\UC}\assign c_{\UC}\assign\widetilde{c}_{\adj}\assign c_{\adj}\assign 0$\;
      $\widetilde{c}_{\samp}\assign(1-c_{\samp})^{-\ell} - 1$\;
      $\widetilde{c}_{\ZZ}\assign(1-c_{\ZZ})^{-\ell-1} - 1$\;
      $\widetilde{c}_{\Root}\assign(1-c_{\Root})^{-\ell-1} - 1$\;
      $K\assign\ceil{(1+\widetilde{c}_{\UC})\cdot(1+\widetilde{c}_{\ZZ})\cdot(1+\widetilde{c}_{\samp})\cdot(1+\widetilde{c}_{\adj})
        \cdot(1+\widetilde{c}_{\Root})\cdot(1+c_{\Sum})\cdot\epsilon^{-\ell-1}}$\;
      $\cQ\assign\varnothing$\;
      $R\assign V(G)$\;
      $n_0\assign n$\;
      \uFor{$i\in[K]$}{%
        \label{alg:equirand:atom:For}
        Find a rational $\xi_i\in\QQ$ satisfying:
        \begin{align*}
          0 & < \xi_i\leq(1-c_{\Root})\cdot\zeta,
          &
          0 & \leq -1 + \zeta - \xi_i + \gamma_i^\ell\cdot\frac{n_{i-1}}{r} \leq c_{\Root}\cdot\zeta,
        \end{align*}
        where
        \begin{align*}
          \gamma_i & \df \widetilde{\epsilon}_i \df (1-c_{\samp})\cdot\widetilde{\xi}_i,
          &
          \widetilde{\xi}_i & \df \xi_i.
        \end{align*}
        \;\label{alg:equirand:atom:xii}
        $r_i\assign\ceil{\gamma_i^\ell\cdot n_{i-1}}$\;\label{alg:equirand:atom:ri}
        $s_0\df n_{i-1}$\;
        \lFor{$j\in[\ell+1]$}{\label{alg:equirand:atom:sj}$s_j\assign\floor{\gamma_i\cdot s_{j-1}}+1$}
        %
        \setcounter{algosplit}{\theAlgoLine}
      }
    \end{algorithm}

    \begin{algorithm}[htbp]
      \caption*{(continued).}
      \DontPrintSemicolon
      \setcounter{AlgoLine}{\thealgosplit}
      \let\oldnl\nl
      \let\nl\relax
      \vspace{-\baselineskip}\InvisibleBegin{
        \global\let\nl\oldnl 
        Run Algorithm~\ref{alg:ext:atom} with parameters:
        \begin{gather*}
          (n,m,\ell,s_0,\ldots,s_\ell,s_{\ell+1},\widetilde{m},G,U)
          \df
          (n,m,\ell,s_0,\ldots,s_\ell,s_{\ell+1},\widetilde{m},G,R).
        \end{gather*}
        \;\label{alg:equirand:atom:call}
        Let $(W,s)$ be the pair returned by the algorithm.\;
        Let $W'$ be a uniformly at random subset of $W$ of size $r_i$.\;
        $\cQ\assign\cP\cup\{(W',r_i)\}$\;
        $R\assign R\setminus W'$\;
        $n_i\assign n_{i-1} - r_i$\;
      }
      $\cP\assign\varnothing$\;
      $r\assign\floor{n/K}$\;
      $b\assign n - K\cdot r$\;
      \For{$(W,s)\in\cQ$}{%
        \label{alg:equirand:atom:secondFor}
        Let $A\subseteq R$ be a set of size $r - s + \One[b > 0]$.\;
        $\cP\assign\cP\cup\{W\cup A\}$\;
        $R\assign R\setminus A$\;
        $b\assign b-1$\;
      }
      \Return{$\cP$}
    \end{algorithm}

    The only other detail we need to consider about correctness is that in step~\ref{thm:equirand:atom:randsubsetcall} of
    item~\ref{thm:equirand:atom}, our argument of existence of the subset $Q_i$ of $\widetilde{Q}_i$ of size exactly $r_i$ that
    is a $(\widetilde{\delta},\xi_i,m)$-atom of $G$ relied on a probabilistic argument via
    Proposition~\ref{prop:randsubset}\ref{prop:randsubset:atom}. Algorithm~\ref{alg:equirand:atom} also uses this result, but
    relies on the actual success of the probabilistic event, which happens with probability at least $1 -
    \widetilde{\rho}_{\samp}$. Since we rely on this success happening a total of $K$ times (and of course, these are
    independent events), the final probability of success of the algorithm is at least
    \begin{gather*}
      (1 - \widetilde{\rho}_{\samp})^K = 1 - \rho_{\samp}.
    \end{gather*}

    \medskip

    We now analyze the time-complexity of Algorithm~\ref{alg:equirand:atom}. Clearly, the loops of
    lines~\ref{alg:equirand:atom:For} and~\ref{alg:equirand:atom:secondFor} execute exactly $K\leq O(\epsilon^{-\ell-1})$ times.
    If we disregard the time-complexity of the computations of $\xi_i$, $\gamma_i$, $\widetilde{\epsilon}_i$ and
    $\widetilde{\xi}_i$ in line~\ref{alg:equirand:atom:xii}, $r_i$ in line~\ref{alg:equirand:atom:ri} and the $s_j$ in
    line~\ref{alg:equirand:atom:sj} (we will account for these later), each iteration of the loop of
    line~\ref{alg:equirand:atom:For} takes time at most $O((\ell+m)\cdot n^{m+1} + n)\leq O((\ell+m)\cdot n^{m+1})$ (the first
    term is due to the call to Algorithm~\ref{alg:ext:atom} in line~\ref{alg:equirand:atom:call}, see
    Proposition~\ref{prop:ext}\ref{prop:ext:atomalg}), which when accounting for the number of executions becomes
    \begin{gather}\label{eq:equirand:atom:For:complexity}
     O\bigl((\ell+m)\cdot\epsilon^{-\ell-1}\cdot n^{m+1}\bigr).
    \end{gather}

    Since each iteration of the loop of line~\ref{alg:equirand:atom:secondFor} clearly takes time $O(n)$, the total
    time-complexity of it is at most
    \begin{gather}\label{eq:equirand:atom:secondFor:complexity}
      O\bigl(\epsilon^{-\ell-1}\cdot n\bigr).
    \end{gather}

    Regarding the time-complexity of the computation of the parameters:
    \begin{align*}
      \widetilde{\delta}, & &
      m, & &
      \zeta, & &
      \widetilde{c}_{\UC}=c_{\UC}=\widetilde{c}_{\adj}=c_{\adj}=0,
      \\
      \widetilde{c}_{\samp}, & &
      \widetilde{c}_{\ZZ}, & &
      \widetilde{c}_{\Root}, & &
      K,
    \end{align*}
    we can give the following respective crude upper bounds:
    \begin{gather*}
      \begin{aligned}
        O\bigl(\len(c_{\atom})\cdot\len(\delta)\bigr),
        & &
        O\Bigl(
        \log(d+1)\cdot\len(c_{\atom})
        + \bigl(\log(d+1) + \len(c_{\atom})\bigr)\cdot\len(\delta)
        \Bigr),
      \end{aligned}
      \\
      \begin{aligned}
        O\bigl(\len(c_{\ZZ})\cdot\len(\epsilon)\bigr),
        & &
        O(1),
      \end{aligned}
      \\
      \begin{aligned}
        O\bigl(\ell^2\cdot\log(\ell+1)\cdot\len(c_{\samp})^2\bigr),
        & &
        O\bigl(\ell^2\cdot\log(\ell+1)\cdot\len(c_{\ZZ})^2\bigr),
      \end{aligned}
      \\
      O\bigl(\ell^2\cdot\log(\ell+1)\cdot\len(c_{\Root})^2\bigr),
      \\
      O\bigl(
      \ell^4\cdot\log(\ell+1)\cdot\len(c_{\ZZ})\cdot\len(c_{\samp})\cdot\len(c_{\Root})\cdot\len(c_{\Sum})\cdot\len(\epsilon)^2
      \bigr),
    \end{gather*}
    which together are bounded by
    \begin{gather}\label{eq:equirand:atom:param:complexity}
      \begin{multlined}
        O\bigl(
        \log(d+1)\cdot\len(c_{\atom})\cdot\len(\delta)
        \\
        +
        \ell^4\cdot\log(\ell+1)\cdot\len(c_{\ZZ})^2\cdot\len(c_{\samp})^2\cdot\len(c_{\Root})^2\cdot\len(c_{\Sum})\cdot\len(\epsilon)^2
        \bigr).
      \end{multlined}
    \end{gather}

    Let us now analyze the time-complexity of the computations of line~\ref{alg:equirand:atom:xii}. The line requires us to find
    the rational $\xi_i$ satisfying
    \begin{align*}
      0 & < \xi_i\leq(1-c_{\Root})\cdot\zeta,
      &
      0 & \leq -1 + \zeta - \xi_i + \gamma_i^\ell\cdot\frac{n_{i-1}}{r} \leq c_{\Root}\cdot\zeta,
    \end{align*}
    where
    \begin{align*}
      \gamma_i & \df \widetilde{\epsilon}_i \df (1-c_{\samp})\cdot\widetilde{\xi}_i,
      &
      \widetilde{\xi}_i & \df \xi_i.
    \end{align*}
    By Lemma~\ref{lem:equicalc}\ref{lem:equicalc:xii}, we know that the middle term $f(\xi_i)$ of the second set of inequalities
    above gives a negative value when evaluated at $\xi_i=0$ and gives a value that is at least $c_{\Root}\cdot\zeta$ when
    evaluated at $(1-c_{\Root})\cdot\zeta$. This means that we can do a simple binary search for the value of $\xi_i$, computing
    the corresponding values of $\gamma_i$, $\widetilde{\epsilon}_i$ and $\widetilde{\xi}_i$ at each step. Note that the
    function $f$ is Lipschitz-continuous with constant at most
    \begin{align*}
      L & \df
      \ell\cdot\frac{n_{i-1}}{r}
      \leq
      \ell\cdot\frac{n}{\floor{n/K}}
      \leq
      \ell\cdot\frac{K}{1 - K/n}
      \\
      & \leq
      O(\ell\cdot K)
      \leq
      O(\ell\cdot\epsilon^{-\ell-1}),
    \end{align*}
    so our binary search has to stop as soon as our interval of search is of length less than $c_{\Root}\cdot\zeta/L$ and since
    the initial length is $(1-c_{\Root})\cdot\zeta$, we conclude that the binary search takes at most
    \begin{align*}
      \log_2\left(\frac{(1-c_{\Root})\cdot\zeta\cdot L}{c_{\Root}\cdot\zeta}\right)
      & \leq
      \log_2\left(\frac{(1-c_{\Root})\cdot\ell}{c_{\Root}\cdot(1-c_{\ZZ})\cdot\epsilon^{\ell+2}}\right)
      +
      O(1)
      \\
      & \leq
      O\bigl(\len(c_{\Root}) + \len(c_{\ZZ}) + \ell\cdot\len(\epsilon)\bigr)
    \end{align*}
    steps. This also in particular means that at all stages the bitlength of the endpoints of the search is at most the above
    plus
    \begin{gather*}
      O\bigl(\len(c_{\samp}) + \len(\zeta)\bigr)
      \leq
      O\bigl(\len(c_{\samp}) + \len(c_{\ZZ}) + \len(\epsilon)\bigr)
    \end{gather*}

    We can give the following crude bound for the time it takes to make the computations of each step of the binary search:
    \begin{gather*}
      O\Bigl(
      \ell^2\cdot\log(\ell+1)\cdot
      \len(c_{\samp})^2\cdot\bigl(\len(c_{\Root}) + \len(c_{\ZZ}) + \ell\cdot\len(\epsilon) + \len(c_{\samp})\bigr)^2
      \cdot\bigl(\log(n+1)\bigr)^2
      \Bigr).
    \end{gather*}

    Accounting for the number of iterations of the binary search and the number of iterations of the loop of
    line~\ref{alg:equirand:atom:For}, the total time-complexity of all executions of line~\ref{alg:equirand:atom:xii} is at most
    \begin{gather}\label{eq:equirand:atom:xii:complexity}
      \begin{aligned}
        \MoveEqLeft
        \begin{multlined}[t]
          O\Bigl(
          \epsilon^{-\ell-1}\cdot
          \bigl(\len(c_{\Root}) + \len(c_{\ZZ}) + \ell\cdot\len(\epsilon)\bigr)\cdot
          \ell^2\cdot\log(\ell+1)
          \\
          \cdot
          \len(c_{\samp})^2\cdot\bigl(\len(c_{\Root}) + \len(c_{\ZZ}) + \ell\cdot\len(\epsilon) + \len(c_{\samp})\bigr)^2
          \cdot\bigl(\log(n+1)\bigr)^2
          \Bigr)
        \end{multlined}
        \\
        & \leq
        O\left(
        \frac{
          \ell^5\cdot\log(\ell+1)\cdot\len(\epsilon)^3\cdot\len(c_{\ZZ})^3\cdot\len(c_{\samp})^4
          \cdot\len(c_{\Root})^3\cdot(\log(n+1))^2
        }{
          \epsilon^{\ell+1}
        }
        \right).
      \end{aligned}
    \end{gather}

    For the calculations of lines~\ref{alg:equirand:atom:ri} and~\ref{alg:equirand:atom:sj}, note that
    \begin{gather*}
      \len(\gamma_i)
      \leq
      O\bigl(\len(\xi_i) + \len(c_{\samp})\bigr)
      \leq
      O\Bigl(
      \len(c_{\Root}) + \len(c_{\ZZ}) + \ell\cdot\len(\epsilon) + \len(c_{\samp})
      \Bigr),
    \end{gather*}
    so taking into account the number of iterations of the loop of line~\ref{alg:equirand:atom:For}, the total time-complexity
    of all executions of these lines are
    \begin{gather*}
      O\Bigl(
      \epsilon^{-\ell-1}\cdot
      \ell\cdot\log(\ell+1)\cdot
      \bigl(\len(c_{\Root}) + \len(c_{\ZZ}) + \ell\cdot\len(\epsilon) + \len(c_{\samp})\bigr)\cdot\log(n+1)
      \Bigr),
      \\
      O\Bigl(
      \epsilon^{-\ell-1}\cdot
      \ell\cdot
      \bigl(\len(c_{\Root}) + \len(c_{\ZZ}) + \ell\cdot\len(\epsilon) + \len(c_{\samp})\bigr)\cdot\log(n+1)
      \Bigr).
    \end{gather*}
    These are clearly dominated by the time-complexity in~\eqref{eq:equirand:atom:xii:complexity}.

    Putting everything together (i.e., equations~\eqref{eq:equirand:atom:For:complexity},
    \eqref{eq:equirand:atom:secondFor:complexity}, \eqref{eq:equirand:atom:param:complexity}
    and~\eqref{eq:equirand:atom:xii:complexity}), the total time-complexity of Algorithm~\ref{alg:equirand:atom} is at most
    \begin{align*}
      \MoveEqLeft
      \begin{multlined}[t]
        O\bigggl(
        (\ell+m)\cdot\epsilon^{-\ell-1}\cdot n^{m+1}
        +
        \epsilon^{-\ell-1}\cdot n
        \\
        +
        \log(d+1)\cdot\len(c_{\atom})\cdot\len(\delta)
        \\
        +
        \ell^4\cdot\log(\ell+1)\cdot\len(c_{\ZZ})^2\cdot\len(c_{\samp})^2\cdot\len(c_{\Root})^2\cdot\len(c_{\Sum})\cdot\len(\epsilon)^2
        \\
        +
        \frac{
          \ell^5\cdot\log(\ell+1)\cdot\len(\epsilon)^3\cdot\len(c_{\ZZ})^3\cdot\len(c_{\samp})^4
          \cdot\len(c_{\Root})^3\cdot(\log(n+1))^2
        }{
          \epsilon^{\ell+1}
        }
        \bigggr)
      \end{multlined}
      \\
      & \leq
      \begin{multlined}[t]
        \frac{\ell}{\epsilon^{\ell+1}}\cdot n^{O(d/(c_{\atom}^2\cdot\delta^2))}
        +
        \log(d+1)\cdot\len(c_{\atom})\cdot\len(\delta)
        \\
        +
        \frac{
          \ell^5\cdot\log(\ell+1)\cdot\len(\epsilon)^3\cdot\len(c_{\ZZ})^3\cdot\len(c_{\samp})^4
          \cdot\len(c_{\Root})^3\cdot\len(c_{\Sum})\cdot(\log(n+1))^2
        }{
          \epsilon^{\ell+1}
        },
      \end{multlined}
    \end{align*}
    where the last inequality follows since
    \begin{gather*}
      m
      =
      \Ceil{C\cdot\frac{d}{c_{\atom}^2\cdot\delta^2}}
      \leq
      O\left(\frac{d}{c_{\atom}^2\cdot\delta^2}\right).
    \end{gather*}
  \item[Item~\ref{thm:equirand:goodalg}.] Algorithm~\ref{alg:equirand:good} is similar to Algorithm~\ref{alg:equirand:atom},
    except that to extract good sets, we call the random-query-oracle model Algorithm~\ref{alg:ext:good}, which itself has a
    probability of failure.

    \begin{algorithm}[htbp]
      \caption{Random-query model algorithm that returns an $\epsilon$-good equipartition of $G$ into
        $K_{\equirand}(\ell,c_{\UC},c_{\ZZ},c_{\samp},c_{\Root},c_{\Sum},\epsilon)$ parts. The success probability is at least
        $(1-\rho_{\UC})\cdot(1-\rho_{\samp})$ and the time-complexity is:
        \begin{multline*}
          O\bigggl(
          \frac{2^\ell\cdot\log(n/\rho_{\UC})\cdot n}{c_{\UC}^2\cdot\epsilon^{\ell+5}}
          \\
          +
          \frac{
            \len(\epsilon)^3\cdot\len(c_{\ZZ})^3\cdot\len(c_{\samp})^4\cdot\len(c_{\UC})^2
            \cdot\len(c_{\Root})^3\cdot\len(c_{\Sum})\cdot\len(\rho_{\UC})
          }{
            \epsilon^{\ell+1}
          }
          \\
          \cdot\ell^5\cdot\log(\ell+1)\cdot(\log(n+1))^2
          \bigggr).
        \end{multline*}
      }
      \label{alg:equirand:good}
      \DontPrintSemicolon
      \KwIn{Numbers $n,d,\ell\in\NN_+$ with $d\leq\ell\leq\log_2(n)$,
        $c_{\UC},c_{\ZZ},c_{\samp},c_{\Root},c_{\Sum},\rho_{\UC}\in(0,1)\cap\QQ$, and
        $\epsilon\in(0,\epsilon_{\equirand}(\ell,c_{\UC},c_{\ZZ},c_{\samp},c_{\Root},c_{\Sum}))$, and a graph $G$ with $\lvert
        G\rvert=n$, $\Lit(G)\leq\ell$ and $\VC(G)\leq d$. We further assume that $n\geq
        n_{\equirand}^{\good}(\ell,c_{\UC},c_{\ZZ},c_{\samp},c_{\Root},c_{\Sum},\epsilon,\rho_{\samp})$ for some choice of
        $\rho_{\samp}\in(0,1)$.}
      \KwOut{An $\epsilon$-good equipartition of $G$ into
        $K_{\equirand}(\ell,c_{\UC},c_{\ZZ},c_{\samp},c_{\Root},c_{\Sum},\epsilon)$ parts. With probability at most $1 -
        (1-\rho_{\UC})\cdot(1-\rho_{\samp})$, the algorithm may return ``\Failure'' or an incorrect partition.}
      $\widetilde{c}_{\UC}\assign (1-c_{\UC})^{-\ell} - 1$\;
      $\widetilde{c}_{\adj}\assign c_{\adj}\assign 0$\;
      $\widetilde{c}_{\samp}\assign(1-c_{\samp})^{-\ell} - 1$\;
      $\widetilde{c}_{\ZZ}\assign(1-c_{\ZZ})^{-\ell-1} - 1$\;
      $\widetilde{c}_{\Root}\assign(1-c_{\Root})^{-\ell-1} - 1$\;
      $K\assign\ceil{(1+\widetilde{c}_{\UC})\cdot(1+\widetilde{c}_{\ZZ})\cdot(1+\widetilde{c}_{\samp})\cdot(1+\widetilde{c}_{\adj})
        \cdot(1+\widetilde{c}_{\Root})\cdot(1+c_{\Sum})\cdot\epsilon^{-\ell-1}}$\;
      $\widetilde{\rho}_{\UC}\assign\rho_{\UC}/K$\;
      $\zeta\assign (1-c_{\ZZ})\cdot\epsilon$\;      
      Let $C$ be the absolute constant of Theorem~\ref{thm:UC}.\;
      $\cQ\assign\varnothing$\;
      $R\assign V(G)$\;
      $n_0\assign n$\;
      \uFor{$i\in[K]$}{%
        \label{alg:equirand:good:For}
        Find a rational $\xi_i\in\QQ$ satisfying:
        \begin{align*}
          0 & < \xi_i\leq(1-c_{\Root})\cdot\zeta,
          &
          0 & \leq -1 + \zeta - \xi_i + \gamma_i^\ell\cdot\frac{n_{i-1}}{r} \leq c_{\Root}\cdot\zeta,
        \end{align*}
        where
        \begin{align*}
          \gamma_i & \df (1-c_{\UC})\cdot\widetilde{\epsilon}_i,
          &
          \widetilde{\epsilon}_i & \df (1-c_{\samp})\cdot\widetilde{\xi}_i,
          &
          \widetilde{\xi}_i & \df \xi_i.
        \end{align*}
        \;\label{alg:equirand:good:xii}
        $r_i\assign\ceil{\gamma_i^\ell\cdot n_{i-1}}$\;\label{alg:equirand:good:ri}
        \setcounter{algosplit}{\theAlgoLine}
      }
    \end{algorithm}

    \begin{algorithm}[htbp]
      \caption*{(continued).}
      \DontPrintSemicolon
      \setcounter{AlgoLine}{\thealgosplit}
      \let\oldnl\nl
      \let\nl\relax
      \vspace{-\baselineskip}\InvisibleBegin{
        \global\let\nl\oldnl 
        $s_0\df n_{i-1}$\;
        \lFor{$j\in[\ell+1]$}{\label{alg:equirand:good:sj}$s_j\assign\floor{\gamma_i\cdot s_{j-1}}+1$}
        $m\assign (C\cdot 4\cdot (d + \ell + 1 +
        \ceil{\log_2(n/\widetilde{\rho}_{\UC})}))/(c_{\UC}^2\cdot\widetilde{\epsilon}_i^2)$\;\label{alg:equirand:good:m}
        $\widetilde{m}\assign\floor{(1 - c_{\UC}/2)\cdot\widetilde{\epsilon}_i\cdot m}$\;\label{alg:equirand:good:wm}
        Run Algorithm~\ref{alg:ext:good} with parameters:
        \begin{gather*}
          (n,m,\ell,s_0,\ldots,s_\ell,s_{\ell+1},\widetilde{m},\cO_G,U)
          \df
          (n,m,\ell,s_0,\ldots,s_\ell,s_{\ell+1},\widetilde{m},\cO_G,R).
        \end{gather*}
        \;\label{alg:equirand:good:call}
        \lIf{The algorithm returned ``\Failure''}{\Return{\Failure.}}
        Let $(W,s)$ be the pair returned by the algorithm.\;
        Let $W'$ be a uniformly at random subset of $W$ of size $r_i$.\;\label{alg:equirand:good:subsample}
        $\cQ\assign\cP\cup\{(W',r_i)\}$\;
        $R\assign R\setminus W'$\;
        $n_i\assign n_{i-1} - r_i$\;
      }
      $\cP\assign\varnothing$\;
      $r\assign\floor{n/K}$\;
      $b\assign n - K\cdot r$\;
      \For{$(W,s)\in\cQ$}{%
        \label{alg:equirand:good:secondFor}
        Let $A\subseteq R$ be a set of size $r - s + \One[b > 0]$.\;
        $\cP\assign\cP\cup\{W\cup A\}$\;
        $R\assign R\setminus A$\;
        $b\assign b-1$\;
      }
      \Return{$\cP$}
    \end{algorithm}

    The proof of correctness of Algorithm~\ref{alg:equirand:good} is analogous to that of Algorithm~\ref{alg:equirand:atom},
    provided all calls to Algorithm~\ref{alg:ext:good} in line~\ref{alg:equirand:good:call} are successful. We also note that
    our particular calculation of the parameter $m$ is
    \begin{gather*}
      m
      \df
      \Ceil{\frac{C\cdot 4\cdot (d + \ell + 1 + \ceil{\log_2(n/\widetilde{\rho}_{\UC})})}{c_{\UC}^2\cdot\widetilde{\epsilon}_i^2}}
      \geq
      \frac{C\cdot 4\cdot (d + \ln((2^{\ell+1}-1)\cdot n/\widetilde{\rho}_{\UC}))}{(c_{\UC}^2\cdot\widetilde{\epsilon}_i^2)},
    \end{gather*}
    so we can indeed apply Proposition~\ref{prop:ext}\ref{prop:ext:goodalg} and use Algorithm~\ref{alg:ext:good}.

    Since we know that the loop of line~\ref{alg:equirand:good:For} runs exactly $K$ times, we know that probability that all
    calls to Algorithm~\ref{alg:equirand:good} are successful is at least
    \begin{gather*}
      (1 - \widetilde{\rho}_{\UC})^K \geq 1 - K\cdot\widetilde{\rho}_{\UC} = 1 - \rho_{\UC}.
    \end{gather*}
    Putting this together with the fact that, just as in Algorithm~\ref{alg:equirand:atom}, each subsampling in
    line~\ref{alg:equirand:good:subsample} is successful with probability at least $1 - \widetilde{\rho}_{\samp}$, the final
    probability of success is at least
    \begin{gather*}
      (1-\rho_{\UC})\cdot(1 - \widetilde{\rho}_{\samp})^K = (1-\rho_{\UC})\cdot(1-\rho_{\samp}).
    \end{gather*}

    \medskip

    We now analyze the time-complexity of Algorithm~\ref{alg:equirand:good}. The loops of lines~\ref{alg:equirand:good:For}
    and~\ref{alg:equirand:good:secondFor} execute exactly $K\leq O(\epsilon^{-\ell-1})$ times and disregarding the
    time-complexity of the computations of $\xi_i$, $\gamma_i$, $\widetilde{\epsilon}_i$ and $\widetilde{\xi}_i$ in
    line~\ref{alg:equirand:good:xii}, $r_i$ in line~\ref{alg:equirand:good:ri}, the $s_j$ in line~\ref{alg:equirand:good:sj},
    $m$ in line~\ref{alg:equirand:good:m} and $\widetilde{m}$ in line~\ref{alg:equirand:good:wm} (we will account for these
    later), each iteration of the loop of line~\ref{alg:equirand:good:For} takes time at most
    \begin{align*}
      O(2^\ell\cdot n\cdot m)
      & \leq
      O\left(
      \frac{2^\ell\cdot(\log(n/\rho_{\UC}) + \ell\cdot\log(1/\epsilon))\cdot n}{c_{\UC}^2\cdot\widetilde{\epsilon}_i^2}
      \right)
      \leq
      O\left(\frac{2^\ell\cdot\log(n/\rho_{\UC})\cdot n}{c_{\UC}^2\cdot\widetilde{\epsilon}_i^2}\right)
      \\
      & \leq
      O\left(\frac{2^\ell\cdot\log(n/\rho_{\UC})\cdot n\cdot K^{2/\ell}}{c_{\UC}^2}\right)
      \leq
      O\left(\frac{2^\ell\cdot\log(n/\rho_{\UC})\cdot n}{c_{\UC}^2\cdot\epsilon^4}\right),
    \end{align*}
    where the second inequality follows since
    \begin{gather*}
      n
      \geq
      n_{\equi}(\ell,c_{\UC},c_{\ZZ},c_{\samp},c_{\adj},c_{\Root},c_{\Sum},\epsilon)
      \geq
      \frac{K\cdot(1-\epsilon)}{\epsilon-\zeta}
      \geq
      \Omega(\epsilon^{-\ell-2}),
    \end{gather*}
    and the third inequality follows from Lemma~\ref{lem:equicalc}\ref{lem:equicalc:xiilowerbound} and $r\df\floor{n/K}$. It is
    also clear that each iteration of the loop of line~\ref{alg:equirand:good:secondFor} takes time at most $O(n)$, so the total
    time-complexity of these loops (excluding the aforementioned computations) is at most
    \begin{gather}\label{eq:equirand:good:loops:complexity}
      O\left(\frac{2^\ell\cdot\log(n/\rho_{\UC})\cdot n}{c_{\UC}^2\cdot\epsilon^{\ell+5}}\right).
    \end{gather}

    Regarding the time-complexity of the computation of the parameters:
    \begin{align*}
      \widetilde{c}_{\UC}, & &
      \widetilde{c}_{\adj} = c_{\adj} = 0, & &
      \widetilde{c}_{\samp}, & &
      \widetilde{c}_{\ZZ},
      \\
      \widetilde{c}_{\Root}, & &
      K, & &
      \widetilde{\rho}_{\UC}, & &
      \zeta,
    \end{align*}
    we can give the following respective crude upper bounds:
    \begin{gather*}
      \begin{aligned}
        O\bigl(\ell^2\cdot\log(\ell+1)\cdot\len(c_{\UC})^2\bigr), & &
        O(1), & &
        O\bigl(\ell^2\cdot\log(\ell+1)\cdot\len(c_{\samp})^2\bigr),
      \end{aligned}
      \\      
      \begin{aligned}
        O\bigl(\ell^2\cdot\log(\ell+1)\cdot\len(c_{\ZZ})^2\bigr),
        & &
        O\bigl(\ell^2\cdot\log(\ell+1)\cdot\len(c_{\Root})^2\bigr),
      \end{aligned}
      \\
      O\bigl(
      \ell^5\cdot\log(\ell+1)\cdot\len(c_{\ZZ})\cdot\len(c_{\samp})\cdot\len(c_{\Root})\cdot\len(c_{\Sum})\cdot\len(\epsilon)^2
      \bigr),
      \\
      \begin{aligned}
        O\bigl(\len(\rho_{\UC})\cdot\ell\cdot\log(1/\epsilon)\bigr),
        & &
        O\bigl(\len(c_{\ZZ})\cdot\len(\epsilon)\bigr),
      \end{aligned}
    \end{gather*}
    which together are upper bounded by
    \begin{gather}\label{eq:equirand:good:param:complexity}
      \begin{multlined}[t]
        O\bigl(
        \ell^5\cdot\log(\ell+1)\cdot\len(c_{\UC})^2\cdot\len(c_{\ZZ})^2\cdot\len(c_{\samp})^2\cdot\len(c_{\Root})^2
        \cdot\len(c_{\Sum})\cdot\len(\epsilon)^2
        \\
        +
        \len(\rho_{\UC})\cdot\ell\cdot\log(1/\epsilon)
        \bigr).
      \end{multlined}
    \end{gather}
      
    The analysis of the time-complexity of the computations of line~\ref{alg:equirand:good:xii} is analogous to the analysis of
    the computations of line~\ref{alg:equirand:atom:xii} in Algorithm~\ref{alg:equirand:atom}
    (see~\eqref{eq:equirand:atom:xii:complexity}), except that now the computation of $\gamma_i$ is slightly more expensive (as
    $\gamma_i=(1-c_{\UC})\cdot\widetilde{\epsilon}_i$), which yields the total time-complexity (already accounting for the
    repetition of the loop of line~\ref{alg:equirand:good:xii}):
    \begin{gather}\label{eq:equirand:good:xii:complexity}
      O\left(
        \frac{
          \ell^5\cdot\log(\ell+1)\cdot\len(\epsilon)^3\cdot\len(c_{\ZZ})^3\cdot\len(c_{\samp})^4\cdot\len(c_{\UC})^2
          \cdot\len(c_{\Root})^3\cdot(\log(n+1))^2
        }{
          \epsilon^{\ell+1}
        }
        \right).
    \end{gather}
    The time-complexities of the calculations of lines~\ref{alg:equirand:good:ri} and~\ref{alg:equirand:good:sj} are completely
    dominated by the above (these are as in Algorithm~\ref{alg:equirand:atom}, accounting for the additional $\len(c_{\UC})$ in
    the bitlength of $\gamma_i$).

    The computation of $m$ in line~\ref{alg:equirand:good:m} has time-complexity at most
    \begin{align*}
      \MoveEqLeft
      O\bigl(
      \log(n+1)\cdot\len(\rho_{\UC})\cdot\ell\cdot\len(\epsilon)\cdot\len(c_{\UC})\cdot\len(\widetilde{\epsilon}_i)
      \bigr)
      \\
      & \leq
      \begin{multlined}[t]
        O\bigl(
        \log(n+1)\cdot\len(\rho_{\UC})\cdot\ell\cdot\len(\epsilon)\cdot\len(c_{\UC})
        \\
        \cdot\bigl(\len(c_{\Root}) + \len(c_{\ZZ}) + \len(c_{\samp}) + \ell\cdot\len(\epsilon)\bigr)
        \bigr)
      \end{multlined}
      \\
      & \leq
      \log(n+1)\cdot\len(\rho_{\UC})\cdot\ell^2\cdot\len(\epsilon)^2\cdot\len(c_{\UC})\cdot\len(c_{\Root})
      \cdot\len(c_{\ZZ})\cdot\len(c_{\samp})
    \end{align*}
    and the computation of $\widetilde{m}$ in line \ref{alg:equirand:good:wm} has time-complexity at most
    \begin{align*}
      \MoveEqLeft
      O\bigl(
      \len(c_{\UC})\cdot\len(\widetilde{\epsilon}_i)\cdot\log(m+1)
      \bigr)
      \\
      & \leq
      \begin{multlined}[t]
        O\biggl(
        \len(c_{\UC})\cdot\len(\widetilde{\epsilon}_i)\cdot
        \Bigl(\log\bigl(\log(n)\bigr) + \log\bigl(\log(1/\rho_{\UC})\bigr) + \log\bigl(\log(1/\epsilon)\bigr)
        \\
        +
        \log(1/c_{\UC}) + \log(1/\widetilde{\epsilon}_i)\Bigr)
        \biggr)
      \end{multlined}
      \\
      & \leq
      O\Bigl(
      \len(c_{\UC})^2\cdot\bigl(\len(c_{\Root}) + \len(c_{\ZZ}) + \len(c_{\samp}) + \ell\cdot\len(\epsilon)\bigr)^2
      \cdot\log(n+1)\cdot\len(\rho_{\UC})
      \Bigr)
      \\
      & \leq
      O\bigl(
      \log(n+1)\cdot\ell^2\cdot\len(\epsilon)^2
      \cdot\len(c_{\UC})^2\cdot\len(c_{\Root})^2\cdot\len(c_{\ZZ})^2\cdot\len(c_{\samp})^2
      \cdot\len(\rho_{\UC})
      \bigr).
    \end{align*}

    Thus, taking also into account that the loop of line~\ref{alg:equirand:good:For} runs $K\leq O(\epsilon^{-\ell-1})$ times,
    the total time-complexity of lines~\ref{alg:equirand:good:m} and~\ref{alg:equirand:good:wm} is at most
    \begin{gather}\label{alg:equirand:good:mwm:complexity}
      O\left(
      \frac{
        \log(n+1)\cdot\ell^2\cdot\len(\epsilon)^2
        \cdot\len(c_{\UC})^2\cdot\len(c_{\Root})^2\cdot\len(c_{\ZZ})^2\cdot\len(c_{\samp})^2
        \cdot\len(\rho_{\UC})
      }{
        \epsilon^{\ell+1}
      }
      \right).
    \end{gather}

    Putting everything together (i.e., equations~\eqref{eq:equirand:good:loops:complexity},
    \eqref{eq:equirand:good:param:complexity}, \eqref{eq:equirand:good:xii:complexity}
    and~\eqref{alg:equirand:good:mwm:complexity}), the total time-complexity of Algorithm~\ref{alg:equirand:good} is at most
    \begin{align*}
      \MoveEqLeft
      \begin{multlined}[t]
        O\bigggl(
        \frac{2^\ell\cdot\log(n/\rho_{\UC})\cdot n}{c_{\UC}^2\cdot\epsilon^{\ell+5}}
        \\
        +
        \ell^5\cdot\log(\ell+1)\cdot\len(c_{\UC})^2\cdot\len(c_{\ZZ})^2\cdot\len(c_{\samp})^2\cdot\len(c_{\Root})^2
        \cdot\len(c_{\Sum})\cdot\len(\epsilon)^2
        \\
        +
        \len(\rho_{\UC})\cdot\ell\cdot\log(1/\epsilon)
        \\
        +
        \frac{
          \ell^5\cdot\log(\ell+1)\cdot\len(\epsilon)^3\cdot\len(c_{\ZZ})^3\cdot\len(c_{\samp})^4\cdot\len(c_{\UC})^2
          \cdot\len(c_{\Root})^3\cdot(\log(n+1))^2
        }{
          \epsilon^{\ell+1}
        }
        \\
        +
        \frac{
        \log(n+1)\cdot\ell^2\cdot\len(\epsilon)^2
        \cdot\len(c_{\UC})^2\cdot\len(c_{\Root})^2\cdot\len(c_{\ZZ})^2\cdot\len(c_{\samp})^2
        \cdot\len(\rho_{\UC})
        }{
          \epsilon^{\ell+1}
        }
        \bigggr)
      \end{multlined}
      \\
      & \leq
      \begin{multlined}[t]
        O\bigggl(
        \frac{2^\ell\cdot\log(n/\rho_{\UC})\cdot n}{c_{\UC}^2\cdot\epsilon^{\ell+5}}
        \\
        +
        \frac{
          \len(\epsilon)^3\cdot\len(c_{\ZZ})^3\cdot\len(c_{\samp})^4\cdot\len(c_{\UC})^2
          \cdot\len(c_{\Root})^3\cdot\len(c_{\Sum})\cdot\len(\rho_{\UC})
        }{
          \epsilon^{\ell+1}
        }
        \\
        \cdot\ell^5\cdot\log(\ell+1)\cdot(\log(n+1))^2
        \bigggr).
        \qedhere
      \end{multlined}
    \end{align*}
  \end{description}
\end{proof}

\begin{discussion}
  Why is there no space-efficient randomized algorithms in Theorem~\ref{thm:equirand}? The reason is very simple: since we use
  randomness to subsample sets, to keep track of which random subset we picked, we would have to either store the set itself or
  the randomness bits used to generate it. Naively, either of these would take space $\Omega(n)$; at this space-complexity
  regime, we could simply run Algorithms~\ref{alg:equirand:good} or~\ref{alg:equirand:atom}. Instead, in the next section, we
  will see how to do ``random set subsampling'' with $O(\log(n))$ bits of randomness; as a consequence, we will get
  derandomizations that are polynomial-time and logarithmic-space.
\end{discussion}

\begin{remark}[Choosing $c_{\atom}$ and $\delta$]
  Similarly to Remark~\ref{rmk:delta}, one can actually optimize the choice of $(c_{\atom},\delta)$ since Lemma~\ref{lem:mono}
  says that a bigger value of $\delta$ leads to stronger atomicity.

  All complexity bounds in Theorem~\ref{thm:equirand} also factor as decreasing functions of the product $c_{\atom}\cdot\delta$,
  so we get the same calculation, which effectively means that in all complexity bounds we can replace $c_{\atom}\cdot\delta$ by
  $(1-o(1))\cdot 2^{-\ell-2}$ by running the algorithms with $c_{\atom}\df 1 - a$ for $a > 0$ small and $\delta\df 2^{-\ell-2}$.
\end{remark}


\section{Equipartitions, without randomness}

\begin{proposition}[A slight improvement over Alon--Nussboim~\protect{\cite[Lemma~2]{AN08}}]\label{prop:highchebyshev}
  Let $h,n\in\NN_+$ with $h$ even, let $p,\delta\in(0,1)$ and let $(\rn{X}_i)_{i=1}^n$ be $\{0,1\}$-valued random variables that
  are $h$-wise independent and each $\rn{X}_i$ is $p$-Bernoulli (i.e., $\PP[\rn{X}_i=1]=p$). Let also
  $\rn{X}\df\sum_{i=1}^n\rn{X}_i$ and suppose that
  \begin{gather}\label{eq:highchebyshev}
    n\geq \frac{(h-1)^2}{p\cdot(1-p)}.
  \end{gather}

  Then
  \begin{gather*}
    \PP\bigl[\lvert\rn{X} - p\cdot n\rvert > \delta\cdot p\cdot n\bigr]
    <
    \frac{1}{\sqrt{2\cdot\pi\cdot h}}
    \cdot
    \left(\frac{2\cdot e\cdot (1-p)}{\delta^2\cdot p\cdot n}\right)^{h/2}.
  \end{gather*}
\end{proposition}

\begin{proof}
  For every $i\in[n]$, let $\rn{Y}_i\df\rn{X}_i - p$ so that $\EE[\rn{Y}_i] = 0$ and let $\rn{Y}\df\sum_{i=1}^n\rn{Y}_i$ so that
  $\rn{Y} = \rn{X} - p\cdot n$. Note that since $(\rn{X}_i)_{i=1}^n$ is $h$-wise independent, so is $(\rn{Y}_i)_{i=1}^n$.

  We now bound the $h$th moment of $\rn{Y}$:
  \begin{align*}
    \EE[\rn{Y}]
    & =
    \EE\left[\left(\sum_{i=1}^n\rn{Y}_i\right)^h\right]
    =
    \sum_{\substack{j\in\{0,\ldots,h\}^n\\\sum_{i=1}^n j_i = h}} \EE\left[\prod_{i=1}^n \rn{Y}_i^{j_i}\right]
    =
    \sum_{\substack{j\in\{0,\ldots,h\}^n\\\sum_{i=1}^n j_i = h}} \prod_{i=1}^n \EE[\rn{Y}_i^{j_i}]
    \\
    & =
    \sum_{\substack{j\in\{0,2,3,\ldots,h\}^n\\\sum_{i=1}^n j_i = h}} \prod_{i=1}^n \EE[\rn{Y}_i^{j_i}]
    =
    \sum_{s=1}^{h/2} \sum_{I\in\binom{[n]}{s}} \sum_{\substack{j\in\{2,3,\ldots,h\}^I\\\sum_{i\in I} j_i = h}}
    \prod_{i\in I} \EE[\rn{Y}_j^{j_i}],
  \end{align*}
  where the third equality follows from $h$-wise independence since the product involves at most $h$ random variables, the
  fourth equality follows since $\EE[\rn{Y}_i]=0$ so whenever $j$ has some entry $j_i=1$, the corresponding term is $0$, and the
  fifth equality follows by keeping track of which entries are non-zero. Now note that since $\rn{X}_i$ is $\{0,1\}$-valued, for
  every $j\geq 2$, we have
  \begin{gather*}
    \EE[\rn{Y}_i^j]
    =
    (1-p)^j\cdot p + (-p)^j\cdot(1-p)
    =
    p\cdot(1-p)\cdot((1-p)^{j-1} + (-1)^j\cdot p^{j-1})
    \leq
    p\cdot(1-p),
  \end{gather*}
  so we get
  \begin{align*}
    \EE[\rn{Y}]
    & =
    \sum_{s=1}^{h/2} \sum_{I\in\binom{[n]}{s}} \sum_{\substack{j\in\{2,\ldots,h\}^I\\\sum_{i\in I} j_i = h}}
    \prod_{i\in I} \EE[\rn{Y}_j^{j_i}]
    \\
    & \leq
    \sum_{s=1}^{h/2}
    \binom{n}{s}\cdot p^s\cdot(1-p)^s
    \cdot \left\lvert\left\{j\in\{2,\ldots,h\}^s \mid \sum_{i=1}^s j_i = h\right\}\right\rvert
    \\
    & =
    \sum_{s=1}^{h/2}
    \binom{n}{s}\cdot p^s\cdot(1-p)^s
    \cdot\left\lvert\left\{j\in[h-1]^s \mid \sum_{i=1}^s j_i = h-s\right\}\right\rvert
    \\
    & =
    \sum_{s=1}^{h/2}
    \binom{n}{s}\cdot p^s\cdot(1-p)^s\cdot\binom{h-s-1}{s-1}
    \leq
    \sum_{s=1}^{h/2}
    \binom{n}{s}\cdot p^s\cdot(1-p)^s\cdot(h-1)^{s-1}.
  \end{align*}

  Let $\Psi_s$ be the term of the last sum above corresponding to $s$ and note that for $s\in[h/2]$, we have
  \begin{gather*}
    \frac{\Psi_s}{\Psi_{s-1}}
    =
    \frac{(n-s+1)}{s}\cdot p\cdot(1-p)\cdot(h-1)
    \geq
    \frac{2\cdot n-h+2}{h}\cdot p\cdot(1-p)\cdot(h-1)
    \geq
    1,
  \end{gather*}
  where the inequality follows from~\eqref{eq:highchebyshev}. Thus, the largest $\Psi_s$ is $\Psi_{h/2}$, from which we get
  \begin{align*}
    \EE[\rn{Y}]
    & \leq
    \frac{h}{2}\cdot\Psi_{s/2}
    =
    \frac{h}{2}\cdot\binom{n}{h/2}\cdot p^{h/2}\cdot(1-p)^{h/2}\cdot(h-1)^{h/2-1}
    \\
    & \leq
    \frac{h}{2}\cdot \frac{n^{h/2}\cdot e^{h/2}}{\sqrt{\pi\cdot h}\cdot(h/2)^{h/2}}
    \cdot p^{h/2}\cdot(1-p)^{h/2}\cdot(h-1)^{h/2-1}
    \\
    & \leq
    \frac{2^{h/2-1}\cdot e^{h/2}\cdot p^{h/2}\cdot(1-p)^{h/2}}{\sqrt{\pi\cdot h}}\cdot n^{h/2},
  \end{align*}
  where the second inequality follows from
  \begin{gather*}
    \binom{n}{h/2} \leq \frac{n^{h/2}}{(h/2)!} \leq \frac{n^{h/2}\cdot e^{h/2}}{\sqrt{\pi\cdot h}\cdot(h/2)^{h/2}},
  \end{gather*}
  using Stirling's Approximation for $(h/2)!$.

  Since $h$ is even, we get
  \begin{align*}
    \PP\bigl[\lvert\rn{X} - p\cdot n\rvert > \delta\cdot p\cdot n\bigr]
    & =
    \PP\bigl[(\rn{X} - p\cdot n)^h >\delta^h\cdot p^h\cdot n^h\bigr]
    =
    \PP\bigl[\rn{Y}^h > \delta^h\cdot p^h\cdot n^h\bigr]
    \\
    & <
    \frac{\EE[\rn{Y}]}{\delta^h\cdot p^h\cdot n^h}
    \leq
    \frac{2^{h/2-1}\cdot e^{h/2}\cdot p^{h/2}\cdot(1-p)^{h/2}}{\sqrt{\pi\cdot h}}\cdot n^{h/2}
    \cdot\frac{1}{\delta^h\cdot p^h\cdot n^h}
    \\
    & =
    \frac{1}{\sqrt{2\cdot\pi\cdot h}}
    \cdot
    \left(\frac{2\cdot e\cdot (1-p)}{\delta^2\cdot p\cdot n}\right)^{h/2},
  \end{align*}
  where the first inequality is Markov's Inequality.
\end{proof}

\begin{definition}[\protect{\cite[Defintion~3.31]{Vad12}}]\label{def:PRG}
  For $h\in\NN_+$, a finite non-empty family $\cF$ of functions of the form $D\to R$ (with same domain $D$ and range $R$) is
  \emph{$h$-wise independent} if for every $D'\subseteq D$ of size $\lvert D'\rvert\leq h$, if $\rn{f}$ is picked uniformly at
  random in $\cF$, then $(\rn{f}(d))_{d\in D'}$ is uniformly distributed in $R^{D'}$.
\end{definition}

\begin{proposition}[\protect{\cite[Construction~3.32 and Proposition~3.33]{Vad12}}]\label{prop:PRG}
  Let $h\in\NN_+$, let $\FF$ be a finite field with $\lvert\FF\rvert\geq h$ and let $\cF$ be the set of all (univariate)
  polynomials of degree at most $h-1$ with coefficients in $\FF$. Then $\cF$ is $h$-wise independent and
  $\lvert\cF\rvert=\lvert\FF\rvert^h$.

  In particular, if $\rn{f}$ is a uniformly at random element of $\cF$, then $(\rn{f}(x))_{x\in\FF}$ is a sequence of $h$-wise
  independent random variables such that each $\rn{f}(x)$ is uniformly distributed in $\FF$.
\end{proposition}

\begin{proof}
  For the first assertion, it suffices to show that for all distinct $x_1,\ldots,x_h\in\FF$, and all $y_1,\ldots,y_h\in\FF$,
  there exists exactly one polynomial $f$ of degree at most $h-1$ such that $f(x_i) = y_i$ for every $i\in[h]$. Existence is by
  Lagrange Interpolation:
  \begin{gather*}
    f(x) \df \sum_{i=1}^h y_i\cdot\prod_{j\in[h]\setminus\{i\}} \frac{x-x_j}{x_i-x_j}
  \end{gather*}
  and uniqueness follows since if $f$ and $\widetilde{f}$ both satisfy this, then $f-\widetilde{f}$ is a polynomial of degree at
  most $h-1$ with $h$ roots, which yields $f = \widetilde{f}$.

  The second assertion follows from the first (and Definition~\ref{def:PRG}).
\end{proof}

\begin{definition}\label{def:PRGbiased}
  For $n,h\in\NN_+$, $p\in[0,1]$ and a finite set $S$, a \emph{pseudorandom generator} of bitlength $n$, independence $h$, bias
  $p$, and seed space $S$ is a function $f\colon S\to\{0,1\}^n$ such that if $\rn{z}$ is picked uniformly at random in $S$, then
  the $n$ coordinates of $f(\rn{z})$ are $h$-wise independent and each coordinate of $f(\rn{z})$ is $p$-Bernoulli.
\end{definition}

The next theorem shows how to algorithmically efficiently compute a pseudorandom generator with a bias that is a dyadic (i.e.,
an element of $\{a/2^n \mid a\in\ZZ,b\in\NN\}$). For a more general version that allows for any rational bias, see
Theorem~\ref{thm:PRGgen} in Appendix~\ref{sec:PRGgen}.

\begin{theorem}\label{thm:PRG}
  Let $n,h,t\in\NN_+$ with $h\leq n\leq 2^t$. Let $\widetilde{p}\in\{0,\ldots,2^t\}$ and let
  $p\df\widetilde{p}/2^t\in\QQ\cap[0,1]$. Then Algorithms~\ref{alg:encodePRG} and~\ref{alg:runPRG} together form a pseudorandom
  generator of bitlength $n$, independence $h$, bias $p$ and seed space $[2^t]^h$ in the sense that:
  \begin{itemize}
  \item Algorithm~\ref{alg:encodePRG} outputs an irreducible polynomial $P\in\FF_2[X]$ of degree $t$.
  \item If $f(Z)_u$ is the output of Algorithm~\ref{alg:runPRG} when given the output $P$ of Algorithm~\ref{alg:encodePRG}, the
    index $u\in[n]$ and $Z\in[2^t]^h$, then $f\colon[2^t]^h\to\{0,1\}^n$ is a pseudorandom generator of bitlength $n$, independence
    $h$, bias $p$ and seed space $[2^t]^h$ (see Definition~\ref{def:PRGbiased}).
  \end{itemize}

  The space-complexities of Algorithms~\ref{alg:encodePRG} and~\ref{alg:runPRG} are
  \begin{align*}
    O(t),
    & &
    O(h\cdot t),
  \end{align*}
  respectively, and the time-complexities are
  \begin{align*}
    O(2^{2\cdot t}\cdot t),
    & &
    O(h\cdot t^2).
  \end{align*}
\end{theorem}

\begin{proof}
  We start by arguing correctness of Algorithm~\ref{alg:encodePRG}.

  \begin{algorithm}[htbp]
    \caption{Computation algorithm of query-into-oracle model with space- and time-complexities
      \begin{align*}
        O(t),
        & &
        O(2^{2\cdot t}\cdot t),
      \end{align*}
      respectively, that returns an irreducible polynomial $P(X)\in\FF_2[X]$ of degree $t$ such that when given to
      Algorithm~\ref{alg:runPRG} produces an oracle for a pseudorandom generator of bitlength $n$, independence $h$, bias $p$
      and seed space $[2^t]^h$. This particular algorithm does not need the values of $n$, $h$ or $p$.}
    \label{alg:encodePRG}
    \DontPrintSemicolon
    \KwIn{A number $t\in\NN_+$.}
    \KwOut{An irreducible polynomial $P(X)\in\FF_2[X]$ of degree $t$.}
    \lIf{$t=1$}{\Return{$X$}}
    Let $d$ be the largest divisor of $t$ smaller than $t$.\;
    $m\assign (2^t-1)/(2^d-1)$\;
    Find the lexicographically first polynomial $P(X)\in\FF_2[X]$ of degree $t$ that divides
    \begin{gather*}
      \widetilde{P}(X) \df \sum_{j=0}^m X^j \in \FF_2[X].
    \end{gather*}
    \;\label{alg:encodePRG:P}
    \Return{$P$}
  \end{algorithm}

  \begin{algorithm}[htbp]
    \caption{Oracle algorithm of query-into-oracle model with space- and time-complexities
      \begin{align*}
        O(h\cdot t),
        & &
        O(h\cdot t^2),
      \end{align*}
      respectively, that when given $P$ computed via Algorithm~\ref{alg:encodePRG}, a seed $Z\in[2^t]^h$ and an index $u\in[n]$,
      outputs a bit $X_u$ such that the corresponding function $Z\mapsto (X_u)_{u\in[n]}$ is a pseudorandom generator of
      bitlength $n$, independence $h$, bias $p$ and seed space $[2^t]^h$.}
    \label{alg:runPRG}
    \DontPrintSemicolon
    \KwIn{An irreducible polynomial $P(X)\in\FF_2[X]$, an integer $\widetilde{p}\in\{0,\ldots,2^t\}$, where $t\df\deg(P)\geq 1$,
      numbers $h,n\in\NN_+$ such that $2^t\geq n\geq h$, a seed $Z\in[2^t]^h$, and an index $u\in[n]$.}
    \KwOut{A bit $X_u\in\{0,1\}$ such that when $Z$ is random with uniform distribution over $[2^t]^n$, then $(X_u)_{u=1}^n$ is a
      collection of $p$-Bernoulli bits that is $h$-wise independent, where $p\df\widetilde{p}/2^t$.}
    \lIf{$\widetilde{p}=0$}{\Return{$0$}}
    \lIf{$\widetilde{p}=2^t$}{\Return{$1$}}
    Write each $Z_j-1\in\{0,\ldots,2^t-1\}$ in base $2$ as
    \begin{gather*}
      Z_j-1 = c^j_{t-1}c^j_{t-2}\cdots c^j_0.
    \end{gather*}
    \;
    Write $u-1\in\{0,\ldots,n-1\}\subseteq\{0,\ldots,2^t-1\}$ in base $2$ as
    \begin{gather*}
      u-1 = b_{t-1}b_{t-2}\cdots b_0.
    \end{gather*}
    \;
    Interpreting the numbers $c^j_\ell$ and $b_\ell$ as elements of $\FF_2$, compute
    \begin{gather*}
      R(X)
      \df
      \left(
      \sum_{j\in[h]}
      \left(
      \left(\sum_{\ell=0}^{t-1} c^j_\ell\cdot X^\ell\right)
      \cdot\left(\sum_{\ell=0}^{t-1} b_\ell\cdot X^\ell\right)^{j-1}
      \right)
      \right)
      \bmod P(X).
    \end{gather*}
    \;\label{alg:runPRG:poly}
    Writing $R(X) = \sum_{\ell=0}^{t-1} r_\ell\cdot X^\ell$ and interpreting the coefficients $r_\ell\in\FF_2$ as
    elements of $\{0,1\}$, let
    \begin{gather*}
      v\df\sum_{\ell=0}^{t-1} r_\ell\cdot 2^\ell\in\{0,\ldots,2^t-1\}.
    \end{gather*}
    \;
    \Return{$\One[v\leq\widetilde{p}]$}
  \end{algorithm}

  It is clear that if $t=1$ then $P(X)=X\in\FF_2[X]$ is an irreducible polynomial of degree $t=1$. Suppose then that
  $t\geq 2$. We know that for any sought irreducible polynomial $P(X)\in\FF_2[X]$ of degree $t$, the quotient
  $\FF_2[X]/(P)$ is isomorphic to the field $\FF_{2^t}$ of size $2^t$.

  In turn, we know that in the algebraic closure $\overline{\FF}_2$ of $\FF_2$, the field $\FF_{2^t}$ is exactly
  given by the $2^t$ distinct roots of the polynomial $R_t(X)\df X^{2^t} - X$ (note that the derivative of $R$ is
  $-1$, so all its roots have multiplicity $1$). Furthermore, the field $\FF_{2^t}$ contains as subfields exactly the
  fields of the form $\FF_{2^f}$, where $f\in\NN_+$ is a divisor of $t$; in fact, $\FF_{2^f}$ consists exactly of the $2^f$
  roots of $R_f(X)\df X^{2^f} - X$, which are themselves roots of $R_t(X)$. Now, the algorithm takes $d\in\NN_+$ as the
  largest divisor of $t$ smaller than $t$ and since
  \begin{gather*}
    R_t(X)
    =
    R_f(X)\cdot\widetilde{P}(X),
  \end{gather*}
  where
  \begin{align*}
    \widetilde{P}(X)
    & \df
    \sum_{j=0}^m X^j,
    &
    m
    & \df
    \frac{2^t-1}{2^d-1}
    =
    \sum_{j=0}^{t/d-1} 2^{j\cdot d},
  \end{align*}
  it follows that the roots of $\widetilde{P}$ in the algebraic closure $\overline{\FF}_2$ are exactly the elements of
  $\FF_{2^t}\setminus\FF_{2^d}$ and since $\FF_{2^d}$ is the largest proper subfield\footnote{This is because the subfields of
  $\FF_{2^t}$ are exactly those of the form $\FF_{2^m}$ with $t$ divisible by $m$; one direction is because $\FF_{2^t}$ must be
  an $\FF_{2^m}$-vector space so $2^t$ has to be a power of $2^m$, hence $m$ must divide $t$; and the other direction is because
  when $m$ divides $t$, we can factor $X^{2^t}-X = (X^{2^m}-X)\cdot\sum_{j=0}^{(2^t-1)/(2^d-1)} X^j$.} of $\FF_{2^t}$, all
  irreducible factors of $\widetilde{P}$ over $\FF_2[X]$ must be of degree $t$, which means that any polynomial in $\FF_2[X]$ of
  degree $t$ that divides $\widetilde{P}$ must be irreducible. Thus, Algorithm~\ref{alg:encodePRG} eventually finds one such
  irreducible $P(X)\in\FF_2[X]$ of degree $t$.

  \smallskip

  We now analyze the space-complexity of Algorithm~\ref{alg:encodePRG}. We need $O(\log(t+1))$ space to compute $d$ and we need
  at most $O(t)$ space to compute $m$.

  Regarding $P$ and $\widetilde{P}$, for the former, we need space at most $O(t)$ to store the current $P$ being tested. For
  $\widetilde{P}$, we don't actually need to store it in memory: we know that all of its coefficients are $1$, so we only need
  to keep track of its degree $m$. Furthermore, when testing divisibility of $\widetilde{P}$ by $P$, when performing long
  division, we only need to keep track of the topmost $t$ coefficients of the remainder as the standard long division algorithm
  progresses (again, because all other coefficients of the remainder will be the ones of $\widetilde{P}$, which are all $1$).
  Thus, the space-complexity of line~\eqref{alg:encodePRG:P} is at most $O(t)$. Thus, the total space-complexity of
  Algorithm~\ref{alg:encodePRG} is at most $O(t)$.

  \smallskip

  For time-complexity of Algorithm~\ref{alg:encodePRG}, note that to compute the largest divisor $d$ that is less than $t$,
  we can instead compute the smallest non-trivial divisor $d'\in\NN$ since $d = t/d'$ (this still only takes space
  $O(\log(t+1))$). We only need to search for $d'$ in $\{2,\ldots,\floor{\sqrt{t}}\}$, so the time-complexity of computing
  $d$ is at most
  \begin{gather*}
    O\Bigl(\sqrt{t}\cdot\bigl(\log(t+1)\bigr)^2\Bigr)
    \leq
    O(t).
  \end{gather*}
  For the computation of $m$, the time-complexity is easily seen to be at most $O(t^2)$. Finally, the time-complexity of
  line~\ref{alg:encodePRG:P} amounts to the time it takes to enumerate all the potential $P$ and the time it takes to test
  divisibility of $\widetilde{P}$ by it. There are a total of $2^t$ many $P$ and to check divisibility of $\widetilde{P}$ by one
  such $P$, we can perform long polynomial division, which takes time proportional to the product of degrees, that is, $O(m\cdot
  t)\leq O(2^t\cdot t)$. Taking into account the enumeration of the at most $2^t$ many $P$, the total time-complexity of
  Algorithm~\ref{alg:encodePRG} is at most $O(2^{2\cdot t}\cdot t)$.

  \medskip

  We now argue correctness of Algorithm~\ref{alg:runPRG}. The result is completely obvious if $\widetilde{p}=0$ (in which case
  $p=0$) or if $\widetilde{p}=2^t$ (in which case $p=1$), so we assume that $\widetilde{p}\in[2^t-1]$. We fix $P$, $h$ and $n$
  of the input of the algorithm, we let $\rn{Z}$ be picked uniformly at random in $[2^t]^h$ and we make the following
  definitions:
  \begin{itemize}
  \item $\rn{X}_u$ as the bit that Algorithm~\ref{alg:runPRG} outputs when further given $\rn{Z}$ and $u\in[n]$.
  \item For a variable of Algorithm~\ref{alg:runPRG}, we will write it in bold face if the variable depends on
    $\rn{Z}$ (which means it will be a random variable that is $\rn{Z}$-measurable) and we will write a subscript $u$ if it
    depends on the value of $u$. Under these conventions, it is easy to check the following dependencies:
    \begin{align*}
      \rn{c}^j_\ell, & &
      (b_\ell)_u, & &
      \rn{R}(X)_u, & &
      \rn{v}_u.
    \end{align*}
  \end{itemize}

  We will now list several claims, culminating in the final claim of correctness, and provide a short proof of these claims
  later:
  \begin{enumerate}[wide, label={\textbf{Claim~\thetheorem.\arabic*.}}, ref={\thetheorem.\arabic*}]
  \item\label{clm:PRG:c} $(\rn{c}^j_\ell)_{(j,\ell)\in[h]\times\{0,\ldots,t-1\}}$ is uniformly distributed
    in the set $\{0,1\}^{[h]\times\{0,\ldots,t-1\}}$. 
  \item\label{clm:PRG:wR} If we define the (random) polynomial
    \begin{gather*}
      \rn{\widetilde{R}}(Y)\in \frac{\FF_2[X]}{(P)}[Y]
    \end{gather*}
    as the polynomial in the variable $Y$ with coefficients in the field $\FF_2[X]/(P)\cong\FF_{2^t}$ given by
    \begin{gather*}
      \rn{\widetilde{R}}(Y)
      \df
      \sum_{j\in[h]}
      \left(
      \left(\sum_{\ell=0}^{t-1} \rn{c}^j_\ell\cdot X^\ell\right)\bmod P(X)
      \right)\cdot Y^{j-1},
    \end{gather*}
    then $\rn{\widetilde{R}}(Y)$ is a uniformly at random polynomial in $(\FF_2[X]/(P))[Y]$ of degree at most $h-1$.
  \item\label{clm:PRG:R} For $\rn{\widetilde{R}}(Y)$ as in Claim~\ref{clm:PRG:wR}, we have
    \begin{gather*}
      \rn{R}(X)_u
      =
      \rn{\widetilde{R}}\left(\sum_{\ell=0}^{t-1} (b_\ell)_u\cdot X^\ell \bmod P(X)\right).
    \end{gather*}
  \item\label{clm:PRG:Runif} $(\rn{R}(X)_u)_{u\in[n]}$ is a sequence of $h$-wise independent random elements of the field
    $\FF_2[X]/(P)\cong\FF_{2^t}$ such that each $\rn{R}(X)_u$ is uniformly distributed in $\FF_2[X]/(P)$.
  \item\label{clm:PRG:v} $(\rn{v}_u)_{u\in[n]}$ is a sequence of $h$-wise independent random elements of
    $\{0,\ldots,2^t-1\}$ such that each $\rn{v}_u$ is uniformly distributed in $\{0,\ldots,2^t-1\}$.
  \item\label{clm:PRG:X} Algorithm~\ref{alg:runPRG} is correct, i.e., $(\rn{X}_u)_{u\in[n]}$ is a sequence of $h$-wise
    independent random variables such that each $\rn{X}_u$ is $p$-Bernoulli.
  \end{enumerate}

  Claim~\ref{clm:PRG:c} follows from the fact that $\rn{Z}$ is uniformly distributed in $[2^t]^h$. Claim~\ref{clm:PRG:wR}
  follows directly from Claim~\ref{clm:PRG:c} and the definition of $\rn{\widetilde{R}}$. Claim~\ref{clm:PRG:R} follows since
  (polynomial) modulus satisfies
  \begin{align*}
    (A+B)\bmod P & = (A\bmod P + B\bmod P)\bmod P,
    \\
    (A\cdot B)\bmod P & = \bigl((A\bmod P)\cdot (B\bmod P)\bigr)\bmod P,
  \end{align*}

  Claim~\ref{clm:PRG:Runif} follows from Claims~\ref{clm:PRG:wR} and~\ref{clm:PRG:R} along with Proposition~\ref{prop:PRG}:
  Claim~\ref{clm:PRG:wR} says that $\rn{\widetilde{R}}(Y)$ is a uniformly at random polynomial of degree at most $h-1$ with
  coefficients in the field $\FF_2[X]/(P)\cong\FF_{2^t}$ and Claim~\ref{clm:PRG:R} says that $\rn{R}(X)_u$ is the evaluation
  of this polynomial at the point
  \begin{gather*}
    \sum_{\ell=0}^{t-1} (b_\ell)_u\cdot X^\ell \bmod P(X)\in\frac{\FF_2[X]}{(P)}\cong\FF_{2^t}
  \end{gather*}
  and since different values of $u$ lead to different values of $((b_\ell)_u)_{\ell\in\{0,\ldots,t-1\}}$, changing $u$ changes
  the point of the field $\FF_2[X]/(P)\cong\FF_{2^t}$ on which we are evaluating the polynomial, so Proposition~\ref{prop:PRG}
  gives us the $h$-wise independence of $(\rn{R}(X))_{u\in[n]}$ and the fact that each $\rn{R}(X)_u$ is uniformly distributed in
  $\FF_2[X]/(P)\cong\FF_{2^t}$.

  Claim~\ref{clm:PRG:v} follows directly from Claim~\ref{clm:PRG:Runif} and the definition of the $\rn{v}_u$.

  Finally, Claim~\ref{clm:PRG:X} follows directly from Claim~\ref{clm:PRG:v} and the fact that $\rn{X}_u =
  \One[\rn{v}_u\leq\widetilde{p}]$ so $\PP[\rn{X}_u=1] = \widetilde{p}/2^t = p$.

  \smallskip

  We now analyze the space-complexity of Algorithm~\ref{alg:runPRG}. Storing the variables
  \begin{align*}
    \bigl((c^j_\ell)_{\ell=0}^{t-1}\bigr)_{j\in[h]}, & &
    (b_\ell)_{\ell=0}^{t-1}, & &
    R, & &
    v
    \intertext{takes space}
    O(h\cdot t), & &
    O(t), & &
    O(t), & &
    O(t),
  \end{align*}
  respectively. Furthermore, it is straightforward to check that all computations require extra space at most proportional to
  the size of the variables being computed (see more details about line~\ref{alg:runPRG:poly} below), so the space-complexity
  of Algorithm~\ref{alg:runPRG} is at most $O(h\cdot t)$.

  \smallskip

  Finally, we analyze the time-complexity of Algorithm~\ref{alg:runPRG}. It is straightforward to check that the most
  time-expensive step is the computation of line~\ref{alg:runPRG:poly}. This can be done efficiently by inductively computing
  for each $\widetilde{h}\in\{0,\ldots,h\}$ the values:
  \begin{gather*}
    \begin{aligned}
      A_{-1}(X) & \df 0, &
      B_0(X) & \df 1,
    \end{aligned}
    \\
    B_{\widetilde{h}}(X)
    \df
    \left(\sum_{\ell=0}^{t-1} b_\ell\cdot X^\ell\right)\cdot B_{\widetilde{h}-1}(X) \bmod P(X),
    \\
    A_{\widetilde{h}}(X)
    \df
    \left(
    A_{\widetilde{h}-1}(X) + \left(\sum_{\ell=0}^{t-1} c^j_\ell\cdot X^\ell\right)\cdot B_{\widetilde{h}}(X)
    \right) \bmod P(X).
  \end{gather*}
  At each step we need to compute two polynomial moduli by $P(X)$ of polynomials in $\FF_2[X]$ of degree at most $2\cdot t$, two
  polynomial multiplications in $\FF_2[X]$ of polynomials of degree at most $t$ each and a polynomial addition in $\FF_2[X]$ of
  two polynomials of degree at most $2\cdot t$. Thus, the time-complexity of line~\ref{alg:runPRG:poly}, hence also of
  Algorithm~\ref{alg:runPRG:poly} is at most $O(h\cdot t^2)$.
\end{proof}

\begin{proposition}\label{prop:PRGsubset}
  Let $s,d,n,m\in\NN_+$, let $G$ be a graph with $\lvert G\rvert=n$, $\VC(G)\leq d$, let $W\subseteq V(G)$ be a
  set with $\lvert W\rvert=s$, let $\widetilde{\delta},\widetilde{\xi},c_{\adj},p\in(0,1]$, let $c_{\samp}\in(0,1/2)$, let
  \begin{gather*}
    \widetilde{\epsilon} \df (1-c_{\samp})\cdot\widetilde{\xi}
  \end{gather*}
  and let $\cH\subseteq\cP(W)$ be a family of subsets of $W$ such that $\lvert H\rvert\leq\widetilde{\epsilon}\cdot s$.

  Suppose $h\in\NN_+$ is an even integer satisfying
  \begin{gather}\label{eq:PRGsubset:highchebyshev}
    \floor{\widetilde{\epsilon}\cdot s}
    \geq
    \begin{dcases*}
      \frac{(h-1)^2}{p\cdot(1-p)}, & if $p < 1$,\\
      0, & if $p=1$.
    \end{dcases*}
  \end{gather}
  For $X\in\{0,1\}^W$, define the set
  \begin{gather*}
    W(X) \df X^{-1}(1) \df \{u\in W \mid X_w=1\}
  \end{gather*}
  and let $E(X)$ be the statement
  \begin{gather*}
    \bigl\lvert\lvert W(X)\rvert - p\cdot s\bigr\rvert
    \leq
    c_{\adj}\cdot c_{\samp}\cdot\widetilde{\xi}\cdot p\cdot s.
  \end{gather*}
  
  Let further $\rn{X}=(\rn{X}_u)_{u\in W}$ be a sequence of $h$-wise independent random variables such that each $\rn{X}_u$ is
  $p$-Bernoulli, let $f\colon S\to\{0,1\}^s$ be a pseudorandom generator of bitlength $s$, independence $h$, bias $p$, and seed
  space $S$ and let $b\colon W\to [s]$ be a bijection. Furthermore, for $z\in S$, we let $X^{f(z),b}\in\{0,1\}^W$ be given by
  $X^{f(z),b}_u\df f(z)_{b(u)}$
  
  Then the following hold:
  \begin{enumerate}
  \item\label{prop:PRGsubset:size} We have
    \begin{gather*}
      \PP\bigl[E(\rn{X})\bigr]
      >
      1 - \frac{1}{\sqrt{2\cdot\pi\cdot h}}
      \cdot
      \left(\frac{2\cdot e\cdot (1-p)}{c_{\adj}^2\cdot c_{\samp}^2\cdot\widetilde{\xi}^2\cdot p\cdot s}\right)^{h/2}.
    \end{gather*}
  \item\label{prop:PRGsubset:conc} We have
    \begin{align*}
      \MoveEqLeft
      \PP\bigl[
        E(\rn{X})\land
        \forall H\in\cH,
        \lvert H\cap W(\rn{X})\rvert
        \leq
        \widetilde{\xi}\cdot\lvert W(\rn{X})\rvert
        \bigr]
      \\
      & >
      \begin{multlined}[t]
        1 -
        \frac{1}{\sqrt{2\cdot\pi\cdot h}}
        \cdot
        \left(\frac{2\cdot e\cdot (1-p)}{c_{\samp}^2\cdot\widetilde{\xi}^2\cdot p}\right)^{h/2}
        \\
        \cdot
        \left(
        \left(\frac{1}{c_{\adj}^2\cdot s}\right)^{h/2}
        +
        \lvert\cH\rvert\cdot\left(\frac{(1-c_{\samp})^2}{(1-c_{\adj})^2\cdot\floor{\widetilde{\epsilon}\cdot s}}\right)^{h/2}
        \right).
      \end{multlined}
    \end{align*}
  \item\label{prop:PRGsubset:good} If $W$ is $\widetilde{\epsilon}$-good in $G$, then
    \begin{align*}
      \MoveEqLeft
      \PP\bigl[
        E(\rn{X})\land
        W(\rn{X})\text{ is $\widetilde{\xi}$-good in $G$}
        \bigr]
      \\
      & >
      \begin{multlined}[t]
        1 -
        \frac{1}{\sqrt{2\cdot\pi\cdot h}}
        \cdot
        \left(\frac{2\cdot e\cdot (1-p)}{c_{\samp}^2\cdot\widetilde{\xi}^2\cdot p}\right)^{h/2}
        \\
        \cdot
        \left(
        \left(\frac{1}{c_{\adj}^2\cdot s}\right)^{h/2}
        +
        n\cdot\left(\frac{(1-c_{\samp})^2}{(1-c_{\adj})^2\cdot\floor{\widetilde{\epsilon}\cdot s}}\right)^{h/2}
        \right).
      \end{multlined}
    \end{align*}

    In particular, if
    \begin{gather*}
      \frac{1}{\sqrt{2\cdot\pi\cdot h}}
      \cdot
      \left(\frac{2\cdot e\cdot (1-p)}{c_{\samp}^2\cdot\widetilde{\xi}^2\cdot p}\right)^{h/2}
      \cdot
      \left(
      \left(\frac{1}{c_{\adj}^2\cdot s}\right)^{h/2}
      +
      n\cdot\left(\frac{(1-c_{\samp})^2}{(1-c_{\adj})^2\cdot\floor{\widetilde{\epsilon}\cdot s}}\right)^{h/2}
      \right)
      \leq
      1,
    \end{gather*}
    then there exists $z\in S$ such that $W(X^{f(z),b})$ is $\widetilde{\xi}$-good in $G$ and $E(X^{f(z),b})$ holds.
  \item\label{prop:PRGsubset:exc} If $W$ is $(\widetilde{\delta},\widetilde{\epsilon})$-excellent in $G$, then
    \begin{align*}
      \MoveEqLeft
      \PP\bigl[
        E(\rn{X})\land
        W(\rn{X})\text{ is $(\widetilde{\delta},\widetilde{\xi})$-excellent in $G$}
        \bigr]
      \\
      & >
      \begin{multlined}[t]
        1 -
        \frac{1}{\sqrt{2\cdot\pi\cdot h}}
        \cdot
        \left(\frac{2\cdot e\cdot (1-p)}{c_{\samp}^2\cdot\widetilde{\xi}^2\cdot p}\right)^{h/2}
        \\
        \cdot
        \left(
        \left(\frac{1}{c_{\adj}^2\cdot s}\right)^{h/2}
        +
        (1+s^d)\cdot\left(\frac{(1-c_{\samp})^2}{(1-c_{\adj})^2\cdot\floor{\widetilde{\epsilon}\cdot s}}\right)^{h/2}
        \right).
      \end{multlined}
    \end{align*}

    In particular, if
    \begin{gather*}
      \frac{1}{\sqrt{2\cdot\pi\cdot h}}
      \cdot
      \left(\frac{2\cdot e\cdot (1-p)}{c_{\samp}^2\cdot\widetilde{\xi}^2\cdot p}\right)^{h/2}
      \cdot
      \left(
      \left(\frac{1}{c_{\adj}^2\cdot s}\right)^{h/2}
      +
      (1+s^d)\cdot\left(\frac{(1-c_{\samp})^2}{(1-c_{\adj})^2\cdot\floor{\widetilde{\epsilon}\cdot s}}\right)^{h/2}
      \right)
      \leq
      1,
    \end{gather*}
    then there exists $z\in S$ such that $W(X^{f(z),b})$ is $(\widetilde{\delta},\widetilde{\xi})$-excellent in $G$ and
    $E(X^{f(z),b})$ holds.
  \item\label{prop:PRGsubset:atom} If $W$ is a $(\widetilde{\delta},\widetilde{\epsilon},m)$-atom of $G$, then
    \begin{align*}
      \MoveEqLeft
      \PP\bigl[
        E(\rn{X})\land
        W(\rn{X})\text{ is a $(\widetilde{\delta},\widetilde{\xi},m)$-atom of $G$}
        \bigr]
      \\
      & >
      \begin{multlined}[t]
        1 -
        \frac{1}{\sqrt{2\cdot\pi\cdot h}}
        \cdot
        \left(\frac{2\cdot e\cdot (1-p)}{c_{\samp}^2\cdot\widetilde{\xi}^2\cdot p}\right)^{h/2}
        \\
        \cdot
        \left(
        \left(\frac{1}{c_{\adj}^2\cdot s}\right)^{h/2}
        +
        n^m\cdot\left(\frac{(1-c_{\samp})^2}{(1-c_{\adj})^2\cdot\floor{\widetilde{\epsilon}\cdot s}}\right)^{h/2}
        \right).
      \end{multlined}
    \end{align*}

    In particular, if
    \begin{gather*}
      \frac{1}{\sqrt{2\cdot\pi\cdot h}}
      \cdot
      \left(\frac{2\cdot e\cdot (1-p)}{c_{\samp}^2\cdot\widetilde{\xi}^2\cdot p}\right)^{h/2}
      \cdot
      \left(
      \left(\frac{1}{c_{\adj}^2\cdot s}\right)^{h/2}
      +
      n^m\cdot\left(\frac{(1-c_{\samp})^2}{(1-c_{\adj})^2\cdot\floor{\widetilde{\epsilon}\cdot s}}\right)^{h/2}
      \right)
      \leq
      1,
    \end{gather*}
    then there exists $z\in S$ such that $W(X^{f(z),b})$ is a $(\widetilde{\delta},\widetilde{\xi},m)$-atom of $G$ and
    $E(X^{f(z),b})$ holds.
  \end{enumerate}
\end{proposition}

\begin{proof}
  The result is trivial if $p=1$, so we assume $p < 1$.
  
  For item~\ref{prop:PRGsubset:size}, we simply note that
  \begin{align*}
    \PP\bigl[E(\rn{X})\bigr]
    & =
    \PP\Bigl[
      \bigl\lvert\lvert W(\rn{X})\rvert - p\cdot s\bigr\rvert
      \leq
      c_{\adj}\cdot c_{\samp}\cdot\widetilde{\xi}\cdot p\cdot s
      \Bigr]
    \\
    & =
    \PP\left[\left\lvert\lvert\sum_{u\in W}\rn{X}_u - p\cdot s\right\rvert
      \leq
      c_{\adj}\cdot c_{\samp}\cdot\widetilde{\xi}\cdot p\cdot s
      \right]
    \\
    & >
    1 - \frac{1}{\sqrt{2\cdot\pi\cdot h}}
    \cdot
    \left(\frac{2\cdot e\cdot (1-p)}{c_{\adj}^2\cdot c_{\samp}^2\cdot\widetilde{\xi}^2\cdot p\cdot s}\right)^{h/2},
  \end{align*}
  where the last inequality follows from Proposition~\ref{prop:highchebyshev} since~\eqref{eq:PRGsubset:highchebyshev} implies
  \begin{gather*}
    s\geq \floor{\widetilde{\epsilon}\cdot s}\geq \frac{(h-1)^2}{p\cdot(1-p)}.    
  \end{gather*}
  
  \medskip

  We now prove item~\ref{prop:PRGsubset:conc}. For $H\in\cH$, since $\lvert H\rvert\leq\widetilde{\epsilon}\cdot s$, we can fix
  a set $H'\subseteq W$ of size $\lvert H'\rvert=\floor{\widetilde{\epsilon}\cdot s}$ with $H\subseteq H'$. We then define the
  event $E_H(\rn{X})$ as
  \begin{gather*}
    \lvert H'\cap W(\rn{X})\rvert
    \leq
    (1+\delta)\cdot p\cdot\floor{\widetilde{\epsilon}\cdot s},
  \end{gather*}
  where
  \begin{gather*}
    \delta
    \df
    \frac{\widetilde{\xi}^2\cdot(1-c_{\adj}\cdot c_{\samp})}{\widetilde{\epsilon}} - 1
    =
    \frac{(1-c_{\adj})\cdot c_{\samp}\cdot\widetilde{\xi}}{1-c_{\samp}}
    \in
    (0,1).
  \end{gather*}

  We claim that the conjunction $E_H(\rn{X})\land E(\rn{X})$ implies
  \begin{gather*}
    \lvert H\cap W(\rn{X})\rvert\leq\widetilde{\xi}\cdot\lvert W(\rn{X})\rvert.
  \end{gather*}
  Indeed, note that within this conjunction, we have
  \begin{align*}
    \lvert H\cap W(\rn{X})\rvert
    & \leq
    \lvert H'\cap W(\rn{X})\rvert
    \leq
    \frac{(1-c_{\adj}\cdot c_{\samp})\cdot\widetilde{\xi}^2}{\widetilde{\epsilon}}\cdot p\cdot\floor{\widetilde{\epsilon}\cdot s}
    \\
    & \leq
    (1-c_{\adj}\cdot c_{\samp})\cdot\widetilde{\xi}^2\cdot p\cdot s
    \leq
    \widetilde{\xi}\cdot\lvert W(\rn{X})\rvert,
  \end{align*}
  where the first inequality follows since $H\subseteq H'$, the second inequality follows from $E_H(\rn{X})$ and the last
  inequality follows since $E(\rn{X})$ implies
  \begin{gather*}
    \lvert W(\rn{X})\rvert \geq (1-c_{\adj}\cdot c_{\samp})\cdot\widetilde{\xi}\cdot p\cdot s.
  \end{gather*}

  On the other hand, note that
  \begin{align*}
    \PP\bigl[E_H(\rn{X})\bigr]
    & =
    \PP\left[
      \sum_{u\in H'} \rn{X}_u
      \leq
      (1+\delta)\cdot p\cdot\lvert H'\rvert
      \right]
    \geq
    \PP\left[
      \left\lvert\sum_{u\in H'} \rn{X}_u - p\cdot\lvert H'\rvert\right\rvert
      \leq
      \delta\cdot p\cdot\lvert H'\rvert
      \right]
    \\
    & >
    1
    -
    \frac{1}{\sqrt{2\cdot\pi\cdot h}}
    \cdot
    \left(\frac{2\cdot e\cdot (1-p)}{\delta^2\cdot p\cdot\lvert H'\rvert}\right)^{h/2}
    \\
    & =
    1
    -
    \frac{1}{\sqrt{2\cdot\pi\cdot h}}
    \cdot
    \left(\frac{
      2\cdot e\cdot (1-p)\cdot(1-c_{\samp})^2
    }{
      (1-c_{\adj})^2\cdot c_{\samp}^2\cdot\widetilde{\xi}^2\cdot p\cdot\floor{\widetilde{\epsilon}\cdot s}
    }
    \right)^{h/2}.
  \end{align*}

  Applying a union bound over all $E_H(\rn{X})$ for $H\in\cH$ and also over the event $E(\rn{X})$ whose probability is bounded
  in item~\ref{prop:PRGsubset:size} yields the result.

  \medskip

  Before we prove the other items, let us note that their final assertions on the existence of $z\in S$ such that
  $E(X^{f(z),b})$ holds and $W(X^{f(z),b})$ is $\widetilde{\xi}$-good, $(\widetilde{\delta},\widetilde{\xi})$-excellent or a
  $(\widetilde{\delta},\widetilde{\xi},m)$-atom follows from the fact that the assumption in each item is what makes the
  corresponding probability be positive, and then the fact that $f$ is a pseudorandom generator of bitlength $s$, independence
  $h$, bias $p$, and seed space $S$, means that we taking $\rn{z}$ uniformly at random in $S$, the random variables
  $\rn{X}^{f(\rn{z}),b}$ satisfy the hypotheses hence with positive probability yield the desired event; this means that at
  least one $z\in S$ must yield the desired event.

  \medskip

  It remains to prove the probability bounds of the remaining items. These proofs are mostly analogous to the proof of
  Proposition~\ref{prop:randsubset}, but we spell them out anyway:

  We now prove item~\ref{prop:PRGsubset:good}. Since $W$ is $\widetilde{\epsilon}$-good in $G$, we know that for every $v\in
  V(G)$, there exists an $\widetilde{\epsilon}$-majority opinion $t_G^{\widetilde{\epsilon}}(v,W)$ of $W$ with respect to $v$ in
  $G$, which means that taking
  \begin{gather*}
    \cH \df \{W\cap N_G^{1-t_G^{\widetilde{\epsilon}}(v,W)}(v) \mid v\in V(G)\},
  \end{gather*}
  we have $\lvert H\rvert\leq\widetilde{\epsilon}\cdot s$ for every $H\in\cH$. On the other hand, it is clear that
  $\lvert\cH\rvert\leq\lvert G\rvert = n$, so applying item~\ref{prop:PRGsubset:conc}, it then suffices to show that if
  \begin{gather*}
    \lvert H\cap W(\rn{X})\rvert
    \leq
    \widetilde{\xi}\cdot\lvert W(\rn{X})\rvert,
  \end{gather*}
  for every $H\in\cH$, then $W(\rn{X})$ is $\widetilde{\epsilon}$-good in $G$. But indeed, for $v\in V(G)$, we note that the set
  $H\df W\cap N_G^{1-t_G^{\widetilde{\epsilon}}(v,W)}$ is in $\cH$, so the above implies
  \begin{gather*}
    \lvert N_G^{1-t_G^{\widetilde{\epsilon}}(v,W)}(v)\cap W(\rn{X})\rvert
    \leq
    \widetilde{\xi}\cdot\lvert W(\rn{X})\rvert,
  \end{gather*}
  i.e., $t_G^{\widetilde{\epsilon}}(v,W)$ is also a $\widetilde{\xi}$-majority opinion of $W(\rn{X})$ with respect to $v$ in
  $G$, hence $W(\rn{X})$ is $\widetilde{\xi}$-good in $G$.

  \medskip

  For item~\ref{prop:PRGsubset:exc}, since $W$ is $(\widetilde{\delta},\widetilde{\epsilon})$-excellent in $G$, for every
  $\widetilde{\delta}$-good set $U$, there exists a $(\widetilde{\delta},\widetilde{\epsilon})$-majority opinion
  $t^{\widetilde{\delta},\widetilde{\epsilon}}_G(W,U)$ of $W$ with respect to $U$ in $G$, which means that taking
  \begin{gather*}
    \cH
    \df
    \{W\cap N_{G,\widetilde{\delta}}^{1-t_G^{\widetilde{\delta},\widetilde{\epsilon}}(W,U)}(U)
    \mid U\text{ is $\widetilde{\delta}$-good in $G$}\},
  \end{gather*}
  where
  \begin{gather*}
    N_{G,\widetilde{\delta}}^b(U)
    \df
    \{v\in V(G) \mid t_G^{\widetilde{\delta}}(v,U) = b\}
  \end{gather*}
  we have $\lvert H\rvert\leq\widetilde{\epsilon}\cdot s$ for every $H\in\cH$.

  On the other hand, for
  \begin{gather*}
    \cH_{G,\widetilde{\delta}}
    \df
    N_{G,\widetilde{\delta}}^1(U)
    \mid U\text{ is $\widetilde{\delta}$-good in $G$}\},
  \end{gather*}
  Lemma~\ref{lem:VCgood} guarantees that $\VC(\cH_{G,\widetilde{\delta}})\leq d$ (as $\widetilde{\delta}<1/(d+1)$ and
  $\VC(G)\leq d$). Thus, we get
  \begin{gather*}
    \lvert\cH\rvert
    \leq
    2\cdot\lvert\cH_{G,\widetilde{\delta}}\rest_W\rvert
    \leq
    2\cdot\sum_{i=0}^d \binom{s}{i}
    \leq
    2\cdot(1+s^d),
  \end{gather*}
  where the first inequality follows from the Sauer--Shelah--Perles Lemma (Lemma~\ref{lem:SSP}).

  Applying item~\ref{prop:PRGsubset:conc}, it then suffices to show that if
  \begin{gather*}
    \lvert H\cap W(\rn{X})\rvert
    \leq
    \widetilde{\xi}\cdot\lvert W(\rn{X})\rvert,
  \end{gather*}
  for every $H\in\cH$, then $W(\rn{X})$ is $(\widetilde{\delta},\widetilde{\xi})$-excellent in $G$. But indeed, for a
  $\widetilde{\delta}$-good set $U$ in $G$, we note that the set $H\df W\cap
  N_{G,\widetilde{\delta}}^{1-t_G^{\widetilde{\delta},\widetilde{\epsilon}}(W,U)}(U)$ is in $\cH$, so the above implies
  \begin{gather*}
    \lvert N_{G,\widetilde{\delta}}^{1-t_G^{\widetilde{\delta},\widetilde{\epsilon}}(W,U)}(U)\cap W(\rn{X})\rvert
    \leq
    \widetilde{\xi}\cdot\lvert W(\rn{X})\rvert,
  \end{gather*}
  i.e., $t_G^{\widetilde{\delta},\widetilde{\epsilon}}(W,U)$ is also a $(\widetilde{\delta},\widetilde{\xi})$-majority opinion
  of $W(\rn{X})$ with respect to $U$ in $G$, so $W(\rn{X})$ is $(\widetilde{\delta},\widetilde{\xi})$-excellent in $G$.

  \medskip

  Finally, for item~\ref{prop:PRGsubset:atom}, since $W$ is a $(\widetilde{\delta},\widetilde{\epsilon},m)$-atom of $G$, we
  know that it is not $(\widetilde{\delta},\widetilde{\epsilon})$-split by any $m$-tuple $x\in V(G)^m$ in $G$, which means that
  there exists $b_x\in\{0,1\}$ such that
  \begin{gather*}
    \lvert N_{G,\widetilde{\delta}}^{b_x}(x)\cap W\rvert\leq\widetilde{\epsilon}\cdot\lvert W\rvert,
  \end{gather*}
  where
  \begin{gather*}
    N_{G,\widetilde{\delta}}^b(x)
    \df
    \{v\in V(G) \mid t_G^{\widetilde{\delta}}(v,x) = b\}.
  \end{gather*}
  This means that if we set
  \begin{gather*}
    \cH \df \{W\cap N_{G,\widetilde{\delta}}^{b_x}(x) \mid x\in V(G)^m\},
  \end{gather*}
  we have $\lvert H\rvert\leq\widetilde{\epsilon}\cdot s$ for every $H\in\cH$.

  It is clear that $\lvert\cH\rvert\leq n^m$, so applying item~\ref{prop:PRGsubset:conc}, it then suffices to show that if
  \begin{gather*}
    \lvert H\cap W(\rn{X})\rvert \leq \xi\cdot r,
  \end{gather*}
  for every $H\in\cH$, then $W(\rn{X})$ is a $(\widetilde{\delta},\widetilde{\xi},m)$-atom of $G$. But indeed, for $x\in
  V(G)^m$, we note that the set $H\df W\cap N_{G,\widetilde{\delta}}^{b_x}(x)$ is in $\cH$, so the above implies
  \begin{gather*}
    \lvert N_{G,\widetilde{\delta}}^{b_x}(x)\cap W(\rn{X})\rvert
    \leq
    \xi\cdot\lvert W(\rn{X})\rvert,
  \end{gather*}
  i.e., $x$ does not $(\widetilde{\delta},\widetilde{\xi})$-split $\rn{\widetilde{W}}$ in $G$, hence $W(\rn{X})$ is a
  $(\widetilde{\delta},\widetilde{\xi},m)$-atom of $G$.
\end{proof}

\begin{lemma}[Calculations for equitable partitions without randomness]\label{lem:equiPRGcalc}
  Let $d,\ell\in\NN_+$, let $c_{\UC}\in[0,1)$, let $c_{\atom},c_{\ZZ},c_{\adj},c_{\Root},c_{\Sum}\in(0,1)$ and let
  $c_{\samp}\in(0,1/2)$. Set $c_{\UC}\df 0$ and as in Lemma~\ref{lem:equicalc}, define the parameters:
  \begin{align*}
    \widetilde{c}_{\UC} & \df \frac{1}{(1-c_{\UC})^\ell} - 1 = 0,
    &
    \widetilde{c}_{\adj} & \df \frac{1}{(1-c_{\adj})^\ell} - 1,
    &
    \widetilde{c}_{\samp} & \df \frac{1}{(1-c_{\samp})^\ell} - 1,
    \\
    \widetilde{c}_{\ZZ} & \df \frac{1}{(1-c_{\ZZ})^{\ell+1}} - 1,
    &
    \widetilde{c}_{\Root} & \df \frac{1}{(1-c_{\Root})^{\ell+1}} - 1,
  \end{align*}
  and for $\epsilon\in(0,1)$, let
  \begin{align*}
    \zeta
    & \df
    (1-c_{\ZZ})\cdot\epsilon.
    \\
    K
    \df
    K_{\equiPRG}(\ell,c_{\ZZ},c_{\samp},c_{\adj},c_{\Root},c_{\Sum},\epsilon)
    & \df
    K_{\equi}(\ell,0,c_{\ZZ},c_{\samp},c_{\adj},c_{\Root},c_{\Sum},\epsilon),
    \\
    \epsilon_{\equiPRG}(\ell,c_{\ZZ},c_{\samp},c_{\adj},c_{\Root},c_{\Sum})
    & \df
    \min\left\{\frac{1}{2}, \epsilon_{\equi}(\ell,0,c_{\ZZ},c_{\samp},c_{\adj},c_{\Root},c_{\Sum})\right\}.
  \end{align*}
  For $n\in\NN_+$ with $n\geq n_{\equi}(\ell,0,c_{\ZZ},c_{\samp},c_{\adj},c_{\Root},c_{\Sum},\epsilon)$ as in
  Lemma~\ref{lem:equicalc}, define the quantities $r$, $\xi_i$, $\widetilde{\xi}_i$, $\widetilde{\epsilon}_i$, $\gamma_i$, $r_i$
  and $n_i$ as in Lemma~\ref{lem:equicalc}.

  Given $\delta\in(0,1)$, let further
  \begin{gather*}
    m \df \Ceil{C\cdot\frac{d}{c_{\atom}^2\cdot\delta^2}},
  \end{gather*}
  where $C$ is the absolute constant of Lemma~\ref{lem:atomconst}.

  Suppose now that $\epsilon < \epsilon_{\equiPRG}(\ell,c_{\ZZ},c_{\samp},c_{\adj},c_{\Root},c_{\Sum})$, let
  \begin{gather*}
    h_{\equiPRG}^{\good}
    \df
    4,
    \\
    \begin{multlined}[t]
      n_{\equiPRG}^{\good}(\ell,c_{\ZZ},c_{\samp},c_{\adj},c_{\Root},c_{\Sum},\epsilon)
      \\
      \df
      \min\bigggl\{n\in\NN_+ \;\bigggm\vert\;
      n\geq n_{\equi}(\ell,0,c_{\ZZ},c_{\samp},c_{\adj},c_{\Root},c_{\Sum},\epsilon)
      \\
      \land
      \Floor{(1-\zeta)^{1+1/\ell}\cdot\left(\frac{r}{n}\right)^{1/\ell}\cdot r}
      \geq
      \frac{2\cdot(h_{\equiPRG}^{\good}-1)^2}{(1-\zeta)\cdot(r/n)}
      \\
      \land
      \frac{2\cdot(1+c_{\adj}\cdot c_{\samp})}{c_{\adj}\cdot(1 - c_{\samp}\cdot(1-c_{\adj}))}
      \leq
      (1-\zeta)^{1+1/\ell}\cdot\left(\frac{r}{n}\right)^{1/\ell}\cdot r
      \\
      \land
      1 \geq
      \frac{1}{\sqrt{2\cdot\pi\cdot h_{\equiPRG}^{\good}}}
      \cdot
      \left(\frac{2\cdot e}{c_{\samp}^2\cdot(1-\zeta)^{1+2/\ell}\cdot (r/n)^{1+2/\ell}}\right)^{h_{\equiPRG}^{\good}/2}
      \\
      \cdot
      \Biggl(
      \left(\frac{1}{c_{\adj}^2\cdot(1-\zeta)\cdot r}\right)^{h_{\equiPRG}^{\good}/2}
      +
      n\cdot\left(\frac{
        (1-c_{\samp})^2
      }{
        (1-c_{\adj})^2\cdot\floor{(1-\zeta)^{1+1/\ell}\cdot(r/n)^{1/\ell}\cdot r}
      }\right)^{h_{\equiPRG}^{\good}/2}
      \Biggr)
      \bigggr\}.
    \end{multlined}
  \end{gather*}
  and let
  \begin{gather*}
    h_{\equiPRG}^{\atom}
    \df
    2\cdot m + 2,
    \\
    \begin{multlined}[t]
      n_{\equiPRG}^{\atom}(\ell,d,c_{\atom},c_{\ZZ},c_{\samp},c_{\adj},c_{\Root},c_{\Sum},\epsilon,\delta)
      \\
      \df
      \min\bigggl\{n\in\NN_+ \;\bigggm\vert\;
      n\geq n_{\equi}(\ell,0,c_{\ZZ},c_{\samp},c_{\adj},c_{\Root},c_{\Sum},\epsilon)
      \\
      \land
      \Floor{(1-\zeta)^{1+1/\ell}\cdot\left(\frac{r}{n}\right)^{1/\ell}\cdot r}
      \geq
      \frac{2\cdot(h_{\equiPRG}^{\atom}-1)^2}{(1-\zeta)\cdot(r/n)}
      \\
      \land
      \frac{2\cdot(1+c_{\adj}\cdot c_{\samp})}{c_{\adj}\cdot(1 - c_{\samp}\cdot(1-c_{\adj}))}
      \leq
      (1-\zeta)^{1+1/\ell}\cdot\left(\frac{r}{n}\right)^{1/\ell}\cdot r.
      \\
      \land
      1 \geq
      \frac{1}{\sqrt{2\cdot\pi\cdot h_{\equiPRG}^{\atom}}}
      \cdot
      \left(\frac{2\cdot e}{c_{\samp}^2\cdot(1-\zeta)^{1+2/\ell}\cdot (r/n)^{1+2/\ell}}\right)^{h_{\equiPRG}^{\atom}/2}
      \\
      \cdot
      \Biggl(
      \left(\frac{1}{c_{\adj}^2\cdot (1-\zeta)\cdot r}\right)^{h_{\equiPRG}^{\atom}/2}
      +
      n^m\cdot\left(\frac{
        (1-c_{\samp})^2
      }{
        (1-c_{\adj})^2\cdot\floor{(1-\zeta)^{1+1/\ell}\cdot(r/n)^{1/\ell}\cdot r}
      }\right)^{h_{\equiPRG}^{\atom}/2}
      \Biggr)
      \bigggr\}.
    \end{multlined}
  \end{gather*}

  For $\epsilon\in(0,\epsilon_{\equiPRG}(\ell,c_{\ZZ},c_{\samp},c_{\adj},c_{\Root},c_{\Sum}))$ and an integer $n\in\NN_+$ with
  \begin{gather*}
    n\geq n_{\equi}(\ell,0,c_{\ZZ},c_{\samp},c_{\adj},c_{\Root},c_{\Sum},\epsilon),
  \end{gather*}
  fix $i\in[K]$ and let $s_0,\ldots,s_\ell$ be defined inductively as
  \begin{align*}
    s_0 & \df n_{i-1}, &
    s_j & \df \floor{\gamma_i\cdot s_{j-1}} + 1.
  \end{align*}
  Fix then $j\in\{0,\ldots,\ell\}$ and set $p\df\gamma_i^{\ell-j}$ and $s\df s_j$. Then the following hold:
  \begin{enumerate}
  \item\label{lem:equiPRGcalc:gammai} If $\epsilon\leq 1/2$, then
    \begin{gather*}
      \gamma_i^\ell\cdot n_{i-1}
      \leq
      p\cdot s
      <
      \gamma_i^\ell\cdot n_{i-1} + \frac{1}{1-\gamma_i}
      \leq
      \gamma_i^\ell\cdot n_{i-1} + 2
      \leq
      r_i + 2.
    \end{gather*}
  \item\label{lem:equiPRGcalc:ps} We have
    \begin{align*}
      p & \geq \gamma_i^\ell \geq \frac{(1-\zeta)\cdot r}{n},
      &
      s & \geq r_i \geq (1-\zeta)\cdot r.
    \end{align*}
  \item\label{lem:equiPRGcalc:h} If $n\geq n_{\equiPRG}^{\good}(\ell,c_{\ZZ},c_{\samp},c_{\adj},c_{\Root},c_{\Sum},\epsilon)$
    and $h\df h_{\equiPRG}^{\good}$ and $j < \ell$, then
    \begin{gather}\label{eq:equiPRGcalc:h}
      \floor{\widetilde{\epsilon}_i\cdot s}
      \geq
      \begin{dcases*}
        \frac{(h-1)^2}{p\cdot(1-p)}, & if $p < 1$,\\
        0, & if $p=1$.
      \end{dcases*}
    \end{gather}

    Similarly, if $n_{\equiPRG}^{\atom}(\ell,d,c_{\atom},c_{\ZZ},c_{\samp},c_{\adj},c_{\Root},c_{\Sum},\epsilon,\delta)$ and
    $h\df h_{\equiPRG}^{\atom}$ and $j < \ell$, then~\eqref{eq:equiPRGcalc:h} holds.
  \item\label{lem:equiPRGcalc:adj} If $n\geq n_{\equiPRG}^{\good}(\ell,c_{\ZZ},c_{\samp},c_{\adj},c_{\Root},c_{\Sum},\epsilon)$, then
    \begin{gather}\label{eq:equiPRGcalc:adj}
      \widetilde{\xi}_i + \frac{c_{\adj}\cdot c_{\samp}\cdot\widetilde{\xi}_i\cdot p\cdot s + 2}{r_i}
      \leq
      \xi_i.
    \end{gather}

    Similarly, if $n_{\equiPRG}^{\atom}(\ell,d,c_{\atom},c_{\ZZ},c_{\samp},c_{\adj},c_{\Root},c_{\Sum},\epsilon,\delta)$,
    then~\eqref{eq:equiPRGcalc:adj} holds.
  \item\label{lem:equiPRGcalc:good} If $n\geq n_{\equiPRG}^{\good}(\ell,c_{\ZZ},c_{\samp},c_{\adj},c_{\Root},c_{\Sum},\epsilon)$
    and $h\df h_{\equiPRG}^{\good}$, then
    \begin{gather*}
      \frac{1}{\sqrt{2\cdot\pi\cdot h}}
      \cdot
      \left(\frac{2\cdot e\cdot (1-p)}{c_{\samp}^2\cdot\widetilde{\xi}_i^2\cdot p}\right)^{h/2}
      \cdot
      \left(
      \left(\frac{1}{c_{\adj}^2\cdot s}\right)^{h/2}
      +
      n\cdot\left(\frac{(1-c_{\samp})^2}{(1-c_{\adj})^2\cdot\floor{\widetilde{\epsilon}_i\cdot s}}\right)^{h/2}
      \right)
      \leq
      1.
    \end{gather*}
  \item\label{lem:equiPRGcalc:atom} If $n\geq
    n_{\equiPRG}^{\atom}(\ell,d,c_{\atom},c_{\ZZ},c_{\samp},c_{\adj},c_{\Root},c_{\Sum},\epsilon,\delta)$ and $h\df
    h_{\equiPRG}^{\atom}$, then
    \begin{gather*}
      \frac{1}{\sqrt{2\cdot\pi\cdot h}}
      \cdot
      \left(\frac{2\cdot e\cdot (1-p)}{c_{\samp}^2\cdot\widetilde{\xi}_i^2\cdot p}\right)^{h/2}
      \cdot
      \left(
      \left(\frac{1}{c_{\adj}^2\cdot s}\right)^{h/2}
      +
      n^m\cdot\left(\frac{(1-c_{\samp})^2}{(1-c_{\adj})^2\cdot\floor{\widetilde{\epsilon}_i\cdot s}}\right)^{h/2}
      \right)
      \leq
      1.
    \end{gather*}
  \item\label{lem:equiPRGcalc:ngood} We have
    \begin{align*}
      \MoveEqLeft
      n_{\equiPRG}^{\good}(\ell,c_{\ZZ},c_{\samp},c_{\adj},c_{\Root},c_{\Sum},\epsilon)
      \\
      & \leq
      \bigl(1 + o_{\epsilon\to 0,c_{\ZZ},c_{\samp},c_{\adj},c_{\Root},c_{\Sum}}(1)\bigr)\cdot
      \sqrt{\frac{2}{\pi}}\cdot\frac{e^2\cdot (1-c_{\samp})^4\cdot K^{4+6/\ell}}{c_{\samp}^4\cdot(1-c_{\adj})^4}.
    \end{align*}
  \item\label{lem:equiPRGcalc:natom} We have
    \begin{align*}
      \MoveEqLeft
      n_{\equiPRG}^{\atom}(\ell,d,c_{\atom},c_{\ZZ},c_{\samp},c_{\adj},c_{\Root},c_{\Sum},\epsilon,\delta)
      \\
      & \leq
      \bigl(1 + o_{\epsilon\to 0,c_{\ZZ},c_{\samp},c_{\adj},c_{\Root},c_{\Sum}}(1)\bigr)\cdot
      \frac{
        2^m\cdot e^{m+1}\cdot (1-c_{\samp})^{2\cdot(m+1)}\cdot K^{(m+1)\cdot(2+3/\ell)}
      }{
        \sqrt{\pi\cdot(m+1)}\cdot c_{\samp}^{2\cdot(m+1)}\cdot(1-c_{\adj})^{2\cdot(m+1)}
      }
      \\
      & \leq
      \bigl(1 + o_{\epsilon\to 0,c_{\ZZ},c_{\samp},c_{\adj},c_{\Root},c_{\Sum}}(1)\bigr)\cdot
      K^{O(d/(c_{\atom}^2\cdot\delta^2))}.
    \end{align*}
  \end{enumerate}
\end{lemma}

\begin{proof}
  For item~\ref{lem:equiPRGcalc:gammai}, a simple induction in $j\in\{0,\ldots,\ell\}$ yields
  \begin{gather}\label{eq:equiPRGcalc:gammai}
    \gamma_i^j\cdot n_{i-1}
    +
    \One[j > 0]
    \leq
    s_j
    \leq
    \gamma_i^j\cdot n_{i-1} + \sum_{t=0}^{j-1} \gamma_i^t
    <
    \gamma_i^j\cdot n_{i-1} + \frac{1}{1-\gamma_i}
    \leq
    \gamma_i^j\cdot n_{i-1} + 2,
  \end{gather}
  where the last inequality follows since $\gamma_i\leq\epsilon\leq 1/2$. Since $s\df s_j$ and $p\df\gamma_i^{\ell-j}$ for one
  of the $j\in\{0,\ldots,\ell\}$, multiplying the above for this $j$ by $p$, we get
  \begin{gather*}
    \gamma_i^\ell\cdot n_{i-1}
    \leq
    p\cdot s
    <
    \gamma_i^\ell\cdot n_{i-1} + 2\cdot p
    \leq
    \gamma_i^\ell\cdot n_{i-1} + 2
    \leq
    r_i + 2,
  \end{gather*}
  where the last inequality follows since $r_i\df\ceil{\gamma_i^\ell\cdot n_{i-1}}$.

  \medskip

  For item~\ref{lem:equiPRGcalc:ps}, we have
  \begin{gather*}
    p = \gamma_i^j \geq \gamma_i^\ell \geq \frac{(1-\zeta)\cdot r}{n},
  \end{gather*}
  where the last inequality follows from the bound
  \begin{gather}\label{eq:equiPRGcalc:xiilowerbound}
    \xi_i \geq \widetilde{\xi}_i \geq \widetilde{\epsilon}_i
    \geq
    \gamma_i \geq \left(\frac{(1 - \zeta)\cdot r}{n}\right)^{1/\ell}
  \end{gather}

  Lemma~\ref{lem:equicalc}\ref{lem:equicalc:xiilowerbound} and
  since~\eqref{eq:equiPRGcalc:gammai} implies
  \begin{gather*}
    s
    \geq
    \gamma_i^j\cdot n_{i-1} + \One[j>0]
    >
    \gamma_i^\ell\cdot n_{i-1},
  \end{gather*}
  and $s$ is an integer, it follows that
  \begin{gather}\label{eq:equiPRGcalc:slowerbound}
    s
    \geq
    \ceil{\gamma_i^\ell\cdot n_{i-1}}
    =
    r_i
    \geq
    (1-\zeta)\cdot r,
  \end{gather}
  where the last inequality follows from Lemma~\ref{lem:equicalc}\ref{lem:equicalc:rilowerbound}.

  \medskip

  Item~\ref{lem:equiPRGcalc:h} follows by using the bounds of item~\ref{lem:equiPRGcalc:ps}, the lower bound on
  $\widetilde{\epsilon}_i$ from Lemma~\ref{lem:equicalc}\ref{lem:equicalc:xiilowerbound}
  (see~\eqref{eq:equiPRGcalc:xiilowerbound} above) and the bound $p\leq\gamma_i\leq 1/2$ (as $j < \ell$) along with the
  condition
  \begin{gather*}
    \Floor{(1-\zeta)^{1+1/\ell}\cdot\left(\frac{r}{n}\right)^{1/\ell}\cdot r}
    \geq
    \frac{2\cdot(h_{\equiPRG}^{\good}-1)^2}{(1-\zeta)\cdot(r/n)}
  \end{gather*}
  in the definition of $n_{\equiPRG}^{\good}(\ell,c_{\ZZ},c_{\samp},c_{\adj},c_{\Root},c_{\Sum},\epsilon)$ (and the
  analogous condition in the definition of
  $n_{\equiPRG}^{\atom}(\ell,d,c_{\atom},c_{\ZZ},c_{\samp},c_{\adj},c_{\Root},c_{\Sum},\epsilon,\delta)$).

  \medskip

  For item~\ref{lem:equiPRGcalc:adj}, since $\widetilde{\xi}_i\df(1-c_{\adj})\cdot\xi_i$, it follows
  that~\eqref{eq:equiPRGcalc:adj} is equivalent to
  \begin{gather*}
    \frac{c_{\adj}\cdot c_{\samp}\cdot(1-c_{\adj})\cdot\xi_i\cdot p\cdot s + 2}{r_i}
    \leq
    c_{\adj}\cdot\xi_i.
  \end{gather*}
  Using the upper bound $p\cdot s < r_i + 2$ of item~\ref{lem:equiPRGcalc:gammai}, it suffices to show
  \begin{gather*}
    c_{\adj}\cdot c_{\samp}\cdot(1-c_{\adj})\cdot\xi_i + \frac{2\cdot(1+c_{\adj}\cdot c_{\samp})}{r_i}
    \leq
    c_{\adj}\cdot\xi_i,
  \end{gather*}
  which in turn is equivalent to
  \begin{gather*}
    \frac{2\cdot(1+c_{\adj}\cdot c_{\samp})}{c_{\adj}\cdot(1 - c_{\samp}\cdot(1-c_{\adj}))}
    \leq
    \xi_i\cdot r_i.
  \end{gather*}
  We now use the lower bounds $r_i\geq(1-\zeta)\cdot r$ and $\xi_i\geq((1-\zeta)\cdot r/n)^{1/\ell}$ from
  Lemma~\ref{lem:equicalc}, items~\ref{lem:equicalc:xiilowerbound} and~\ref{lem:equicalc:rilowerbound} to conclude that it
  suffices to show
  \begin{gather*}
    \frac{2\cdot(1+c_{\adj}\cdot c_{\samp})}{c_{\adj}\cdot(1 - c_{\samp}\cdot(1-c_{\adj}))}
    \leq
    (1-\zeta)^{1+1/\ell}\cdot\left(\frac{r}{n}\right)^{1/\ell}\cdot r.
  \end{gather*}
  But this is exactly one of the conditions in the definitions of
  $n_{\equiPRG}^{\good}(\ell,c_{\ZZ},c_{\samp},c_{\adj},c_{\Root},c_{\Sum},\epsilon)$ and
  $n_{\equiPRG}^{\atom}(\ell,d,c_{\atom},c_{\ZZ},c_{\samp},c_{\adj},c_{\Root},c_{\Sum},\epsilon,\delta)$.

  \medskip

  We now prove item~\ref{lem:equiPRGcalc:good}. Using the bounds
  \begin{align}\label{eq:equiPRGcalc:lowerbounds}
    \widetilde{\xi}_i & \geq \widetilde{\epsilon}_i \geq \left(\frac{(1-\zeta)\cdot r}{n}\right)^{1/\ell},
    &
    p & \geq \widetilde{\epsilon}_i^\ell,
    &
    s & \geq (1-\zeta)\cdot r,
  \end{align}
  from Lemma~\ref{lem:equicalc}, items~\ref{lem:equicalc:xiilowerbound} and~\ref{lem:equicalc:rilowerbound}
  and~\eqref{eq:equiPRGcalc:slowerbound}, it suffices to show
  \begin{multline*}
    \frac{1}{\sqrt{2\cdot\pi\cdot h}}
    \cdot
    \left(\frac{2\cdot e\cdot (1-p)}{c_{\samp}^2\cdot(1-\zeta)^{1+2/\ell}\cdot(r/n)^{1+2/\ell}}\right)^{h/2}
    \\
    \cdot
    \left(
    \left(\frac{1}{c_{\adj}^2\cdot(1-\zeta)\cdot r}\right)^{h/2}
    +
    n\cdot\left(\frac{(1-c_{\samp})^2}{(1-c_{\adj})^2\cdot\floor{(1-\zeta)^{1+1/\ell}\cdot(r/n)^{1/\ell}\cdot r}}\right)^{h/2}
    \right)
    \leq
    1,
  \end{multline*}
  which is exactly one of the conditions in the definition of
  $n_{\equiPRG}^{\good}(\ell,c_{\ZZ},c_{\samp},c_{\adj},c_{\Root},c_{\Sum},\epsilon)$.

  \medskip

  For item~\ref{lem:equiPRGcalc:atom}, using the bounds~\eqref{eq:equiPRGcalc:lowerbounds}, it suffices to show
  \begin{multline*}
    \frac{1}{\sqrt{2\cdot\pi\cdot h_{\equiPRG}^{\atom}}}
    \cdot
    \left(\frac{2\cdot e}{c_{\samp}^2\cdot(1-\zeta)^{1+2/\ell}\cdot (r/n)^{1+2/\ell}}\right)^{h_{\equiPRG}^{\atom}/2}
    \\
    \cdot
    \left(
    \left(\frac{1}{c_{\adj}^2\cdot (1-\zeta)\cdot r}\right)^{h_{\equiPRG}^{\atom}/2}
    +
    n^m\cdot\left(\frac{
      (1-c_{\samp})^2
    }{
      (1-c_{\adj})^2\cdot\floor{(1-\zeta)^{1+1/\ell}\cdot(r/n)^{1/\ell}\cdot r}
    }\right)^{h_{\equiPRG}^{\atom}/2}
    \right)
    \leq
    1,
  \end{multline*}
  which is one of the conditions in the definition of
  $n_{\equiPRG}^{\atom}(\ell,d,c_{\atom},c_{\ZZ},c_{\samp},c_{\adj},c_{\Root},c_{\Sum},\epsilon,\delta)$.

  \medskip

  We now prove item~\ref{lem:equiPRGcalc:ngood}. For this, we note that since
  \begin{gather*}
    r
    =
    \Floor{\frac{n}{K}}
    =
    \bigl(1 + o_{c_{\ZZ},c_{\samp},c_{\adj},c_{\Root},c_{\Sum}}(1)\bigr)\cdot\frac{n}{K}
  \end{gather*}
  and using Lemma~\ref{lem:equicalc}\ref{lem:equicalc:nequi}, the conditions in the definition of
  $n_{\equiPRG}^{\good}(\ell,c_{\ZZ},c_{\samp},c_{\adj},c_{\Root},c_{\Sum},\epsilon)$ are asymptotically equivalent to
  \begin{gather*}
    n
    \geq
    \bigl(1+o_{\epsilon\to 0,c_{\ZZ},c_{\samp},c_{\adj},c_{\Root},c_{\Sum}}(1)\bigr)
    \cdot
    \max\left\{\frac{1-c_{\ZZ}}{c_{\ZZ}}, \frac{1+c_{\Sum}}{(1-c_{\Root})\cdot c_{\Sum}}\right\}
    \cdot\frac{K}{(1-c_{\ZZ})\cdot\epsilon},
    \\
    \bigl(1 + o_{\epsilon\to 0,c_{\ZZ},c_{\samp},c_{\adj},c_{\Root},c_{\Sum}}(1)\bigr)\cdot
    \frac{n}{K^{1 + 1/\ell}}
    \geq
    18\cdot K,
    \\
    \bigl(1 + o_{n\to\infty,\epsilon\to 0,c_{\ZZ},c_{\samp},c_{\adj},c_{\Root},c_{\Sum}}(1)\bigr)\cdot
    \frac{2\cdot(1+c_{\adj}\cdot c_{\samp})}{c_{\adj}\cdot(1 - c_{\samp}\cdot(1-c_{\adj}))}
    \leq
    \frac{n}{K^{1+1/\ell}},
    \\
    \begin{multlined}[t]
      1
      \geq
      \bigl(1 + o_{n\to\infty,\epsilon\to 0,c_{\ZZ},c_{\samp},c_{\adj},c_{\Root},c_{\Sum}}(1)\bigr)\cdot
      \frac{1}{\sqrt{8\cdot\pi}}
      \cdot
      \left(\frac{2\cdot e\cdot K^{1+2/\ell}}{c_{\samp}^2}\right)^2
      \cdot
      \Biggl(
      \left(\frac{K}{c_{\adj}^2\cdot n}\right)^2
      \\
      +
      n\cdot\left(\frac{
        (1-c_{\samp})^2\cdot K^{1+1/\ell}
      }{
        (1-c_{\adj})^2\cdot n
      }\right)^2
      \Biggr).
    \end{multlined}
  \end{gather*}

  It is straightforward to check that the asymptotically harshest requirement is the last one and deduce that
  \begin{align*}
    \MoveEqLeft
    n_{\equiPRG}^{\good}(\ell,c_{\ZZ},c_{\samp},c_{\adj},c_{\Root},c_{\Sum},\epsilon)
    \\
    & \leq
    \bigl(1 + o_{\epsilon\to 0,c_{\ZZ},c_{\samp},c_{\adj},c_{\Root},c_{\Sum}}(1)\bigr)\cdot
    \sqrt{\frac{2}{\pi}}\cdot\frac{e^2\cdot (1-c_{\samp})^4\cdot K^{4+6/\ell}}{c_{\samp}^4\cdot(1-c_{\adj})^4}.
  \end{align*}

  \medskip

  Finally, we prove item~\ref{lem:equiPRGcalc:natom}. Using again the expression for the asymptotic behavior of $r$, the
  conditions in the definition of
  $n_{\equiPRG}^{\atom}(\ell,d,c_{\atom},c_{\ZZ},c_{\samp},c_{\adj},c_{\Root},c_{\Sum},\epsilon,\delta)$ are asymptotically
  equivalent to
  \begin{gather*}
    n
    \geq
    \bigl(1+o_{\epsilon\to 0,c_{\ZZ},c_{\samp},c_{\adj},c_{\Root},c_{\Sum}}(1)\bigr)
    \cdot
    \max\left\{\frac{1-c_{\ZZ}}{c_{\ZZ}}, \frac{1+c_{\Sum}}{(1-c_{\Root})\cdot c_{\Sum}}\right\}
    \cdot\frac{K}{(1-c_{\ZZ})\cdot\epsilon},
    \\
    \bigl(1 + o_{\epsilon\to 0,c_{\ZZ},c_{\samp},c_{\adj},c_{\Root},c_{\Sum}}(1)\bigr)\cdot
    \frac{n}{K^{1 + 1/\ell}}
    \geq
    2\cdot(2\cdot m+1)^2\cdot K,
    \\
    \bigl(1 + o_{n\to\infty,\epsilon\to 0,c_{\ZZ},c_{\samp},c_{\adj},c_{\Root},c_{\Sum}}(1)\bigr)\cdot
    \frac{2\cdot(1+c_{\adj}\cdot c_{\samp})}{c_{\adj}\cdot(1 - c_{\samp}\cdot(1-c_{\adj}))}
    \leq
    \frac{n}{K^{1+1/\ell}},
    \\
    \begin{multlined}[t]
      1
      \geq
      \bigl(1 + o_{n\to\infty,\epsilon\to 0,c_{\ZZ},c_{\samp},c_{\adj},c_{\Root},c_{\Sum}}(1)\bigr)\cdot
      \frac{1}{\sqrt{2\cdot\pi\cdot(2\cdot m+2)}}
      \cdot
      \left(\frac{2\cdot e\cdot K^{1+2/\ell}}{c_{\samp}^2}\right)^{m+1}
      \\
      \cdot
      \left(
      \left(\frac{K}{c_{\adj}^2\cdot n}\right)^{m+1}
      +
      n^m\cdot\left(\frac{
        (1-c_{\samp})^2\cdot K^{1+1/\ell}
      }{
        (1-c_{\adj})^2\cdot n
      }\right)^{m+1}
      \right).
    \end{multlined}
  \end{gather*}

  It is straightforward to check that the asymptotically harshest requirement is the last one and deduce that
  \begin{align*}
    \MoveEqLeft
    n_{\equiPRG}^{\atom}(\ell,d,c_{\atom},c_{\ZZ},c_{\samp},c_{\adj},c_{\Root},c_{\Sum},\epsilon,\delta)
    \\
    & \leq
    \bigl(1 + o_{\epsilon\to 0,c_{\ZZ},c_{\samp},c_{\adj},c_{\Root},c_{\Sum}}(1)\bigr)\cdot
    \frac{
      2^m\cdot e^{m+1}\cdot (1-c_{\samp})^{2\cdot(m+1)}\cdot K^{(m+1)\cdot(2+3/\ell)}
    }{
      \sqrt{\pi\cdot(m+1)}\cdot c_{\samp}^{2\cdot(m+1)}\cdot(1-c_{\adj})^{2\cdot(m+1)}
    }
    \\
    & \leq
    \bigl(1 + o_{\epsilon\to 0,c_{\ZZ},c_{\samp},c_{\adj},c_{\Root},c_{\Sum}}(1)\bigr)\cdot
    K^{O(d/(c_{\atom}^2\cdot\delta^2))}.
    \qedhere
  \end{align*}
\end{proof}

\begin{theorem}[Equitable partitions without randomness]\label{thm:equiPRG}
  Let $\ell,d\in\NN_+$ with $d\leq\ell$, let $\delta,\epsilon\in(0,1)$, let
  $c_{\atom},c_{\ZZ},c_{\adj},c_{\Root},c_{\Sum}\in(0,1)\cap\QQ$, let $c_{\samp}\in(0,1/2)\in\QQ$ and let $G$ be a graph with
  $\lvert G\rvert=n$, $\Lit(G)\leq\ell$ and $\VC(G)\leq d$. Let also
  \begin{align*}
    \MoveEqLeft
    \epsilon_{\equiPRG}(\ell,c_{\ZZ},c_{\samp},c_{\adj},c_{\Root},c_{\Sum})
    \\
    & \df
    \min\left\{\frac{1}{2}, \epsilon_{\equi}(\ell,0,c_{\ZZ},c_{\samp},c_{\adj},c_{\Root},c_{\Sum})\right\}
    \\
    & \geq
    \begin{dcases*}
      \bigl(1+o_{c_{\Sum}\to 0,c_{\UC}\to 0,c_{\samp}\to 0,c_{\adj}\to 0,c_{\ZZ},c_{\Root}}(1)\bigr)\cdot
      \frac{2\cdot(1-c_{\ZZ})\cdot\ln(1/c_{\Sum})}{c_{\Sum}},
      & if $\ell=1$,
      \\
      \bigl(1 + o_{c_{\UC}\to 0, c_{\samp}\to 0, c_{\adj}\to 0}(1)\bigr)\cdot
      \frac{
        (1+c_{\Sum})^{\ell-1}\cdot(1-c_{\ZZ})
      }{
        (1-c_{\Root})^\ell\cdot(1-1/\ell)^\ell\cdot c_{\Sum}^\ell
      },
      & if $\ell\geq 2$.
    \end{dcases*}
  \end{align*}
  be as in Lemmas~\ref{lem:equicalc} and~\ref{lem:equiPRGcalc} and for
  $\epsilon\in(0,\epsilon_{\equiPRG}(\ell,c_{\ZZ},c_{\samp},c_{\adj},c_{\Root},c_{\Sum}))$, let
  \begin{align*}
    K
    & \df
    K_{\equiPRG}(\ell,c_{\ZZ},c_{\samp},c_{\adj},c_{\Root},c_{\Sum},\epsilon)
    \df
    K_{\equi}(\ell,0,c_{\ZZ},c_{\samp},c_{\adj},c_{\Root},c_{\Sum},\epsilon)
    \\
    & \df
    \ceil{
      (1+\widetilde{c}_{\ZZ})\cdot(1+\widetilde{c}_{\samp})\cdot(1+\widetilde{c}_{\adj})
      \cdot(1+\widetilde{c}_{\Root})\cdot(1+c_{\Sum})\cdot\epsilon^{-\ell-1}
    },
    \\
    \MoveEqLeft
    n_{\equiPRG}^{\good}(\ell,c_{\ZZ},c_{\samp},c_{\adj},c_{\Root},c_{\Sum},\epsilon)
    \\
    & \leq
    \bigl(1 + o_{\epsilon\to 0,c_{\ZZ},c_{\samp},c_{\adj},c_{\Root},c_{\Sum}}(1)\bigr)\cdot
    \sqrt{\frac{2}{\pi}}\cdot\frac{e^2\cdot (1-c_{\samp})^4\cdot K^{4+6/\ell}}{c_{\samp}^4\cdot(1-c_{\adj})^4}.
    \\
    \MoveEqLeft
    n_{\equiPRG}^{\atom}(\ell,d,c_{\atom},c_{\ZZ},c_{\samp},c_{\adj},c_{\Root},c_{\Sum},\epsilon,\delta)
    \\
    & \leq
    \bigl(1 + o_{\epsilon\to 0,c_{\ZZ},c_{\samp},c_{\adj},c_{\Root},c_{\Sum}}(1)\bigr)\cdot
    K^{O(d/(c_{\atom}^2\cdot\delta^2))}.
  \end{align*}
  also be as in Lemmas~\ref{lem:equicalc} and~\ref{lem:equiPRGcalc}. Let $n\in\NN_+$. For the items below about good
  equipartitions, we assume $n\geq n_{\equiPRG}^{\good}(\ell,c_{\ZZ},c_{\samp},c_{\adj},c_{\Root},c_{\Sum},\epsilon)$ and for
  the items below about atomic equipartitions, we assume $n\geq
  n_{\equiPRG}^{\atom}(\ell,d,c_{\atom},c_{\ZZ},c_{\samp},c_{\adj},c_{\Root},c_{\Sum},\epsilon,\delta)$. Then the following
  hold:
  \begin{enumerate}
  \item\label{thm:equiPRG:good} In the deterministic model, Algorithm~\ref{alg:equiPRG:good} computes an $\epsilon$-good
    equipartition of $G$ into $K$ parts in time
    \begin{multline*}
      2^{O(\ell\cdot(\len(c_{\Root}) + \len(c_{\ZZ}) + \ell\cdot\len(\epsilon) + \len(c_{\samp}) + \len(c_{\adj})))}\cdot n^6
      \\
      +
      O\Bigl(
      \ell^4\cdot\log(\ell+1)\cdot\len(c_{\ZZ})^2\cdot\len(c_{\samp})^2\cdot\len(c_{\Root})^2\cdot\len(c_{\Sum})\cdot\len(\epsilon)^2
      \cdot\bigl(\log(n+1)\bigr)^2
      \Bigr).
    \end{multline*}
  \item\label{thm:equiPRG:atom} If $(1+c_{\atom})\cdot\delta < 1/2^{\ell+1}$, then in the deterministic model,
    Algorithm~\ref{alg:equiPRG:atom} computes a $(\delta,\epsilon)$-atomic equipartition of $G$ into $K$ parts in time
    \begin{multline*}
      2^{O(d\cdot\ell\cdot(\len(c_{\Root})+\len(c_{\ZZ})+\ell\cdot\len(\epsilon)+\len(c_{\samp}) + \len(c_{\adj})))/(c_{\atom}^2\cdot\delta^2)}
      \cdot n^{O(d/c_{\atom}^2\cdot\delta^2)}
      \\
      +
      O\bigl(
      \log(d+1)\cdot\len(c_{\atom})\cdot\len(\delta)
      \\
      +
      \ell^4\cdot\log(\ell+1)\cdot\len(c_{\ZZ})^2\cdot\len(c_{\samp})^2\cdot\len(c_{\Root})^2\cdot\len(c_{\Sum})\cdot\len(\epsilon)^2
      \cdot\bigl(\log(n+1)\bigr)^2
      \Bigr),
    \end{multline*}
  \item\label{thm:equiPRG:goodspace} In the query-into-oracle model, Algorithms~\ref{alg:equiPRG:goodspace}
    and~\ref{alg:equiPRG:atomspaceoracle} compute an $\epsilon$-good equipartition into $K$ parts. The space-complexities of
    Algorithms~\ref{alg:equiPRG:goodspace} and~\ref{alg:equiPRG:atomspaceoracle} are
    \begin{gather*}
      \begin{multlined}[t]
        O\bigggl(
        \epsilon^{-\ell-1}\cdot
        \left(\log\left(\frac{n}{\epsilon^\ell}\right)
        +
        \ell
        \cdot\bigl(\len(c_{\Root}) + \len(c_{\ZZ}) + \ell\cdot\len(\epsilon) + \len(c_{\samp}) + \len(c_{\adj})\bigr)
        \right)
        \\
        + \ell\cdot\len(c_{\adj})
        \bigggr),
      \end{multlined}
      \\
      \begin{multlined}[t]
        O\bigggl(
        \epsilon^{-\ell-1}\cdot
        \left(
        \log\left(\frac{n}{\epsilon^\ell}\right)\right)
        \\
        +
        \ell\cdot\bigl(\len(c_{\Root}) + \len(c_{\ZZ}) + \ell\cdot\len(\epsilon)
        + \len(c_{\samp}) + \len(c_{\adj})\bigr)
        \bigggr),
      \end{multlined}
    \end{gather*}
    respectively, and the time-complexities are
    \begin{gather*}
      \begin{multlined}[t]
        2^{O(\ell\cdot(\len(c_{\Root}) + \len(c_{\ZZ}) + \ell\cdot\len(\epsilon) + \len(c_{\samp}) + \len(c_{\adj})))}
        \cdot n^{O(\epsilon^{-\ell-1})}
        \\
        +
        \ell^4\cdot\log(\ell+1)\cdot\len(c_{\ZZ})^2\cdot\len(c_{\samp})^2\cdot\len(c_{\Root})^2\cdot\len(c_{\Sum})\cdot\len(\epsilon)^2
        \cdot\bigl(\log(n+1)\bigr)^2,
      \end{multlined}
      \\
      (2\cdot n)^{O(\epsilon^{-\ell-1})}
      \cdot\bigggl(
      \epsilon^{-2\cdot\ell-2}
      +
      \ell^2\cdot\bigl(\len(c_{\Root}) + \len(c_{\ZZ}) + \ell\cdot\len(\epsilon)
      + \len(c_{\samp}) + \len(c_{\adj})\bigr)^2
      \bigggr).
    \end{gather*}
  \item\label{thm:equiPRG:atomspace} If $(1+c_{\atom})\cdot\delta < 1/2^{\ell+1}$, then in the query-into-oracle model,
    Algorithms~\ref{alg:equiPRG:atomspace} and~\ref{alg:equiPRG:atomspaceoracle} compute a $(\delta,\epsilon)$-atomic
    equipartition into $K$ parts. The space-complexities of Algorithms~\ref{alg:equiPRG:atomspace}
    and~\ref{alg:equiPRG:atomspaceoracle} are
    \begin{gather*}
      \begin{multlined}[t]
        O\bigggl(
        \epsilon^{-\ell-1}\cdot
        \Biggl(
        \Biggl(\log\left(\frac{n}{\epsilon^\ell}\right)
        \\
        +
        \frac{\ell\cdot d}{c_{\atom}^2\cdot\delta^2}
        \cdot\bigl(\len(c_{\Root}) + \len(c_{\ZZ}) + \ell\cdot\len(\epsilon) + \len(c_{\samp}) + \len(c_{\adj})\bigr)
        \Biggr)
        \Biggr)
        \\
        +
        \len(c_{\atom}) + \len(\delta)
        + \ell\cdot\len(c_{\adj})
        \bigggr),
      \end{multlined}
      \\
      \begin{multlined}[t]
        O\bigggl(
        \epsilon^{-\ell-1}\cdot
        \left(
        \log\left(\frac{n}{c_{\atom}\cdot\delta\cdot\epsilon^\ell}\right)\right)
        \\
        +
        \frac{\ell\cdot d}{c_{\atom}^2\cdot\delta^2}\cdot\bigl(\len(c_{\Root}) + \len(c_{\ZZ}) + \ell\cdot\len(\epsilon)
        + \len(c_{\samp}) + \len(c_{\adj})\bigr)
        \bigggr),
      \end{multlined}
    \end{gather*}
    respectively, and the time-complexities are
    \begin{gather*}
      \begin{multlined}[t]
        2^{O(\ell\cdot d\cdot(\len(c_{\Root}) + \len(c_{\ZZ}) + \ell\cdot\len(\epsilon) + \len(c_{\samp}) + \len(c_{\adj})))/(c_{\atom}^2\cdot\delta^2)}
        \cdot n^{O(\epsilon^{-\ell-1}/(c_{\atom}^2\cdot\delta^2))}
        \\
        +
        \log(d+1)\cdot\len(c_{\atom})\cdot\len(\delta)
        \\
        +
        \ell^4\cdot\log(\ell+1)\cdot\len(c_{\ZZ})^2\cdot\len(c_{\samp})^2\cdot\len(c_{\Root})^2\cdot\len(c_{\Sum})\cdot\len(\epsilon)^2
        \cdot\bigl(\log(n+1)\bigr)^2,
      \end{multlined}
      \\
      \begin{multlined}[t]
        (2\cdot n)^{O(\epsilon^{-\ell-1})}
        \cdot\bigggl(
        \frac{d}{c_{\atom}^2\cdot\delta^2\cdot\epsilon^{2\cdot\ell+2}}
        \\
        +
        \frac{\ell^2\cdot d}{c_{\atom}^2\cdot\delta^2}\cdot\bigl(\len(c_{\Root}) + \len(c_{\ZZ}) + \ell\cdot\len(\epsilon)
        + \len(c_{\samp}) + \len(c_{\adj})\bigr)^2
        \bigggr).
      \end{multlined}
    \end{gather*}
  \end{enumerate}

  Furthermore, for
  \begin{align*}
    \epsilon'
    & \df
    \begin{dcases*}
      \frac{1 - \sqrt{1 - 8\cdot\epsilon\cdot(1-\epsilon)}}{2}, & if $\epsilon < (2-\sqrt{2})/4$,\\
      1/2, & otherwise,
    \end{dcases*}
    \\
    & =
    2\cdot\epsilon + O_{\epsilon\to 0}(\epsilon^2),
    \\
    \epsilon''
    & \df
    \frac{\epsilon'}{1-\epsilon}
    =
    \frac{1 - \sqrt{1 - 8\cdot\epsilon\cdot(1-\epsilon)}}{2\cdot(1-\epsilon)}
    =
    2\cdot\epsilon + O_{\epsilon\to 0}(\epsilon^2),
  \end{align*}
  all $(\delta,\epsilon)$-atomic partitions are $(\delta,\epsilon)$-excellent and $\epsilon$-good; and all $\epsilon$-good
  partitions of $G$ are in particular totally $\epsilon'$-homogeneous and $(\epsilon,\epsilon'')$-excellent.
\end{theorem}

\begin{proof}
  Just as in Theorems~\ref{thm:nonequi} and~\ref{thm:equirand}, the final assertions on how to translate atomic into excellent,
  into good, into totally homogeneous and back into excellent follow from Lemmas~\ref{lem:atom->exc}
  and~\ref{lem:good->hom}\ref{lem:good->hom:partition}.

  The idea of the proof is the same as that of Theorem~\ref{thm:equirand}, except that since we do not have access to
  randomness, we will use the pseudorandom generator from Theorem~\ref{thm:PRG}. We will ensure that the seed space size is
  polynomial in $n$, which means that we will be able to enumerate all the seeds and test whether they generate sets satisfying
  the conditions of Proposition~\ref{prop:PRGsubset}. One small technicality is that the pseudorandom generator of
  Theorem~\ref{thm:PRG} has to have a bias that is dyadic (i.e., an element of $\{a/2^n \mid a\in\ZZ, n\in\NN\}$), we need to
  ensure that the $\gamma_i$ are dyadics (see also Discussion~\ref{dsc:dyadic} below); on a small technical note: it will be
  convenient to assume $V(G)=[n]$.

  Again, we will prove the items slightly out of order for simplicity sake:
  \begin{description}[wide, itemsep={3ex}]
  \item[Item~\ref{thm:equiPRG:atom}.] We start with the proof of correctness of Algorithm~\ref{alg:equiPRG:atom}.

    \begin{algorithm}[htbp]
      \caption{Deterministic model algorithm that returns a $(\delta,\epsilon)$-atomic equipartition of $G$ into
        $K_{\equiPRG}(\ell,c_{\ZZ},c_{\samp},c_{\adj},c_{\Root},c_{\Sum},\epsilon)$ parts in time
        \begin{multline*}
          2^{O(d\cdot\ell\cdot(\len(c_{\Root})+\len(c_{\ZZ})+\ell\cdot\len(\epsilon)+\len(c_{\samp}) + \len(c_{\adj})))/(c_{\atom}^2\cdot\delta^2)}
          \cdot n^{O(d/c_{\atom}^2\cdot\delta^2)}
          \\
          +
          O\bigl(
          \log(d+1)\cdot\len(c_{\atom})\cdot\len(\delta)
          \\
          +
          \ell^4\cdot\log(\ell+1)\cdot\len(c_{\ZZ})^2\cdot\len(c_{\samp})^2\cdot\len(c_{\Root})^2\cdot\len(c_{\Sum})\cdot\len(\epsilon)^2
          \cdot\bigl(\log(n+1)\bigr)^2
          \Bigr),
        \end{multline*}
      }
      \label{alg:equiPRG:atom}
      \DontPrintSemicolon
      \KwIn{Numbers $n,d,\ell\in\NN_+$ with $d\leq\ell\leq\log_2(n)$,
        $c_{\atom},c_{\ZZ},c_{\adj},c_{\Root},c_{\Sum},\delta\in(0,1)\cap\QQ$ with $(1+c_{\atom})\cdot\delta < 1/2^{\ell+1}$,
        $c_{\samp}\in(0,1/2)\cap\QQ$, and
        $\epsilon\in(0,\epsilon_{\equiPRG}(\ell,c_{\ZZ},c_{\samp},c_{\adj},c_{\Root},c_{\Sum}))$, and a graph $G$ with
        $V(G)=[n]$, $\Lit(G)\leq\ell$ and $\VC(G)\leq d$. We further assume that $n\geq
        n_{\equiPRG}^{\atom}(\ell,d,c_{\ZZ},c_{\samp},c_{\adj},c_{\Root},c_{\Sum},\epsilon,\delta)$.}
      \KwOut{A $(\delta,\epsilon)$-atomic equipartition of $G$ into
        $K_{\equiPRG}(\ell,c_{\ZZ},c_{\samp},c_{\adj},c_{\Root},c_{\Sum},\epsilon)$ parts.}
      $\widetilde{\delta}\assign(1+c_{\atom})\cdot\delta$\;
      Let $m\assign\ceil{C\cdot d/(c_{\atom}^2\cdot\delta^2)}$, where $C$ is the absolute constant of
      Lemma~\ref{lem:atomconst}.\;\label{alg:equiPRG:atom:m}
      $\widetilde{m}\assign\floor{\widetilde{\delta}\cdot m}$\;
      $\zeta\assign (1-c_{\ZZ})\cdot\epsilon$\;
      $\widetilde{c}_{\UC}\assign c_{\UC}\assign 0$\;
      $\widetilde{c}_{\adj}\assign (1-c_{\adj})^{-\ell} - 1$\;
      $\widetilde{c}_{\samp}\assign(1-c_{\samp})^{-\ell} - 1$\;
      $\widetilde{c}_{\ZZ}\assign(1-c_{\ZZ})^{-\ell-1} - 1$\;
      $\widetilde{c}_{\Root}\assign(1-c_{\Root})^{-\ell-1} - 1$\;
      $K\assign\ceil{(1+\widetilde{c}_{\UC})\cdot(1+\widetilde{c}_{\ZZ})\cdot(1+\widetilde{c}_{\samp})\cdot(1+\widetilde{c}_{\adj})
        \cdot(1+\widetilde{c}_{\Root})\cdot(1+c_{\Sum})\cdot\epsilon^{-\ell-1}}$\;
      $r\assign\floor{n/K}$\;
      $h\assign h_{\equiPRG}^{\atom}\df 2\cdot m + 2$\;
      \For{$x\in V(G)^m$}{%
        \label{alg:equiPRG:atom:prexFor}
        \For{$v\in V(G)$}{%
          \label{alg:equiPRG:atom:prevFor}
          $d_G(v,x)\assign\lvert\{i\in[m]\mid x_i\in N_G(v)\}\rvert$\;
          \lIf{$m-d_G(v,x)\leq\widetilde{m}$}{$t_G^{\widetilde{\delta}}(v,x)\assign 1$}
          \lElseIf{$d_G(v,x)\leq\widetilde{m}$}{$t_G^{\widetilde{\delta}}(v,x)\assign 0$}
        }
      }
      $\cQ\assign\varnothing$\;
      $R\assign V(G)$\;
      $n_0\assign n$\;
      \setcounter{algosplit}{\theAlgoLine}
    \end{algorithm}

    \begin{algorithm}[htbp]
      \caption*{(continued).}
      \DontPrintSemicolon
      \setcounter{AlgoLine}{\thealgosplit}
      \uFor{$i\in[K]$}{%
        \label{alg:equiPRG:atom:For}
        Find a rational $\xi_i\in\QQ$ satisfying:
        \begin{align*}
          0 & < \xi_i\leq(1-c_{\Root})\cdot\zeta,
          &
          0 & \leq -1 + \zeta - \xi_i + \gamma_i^\ell\cdot\frac{n_{i-1}}{r} \leq c_{\Root}\cdot\zeta,
        \end{align*}
        where
        \begin{align*}
          \gamma_i & \df \widetilde{\epsilon}_i \df (1-c_{\samp})\cdot\widetilde{\xi}_i,
          &
          \widetilde{\xi}_i & \df (1-c_{\adj})\cdot\xi_i,
        \end{align*}
        and $\gamma_i$ is a dyadic.\;\label{alg:equiPRG:atom:xii}
        $r_i\assign\ceil{\gamma_i^\ell\cdot n_{i-1}}$\;\label{alg:equiPRG:atom:ri}
        $s_0\df n_{i-1}$\;
        \lFor{$j\in[\ell+1]$}{\label{alg:equiPRG:atom:sj}$s_j\assign\floor{\gamma_i\cdot s_{j-1}}+1$}
        Run Algorithm~\ref{alg:ext:atom} with parameters:
        \begin{gather*}
          (n,m,\ell,s_0,\ldots,s_\ell,s_{\ell+1},\widetilde{m},G,U)
          \df
          (n,m,\ell,s_0,\ldots,s_\ell,s_{\ell+1},\widetilde{m},G,R).
        \end{gather*}
        \;\label{alg:equiPRG:atom:call}
        Let $(W,s)$ be the pair returned by the algorithm and let $j\in\{0,\ldots,\ell\}$ be such that $s=s_j$.\;
        Let $p\df\gamma_i^{\ell-j}$ and write it as $\widetilde{p}/2^t$ for $\widetilde{p},t\in\NN_+$ with
        $t\geq\ceil{\log_2(n)}$ smallest possible.\;
        \label{alg:equiPRG:atom:p}
        Run Algorithm~\ref{alg:encodePRG} with parameter $t\df t$ and let $P$ be the returned polynomial.\;
        \For{$Z\in[2^t]^h$}{%
          \label{alg:equiPRG:atom:ZFor}
          \For{$u\in[n]$}{%
            \label{alg:equiPRG:atom:uFor}
            Run Algorithm~\ref{alg:runPRG} with parameters:
            \begin{gather*}
              (P,\widetilde{p},h,n,Z,u)\df(P,\widetilde{p},h,n,Z,u).
            \end{gather*}
            \;
            Let $X_u$ be the bit returned by the algorithm.\;
          }
          $W(X)\assign\{u\in W \mid X_u=1\}$\;
          \lIf{$\lvert\lvert W(X)\rvert - p\cdot s\rvert > c_{\adj}\cdot c_{\samp}\cdot\widetilde{\xi}_i\cdot p\cdot s$}{
            \Continue
          }
          \uFor(\tcp*[f]{Check if $W(X)$ is a $(\widetilde{\delta},\widetilde{\xi}_i,m)$-atom}){$x\in V(G)^m$}{%
            \label{alg:equiPRG:atom:xFor}
            $k^0\assign\lvert\{u\in W(X) \mid t_G^{\widetilde{\delta}}(u,x) = 0\}\rvert$\;
            $k^1\assign\lvert\{u\in W(X) \mid t_G^{\widetilde{\delta}}(u,x) = 1\}\rvert$\;
            \If(\tcp*[f]{$W(X)$ is $(\widetilde{\delta},\widetilde{\xi}_i)$-split by $x$}){%
              $k^0 > \widetilde{\xi}_i\cdot\lvert W(X)\rvert$ and $k^1 > \widetilde{\xi}_i\cdot\lvert W(X)\rvert$
            }{%
              \Break\tcp*[r]{$W(X)$ is not a $(\widetilde{\delta},\widetilde{\xi}_i,m)$-atom}
            }
          }
          \Else(\tcp*[f]{Executes only if ``for'' loop of line~\ref{alg:equiPRG:atom:xFor} didn't break}){%
            \Break\tcp*[r]{$W(X)$ is a $(\widetilde{\delta},\widetilde{\xi}_i,m)$-atom}
          }
        }
        \setcounter{algosplit}{\theAlgoLine}
      }
    \end{algorithm}

    \begin{algorithm}[htbp]
      \caption*{(continued).}
      \DontPrintSemicolon
      \setcounter{AlgoLine}{\thealgosplit}
      \let\oldnl\nl
      \let\nl\relax
      \vspace{-\baselineskip}\InvisibleBegin{
        \global\let\nl\oldnl 
        Let $A_+\subseteq R\setminus W(X)$ and $A_-\subseteq W(X)$ be sets of sizes
        \begin{align*}
          \lvert A_+\rvert & = \max\{0, r_i - \lvert W(X)\rvert\}, &
          \lvert A_-\rvert & = \max\{0, \lvert W(X)\rvert - r_i\}.
        \end{align*}
        \;\label{alg:equiPRG:atom:A+A-}
        $W'\assign A_+\cup W(X)\setminus A_-$\;
        $\cQ\assign\cP\cup\{(W',r_i)\}$\;
        $R\assign R\setminus W'$\;
        $n_i\assign n_{i-1} - r_i$\;
      }
      $\cP\assign\varnothing$\;
      $b\assign n - K\cdot r$\;
      \For{$(W,s)\in\cQ$}{%
        \label{alg:equiPRG:atom:secondFor}
        Let $A\subseteq R$ be a set of size $r - s + \One[b > 0]$.\;
        $\cP\assign\cP\cup\{W\cup A\}$\;
        $R\assign R\setminus A$\;
        $b\assign b-1$\;
      }
      \Return{$\cP$}
    \end{algorithm}

    First, we note that the block of line~\ref{alg:equiPRG:atom:prexFor} correctly computes all the
    $\widetilde{\delta}$-majority opinions $t_G^{\widetilde{\delta}}(v,x)$ of all $m$-tuples $x\in V(G)^m$ with respect to all
    vertices $v\in V(G)$ in $G$ when they exist (and leaves $t_G^{\widetilde{\delta}}(v,x)$ undefined when there is no
    $\widetilde{\delta}$-majority opinion; recalling also that since $\widetilde{\delta} < 1/2^{\ell+1}\leq 1/2$ there can only
    be at most one such $\widetilde{\delta}$-majority opinion). This is the same argument as in
    Proposition~\ref{prop:ext}\ref{prop:ext:atomalg}: since $\widetilde{m} = \floor{\widetilde{\delta}\cdot m}$, for
    $d_G(v,x)\df\lvert\{i\in[m]\mid x_i\in N_G(v)\}\rvert\in\NN$, we have
    \begin{gather*}
      d_G(v,x) \geq (1-\widetilde{\delta})\cdot m
      \iff
      m-d_G(v,x) \leq \widetilde{\delta}\cdot m
      \iff
      m-d_G(v,x) \leq \widetilde{m},
      \\
      d_G(v,x) \leq \widetilde{\delta}\cdot m
      \iff
      d_G(v,x) \leq \widetilde{m}.
    \end{gather*}

    Now, since all the parameters $\zeta$, $K$, $\xi_i$, $\widetilde{\xi}_i$, $\widetilde{\epsilon}_i$, $\gamma_i$ of
    Algorithm~\ref{alg:equiPRG:atom} are being picked according to Lemma~\ref{lem:equicalc}, item~\ref{lem:equicalc:xii} of the
    lemma ensures that the rationals $\xi_i\in\QQ$ in line~\ref{alg:equiPRG:atom:xii} indeed exist satisfying
    \begin{align}\label{eq:equiPRG:atom:xii}
      0 & < \xi_i\leq\zeta,
      &
      -1 + \zeta - \xi_i + \gamma_i^\ell\cdot\frac{n_{i-1}}{r} & = 0,
    \end{align}
    where
    \begin{align*}
      \gamma_i & \df \widetilde{\epsilon}_i \df (1-c_{\samp})\cdot\widetilde{\xi}_i,
      &
      \widetilde{\xi}_i & \df (1-c_{\adj})\cdot\xi_i,
    \end{align*}
    with $\gamma_i$ dyadic.

    Let us focus on the $i$th iteration of the loop of line~\ref{alg:equiPRG:atom:For}: by
    Proposition~\ref{prop:ext}\ref{prop:ext:atomalg}, we know that the pair $(W,s)$ returned by Algorithm~\ref{alg:ext:atom} in
    line~\ref{alg:equiPRG:atom:call} is such that $\lvert W\rvert = s\in\{s_0,\ldots,s_\ell\}$ and $W$ is a
    $(\widetilde{\delta},\widetilde{\epsilon}_i,m)$-atom of $G$. Then by Theorem~\ref{thm:PRG}, Algorithms~\ref{alg:encodePRG}
    and~\ref{alg:runPRG} together are responsible for making the variables $Z\in[2^t]^h$ and $X\in\{0,1\}^n$ be such that the
    function $Z\mapsto X$ is a pseudorandom generator of bitlength $n$, independence $h$, bias
    $p\df\widetilde{p}/2^t=\gamma_i^{\ell-j}$ and seed space $[2^t]^h$.

    In turn, by Proposition~\ref{prop:PRGsubset}\ref{prop:PRGsubset:atom}, we know that if
    \begin{gather}\label{eq:equiPRG:atom:conditions}
      \floor{\widetilde{\epsilon}_i\cdot s}
      \geq
      \begin{dcases*}
        \frac{(h-1)^2}{p\cdot(1-p)}, & if $p < 1$,\\
        0, & if $p=1$,
      \end{dcases*}
      \\
      \frac{1}{\sqrt{2\cdot\pi\cdot h}}
      \cdot
      \left(\frac{2\cdot e\cdot (1-p)}{c_{\samp}^2\cdot\widetilde{\xi}_i^2\cdot p}\right)^{h/2}
      \cdot
      \left(
      \left(\frac{1}{c_{\adj}^2\cdot s}\right)^{h/2}
      +
      n^m\cdot\left(\frac{(1-c_{\samp})^2}{(1-c_{\adj})^2\cdot\floor{\widetilde{\epsilon}_i\cdot s}}\right)^{h/2}
      \right)
      \leq
      1,
    \end{gather}
    then there must exist $Z\in[2^t]^h$ such that $W(X)\df\{u\in W \mid X_u=1\}$ is a
    $(\widetilde{\delta},\widetilde{\xi}_i,m)$-atom of $G$ and
    \begin{gather}\label{eq:equiPRG:atom:WXsize}
      \bigl\lvert\lvert W(X)\rvert - p\cdot s\bigr\rvert
      \leq
      c_{\adj}\cdot c_{\samp}\cdot\widetilde{\xi}_i\cdot p\cdot s.
    \end{gather}
    Since~\eqref{eq:equiPRG:atom:conditions} is guaranteed to hold by Lemma~\ref{lem:equiPRGcalc}, items~\ref{lem:equiPRGcalc:h}
    and~\ref{lem:equiPRGcalc:atom}, we are guaranteed the existence of such $Z$. It is then straightforward to check that the
    result of the loop of line~\eqref{alg:equiPRG:atom:ZFor} is that when it is over, the $Z$ variable (and $W(X)$ by
    association) has the desired properties.

    Now the final part of the loop of line~\ref{alg:equiPRG:atom:For} that starts on line~\ref{alg:equiPRG:atom:A+A-} simply
    takes $W(X)$ and either adds or removes vertices (as few as possible) to make it size $r_i$ (note that this is always
    possible since the remainder $R$ is guaranteed to have size $n_{i-1}\geq r\geq r_i$ by Lemma~\ref{lem:equicalc},
    items~\ref{lem:equicalc:eval} and~\ref{lem:equicalc:ri}).

    Since $W(X)$ is a $(\widetilde{\delta},\widetilde{\xi}_i,m)$-atom of $G$ and we have the
    bound~\eqref{eq:equiPRG:atom:WXsize} on the size of $W(X)$, by Corollary~\ref{cor:cont}\ref{cor:cont:atom}, it follows that
    the set $W'\df A_+\cup W(X)\setminus A_-$ of size $r_i$ that is included in $\cQ$ is a
    $(\widetilde{\delta},\epsilon_{\symdiff},m)$-atom of $G$, where
    \begin{align*}
      \epsilon_{\symdiff}
      & \df
      \widetilde{\xi}_i + \frac{\lvert W(X)\symdiff W'\rvert}{r_i}
      \leq
      \widetilde{\xi}_i + \frac{\lvert p\cdot s - r_i\rvert + c_{\adj}\cdot c_{\samp}\cdot\widetilde{\xi}_i\cdot p\cdot s}{r_i}
      \\
      & \leq
      \widetilde{\xi}_i + \frac{c_{\adj}\cdot c_{\samp}\cdot\widetilde{\xi}_i\cdot p\cdot s + 2}{r_i}
      \leq
      \xi_i,
    \end{align*}
    where the second inequality follows since Lemma~\ref{lem:equiPRGcalc}\ref{lem:equiPRGcalc:gammai} implies
    $r_i\df\ceil{\gamma_i^\ell\cdot n_{i-1}}\leq p\cdot s\leq r_i+2$ and the last inequality follows from
    Lemma~\ref{lem:equiPRGcalc}\ref{lem:equiPRGcalc:adj}. Therefore $W'$ is a $(\widetilde{\delta},\xi_i,m)$-atom of $G$ of size
    $r_i$.

    This means that when the loop of line~\ref{alg:equiPRG:atom:For} concludes, the variable $\cQ$ holds a set of $K$ pairs of the
    form $(W'[i],r_i)$, where $(W'[i],r_i)$ is the pair added in the $i$th iteration and is such that $W'[i]$ is a
    $(\widetilde{\delta},\xi_i,m)$-atom of $G$ of size $r_i$ (and the $W'[i]$ are pairwise disjoint).
    
    Finally, the final block of line~\ref{alg:equiPRG:atom:secondFor} takes these sets $W'[i]$ and adds vertices of the
    remainder $R\df V(G)\setminus\bigcup_{i\in[K]} W'[i]$ so as to make them all have size either $r$ or $r+1$ and remain
    pairwise disjoint. This is possible since $r = \floor{n/K}$ and $r_i\leq r$ (by
    Lemma~\ref{lem:equicalc}\ref{lem:equicalc:ri}). It is then clear from the sizes of the sets that the resulting $\cP$ is an
    equipartition of $V(G)$ (i.e., there are no vertices left as we correctly dealt with the potential non-divisibility of $n$
    by $r\df\floor{n/K}$).

    We now argue that the parts of $\cP$ are indeed $(\widetilde{\delta},\epsilon,m)$-atoms of $G$. Fix one such part $P$ and
    assume it came from $W'[i]$, which we know is a $(\widetilde{\delta},\xi_i,m)$-atom of $G$ of size $r_i$. We can let
    $P'\subseteq P$ be such that $W'[i]\subseteq P'$ and $\lvert P'\rvert=r$ so that $\lvert P\setminus P'\rvert\leq 1$. Then by
    Lemma~\ref{lem:upcont}\ref{lem:upcont:atom}, we know that $P'$ is a $(\widetilde{\delta},\epsilon_+,m)$-atom of $G$,
    where
    \begin{gather*}
      \epsilon_+
      \df
      (\xi_i-1)\cdot\frac{\lvert W'[i]\rvert}{\lvert P'\rvert} + 1
      =
      (\xi_i-1)\cdot\frac{r_i}{r} + 1
      \leq
      \zeta,
    \end{gather*}
    where the last inequality follows from Lemma~\ref{lem:equicalc}\ref{lem:equicalc:cont}, so $P'$ is a
    $(\widetilde{\delta},\zeta,m)$-atom of $G$.

    On the other hand, since $\lvert P\setminus P'\rvert\leq 1$, using Lemma~\ref{lem:upcont}\ref{lem:upcont:atom} again, we
    know that $P$ is a $(\widetilde{\delta},\epsilon_{++},m)$-atom of $G$, where
    \begin{gather*}
      \epsilon_{++}
      \df
      (\zeta-1)\cdot\frac{r}{r+1} + 1
      \leq
      \epsilon,
    \end{gather*}
    where the last inequality follows from Lemma~\ref{lem:equicalc:ZZ}.

    Thus, $\cP$ is a equipartition into $(\widetilde{\delta},\epsilon,m)$-atoms of $G$ with $K$ parts. Finally, since
    $\widetilde{\delta}=(1+c_{\atom})\cdot\delta$, by Lemma~\ref{lem:atomconst} and our choice of $m$, it follows that these
    sets are $(\delta,\epsilon)$-atoms of $G$, concluding the proof of correctness of Algorithm~\ref{alg:equiPRG:atom}.

    \medskip

    We now analyze the time-complexity of Algorithm~\ref{alg:equiPRG:atom}. The computations of the parameters
    \begin{align*}
      \widetilde{\delta}, & &
      m, & &
      \widetilde{m}, & &
      \zeta, & &
      c_{\UC}=\widetilde{c}_{\UC}=0,
      \widetilde{c}_{\adj},
      \\
      \widetilde{c}_{\samp}, & &
      \widetilde{c}_{\ZZ}, & &
      \widetilde{c}_{\Root}, & &
      K, & &
      h, & &
      r,
    \end{align*}
    are easily seen to take time at most
    \begin{gather*}
      \begin{aligned}
        O\bigl(\len(c_{\atom})\cdot\len(\delta)\bigr),
        & &
        O\Bigl(
        \log(d+1)\cdot\len(c_{\atom})
        + \bigl(\log(d+1) + \len(c_{\atom})\bigr)\cdot\len(\delta)
        \Bigr),
      \end{aligned}
      \\
      \begin{aligned}
        O\bigl(\len(c_{\ZZ})\cdot\len(\epsilon)\bigr),
        & &
        O(1),
      \end{aligned}
      \\
      \begin{aligned}
        O\bigl(\ell^2\cdot\log(\ell+1)\cdot\len(c_{\adj})^2\bigr),
        &
        O\bigl(\ell^2\cdot\log(\ell+1)\cdot\len(c_{\samp})^2\bigr),
      \end{aligned}
      \\      
      \begin{aligned}
        O\bigl(\ell^2\cdot\log(\ell+1)\cdot\len(c_{\ZZ})^2\bigr),
        & &
        O\bigl(\ell^2\cdot\log(\ell+1)\cdot\len(c_{\Root})^2\bigr),
      \end{aligned}
      \\
      O\bigl(
      \ell^4\cdot\log(\ell+1)\cdot\len(c_{\ZZ})\cdot\len(c_{\samp})\cdot\len(c_{\Root})\cdot\len(c_{\Sum})\cdot\len(\epsilon)^2
      \bigr),
      \\
      O\bigl(\log(d+1) + \len(c_{\atom}) + \len(\delta)\bigr),
      \\
      O\Bigl(\bigl(\log(n+1)\bigr)^2\Bigr),
    \end{gather*}
    respectively. These together are bounded by
    \begin{gather}\label{eq:equiPRG:atom:param:complexity}
      \begin{multlined}
        O\Bigl(
        \log(d+1)\cdot\len(c_{\atom})\cdot\len(\delta)
        \\
        +
        \ell^4\cdot\log(\ell+1)\cdot\len(c_{\ZZ})^2\cdot\len(c_{\samp})^2\cdot\len(c_{\Root})^2\cdot\len(c_{\Sum})\cdot\len(\epsilon)^2
        \cdot\bigl(\log(n+1)\bigr)^2
        \Bigr).
      \end{multlined}
    \end{gather}

    Clearly, the loops of lines~\ref{alg:equiPRG:atom:For} and~\ref{alg:equirand:atom:secondFor} execute exactly $K\leq
    O(\epsilon^{-\ell-1})$ times. Let us analyze the time-complexity of the computations of line~\ref{alg:equirand:atom:xii}.
    The line requires us to find the rational satisfying
    \begin{align*}
      0 & < \xi_i\leq(1-c_{\Root})\cdot\zeta,
      &
      0 & \leq -1 + \zeta - \xi_i + \gamma_i^\ell\cdot\frac{n_{i-1}}{r} \leq c_{\Root}\cdot\zeta,
    \end{align*}
    where
    \begin{align*}
      \gamma_i & \df \widetilde{\epsilon}_i \df (1-c_{\samp})\cdot\widetilde{\xi}_i,
      &
      \widetilde{\xi}_i & \df (1-c_{\adj})\cdot\xi_i,
    \end{align*}
    and such that $\gamma_i$ is a dyadic. The argument is a slight extension of how we found $\xi_i$ in
    Algorithm~\ref{alg:equirand:atom}: if we let $f(\xi_i)$ be the middle term of the second set of inequalities above, then
    Lemma~\ref{lem:equicalc}\ref{lem:equicalc:xii} says that $f(0) < 0$ and $f((1-c_{\Root})\cdot\zeta)\geq
    c_{\Root}\cdot\zeta$, so to find the desired $\xi_i$, we proceed as follows: first, we perform a binary search to find an
    auxiliary $\xi'_i$ satisfying instead
    \begin{align*}
      0 & < \xi'_i\leq(1-c_{\Root})\cdot\zeta,
      &
      \frac{1}{3}\cdot c_{\Root}\cdot\zeta
      & \leq
      -1 + \zeta - \xi'_i + (\gamma'_i)^\ell\cdot\frac{n_{i-1}}{r}
      \leq
      \frac{2}{3}\cdot c_{\Root}\cdot\zeta,
    \end{align*}
    where $\gamma'_i$ is defined in terms of $\xi'_i$ in the same manner as $\gamma_i$ is defined in terms of $\xi_i$. Since the
    function $f$ is Lipschitz-continuous with constant at most
    \begin{align*}
      L & \df
      \ell\cdot\frac{n_{i-1}}{r}
      \leq
      \ell\cdot\frac{n}{\floor{n/K}}
      \leq
      \ell\cdot\frac{K}{1 - K/n}
      \\
      & \leq
      O(\ell\cdot K)
      \leq
      O(\ell\cdot\epsilon^{-\ell-1}),
    \end{align*}
    this initial binary search is guaranteed to take at most
    \begin{align*}
      \log_2\left(\frac{(1-c_{\Root})\cdot\zeta\cdot L}{c_{\Root}\cdot\zeta/3}\right)
      & \leq
      \log_2\left(\frac{3\cdot (1-c_{\Root})\cdot\ell}{c_{\Root}\cdot(1-c_{\ZZ})\cdot\epsilon^{\ell+2}}\right)
      +
      O(1)
      \\
      & \leq
      O\bigl(\len(c_{\Root}) + \len(c_{\ZZ}) + \ell\cdot\len(\epsilon)\bigr)
    \end{align*}
    steps. This also implies that at all stages of the binary search, the bitlength of the endpoints is at most
    \begin{align*}
      \MoveEqLeft
      O\bigl(\len(c_{\Root}) + \len(c_{\ZZ}) + \ell\cdot\len(\epsilon) + \len(c_{\samp}) + \len(\zeta)\bigr)
      \\
      & \leq
      O\bigl(\len(c_{\Root}) + \len(c_{\ZZ}) + \ell\cdot\len(\epsilon) + \len(c_{\samp}) + \len(c_{\adj})\bigr).
    \end{align*}
    After this is done, we take
    \begin{gather*}
      \xi''_i \df \xi'_i + \frac{c_{\Root}\cdot\zeta}{3\cdot L}
    \end{gather*}
    so that Lipschitz-continuity implies that for every $\xi_i\in[\xi'_i,\xi''_i]$, we have $\lvert f(\xi_i)-f(\xi_i')\rvert\leq
    c_{\Root}\cdot\zeta/3$, hence $0\leq f(\xi_i)\leq c_{\Root}\cdot\zeta$. Translating the condition $\xi_i\in[\xi'_i,\xi''_i]$
    to $\gamma_i$ yields exactly
    \begin{gather*}
      \gamma'_i \leq \gamma_i \leq \gamma''_i
    \end{gather*}
    where $\gamma''_i$ is defined in terms of $\xi''_i$ in the same manner as $\gamma_i$ is defined in terms of $\xi_i$ and
    since $\gamma'_i < \gamma''_i$, we can find a dyadic $\gamma_i$ satisfying the above of bitlength at most
    \begin{gather*}
      O\bigl(\len(\gamma'_i) + \len(\gamma''_i - \gamma'_i)\bigr)
      \leq
      O\bigl(\len(c_{\Root}) + \len(c_{\ZZ}) + \ell\cdot\len(\epsilon) + \len(c_{\samp}) + \len(c_{\adj})\bigr).
    \end{gather*}
    Accounting for the fact that the computation of $\xi''_i$ takes time at most
    \begin{align*}
      \MoveEqLeft
      O\Bigl(
      \len(\xi'_i) + \len(c_{\Root})\cdot\bigl(\len(c_{\ZZ})+\len(\epsilon)\bigr)\cdot\log(\ell+1)\cdot\bigl(\log(n+1)\bigr)^2
      \Bigr)
      \\
      & \leq
      O\Bigl(\len(c_{\Root})\cdot\len(c_{\ZZ})\cdot\len(c_{\samp})\cdot\len(c_{\adj})
      \cdot\ell\cdot\len(\epsilon)\cdot\bigl(\log(n+1)\bigr)^2\Bigr),
    \end{align*}
    the time-complexity of a single execution of line~\ref{alg:equiPRG:atom:xii} is at most
    \begin{gather}\label{alg:equiPRG:atom:xii:single}
      O\Bigl(
      \len(c_{\Root})\cdot\len(c_{\ZZ})\cdot\len(c_{\samp})\cdot\len(c_{\adj})
      \cdot\ell\cdot\len(\epsilon)\cdot\bigl(\log(n+1)\bigr)^2
      \Bigr).
    \end{gather}

    It is straightforward to see that the time-complexity of the computations of $r_i$, the $s_j$ and $p$ in
    lines~\ref{alg:equiPRG:atom:ri}, \ref{alg:ext:atom} and~\ref{alg:equiPRG:atom:p} are completely dominated
    by~\eqref{alg:equiPRG:atom:xii:single}.

    By Proposition~\ref{prop:ext}\ref{prop:ext:atomalg}, the call to Algorithm~\ref{alg:ext:atom} in
    line~\ref{alg:equiPRG:atom:call} costs $O((\ell+m)\cdot n^{m+1})$.

    Note further that from our calculation of $\gamma_i$, the variable $t$ computed in line~\ref{alg:equiPRG:atom:p}
    \begin{align*}
      t
      & \leq
      \len(p)
      \leq
      \max\{\log_2(n)+1, \ell\cdot\len(\gamma_i)\}
      \\
      & \leq
      \log_2(n) + O\Bigl(\ell\cdot\bigl(\len(c_{\Root}) + \len(c_{\ZZ}) + \ell\cdot\len(\epsilon)
      + \len(c_{\samp}) + \len(c_{\adj})\bigr)\Bigr).
    \end{align*}
    This in particular means that the inner loop of line~\ref{alg:equiPRG:atom:ZFor} executes at most
    \begin{gather*}
      (2^t)^h
      =
      (2^t)^{2\cdot m + 2}
      \leq
      2^{O(m\cdot\ell\cdot(\len(c_{\Root}) + \len(c_{\ZZ}) + \ell\cdot\len(\epsilon) + \len(c_{\samp}) + \len(c_{\adj})))}\cdot n^{2\cdot m+2}.
    \end{gather*}
    times. The further inner loop of line~\ref{alg:equiPRG:atom:uFor} executes $n$ times, with each iteration costing $O(h\cdot
    t^2)$ time, which amounts to a total time of $O(n\cdot h\cdot t^2)$. Finally, the remaining part of the loop of
    line~\ref{alg:equiPRG:atom:ZFor}, which includes the loop of line~\ref{alg:equiPRG:atom:xFor}, is easily seen to take time
    at most $O(n^{m+1})$. Thus, the time-complexity of the loop of line~\ref{alg:equiPRG:atom:ZFor}, accounting for all the
    $(2^t)^h$ times it executes (within a single iteration of the loop of line~\ref{alg:equiPRG:atom:uFor}), is at most
    \begin{gather}\label{eq:equiPRG:atom:ZFor:complexity}
      2^{O(m\cdot\ell\cdot(\len(c_{\Root}) + \len(c_{\ZZ}) + \ell\cdot\len(\epsilon) + \len(c_{\samp}) + \len(c_{\adj})))}\cdot n^{3\cdot m+3}.
    \end{gather}
    This completely dominates the time-complexity of the call to Algorithm~\ref{alg:ext:atom} in
    line~\ref{alg:equiPRG:atom:call}.

    Finally, since the computations of the remainder of the loop of line~\ref{alg:equiPRG:atom:uFor} starting in
    line~\ref{alg:equiPRG:atom:A+A-} take time at most $O(n)$, it follows that the total time-complexity of the loop of
    line~\ref{alg:equiPRG:atom:For}, accounting for all $K\leq O(\epsilon^{-\ell-1})$ times it executes, is at most
    \begin{gather}\label{eq:equiPRG:atom:For:complexity}
      \begin{aligned}
        \MoveEqLeft
        \begin{multlined}[t]
          O\biggl(
          \epsilon^{-\ell-1}\cdot
          \Bigl(
          \len(c_{\Root})\cdot\len(c_{\ZZ})\cdot\len(c_{\samp})
          \cdot\ell\cdot\len(\epsilon)\cdot\bigl(\log(n+1)\bigr)^2
          \\
          +
          2^{O(m\cdot\ell\cdot(\len(c_{\Root}) + \len(c_{\ZZ}) + \ell\cdot\len(\epsilon) + \len(c_{\samp}) + \len(c_{\adj})))}\cdot n^{3\cdot m+3}.
          +
          n
          \Bigr)
          \biggr)
        \end{multlined}
        \\
        & \leq
        2^{O(m\cdot\ell\cdot(\len(c_{\Root}) + \len(c_{\ZZ}) + \ell\cdot\len(\epsilon) + \len(c_{\samp}) + \len(c_{\adj})))}\cdot n^{3\cdot m+3}.
      \end{aligned}
    \end{gather}

    The remainder of algorithm after the loop of line~\ref{alg:equiPRG:atom:uFor} is concluded is easily seen to take time at
    most $O(\epsilon^{-\ell-1}\cdot n)$, which is dominated by~\eqref{eq:equiPRG:atom:For:complexity}. Thus, the final
    time-complexity of Algorithm~\ref{alg:equiPRG:atom} is at most the sum of~\eqref{eq:equiPRG:atom:param:complexity}
    and~\eqref{eq:equiPRG:atom:For:complexity}:
    \begin{align*}
      \MoveEqLeft
      \begin{multlined}[t]
        O\Bigl(
        \log(d+1)\cdot\len(c_{\atom})\cdot\len(\delta)
        \\
        +
        \ell^4\cdot\log(\ell+1)\cdot\len(c_{\ZZ})^2\cdot\len(c_{\samp})^2\cdot\len(c_{\Root})^2\cdot\len(c_{\Sum})\cdot\len(\epsilon)^2
        \cdot\bigl(\log(n+1)\bigr)^2
        \Bigr)
        \\
        +
        2^{O(m\cdot\ell\cdot(\len(c_{\Root}) + \len(c_{\ZZ}) + \ell\cdot\len(\epsilon) + \len(c_{\samp}) + \len(c_{\adj})))}\cdot n^{3\cdot m+3}
      \end{multlined}
      \\
      & \leq
      \begin{multlined}[t]
        2^{O(d\cdot\ell\cdot(\len(c_{\Root})+\len(c_{\ZZ})+\ell\cdot\len(\epsilon)+\len(c_{\samp}) + \len(c_{\adj})))/(c_{\atom}^2\cdot\delta^2)}
        \cdot n^{O(d/c_{\atom}^2\cdot\delta^2)}
        \\
        +
        O\bigl(
        \log(d+1)\cdot\len(c_{\atom})\cdot\len(\delta)
        \\
        +
        \ell^4\cdot\log(\ell+1)\cdot\len(c_{\ZZ})^2\cdot\len(c_{\samp})^2\cdot\len(c_{\Root})^2\cdot\len(c_{\Sum})\cdot\len(\epsilon)^2
        \cdot\bigl(\log(n+1)\bigr)^2
        \Bigr),
      \end{multlined}
    \end{align*}
    where the last inequality follows since
    \begin{gather*}
      m
      =
      \Ceil{C\cdot\frac{d}{c_{\atom}^2\cdot\delta^2}}
      \leq
      O\left(\frac{d}{c_{\atom}^2\cdot\delta^2}\right).
    \end{gather*}
  \item[Item~\ref{thm:equiPRG:good}.] Algorithm~\ref{alg:equiPRG:good} is essentially the same as
    Algorithm~\ref{alg:equiPRG:atom}, except that we use $m\df 1$ as $(\delta,\epsilon,1)$-atomicity is equivalent to
    $\epsilon$-goodness (provided $\delta<1$), as such, we omit its proof of correctness as it is completely analogous to that
    of Algorithm~\ref{alg:equiPRG:atom} (but using the results for goodness instead of atomicity of all invoked lemmas and
    propositions).

    \begin{algorithm}[htbp]
      \caption{Deterministic model algorithm that returns an $\epsilon$-good equipartition of $G$ into
        $K_{\equiPRG}(\ell,c_{\ZZ},c_{\samp},c_{\adj},c_{\Root},c_{\Sum},\epsilon)$ parts in time
        \begin{multline*}
          2^{O(\ell\cdot(\len(c_{\Root}) + \len(c_{\ZZ}) + \ell\cdot\len(\epsilon) + \len(c_{\samp}) + \len(c_{\adj})))}\cdot n^6
          \\
          +
          O\bigl(
          \ell^4\cdot\log(\ell+1)\cdot\len(c_{\ZZ})^2\cdot\len(c_{\samp})^2\cdot\len(c_{\Root})^2\cdot\len(c_{\Sum})\cdot\len(\epsilon)^2
          \cdot\bigl(\log(n+1)\bigr)^2
          \bigr).
        \end{multline*}
      }
      \label{alg:equiPRG:good}
      \DontPrintSemicolon
      \KwIn{Numbers $n,\ell\in\NN_+$ with $\ell\leq\log_2(n)$, $c_{\ZZ},c_{\adj},c_{\Root},c_{\Sum}\in(0,1)\cap\QQ$,
        $c_{\samp}\in(0,1/2)\cap\QQ$, and
        $\epsilon\in(0,\epsilon_{\equiPRG}(\ell,c_{\ZZ},c_{\samp},c_{\adj},c_{\Root},c_{\Sum}))$, and a graph $G$ with
        $V(G)=[n]$ and $\Lit(G)\leq\ell$. We further assume that $n\geq
        n_{\equiPRG}^{\good}(\ell,c_{\ZZ},c_{\samp},c_{\adj},c_{\Root},c_{\Sum},\epsilon)$.}
      \KwOut{An $\epsilon$-good equipartition of $G$ into
        $K_{\equiPRG}(\ell,c_{\ZZ},c_{\samp},c_{\adj},c_{\Root},c_{\Sum},\epsilon)$ parts.}
      $m\assign 1$\;
      $\widetilde{m}\assign 1$\tcp*[r]{$\widetilde{m} = \floor{\widetilde{\delta}\cdot m}$ for $\widetilde{\delta} < 1$}
      $\zeta\assign (1-c_{\ZZ})\cdot\epsilon$\;
      $\widetilde{c}_{\UC}\assign c_{\UC}\assign 0$\;
      $\widetilde{c}_{\adj}\assign (1-c_{\adj})^{-\ell} - 1$\;
      $\widetilde{c}_{\samp}\assign(1-c_{\samp})^{-\ell} - 1$\;
      $\widetilde{c}_{\ZZ}\assign(1-c_{\ZZ})^{-\ell-1} - 1$\;
      $\widetilde{c}_{\Root}\assign(1-c_{\Root})^{-\ell-1} - 1$\;
      $K\assign\ceil{(1+\widetilde{c}_{\UC})\cdot(1+\widetilde{c}_{\ZZ})\cdot(1+\widetilde{c}_{\samp})\cdot(1+\widetilde{c}_{\adj})
        \cdot(1+\widetilde{c}_{\Root})\cdot(1+c_{\Sum})\cdot\epsilon^{-\ell-1}}$\;
      $r\assign\floor{n/K}$\;
      $h\assign h_{\equiPRG}^{\good}\df 4$\;
      $\cQ\assign\varnothing$\;
      $R\assign V(G)$\;
      $n_0\assign n$\;
      \uFor{$i\in[K]$}{%
        \label{alg:equiPRG:good:For}
        Find a rational $\xi_i\in\QQ$ satisfying:
        \begin{align*}
          0 & < \xi_i\leq(1-c_{\Root})\cdot\zeta,
          &
          0 & \leq -1 + \zeta - \xi_i + \gamma_i^\ell\cdot\frac{n_{i-1}}{r} \leq c_{\Root}\cdot\zeta,
        \end{align*}
        where
        \begin{align*}
          \gamma_i & \df \widetilde{\epsilon}_i \df (1-c_{\samp})\cdot\widetilde{\xi}_i,
          &
          \widetilde{\xi}_i & \df \xi_i,
        \end{align*}
        and $\gamma_i$ is a dyadic.\;\label{alg:equiPRG:good:xii}
        $r_i\assign\ceil{\gamma_i^\ell\cdot n_{i-1}}$\;\label{alg:equiPRG:good:ri}
        $s_0\df n_{i-1}$\;
        \lFor{$j\in[\ell+1]$}{\label{alg:equiPRG:good:sj}$s_j\assign\floor{\gamma_i\cdot s_{j-1}}+1$}
        %
        \setcounter{algosplit}{\theAlgoLine}
      }
    \end{algorithm}

    \begin{algorithm}[htbp]
      \caption*{(continued).}
      \DontPrintSemicolon
      \setcounter{AlgoLine}{\thealgosplit}
      \let\oldnl\nl
      \let\nl\relax
      \vspace{-\baselineskip}\InvisibleBegin{
        \global\let\nl\oldnl 
        Run Algorithm~\ref{alg:ext:good} with parameters:
        \begin{gather*}
          (n,m,\ell,s_0,\ldots,s_\ell,s_{\ell+1},\widetilde{m},G,U)
          \df
          (n,m,\ell,s_0,\ldots,s_\ell,s_{\ell+1},\widetilde{m},G,R).
        \end{gather*}
        \;\label{alg:equiPRG:good:call}
        Let $(W,s)$ be the pair returned by the algorithm and let $j\in\{0,\ldots,\ell\}$ be such that $s=s_j$.\;
        Let $p\df\gamma_i^{\ell-j}$ and write it as $\widetilde{p}/2^t$ for $\widetilde{p},t\in\NN_+$ with
        $t\geq\ceil{\log_2(n)}$ smallest possible.\;
        Run Algorithm~\ref{alg:encodePRG} with parameter $t\df t$ and let $P$ be the returned polynomial.\;
        \For{$Z\in[2^t]^h$}{%
          \For{$u\in[n]$}{%
            Run Algorithm~\ref{alg:runPRG} with parameters:
            \begin{gather*}
              (P,\widetilde{p},h,n,Z,u)\df(P,\widetilde{p},h,n,Z,u).
            \end{gather*}
            \;
            Let $X_u$ be the bit returned by the algorithm.\;
          }
          $W(X)\assign\{u\in W \mid X_u=1\}$\;
          \lIf{$\lvert\lvert W(X)\rvert - p\cdot s\rvert > c_{\adj}\cdot c_{\samp}\cdot\widetilde{\xi}_i\cdot p\cdot s$}{%
            \Continue
          }
          \uFor(\tcp*[f]{Check if $W(X)$ is $\widetilde{\xi}_i$-good}){$x\in V(G)$}{%
            \label{alg:equiPRG:good:xFor}
            \If(\tcp*[f]{$W(X)$ is not $\widetilde{\xi}_i$-good}){%
              $\widetilde{\xi}_i < \lvert W(X)\cap N_G(x)\rvert/\lvert W(X)\rvert < 1 - \widetilde{\xi}_i$
            }{%
              \Break
            }
          }
          \Else(\tcp*[f]{Executes only if ``for'' loop of line~\ref{alg:equiPRG:good:xFor} didn't break}){%
            \Break\tcp*[r]{$W(X)$ is $\widetilde{\xi}_i$-good}
          }
        }
        Let $A_+\subseteq R\setminus W(X)$ and $A_-\subseteq W(X)$ be sets of sizes
        \begin{align*}
          \lvert A_+\rvert & = \max\{0, r_i - \lvert W(X)\rvert\}, &
          \lvert A_-\rvert & = \max\{0, \lvert W(X)\rvert - r_i\}.
        \end{align*}
        \;
        $W'\assign A_+\cup W(X)\setminus A_-$\;
        $\cQ\assign\cP\cup\{(W',r_i)\}$\;
        $R\assign R\setminus W'$\;
        $n_i\assign n_{i-1} - r_i$\;
      }
      $\cP\assign\varnothing$\;
      $b\assign n - K\cdot r$\;
      \For{$(W,s)\in\cQ$}{%
        \label{alg:equiPRG:good:secondFor}
        Let $A\subseteq R$ be a set of size $r - s + \One[b > 0]$.\;
        $\cP\assign\cP\cup\{W\cup A\}$\;
        $R\assign R\setminus A$\;
        $b\assign b-1$\;
      }
      \Return{$\cP$}
    \end{algorithm}

    Furthermore, the time-complexity analysis of Algorithm~\ref{alg:equiPRG:good} is also completely analogous to that of
    Algorithm~\ref{alg:equiPRG:atom} and yields the bound
    \begin{align*}
      \MoveEqLeft
      \begin{multlined}[t]
        O\Bigl(
        \ell^4\cdot\log(\ell+1)\cdot\len(c_{\ZZ})^2\cdot\len(c_{\samp})^2\cdot\len(c_{\Root})^2\cdot\len(c_{\Sum})\cdot\len(\epsilon)^2
        \cdot\bigl(\log(n+1)\bigr)^2
        \Bigr)
        \\
        +
        2^{O(m\cdot\ell\cdot(\len(c_{\Root}) + \len(c_{\ZZ}) + \ell\cdot\len(\epsilon) + \len(c_{\samp}) + \len(c_{\adj})))}\cdot n^{3\cdot m+3}.
      \end{multlined}
      \\
      & \leq
      \begin{multlined}[t]
        2^{O(\ell\cdot(\len(c_{\Root}) + \len(c_{\ZZ}) + \ell\cdot\len(\epsilon) + \len(c_{\samp}) + \len(c_{\adj})))}\cdot n^6
        \\
        +
        O\Bigl(
        \ell^4\cdot\log(\ell+1)\cdot\len(c_{\ZZ})^2\cdot\len(c_{\samp})^2\cdot\len(c_{\Root})^2\cdot\len(c_{\Sum})\cdot\len(\epsilon)^2
        \cdot\bigl(\log(n+1)\bigr)^2
        \Bigr).
      \end{multlined}
    \end{align*}
  \item[Item~\ref{thm:equiPRG:atomspace}.] Similarly to Theorem~\ref{thm:nonequi}\ref{thm:nonequi:atomspace}, the idea behind
    Algorithms~\ref{alg:equiPRG:atomspace} and~\ref{alg:equiPRG:atomspaceoracle} is that they are essentially mimicking
    Algorithm~\ref{alg:equiPRG:atom} and only keeping in memory the variables that are strictly necessary to determine which
    part each vertex belongs to. Furthermore, they also use an auxiliary algorithm,
    Algorithm~\ref{alg:equiPRG:atomspaceauxiliary} that encodes an oracle of the set of vertices that have not yet been covered
    by the growing partition. The proof of correctness is long but straightforward, so we omit it.

    \begin{algorithm}[htbp]
      \caption{Computation algorithm of query-into-oracle model with space- and time-complexities
        \begin{gather*}
          \begin{multlined}[t]
            O\bigggl(
            \epsilon^{-\ell-1}\cdot
            \Biggl(\log\left(\frac{n}{\epsilon^\ell}\right)
            \\
            +
            \frac{\ell\cdot d}{c_{\atom}^2\cdot\delta^2}
            \cdot\bigl(\len(c_{\Root}) + \len(c_{\ZZ}) + \ell\cdot\len(\epsilon) + \len(c_{\samp}) + \len(c_{\adj})\bigr)
            \Biggr)
            \\
            +
            \len(c_{\atom}) + \len(\delta)
            + \ell\cdot\len(c_{\adj})
            \bigggr).
          \end{multlined}
          \\
          \begin{multlined}[t]
            2^{O(\ell\cdot d\cdot(\len(c_{\Root}) + \len(c_{\ZZ}) + \ell\cdot\len(\epsilon) + \len(c_{\samp}) + \len(c_{\adj})))/(c_{\atom}^2\cdot\delta^2)}
            \cdot n^{O(\epsilon^{-\ell-1}/(c_{\atom}^2\cdot\delta^2))}
            \\
            +
            \log(d+1)\cdot\len(c_{\atom})\cdot\len(\delta)
            \\
            +
            \ell^4\cdot\log(\ell+1)\cdot\len(c_{\ZZ})^2\cdot\len(c_{\samp})^2\cdot\len(c_{\Root})^2\cdot\len(c_{\Sum})\cdot\len(\epsilon)^2
            \cdot\bigl(\log(n+1)\bigr)^2.     
          \end{multlined}
        \end{gather*}
        respectively, that returns
        $(m,\widetilde{m},h,r,(\sigma^k,x^k,s^k,\widetilde{s}^k,r_k,P_k,\widetilde{p}_k,Z^k)_{k=1}^K)$ such that when given to
        Algorithm~\ref{alg:equiPRG:atomspaceoracle} produces an oracle for a $(\delta,\epsilon)$-atomic partition of $G$ into
        $K\df K_{\equiPRG}(\ell,c_{\ZZ},c_{\samp},c_{\adj},c_{\Root},c_{\Sum},\epsilon)$ parts.}
      \label{alg:equiPRG:atomspace}
      \DontPrintSemicolon
      \KwIn{Numbers $n,d,\ell\in\NN_+$ with $d\leq\ell\leq\log_2(n)$,
        $c_{\atom},c_{\ZZ},c_{\adj},c_{\Root},c_{\Sum},\delta\in(0,1)\cap\QQ$ with $(1+c_{\atom})\cdot\delta < 1/2^{\ell+1}$,
        $c_{\samp}\in(0,1/2)\cap\QQ$, and
        $\epsilon\in(0,\epsilon_{\equiPRG}(\ell,c_{\ZZ},c_{\samp},c_{\adj},c_{\Root},c_{\Sum}))$, and a query-oracle $\cO_G$ for
        a graph $G$ with $\lvert G\rvert=n$, $\Lit(G)\leq\ell$ and $\VC(G)\leq d$. We further assume that $n\geq
        n_{\equiPRG}^{\atom}(\ell,d,c_{\ZZ},c_{\samp},c_{\adj},c_{\Root},c_{\Sum},\epsilon,\delta)$.}
      \KwOut{A tuple $(m,\widetilde{m},h,r,(\sigma^k,x^k,s^k,\widetilde{s}^k,r_k,P_k,\widetilde{p}_k,Z^k)_{k=1}^K)$, with $K=
        K_{\equiPRG}(\ell,c_{\ZZ},c_{\samp},c_{\adj},c_{\Root},c_{\Sum},\epsilon)$, where $m,r,h\in\NN_+$, $\widetilde{m}\in\NN$
        and for each $k\in[K]$, we have $\sigma^k\in\{0,1\}^{<\ell+1}$, $x^k\in (V(G)^m)^{\lvert\sigma\rvert}$,
        $s^k,\widetilde{s}^k,r_k\in[n]$, $P_k\in\FF_2[X]$ is an irreducible polynomial, $\widetilde{p}_k\in\{0,\ldots,2^{t_k}\}$
        and $Z^k\in[2^{t_k}]^h$, where $t_k\df\deg(P_k)$ such that when given to Algorithm~\ref{alg:equiPRG:atomspaceoracle}
        produces an oracle for a $(\delta,\epsilon)$-atomic partition of $G$ into $K$ parts.}
      $\widetilde{\delta}\assign(1+c_{\atom})\cdot\delta$\;
      Let $m\assign\ceil{C\cdot d/(c_{\atom}^2\cdot\delta^2)}$, where $C$ is the absolute constant of
      Lemma~\ref{lem:atomconst}.\;\label{alg:equiPRG:atomspace:m}
      $\widetilde{m}\assign\floor{\widetilde{\delta}\cdot m}$\;
      $\zeta\assign (1-c_{\ZZ})\cdot\epsilon$\;
      $\widetilde{c}_{\UC}\assign c_{\UC}\assign 0$\;
      $\widetilde{c}_{\adj}\assign (1-c_{\adj})^{-\ell} - 1$\;
      $\widetilde{c}_{\samp}\assign(1-c_{\samp})^{-\ell} - 1$\;
      $\widetilde{c}_{\ZZ}\assign(1-c_{\ZZ})^{-\ell-1} - 1$\;
      $\widetilde{c}_{\Root}\assign(1-c_{\Root})^{-\ell-1} - 1$\;
      \setcounter{algosplit}{\theAlgoLine}
    \end{algorithm}

    \begin{algorithm}[htbp]
      \caption*{(continued).}
      \DontPrintSemicolon
      \setcounter{AlgoLine}{\thealgosplit}
      $K\assign\ceil{(1+\widetilde{c}_{\UC})\cdot(1+\widetilde{c}_{\ZZ})\cdot(1+\widetilde{c}_{\samp})\cdot(1+\widetilde{c}_{\adj})
        \cdot(1+\widetilde{c}_{\Root})\cdot(1+c_{\Sum})\cdot\epsilon^{-\ell-1}}$\;
      $r\assign\floor{n/K}$\;
      $h\assign h_{\equiPRG}^{\atom}\df 2\cdot m + 2$\;
      $\widetilde{n}\assign n$\;
      \uFor{$i\in[K]$}{%
        \label{alg:equiPRG:atomspace:For}
        Find a rational $\xi\in\QQ$ satisfying:
        \begin{align*}
          0 & < \xi\leq(1-c_{\Root})\cdot\zeta,
          &
          0 & \leq -1 + \zeta - \xi + \gamma^\ell\cdot\frac{\widetilde{n}}{r} \leq c_{\Root}\cdot\zeta,
        \end{align*}
        where
        \begin{align*}
          \gamma & \df \widetilde{\epsilon} \df (1-c_{\samp})\cdot\widetilde{\xi},
          &
          \widetilde{\xi} & \df \xi,
        \end{align*}
        and $\gamma$ is a dyadic.\;\label{alg:equiPRG:atomspace:xii}
        $r_i\assign\ceil{\gamma^\ell\cdot\widetilde{n}}$\;\label{alg:equiPRG:atomspace:ri}
        $s_0\df\widetilde{n}$\;
        \lFor{$j\in[\ell+1]$}{\label{alg:equiPRG:atomspace:sj}$s_j\assign\floor{\gamma_i\cdot s_{j-1}}+1$}
        Let $\cO$ be Algorithm~\ref{alg:equiPRG:atomspaceauxiliary} with parameters:
        \begin{multline*}
          \Bigl(n,\cO_G,
          \bigl(m,\widetilde{m},h,r,(\sigma^k,x^k,s^k,\widetilde{s}^k,r_k,P_k,\widetilde{p}_k,Z^k)_{k=1}^K\bigr),
          u\Bigr)
          \\
          =
          \Bigl(n,\cO_G,
          \bigl(m,\widetilde{m},h,r,(\sigma^k,x^k,s^k,\widetilde{s}^k,r_k,P_k,\widetilde{p}_k,Z^k)_{k=1}^{i-1}\bigr),
          \place\Bigr).
        \end{multline*}
        \;
        Run Algorithm~\ref{alg:ext:atomspace} with parameters:
        \begin{gather*}
          (n,m,\ell,s_0,\ldots,s_\ell,s_{\ell+1},\widetilde{m},\cO_G,\cO_U)
          \df
          (n,m,\ell,s_0,\ldots,s_\ell,s_{\ell+1},\widetilde{m},\cO_G,\cO).
        \end{gather*}
        \;\label{alg:equiPRG:atomspace:call}
        Let $(\sigma,(x_j)_{j=0}^{\lvert\sigma\rvert-1})$ be the tuple returned by the algorithm.\tcp*[r]{Encodes set $W$}
        $\sigma^i\assign\sigma$\;
        $x^i\assign (x_j)_{j=0}^{\lvert\sigma\rvert-1}$\;
        $s^i\assign s_{\lvert\sigma^i\rvert}$\tcp*[r]{$s^i\df\lvert W\rvert$}
        Let $p\df\gamma_i^{\ell-\lvert\sigma\rvert}$ and write it as $\widetilde{p}_i/2^t$ for $\widetilde{p}_i,t\in\NN_+$ with
        $t\geq\ceil{\log_2(n)}$ smallest possible.\;
        \setcounter{algosplit}{\theAlgoLine}
      }
    \end{algorithm}

    \begin{algorithm}[htbp]
      \caption*{(continued).}
      \DontPrintSemicolon
      \setcounter{AlgoLine}{\thealgosplit}
      \let\oldnl\nl
      \let\nl\relax
      \vspace{-\baselineskip}\uInvisibleBegin{
        \global\let\nl\oldnl 
        Run Algorithm~\ref{alg:encodePRG} with parameter $t\df t$ and let $P_i$ be the returned polynomial.\;
        \uFor{$Z^i\in[2^t]^h$}{%
          \tcp*[l]{Use Algorithm~\ref{alg:runPRG} to implicitly generate pseudorandom $X$ from $Z^i$}
          $\widetilde{s}^i\assign 0$\tcp*[r]{Size of $W(X)\df\{u\in W\mid X_u=1\}$}
          \For(\tcp*[f]{Compute $\widetilde{s}^i\df\lvert W(X)\rvert$}){$u\in[n]$}{%
            Run Algorithm~\ref{alg:ext:atomspaceoracle} with parameters:
            \begin{gather*}
              (n,m,\widetilde{m},\cO_G,\cO_U,(\sigma,(x_j)_{j=0}^{\lvert\sigma\rvert-1}),u)
              \df
              (n,m,\widetilde{m},\cO_G,\cO,(\sigma,(x_j)_{j=0}^{\lvert\sigma\rvert-1}),\place)
            \end{gather*}
            \;
            \lIf(\tcp*[f]{$u\notin W$}){The algorithm returned $0$}{\Continue}
            Run Algorithm~\ref{alg:runPRG} with parameters:
            \begin{gather*}
              (P,\widetilde{p},h,n,Z,u)\df(P_i,\widetilde{p}_i,h,n,Z^i,u).
            \end{gather*}
            \;
            \lIf(\tcp*[f]{$u\in W(X)$}){The algorithm returned $1$}{$\widetilde{s}^i\assign\widetilde{s}^i+1$}
          }
          \lIf{$\lvert\widetilde{s}^i - p\cdot s^i\rvert > c_{\adj}\cdot c_{\samp}\cdot\widetilde{\xi}\cdot p\cdot s^i$}{%
            \Continue
          }
          \uFor(\tcp*[f]{Check if $W(X)$ is a $(\widetilde{\delta},\widetilde{\xi},m)$-atom of $G$}){$z\in[n]^m$}{%
            \label{alg:equiPRG:atomspace:xFor}
            $k^0\assign 0$\tcp*[r]{Will count $\lvert\{u\in W(X) \mid t_G^{\widetilde{\delta}}(u,z) = 0\}\rvert$}
            $k^1\assign 0$\tcp*[r]{Will count $\lvert\{u\in W(X) \mid t_G^{\widetilde{\delta}}(u,z) = 1\}\rvert$}
            \uFor(\tcp*[f]{Compute $k^0$ and $k^1$}){$u\in [n]$}{%
              Run Algorithm~\ref{alg:ext:atomspaceoracle} with parameters:
              \begin{gather*}
                (n,m,\widetilde{m},\cO_G,\cO_U,(\sigma,(x_j)_{j=0}^{\lvert\sigma\rvert-1}),u)
                \df
                (n,m,\widetilde{m},\cO_G,\cO,(\sigma,(x_j)_{j=0}^{\lvert\sigma\rvert-1}),\place)
              \end{gather*}
              \;
              \lIf(\tcp*[f]{$u\notin W$}){The algorithm returned $0$}{\Continue}
              Run Algorithm~\ref{alg:runPRG} with parameters:
              \begin{gather*}
                (P,\widetilde{p},h,n,Z,u)\df(P_i,\widetilde{p}_i,h,n,Z^i,u).
              \end{gather*}
              \;
              \lIf(\tcp*[f]{$u\notin W(X)$}){The algorithm returned $0$}{\Continue}
              %
              \setcounter{algosplit}{\theAlgoLine}
            }
          }
        }
      }
    \end{algorithm}

    \begin{algorithm}[htbp]
      \caption*{(continued).}
      \DontPrintSemicolon
      \setcounter{AlgoLine}{\thealgosplit}
      \let\oldnl\nl
      \let\nl\relax
      \vspace{-\baselineskip}\InvisibleBegin{
        \vspace{-\baselineskip}\InvisibleBegin{
          \vspace{-\baselineskip}\InvisibleBegin{
            \vspace{-\baselineskip}\InvisibleBegin{
              \global\let\nl\oldnl 
              $a\assign 0$\tcp*[r]{Will count $\lvert\{j\in[m] \mid z_j\in N_G(u)\}\rvert$}
              \For{$j\in[m]$}{%
                \lIf{$\cO_G(u,z_j)=1$}{$a\assign a+1$}
              }
              \lIf(\tcp*[f]{$t_G^{\widetilde{\delta}}(u,z)=1$}){%
                $a\geq m - \widetilde{m}$
              }{%
                \label{alg:equiPRG:atomspace:k1update}
                $k^1\assign k^1+1$
              }
              \lElseIf(\tcp*[f]{$t_G^{\widetilde{\delta}}(u,z)=0$}){%
                $a\leq\widetilde{m}$
              }{%
                \label{alg:equiPRG:atomspace:k0update}
                $k^0\assign k^0+1$
              }
            }
            \If(\tcp*[f]{$W(X)$ is $(\widetilde{\delta},\widetilde{\xi}_i)$-split by $z$}){%
              $k^0 > \widetilde{\xi}\cdot\widetilde{s}^i$ and $k^1 > \widetilde{\xi}\cdot\widetilde{s}^i$}{%
              \Break\tcp*[r]{$W(X)$ is not a $(\widetilde{\delta},\widetilde{\xi}_i,m)$-atom}
            }
          }
          \Else(\tcp*[f]{Executes only if ``for'' loop of line~\ref{alg:equiPRG:atomspace:xFor} didn't break}){%
            \Break\tcp*[r]{$W(X)$ is a $(\widetilde{\delta},\widetilde{\xi}_i,m)$-atom}
          }
        }
        $\widetilde{n}\assign \widetilde{n} - r_k$\;
      }
      \Return{$(m,\widetilde{m},h,r,(\sigma^k,x^k,s^k,\widetilde{s}^k,r_k,P_k,\widetilde{p}_k,Z^k)_{k=1}^K)$}
    \end{algorithm}

    \begin{algorithm}[htbp]
      \caption{Oracle algorithm of query-into-oracle model with space- and time-complexities
        \begin{gather*}
          O\Bigl(K\cdot\bigl(\log(K+1) + \log(n+1) + \log(m+1) + h\cdot t_*\bigr)\Bigr),
          \\
          O\bigl((K^2\cdot m\cdot n^2 + K\cdot h\cdot t_*^2\cdot n)\cdot 2^K\cdot n^{2\cdot K - 1}\bigr),
        \end{gather*}
        respectively, where $t_*$ is such that $\deg(P_k)\leq t_*$ for every $k\in[K]$. When given
        $(m,\widetilde{m},h,r,(\sigma^k,x^k,s^k,\widetilde{s}^k,r_k,P_k,\widetilde{p}_k,Z^k)_{k=1}^K)$ computed via
        Algorithm~\ref{alg:equiPRG:atomspace}, produces an oracle for a $(\delta,\epsilon)$-atomic partition of $G$ into exactly
        $K\df K_{\equiPRG}(\ell,c_{\ZZ},c_{\samp},c_{\adj},c_{\Root},c_{\Sum},\epsilon)$ parts. If instead the parameters
        $(m,\widetilde{m},h,r,(\sigma^k,x^k,s^k,\widetilde{s}^k,r_k,P_k,\widetilde{p}_k,Z^k)_{k=1}^K)$ is computed via
        Algorithm~\ref{alg:equiPRG:goodspace}, then the oracle produced is for an $\epsilon$-good partition of $G$ with $K\df
        K_{\equiPRG}(\ell,c_{\ZZ},c_{\samp},c_{\adj},c_{\Root},c_{\Sum},\epsilon)$ parts.}
      \label{alg:equiPRG:atomspaceoracle}
      \DontPrintSemicolon
      \KwIn{A number $n\in\NN_+$, a query-oracle $\cO_G$ for a graph $G$ with $V(G)=[n]$, a tuple
        $(m,\widetilde{m},h,r,(\sigma^k,x^k,s^k,\widetilde{s}^k,r_k,P_k,\widetilde{p}_k,Z^k)_{k=1}^K)$ provided by
        Algorithm~\ref{alg:equiPRG:atomspace} or Algorithm~\ref{alg:equiPRG:goodspace} ran with the same parameters (and a
        choice of the other parameters required by it and satisfying the required hypotheses), a vertex $u\in V(G)$ and an index
        $k\in[K]$.}
      \KwOut{A bit $b\in\{0,1\}$ such that if $Q_k$ is the set of $u\in V(G)$ for which the algorithm returns $1$ for $k\in[K]$,
        then $\cP = \{Q_k \mid k\in[K]\}$ is a partition of $G$. If the input came from Algorithm~\ref{alg:equiPRG:atomspace},
        then $\cP$ is guaranteed to be $(\delta,\epsilon)$-atomic, where $\delta$ and $\epsilon$ are the parameters used by
        Algorithm~\ref{alg:equiPRG:atomspace}. If the input came from Algorithm~\ref{alg:equiPRG:goodspace}, then $\cP$ is
        guaranteed to be $\epsilon$-good, where $\epsilon$ is the parameter used by Algorithm~\ref{alg:equiPRG:goodspace}.}
      Run Algorithm~\ref{alg:equiPRG:atomspaceauxiliary} with parameters:
      \begin{multline*}
        \Bigl(n,\cO_G,
        \bigl(m,\widetilde{m},h,r,(\sigma^k,x^k,s^k,\widetilde{s}^k,r_k,P_k,\widetilde{p}_k,Z^k)_{k=1}^K\bigr),
        u\Bigr)
        \\
        =
        \Bigl(n,\cO_G,
        \bigl(m,\widetilde{m},h,r,(\sigma^{\widetilde{k}},x^{\widetilde{k}},s^{\widetilde{k}},\widetilde{s}^{\widetilde{k}},
        r_{\widetilde{k}},P_{\widetilde{k}},\widetilde{p}_{\widetilde{k}},Z^{\widetilde{k}})_{{\widetilde{k}}=1}^{k-1}\bigr),
        u\Bigr).
      \end{multline*}
      \;
      \lIf(\tcc*[f]{$u$ is in the main portion of the first $k-1$ parts}){The algorithm returned $0$}{\Return{$0$}}
      Run Algorithm~\ref{alg:equiPRG:atomspaceauxiliary} with parameters:
      \begin{multline*}
        \Bigl(n,\cO_G,
        \bigl(m,\widetilde{m},h,r,(\sigma^k,x^k,s^k,\widetilde{s}^k,r_k,P_k,\widetilde{p}_k,Z^k)_{k=1}^K\bigr),
        u\Bigr)
        \\
        =
        \Bigl(n,\cO_G,
        \bigl(m,\widetilde{m},h,r,(\sigma^{\widetilde{k}},x^{\widetilde{k}},s^{\widetilde{k}},\widetilde{s}^{\widetilde{k}},
        r_{\widetilde{k}},P_{\widetilde{k}},\widetilde{p}_{\widetilde{k}},Z^{\widetilde{k}})_{{\widetilde{k}}=1}^k\bigr),
        u\Bigr).
      \end{multline*}
      \;
      \lIf(\tcc*[f]{$u$ is in the main portion of the $k$th part}){The algorithm returned $0$}{\Return{$1$}}
      \setcounter{algosplit}{\theAlgoLine}
    \end{algorithm}

    \begin{algorithm}[htbp]
      \caption*{(continued).}
      \DontPrintSemicolon
      \setcounter{AlgoLine}{\thealgosplit}
      \let\oldnl\nl
      \let\nl\relax
      \global\let\nl\oldnl 
      \tcc*[l]{At this point, we know $u$ is not in the main portion of the $k$th part, but could be an absorbed vertex of
        remainder.}
      Run Algorithm~\ref{alg:equiPRG:atomspaceauxiliary} with parameters:
      \begin{multline*}
        \Bigl(n,\cO_G,
        \bigl(m,\widetilde{m},h,r,(\sigma^k,x^k,s^k,\widetilde{s}^k,r_k,P_k,\widetilde{p}_k,Z^k)_{k=1}^K\bigr),
        u\Bigr)
        \\
        =
        \Bigl(n,\cO_G,
        \bigl(m,\widetilde{m},h,r,(\sigma^{\widetilde{k}},x^{\widetilde{k}},s^{\widetilde{k}},\widetilde{s}^{\widetilde{k}},
        r_{\widetilde{k}},P_{\widetilde{k}},\widetilde{p}_{\widetilde{k}},Z^{\widetilde{k}})_{{\widetilde{k}}=1}^K\bigr),
        u\Bigr).
      \end{multline*}
      \;
      \lIf(\tcp*[f]{$u$ is not in the remainder}){The algorithm returned $0$}{\Return{$0$}}
      \tcp*[l]{Now $u$ is in the remainder; need to find which part absorbed it}
      $i\assign 1$\tcp*[r]{Index of part absorbing current vertex}
      $a\assign 0$\tcp*[r]{Number of vertices of remainder absorbed by current part so far}
      \For(\tcp*[f]{Vertex that potentially is in the remainder}){$w\in[u-1]$}{%
        Run Algorithm~\ref{alg:equiPRG:atomspaceauxiliary} with parameters:
        \begin{multline*}
          \Bigl(n,\cO_G,
          \bigl(m,\widetilde{m},h,r,(\sigma^k,x^k,s^k,\widetilde{s}^k,r_k,P_k,\widetilde{p}_k,Z^k)_{k=1}^K\bigr),
          u\Bigr)
          \\
          =
          \Bigl(n,\cO_G,
          \bigl(m,\widetilde{m},h,r,(\sigma^{\widetilde{k}},x^{\widetilde{k}},s^{\widetilde{k}},\widetilde{s}^{\widetilde{k}},
          r_{\widetilde{k}},P_{\widetilde{k}},\widetilde{p}_{\widetilde{k}},Z^{\widetilde{k}})_{{\widetilde{k}}=1}^K\bigr),
          w\Bigr).
        \end{multline*}
        \;
        \If(\tcp*[f]{$w$ is in remainder and is absorbed by $i$th part}){The algorithm returned $1$}{%
          $a\assign a+1$\;
          \If(\tcp*[f]{$i$th part already absorbed all its remainder vertices}){$a\geq r - r_k$}{%
            $i\assign i+1$\;
            \lIf(\tcp*[f]{$u$ is absorbed by a part after the $k$th}){$i > k$}{\Return{$0$}}
            $a\assign 0$\;
          }
        }
      }
      \Return{$1$}\tcp*[f]{$u$ is absorbed by $k$th part}
    \end{algorithm}

    \begin{algorithm}[htbp]
      \caption{Auxiliary oracle algorithm such that when given the parameters
        $(m,\widetilde{m},h,r,(\sigma^k,x^k,s^k,\widetilde{s}^k,r_k,P_k,\widetilde{p}_k,Z^k)_{k=1}^K)$ provided within the
        execution of one of Algorithms~\ref{alg:equiPRG:atomspace} or~\ref{alg:equiPRG:atomspace}, produces an oracle for the
        set of vertices that have not yet been covered by the main portion of the parts of the partition. The space- and
        time-complexities are
        \begin{gather*}
          O\Bigl(K\cdot\bigl(\log(K+1) + \log(n+1) + \log(m+1) + h\cdot t_*\bigr)\Bigr).
          \\
          O\bigl(
          1 + (K^2\cdot m\cdot n^2 + K\cdot h\cdot t_*^2\cdot n)\cdot 2^K\cdot n^{2\cdot K - 2}
          \bigr)
        \end{gather*}
        respectively, where $t_*$ is such that $\deg(P_k)\leq t_*$ for every $k\in[K]$.}
      \label{alg:equiPRG:atomspaceauxiliary}
      \DontPrintSemicolon
      \KwIn{A number $n\in\NN_+$, a query-oracle $\cO_G$ for a graph $G$ with $V(G)=[n]$, a tuple
        $(m,\widetilde{m},h,r,(\sigma^k,x^k,s^k,\widetilde{s}^k,r_k,P_k,\widetilde{p}_k,Z^k)_{k=1}^K)$ provided by one of
        Algorithms~\ref{alg:equiPRG:atomspace} or~\ref{alg:equiPRG:atomspace} ran with the same parameters (and a choice of the
        other parameters required by it and satisfying the required hypotheses), and a vertex $u\in V(G)$.}
      \KwOut{A bit $b\in\{0,1\}$ such that if $W$ is the set of $u\in V(G)$ for which the algorithm returns $1$, then the main
        part of the sets defined by $(\sigma^k,x^k,s^k,\widetilde{s}^k,r_k,P_k,\widetilde{p}_k,Z^k)_{k=1}^K$ cover exactly
        $V(G)\setminus W$ (so $W$ is the remainder set).}
      \lIf{$K=0$}{\Return{$1$}}
      \uFor(\tcp*[f]{Check if $u$ is in $\widetilde{k}$th part}){$\widetilde{k}\in[K]$}{%
        \label{alg:equiPRG:atomspaceauxiliary:wkFor}
        Let $\cO$ be Algorithm~\ref{alg:equiPRG:atomspaceauxiliary} with parameters:
        \begin{multline*}
          \Bigl(n,\cO_G,
          \bigl(m,\widetilde{m},h,r,(\sigma^k,x^k,s^k,\widetilde{s}^k,r_k,P_k,\widetilde{p}_k,Z^k)_{k=1}^K\bigr),
          u\Bigr)
          \\
          =
          \Bigl(n,\cO_G,
          \bigl(m,\widetilde{m},h,r,(\sigma^k,x^k,s^k,\widetilde{s}^k,r_k,P_k,\widetilde{p}_k,Z^k)_{k=1}^{\widetilde{k}-1}\bigr),
          \place\Bigr).
        \end{multline*}
        \;
        Run Algorithm~\ref{alg:ext:atomspaceoracle} with parameters:
        \begin{gather*}
          \bigl(n,m,\widetilde{m},\cO_G,\cO_U,(\sigma,(x_j)_{j=0}^{\lvert\sigma\rvert-1}),u\bigr)
          \df
          \bigl(n,m,\widetilde{m},\cO_G,\cO,(\sigma^{\widetilde{k}},x^{\widetilde{k}}),u\bigr).
        \end{gather*}
        \;
        \uIf(\tcp*[f]{$u$ is in $W$ portion of $\widetilde{k}$th part}){The algorithm returned $1$}{%
          Run Algorithm~\ref{alg:runPRG} with parameters:
          \begin{gather*}
            (P,\widetilde{p},h,n,Z,u)\df(P_{\widetilde{k}},\widetilde{p}_{\widetilde{k}},h,n,Z^{\widetilde{k}},u).
          \end{gather*}
          \;
          \uIf(\tcp*[f]{$u$ is in $W(X)$ portion of $\widetilde{k}$th part}){The algorithm returned $1$}{%
            \If(\tcp*[f]{Main portion of $\widetilde{k}$th part contains $W(X)$}){%
              $r_{\widetilde{k}}\geq\widetilde{s}^{\widetilde{k}}$
            }{%
              \Return{$0$}\tcp*[r]{$u$ in main portion of $\widetilde{k}$th part}
            }
            \setcounter{algosplit}{\theAlgoLine}
          }
        }
      }
    \end{algorithm}

    \begin{algorithm}[htbp]
      \caption*{(continued).}
      \DontPrintSemicolon
      \setcounter{AlgoLine}{\thealgosplit}
      \let\oldnl\nl
      \let\nl\relax
      \vspace{-\baselineskip}\uInvisibleBegin{
        \vspace{-\baselineskip}\InvisibleBegin{
          \vspace{-\baselineskip}\InvisibleBegin{
            \global\let\nl\oldnl 
            \tcc*[l]{At this point, $u$ is in $W(X)$ and main portion of $\widetilde{k}$th part is first $r_{\widetilde{k}}$
              vertices of $W(X)$; let us check if $u$ is among them}
            $a\assign 0$\tcp*[r]{Vertices of $W(X)$ seen so far}
            \uFor{$w\in[u-1]$}{%
              \label{alg:equiPRG:atomspaceauxiliary:wFor}
              Run Algorithm~\ref{alg:ext:atomspaceoracle} with parameters:
              \begin{gather*}
                \bigl(n,m,\widetilde{m},\cO_G,\cO_U,(\sigma,(x_j)_{j=0}^{\lvert\sigma\rvert-1}),u\bigr)
                \df
                \bigl(n,m,\widetilde{m},\cO_G,\cO,(\sigma^{\widetilde{k}},x^{\widetilde{k}}),w\bigr).
              \end{gather*}
              \;
              \lIf(\tcp*[f]{$w$ is not in $W$ portion of $\widetilde{k}$th part}){The algorithm returned $0$}{\Continue}
              Run Algorithm~\ref{alg:runPRG} with parameters:
              \begin{gather*}
                (P,\widetilde{p},h,n,Z,u)\df(P_{\widetilde{k}},\widetilde{p}_{\widetilde{k}},h,n,Z^{\widetilde{k}},w).
              \end{gather*}
              \;
              \lIf(\tcp*[f]{$w$ is not in $W(X)$ portion of $\widetilde{k}$th part}){The algorithm returned $0$}{\Continue}
              $a\assign a+1$\tcp*[r]{$w$ in $W(X)$ portion of $\widetilde{k}$th part, count it}
              \lIf(\tcp*[f]{All vertices of $\widetilde{k}$th part already seen}){$a\geq r_{\widetilde{k}}$}{\Break}
            }
            \Else(\tcp*[f]{Executes only if ``for'' loop of line~\ref{alg:equiPRG:atomspaceauxiliary:wFor} didn't break}){%
              \Return{$0$}\tcp*[r]{$u$ in main portion of $\widetilde{k}$th part (first $r_{\widetilde{k}}$ vertices)}              
            }
            \Continue\tcp*[r]{$w$ not in $\widetilde{k}$th part, continues loop of line~\ref{alg:equiPRG:atomspaceauxiliary:wkFor}}
          }
        }
        \tcc*[l]{At this point, $u$ is not in $W(X)$ portion of $\widetilde{k}$th part, but if
          $r_{\widetilde{k}} > \widetilde{s}^{\widetilde{k}}$, then $u$ could be an added vertex}
        \If(\tcp*[f]{There are no added vertices}){$r_{\widetilde{k}}\leq\widetilde{s}^{\widetilde{k}}$}{%
          \Continue\tcp*[r]{$w$ not in $\widetilde{k}$th part, continues loop of line~\ref{alg:equiPRG:atomspaceauxiliary:wkFor}}
        }
        $a\assign 0$\tcp*[r]{Vertices added to $\widetilde{k}$th part seen so far}
        \uFor{$w\in[u-1]$}{%
          \label{alg:equiPRG:atomspaceauxiliary:secondwFor}
          Run Algorithm~\ref{alg:ext:atomspaceoracle} with parameters:
          \begin{gather*}
            \bigl(n,m,\widetilde{m},\cO_G,\cO_U,(\sigma,(x_j)_{j=0}^{\lvert\sigma\rvert-1}),u\bigr)
            \df
            \bigl(n,m,\widetilde{m},\cO_G,\cO,(\sigma^{\widetilde{k}},x^{\widetilde{k}}),w\bigr).
          \end{gather*}
          \;
          \setcounter{algosplit}{\theAlgoLine}
        }
      }
    \end{algorithm}

    \begin{algorithm}[htbp]
      \caption*{(continued).}
      \DontPrintSemicolon
      \setcounter{AlgoLine}{\thealgosplit}
      \let\oldnl\nl
      \let\nl\relax
      \vspace{-\baselineskip}\InvisibleBegin{
        \vspace{-\baselineskip}\InvisibleBegin{
          \global\let\nl\oldnl 
          \If(\tcp*[f]{$w$ is in $W$ portion of $\widetilde{k}$th part}){The algorithm returned $1$}{%
            Run Algorithm~\ref{alg:runPRG} with parameters:
            \begin{gather*}
              (P,\widetilde{p},h,n,Z,u)\df(P_{\widetilde{k}},\widetilde{p}_{\widetilde{k}},h,n,Z^{\widetilde{k}},w).
            \end{gather*}
            \;
            \If(\tcp*[f]{$w$ is in $W(X)$ portion of $\widetilde{k}$th part}){The algorithm returned $1$}{%
              \Continue
            }
          }
          $a\assign a+1$\tcp*[r]{$w$ is not in $W(X)$ portion of $\widetilde{k}$th part, count it}
          \lIf(\tcp*[f]{All vertices added to $\widetilde{k}$th part already seen}){%
            $a\geq r_{\widetilde{k}} - \widetilde{s}^{\widetilde{k}}$
          }{\Break}
        }
        \Else(\tcp*[f]{Executes only if ``for'' loop of line~\ref{alg:equiPRG:atomspaceauxiliary:secondwFor} didn't break}){%
          \Return{$0$}\tcp*[r]{$u$ is an added vertex of $\widetilde{k}$th part}
        }
      }
      \Return{$1$}
    \end{algorithm}

    We now analyze the space- and time-complexities of the algorithms. We start with
    Algorithm~\ref{alg:equiPRG:atomspaceauxiliary}. Since the algorithm passes itself as an oracle parameter to
    Algorithm~\ref{alg:ext:atomspaceoracle}, our analysis is recursive:
    \begin{claim}\label{clm:equiPRG:atomspaceauxiliary}
      For an input $(n,\cO_G,(m,\widetilde{m},h,r,(\sigma^k,x^k,s^k,\widetilde{s}^k,r_k,P_k,\widetilde{p}_k,Z^k)_{k=1}^K),u)$
      and for $i\in\{0,\ldots,K\}$, let $\cO_i$ be Algorithm~\ref{alg:nonequi:atomspaceauxiliary} when ran with the parameters
      \begin{gather*}
      \bigl(n,\cO_G,(m,\widetilde{m},h,r,(\sigma^k,x^k,s^k,\widetilde{s}^k,r_k,P_k,\widetilde{p}_k,Z^k)_{k=1}^i),u\bigr)
      \end{gather*}
      and let $S(i)$ and $T(i)$ its space- and time-complexities. Suppose further that $\deg(P_k)\leq t_*$ for every $k\in[K]$.
      Then we have
      \begin{align*}
        S(K)
        & \leq
        O\Bigl(K\cdot\bigl(\log(K+1) + \log(n+1) + \log(m+1) + h\cdot t_*\bigr)\Bigr).
        \\
        T(i)
        & \leq
        O\bigl(
        1 + (i^2\cdot m\cdot n^2 + i\cdot h\cdot t_*^2\cdot n)\cdot 2^i\cdot n^{2\cdot i - 2}
        \bigr)
      \end{align*}
      for every $i\in\{0,\ldots,K\}$.
    \end{claim}

    \begin{proofof}{Claim~\ref{clm:equiPRG:atomspaceauxiliary}}
      Let $\widetilde{S}(i)$, $\widetilde{T}(i)$ and $\widetilde{U}(i)$ be the space-, time- and $U$-oracle-complexities of
      Algorithm~\ref{alg:ext:atomspaceoracle} when ran with $(n,m,\widetilde{m},\cO_G,\cO_i,(\sigma^i,x^i),u)$. By
      Proposition~\ref{prop:ext}\ref{prop:ext:atomspace}, we have
      \begin{align*}
        \widetilde{S}(i) & = O\bigl(\log(\lvert\sigma\rvert+1) + \log(m+1) + \log(n+1)\bigr)
        &
        \widetilde{T}(i) & = O\bigl((i+1)\cdot m\cdot n\bigr),
        &
        \widetilde{U}(i) & \leq n.
      \end{align*}

      We note that in the iteration of the loop of line~\ref{alg:equiPRG:atomspaceauxiliary:wkFor} corresponding to
      $\widetilde{k}\in[K]$, Algorithm~\ref{alg:ext:atomspaceoracle} gets called at most $2\cdot n$ times with input
      $(n,m,\widetilde{m},\cO_G,\cO_i,(\sigma^{\widetilde{k}},x^{\widetilde{k}}),u)$. Furthermore, Algorithm~\ref{alg:runPRG}
      gets called with input $(P_{\widetilde{k}},\widetilde{p}_k,h,n,Z^{\widetilde{k}},u)$ at most $2\cdot n$ times. This means
      that we get the following recursive formulas for $i\in[K]$:
      \begin{align*}
        S(i)
        & =
        \log_2(i+1) + \log_2(n+1)
        + \max_{j\in\{0,\ldots,i-1\}} \Bigl(\widetilde{S}(j) + S(j) +  O\bigl(h\cdot\deg(P_j)\bigr)\Bigr)
        \\
        & \leq
        \log_2(i+1) + \log_2(n+1) + O\bigl(\log(i) + \log(m+1) + \log(n+1) + h\cdot t_*\bigr)
        + \max_{j\in\{0,\ldots,i-1\}} S(j)
        \\
        & \leq
        O\bigl(\log(i+1) + \log(n+1) + \log(m+1) + h\cdot t_*\bigr) + \max_{j\in\{0,\ldots,i-1\}} S(j),
        \\
        T(i)
        & \leq
        \sum_{j=0}^{i-1} 2\cdot n
        \cdot\Bigl(\widetilde{T}(j) + \widetilde{U}(j)\cdot T(j) + O\bigl(h\cdot\deg(P_j)^2\bigr)\Bigr)
        \\
        & \leq
        \sum_{j=0}^{i-1} 2\cdot n
        \cdot\Bigl(O\bigl((j+1)\cdot m\cdot n\bigr) + n\cdot T(j) + O\bigl(h\cdot t_*^2\bigr)\Bigr)
        \\
        & \leq
        O\left(\binom{i+1}{2}\cdot m\cdot n^2 + j\cdot h\cdot t_*^2\cdot n\right) + 2\cdot n^2\cdot\sum_{j=0}^{i-1} T(j).
      \end{align*}

      Therefore, we conclude that
      \begin{gather*}
        S(K)
        \leq
        O\Bigl(K\cdot\bigl(\log(K+1) + \log(n+1) + \log(m+1) + h\cdot t_*\bigr)\Bigr).
      \end{gather*}

      For $T$, it is straightforward to prove by induction in $i$ that
      \begin{gather*}
        T(i)
        \leq
        O\left(\binom{i+1}{2}\cdot m\cdot n^2\right)
        +
        \sum_{j=1}^{i-1}\left(\binom{j+1}{2}\cdot m\cdot n^2 + j\cdot h\cdot t_*^2\cdot n\right)
        \cdot(2\cdot n^2)\cdot(2\cdot n^2 + 1)^{i-j-1},
      \end{gather*}
      from which we get
      \begin{gather*}
        T(i)
        \leq
        O\bigl(
        1 + (i^2\cdot m\cdot n^2 + i\cdot h\cdot t_*^2\cdot n)\cdot 2^i\cdot n^{2\cdot i - 2}
        \bigr)
      \end{gather*}
      for every $i\in\{0,\ldots,K\}$.
    \end{proofof}

    We now proceed to the analysis of oracle Algorithm~\ref{alg:equiPRG:atomspaceoracle}. Note that it calls the auxiliary
    algorithm, Algorithm~\ref{alg:equiPRG:atomspaceauxiliary}, at most $u+2\leq n+2$ times. It is also clear that it needs
    internal memory to store the variables $i\in[K]$, $a\in\{0,\ldots,n\}$ and $w\in[n]$, which takes space
    \begin{gather*}
      O\bigl(\log(K+1) + \log(n+1)\bigr).
    \end{gather*}

    Taking into account our earlier analysis of the auxiliary Algorithm~\ref{alg:equiPRG:atomspaceauxiliary} in
    Claim~\ref{clm:equiPRG:atomspaceauxiliary}, we conclude that the space- and time-complexities of
    Algorithm~\ref{alg:nonequi:atomspaceoracle} are at most
    \begin{gather*}
      O\Bigl(K\cdot\bigl(\log(K+1) + \log(n+1) + \log(m+1) + h\cdot t_*\bigr)\Bigr),
      \\
      O\bigl((K^2\cdot m\cdot n^2 + K\cdot h\cdot t_*^2\cdot n)\cdot 2^K\cdot n^{2\cdot K - 1}\bigr),
    \end{gather*}
    respectively.

    Taking into account the fact that this will be executed with parameters from the computation
    Algorithm~\ref{alg:equiPRG:atomspace}, we know that
    \begin{align*}
      K
      & \leq
      O(\epsilon^{-\ell-1})
      \\
      m
      & \leq
      O\left(\frac{d}{c_{\atom}^2\cdot\delta^2}\right)
      \leq
      O\left(\frac{\ell}{c_{\atom}^2\cdot\delta^2}\right),
      \\
      h
      & \df
      h_{\equiPRG}^{\atom}
      =
      2\cdot m + 2
      \leq
      O\left(\frac{\ell}{c_{\atom}^2\cdot\delta^2}\right),
      \\
      t_*
      & \df
      \max\bigl(\deg(P_i)\bigr)
      \leq
      \log_2(n) + O\Bigl(\ell\cdot\bigl(\len(c_{\Root}) + \len(c_{\ZZ}) + \ell\cdot\len(\epsilon)
      + \len(c_{\samp}) + \len(c_{\adj})\bigr)\Bigr).
    \end{align*}
    (the bound on the degrees of the $P_i$ comes since in Algorithm~\ref{alg:equiPRG:atomspace}, these are calculated exactly as
    in Algorithm~\ref{alg:equiPRG:atom}) so the space- and time-complexities become
    \begin{gather*}
      \begin{multlined}[t]
      O\bigggl(
      \epsilon^{-\ell-1}\cdot
      \left(
      \log\left(\frac{n}{c_{\atom}\cdot\delta\cdot\epsilon^\ell}\right)\right)
      \\
      +
      \frac{\ell\cdot d}{c_{\atom}^2\cdot\delta^2}\cdot\bigl(\len(c_{\Root}) + \len(c_{\ZZ}) + \ell\cdot\len(\epsilon)
      + \len(c_{\samp}) + \len(c_{\adj})\bigr)
      \bigggr),
      \end{multlined}
      \\
      \begin{multlined}[t]
        (2\cdot n)^{O(\epsilon^{-\ell-1})}
        \cdot\bigggl(
        \frac{d}{c_{\atom}^2\cdot\delta^2\cdot\epsilon^{2\cdot\ell+2}}
        \\
        +
        \frac{\ell^2\cdot d}{c_{\atom}^2\cdot\delta^2}\cdot\bigl(\len(c_{\Root}) + \len(c_{\ZZ}) + \ell\cdot\len(\epsilon)
        + \len(c_{\samp}) + \len(c_{\adj})\bigr)^2
        \bigggr).
      \end{multlined}
    \end{gather*}

    \smallskip

    Finally, we analyze the computation Algorithm~\ref{alg:equiPRG:atomspace}. Let us first mention that the calculations of
    line~\ref{alg:equiPRG:atomspace:xii} are performed just like the calculations of line~\ref{alg:equiPRG:atom:xii} of
    Algorithm~\ref{alg:equiPRG:atom}. We know that
    Algorithm~\ref{alg:equiPRG:atomspace} needs space to store/compute the following variables:
    \begin{gather*}
      \begin{aligned}
        \widetilde{\delta}, & &
        m, & &
        \widetilde{m},
        \\
        \zeta, & &
        \widetilde{c}_{\adj}, & &
        \widetilde{c}_{\samp},
        \\
        \widetilde{c}_{\ZZ}, & &
        \widetilde{c}_{\Root}, & &
        K,
        \\
        r, & &
        h, & &
        \widetilde{n},
        \\
        & &
        i, & &
        \xi,\gamma
        \\
        & &
        s_0,\ldots,s_{\ell+1}, & &
        k^0,
        \\
        & &
        k^1, & &
        a,
      \end{aligned}
      \\
      (\sigma^i,x^i,s^i,\widetilde{s}^i,r_i,P_i,\widetilde{p}_i,Z^i),
    \end{gather*}
    These require the following space, respectively:
    \begin{gather*}
      \begin{aligned}
        O\bigl(\len(c_{\atom}) + \len(\delta)\bigr), & &
        O\bigl(\log(d) + \len(c_{\atom}) + \len(\delta)\bigr), & &
        O\bigl(\log(d) + \len(c_{\atom}) + \len(\delta)\bigr), & &
        \\
        O\bigl(\len(c_{\ZZ}) + \len(\epsilon)\bigr), & &
        O\bigl(\ell\cdot\len(c_{\adj})\bigr), & &
        O\bigl(\ell\cdot\len(c_{\samp})\bigr),
        \\
        O\bigl(\ell\cdot\len(c_{\ZZ})\bigr), & &
        O\bigl(\ell\cdot\len(c_{\Root})\bigr), & &
        O\bigl(\ell\cdot\log(\epsilon)\bigr),
        \\
        \log(n+1), & &
        O\bigl(\log(d) + \len(c_{\atom}) + \len(\delta)\bigr), & &
        \log(n+1),
      \end{aligned}
      \\      
      \begin{aligned}
        O\bigl(\ell\cdot\log(\epsilon)\bigr), & &
        O\bigl(\len(c_{\Root}) + \len(c_{\ZZ}) + \ell\cdot\len(\epsilon) + \len(c_{\samp}) + \len(c_{\adj})\bigr),
        \\
        \ell\cdot\log(n+1), & &
        \log(n+1),
        \\
        \log(n+1), & &
        O\bigl(\log(d) + \len(c_{\atom}) + \len(\delta)\bigr),
      \end{aligned}
      \\
      \begin{multlined}[t]
        O\biggl(\epsilon^{-\ell-1}\cdot\Bigl(\ell\cdot m\cdot\log(n+1)
        + \log(n)
        \\
        + \ell\cdot\bigl(\len(c_{\Root}) + \len(c_{\ZZ}) + \ell\cdot\len(\epsilon) + \len(c_{\samp}) + \len(c_{\adj})\bigr)
        \Bigr)
        \biggr).
      \end{multlined}
    \end{gather*}
    Together, these can be upper bounded by:
    \begin{multline*}
      O\bigggl(
      \epsilon^{-\ell-1}\cdot
      \left(
      \frac{d\cdot\ell\cdot\log(n+1)}{c_{\atom}^2\cdot\delta^2}
      +
      \ell\cdot\bigl(\len(c_{\Root}) + \len(c_{\ZZ}) + \ell\cdot\len(\epsilon) + \len(c_{\samp}) + \len(c_{\adj})\bigr)
      \right)
      \\
      +
      \len(c_{\atom}) + \len(\delta)
      + \ell\cdot\len(c_{\adj})
      \bigggr)
    \end{multline*}

    Since the algorithm calls Algorithm~\ref{alg:ext:atomspace} with Algorithm~\ref{alg:equiPRG:atomspaceauxiliary} as an
    oracle, we need to account for the space-complexity of these algorithms, which are
    \begin{gather*}
      O\bigl(\ell\cdot m\cdot\log(n+1)\bigr)
      \leq
      O\left(\frac{\ell\cdot d\cdot\log(n+1)}{c_{\atom}^2\cdot\delta^2}\right)
      \\
      \begin{aligned}
        \MoveEqLeft
        O\Bigl(K\cdot\bigl(\log(K+1) + \log(n+1) + \log(m+1) + h\cdot t_*\bigr)\Bigr)
        \\
        & \leq
        \begin{multlined}[t]
          O\bigggl(
          \epsilon^{-\ell-1}\cdot
          \Biggl(\log\left(\frac{n}{\epsilon^\ell}\right)
          \\
          +
          \frac{\ell\cdot d}{c_{\atom}^2\cdot\delta^2}
          \cdot\bigl(\len(c_{\Root}) + \len(c_{\ZZ}) + \ell\cdot\len(\epsilon) + \len(c_{\samp}) + \len(c_{\adj})\bigr)
          \Biggr)
          \bigggr).
        \end{multlined}
      \end{aligned}
    \end{gather*}

    Thus, the total space-complexity of Algorithm~\ref{alg:equiPRG:atomspace} is at most
    \begin{multline*}
      O\bigggl(
      \epsilon^{-\ell-1}\cdot
      \left(\log\left(\frac{n}{\epsilon^\ell}\right)
      +
      \frac{\ell\cdot d}{c_{\atom}^2\cdot\delta^2}
      \cdot\bigl(\len(c_{\Root}) + \len(c_{\ZZ}) + \ell\cdot\len(\epsilon) + \len(c_{\samp}) + \len(c_{\adj})\bigr)
      \right)
      \\
      +
      \len(c_{\atom}) + \len(\delta)
      + \ell\cdot\len(c_{\adj})
      \bigggr).
    \end{multline*}

    \smallskip

    Finally, we analyze the time-complexity of Algorithm~\ref{alg:equiPRG:atomspace}. The analysis of the computation of
    parameters is the same as the one done for Algorithm~\ref{alg:equiPRG:atom}. The computation of the initial parameters
    before the loop of line~\ref{alg:equiPRG:atomspace:For} takes time
    \begin{gather}\label{eq:equiPRG:atomspace:param:complexity}
      \begin{multlined}
        O\Bigl(
        \log(d+1)\cdot\len(c_{\atom})\cdot\len(\delta)
        \\
        +
        \ell^4\cdot\log(\ell+1)\cdot\len(c_{\ZZ})^2\cdot\len(c_{\samp})^2\cdot\len(c_{\Root})^2\cdot\len(c_{\Sum})\cdot\len(\epsilon)^2
        \cdot\bigl(\log(n+1)\bigr)^2
        \Bigr).
      \end{multlined}
    \end{gather}

    Let us now focus on the $i$iteration of the loop of line~\ref{alg:equiPRG:atomspace:For}. The computation of the parameters
    $\xi$, $\widetilde{\xi}$, $\widetilde{\epsilon}$ and $\gamma$ in line~\ref{alg:equiPRG:atomspace:xii} takes time
    \begin{gather}\label{alg:equiPRG:atomspace:xii:single}
      O\Bigl(
      \len(c_{\Root})\cdot\len(c_{\ZZ})\cdot\len(c_{\samp})\cdot\len(c_{\adj})
      \cdot\ell\cdot\len(\epsilon)\cdot\bigl(\log(n+1)\bigr)^2
      \Bigr).
    \end{gather}

    If $\cO_i$ is Algorithm~\ref{alg:equiPRG:atomspaceauxiliary} as in Claim~\ref{clm:equiPRG:atomspaceauxiliary}, then in the
    $i$th iteration of the loop of line~\ref{alg:equiPRG:atomspace:For}, there are at most:
    \begin{itemize}
    \item One call to Algorithm~\ref{alg:ext:atomspace} using oracle $\cO_{i-1}$.
    \item One call to Algorithm~\ref{alg:encodePRG}.
    \item $2^{t\cdot h}\cdot (n + n^{m+1})$ calls to Algorithm~\ref{alg:ext:atomspaceoracle}.
    \item $2^{t\cdot h}\cdot (n + n^{m+1})$ calls to Algorithm~\ref{alg:runPRG}.
    \end{itemize}
    It is also straightforward to check that the remainder of the computations of the $i$th iteration of the loop of
    line~\ref{alg:equiPRG:atomspace:For} are completely dominated by these.

    Using Proposition~\ref{prop:ext}\ref{prop:ext:atomspace}, Theorem~\ref{thm:PRG} and
    Claim~\ref{clm:equiPRG:atomspaceauxiliary}, the above give the following time-complexity bounds:
    \begin{gather*}
      O\bigl(
      2^\ell\cdot m\cdot n^m
      +
      (2^{\ell+1} - 1)\cdot n^{m+1}
      \cdot (1+(i-1)^2\cdot m\cdot n^2 + (i-1)\cdot h\cdot t_*^2\cdot n)\cdot 2^{i-1}\cdot n^{2\cdot i - 4}
      \bigr)
      \\
      O(2^{2\cdot t}\cdot t),
      \\
      O\Bigl(
      2^{t\cdot h}\cdot n^{m+1}
      \cdot\bigl(
      \ell\cdot m\cdot n +
      n\cdot(1+(i-1)^2\cdot m\cdot n^2 + (i-1)\cdot h\cdot t_*^2\cdot n)\cdot 2^{i-1}\cdot n^{2\cdot i - 4}
      \bigr)
      \Bigr),
      \\
      O\bigl(
      2^{t\cdot h}\cdot n^{m+1}
      \cdot h\cdot t^2
      \bigr).
    \end{gather*}
    Recalling that
    \begin{gather*}
      t
      \leq
      t_*
      \leq
      \log_2(n) + O\Bigl(\ell\cdot\bigl(\len(c_{\Root}) + \len(c_{\ZZ}) + \ell\cdot\len(\epsilon)
      + \len(c_{\samp}) + \len(c_{\adj})\bigr)\Bigr),
    \end{gather*}
    and $h=2\cdot m + 2$, the above can be bounded together by
    \begin{gather*}
      t_*^2\cdot n^{O(m + K)}\cdot 2^{t\cdot h}
      \leq
      2^{O(\ell\cdot d\cdot(\len(c_{\Root}) + \len(c_{\ZZ}) + \ell\cdot\len(\epsilon) + \len(c_{\samp}) + \len(c_{\adj})))/(c_{\atom}^2\cdot\delta^2)}
      \cdot n^{O(d\cdot\epsilon^{-\ell-1}/(c_{\atom}^2\cdot\delta^2))}.
    \end{gather*}
    Taking into account all $K\leq O(\epsilon^{-\ell-1})$ executions of the loop does not change the bound above.

    Then the final time-complexity of Algorithm~\ref{alg:equiPRG:atomspace} is bounded by:
    \begin{multline*}
      2^{O(\ell\cdot d\cdot(\len(c_{\Root}) + \len(c_{\ZZ}) + \ell\cdot\len(\epsilon) + \len(c_{\samp}) + \len(c_{\adj})))/(c_{\atom}^2\cdot\delta^2)}
      \cdot n^{O(d\cdot\epsilon^{-\ell-1}/(c_{\atom}^2\cdot\delta^2))}
      \\
      +
      \log(d+1)\cdot\len(c_{\atom})\cdot\len(\delta)
      \\
      +
      \ell^4\cdot\log(\ell+1)\cdot\len(c_{\ZZ})^2\cdot\len(c_{\samp})^2\cdot\len(c_{\Root})^2\cdot\len(c_{\Sum})\cdot\len(\epsilon)^2
      \cdot\bigl(\log(n+1)\bigr)^2.     
    \end{multline*}
  \item[Item~\ref{thm:equiPRG:goodspace}.] Algorithm~\ref{alg:equiPRG:goodspace} is essentially
    Algorithm~\ref{alg:equiPRG:atomspace} with $m\df 1$ and $\delta\df 1/2^{\ell+2}$, and we use the fact that
    $\epsilon$-goodness is equivalent to $(\delta,\epsilon,1)$-atomicity provided $\delta < 1$.

    \begin{algorithm}[htbp]
      \caption{Computation algorithm of query-into-oracle model with space- and time-complexities
        \begin{gather*}
          \begin{multlined}[t]
            O\bigggl(
            \epsilon^{-\ell-1}\cdot
            \left(\log\left(\frac{n}{\epsilon^\ell}\right)
            +
            \ell
            \cdot\bigl(\len(c_{\Root}) + \len(c_{\ZZ}) + \ell\cdot\len(\epsilon) + \len(c_{\samp}) + \len(c_{\adj})\bigr)
            \right)
            \\
            + \ell\cdot\len(c_{\adj})
            \bigggr)
          \end{multlined}
          \\
          \begin{multlined}[t]
            2^{O(\ell\cdot(\len(c_{\Root}) + \len(c_{\ZZ}) + \ell\cdot\len(\epsilon) + \len(c_{\samp}) + \len(c_{\adj})))}
            \cdot n^{O(\epsilon^{-\ell-1})}
            \\
            +
            \ell^4\cdot\log(\ell+1)\cdot\len(c_{\ZZ})^2\cdot\len(c_{\samp})^2\cdot\len(c_{\Root})^2\cdot\len(c_{\Sum})\cdot\len(\epsilon)^2
            \cdot\bigl(\log(n+1)\bigr)^2.
          \end{multlined}
        \end{gather*}
        respectively, that returns
        $(m,\widetilde{m},h,r,(\sigma^k,x^k,s^k,\widetilde{s}^k,r_k,P_k,\widetilde{p}_k,Z^k)_{k=1}^K)$ such that when given to
        Algorithm~\ref{alg:equiPRG:atomspaceoracle} produces an oracle for an $\epsilon$-good partition of $G$ into $K\df
        K_{\equiPRG}(\ell,c_{\ZZ},c_{\samp},c_{\adj},c_{\Root},c_{\Sum},\epsilon)$ parts.}
      \label{alg:equiPRG:goodspace}
      \DontPrintSemicolon
      \KwIn{Numbers $n,\ell\in\NN_+$ with $\ell\leq\log_2(n)$, $c_{\ZZ},c_{\adj},c_{\Root},c_{\Sum}\in(0,1)\cap\QQ$,
        $c_{\samp}\in(0,1/2)\cap\QQ$, and
        $\epsilon\in(0,\epsilon_{\equiPRG}(\ell,c_{\ZZ},c_{\samp},c_{\adj},c_{\Root},c_{\Sum}))$, and a query-oracle $\cO_G$ for
        a graph $G$ with $\lvert G\rvert=n$, $\Lit(G)\leq\ell$ and $\VC(G)\leq d$. We further assume that $n\geq
        n_{\equiPRG}^{\good}(\ell,c_{\ZZ},c_{\samp},c_{\adj},c_{\Root},c_{\Sum},\epsilon)$.}
      \KwOut{A tuple $(m,\widetilde{m},h,r,(\sigma^k,x^k,s^k,\widetilde{s}^k,r_k,P_k,\widetilde{p}_k,Z^k)_{k=1}^K)$, with $K=
        K_{\equiPRG}(\ell,c_{\ZZ},c_{\samp},c_{\adj},c_{\Root},c_{\Sum},\epsilon)$, where $m,r,h\in\NN_+$, $\widetilde{m}\in\NN$
        and for each $k\in[K]$, we have $\sigma^k\in\{0,1\}^{<\ell+1}$, $x^k\in (V(G)^m)^{\lvert\sigma\rvert}$,
        $s^k,\widetilde{s}^k,r_k\in[n]$, $P_k\in\FF_2[X]$ is an irreducible polynomial, $\widetilde{p}_k\in\{0,\ldots,2^{t_k}\}$
        and $Z^k\in[2^{t_k}]^h$, where $t_k\df\deg(P_k)$ such that when given to Algorithm~\ref{alg:equiPRG:atomspaceoracle}
        produces an oracle for an $\epsilon$-good partition of $G$ into $K$ parts.}
      Let $m\assign 1$\;
      $\widetilde{m}\assign 1$\tcp*[r]{$\widetilde{m} = \floor{\widetilde{\delta}\cdot m}$ for $\widetilde{\delta} < 1$}
      $\zeta\assign (1-c_{\ZZ})\cdot\epsilon$\;
      $\widetilde{c}_{\UC}\assign c_{\UC}\assign 0$\;
      $\widetilde{c}_{\adj}\assign (1-c_{\adj})^{-\ell} - 1$\;
      $\widetilde{c}_{\samp}\assign(1-c_{\samp})^{-\ell} - 1$\;
      $\widetilde{c}_{\ZZ}\assign(1-c_{\ZZ})^{-\ell-1} - 1$\;
      $\widetilde{c}_{\Root}\assign(1-c_{\Root})^{-\ell-1} - 1$\;
      \setcounter{algosplit}{\theAlgoLine}
    \end{algorithm}

    \begin{algorithm}[htbp]
      \caption*{(continued).}
      \DontPrintSemicolon
      \setcounter{AlgoLine}{\thealgosplit}
      $K\assign\ceil{(1+\widetilde{c}_{\UC})\cdot(1+\widetilde{c}_{\ZZ})\cdot(1+\widetilde{c}_{\samp})\cdot(1+\widetilde{c}_{\adj})
        \cdot(1+\widetilde{c}_{\Root})\cdot(1+c_{\Sum})\cdot\epsilon^{-\ell-1}}$\;
      $r\assign\floor{n/K}$\;
      $h\assign h_{\equiPRG}^{\good}\df 4$\;
      $\widetilde{n}\assign n$\;
      \uFor{$i\in[K]$}{%
        \label{alg:equiPRG:goodspace:For}
        Find a rational $\xi\in\QQ$ satisfying:
        \begin{align*}
          0 & < \xi\leq(1-c_{\Root})\cdot\zeta,
          &
          0 & \leq -1 + \zeta - \xi + \gamma^\ell\cdot\frac{\widetilde{n}}{r} \leq c_{\Root}\cdot\zeta,
        \end{align*}
        where
        \begin{align*}
          \gamma & \df \widetilde{\epsilon} \df (1-c_{\samp})\cdot\widetilde{\xi},
          &
          \widetilde{\xi} & \df \xi,
        \end{align*}
        and $\gamma$ is a dyadic.\;\label{alg:equiPRG:goodspace:xii}
        $r_i\assign\ceil{\gamma^\ell\cdot\widetilde{n}}$\;\label{alg:equiPRG:goodspace:ri}
        $s_0\df\widetilde{n}$\;
        \lFor{$j\in[\ell+1]$}{\label{alg:equiPRG:goodspace:sj}$s_j\assign\floor{\gamma_i\cdot s_{j-1}}+1$}
        Let $\cO$ be Algorithm~\ref{alg:equiPRG:atomspaceauxiliary} with parameters:
        \begin{multline*}
          \Bigl(n,\cO_G,
          \bigl(m,\widetilde{m},h,r,(\sigma^k,x^k,s^k,\widetilde{s}^k,r_k,P_k,\widetilde{p}_k,Z^k)_{k=1}^K\bigr),
          u\Bigr)
          \\
          =
          \Bigl(n,\cO_G,
          \bigl(m,\widetilde{m},h,r,(\sigma^k,x^k,s^k,\widetilde{s}^k,r_k,P_k,\widetilde{p}_k,Z^k)_{k=1}^{i-1}\bigr),
          \place\Bigr).
        \end{multline*}
        \;
        Run Algorithm~\ref{alg:ext:atomspace} with parameters:
        \begin{gather*}
          (n,m,\ell,s_0,\ldots,s_\ell,s_{\ell+1},\widetilde{m},\cO_G,\cO_U)
          \df
          (n,m,\ell,s_0,\ldots,s_\ell,s_{\ell+1},\widetilde{m},\cO_G,\cO).
        \end{gather*}
        \;\label{alg:equiPRG:goodspace:call}
        Let $(\sigma,(x_j)_{j=0}^{\lvert\sigma\rvert-1})$ be the tuple returned by the algorithm.\tcp*[r]{Encodes set $W$}
        $\sigma^i\assign\sigma$\;
        $x^i\assign (x_j)_{j=0}^{\lvert\sigma\rvert-1}$\;
        $s^i\assign s_{\lvert\sigma^i\rvert}$\tcp*[r]{$s^i\df\lvert W\rvert$}
        Let $p\df\gamma_i^{\ell-\lvert\sigma\rvert}$ and write it as $\widetilde{p}_i/2^t$ with $t\geq\ceil{\log_2(n)}$.\;
        \setcounter{algosplit}{\theAlgoLine}
      }
    \end{algorithm}

    \begin{algorithm}[htbp]
      \caption*{(continued).}
      \DontPrintSemicolon
      \setcounter{AlgoLine}{\thealgosplit}
      \let\oldnl\nl
      \let\nl\relax
      \vspace{-\baselineskip}\uInvisibleBegin{
        \global\let\nl\oldnl 
        Run Algorithm~\ref{alg:encodePRG} with parameter $t\df t$ and let $P_i$ be the returned polynomial.\;
        \uFor{$Z^i\in[2^t]^h$}{%
          \tcp*[l]{Use Algorithm~\ref{alg:runPRG} to implicitly generate pseudorandom $X$ from $Z^i$}
          $\widetilde{s}^i\assign 0$\tcp*[r]{Size of $W(X)\df\{u\in W\mid X_u=1\}$}
          \For(\tcp*[f]{Compute $\widetilde{s}^i\df\lvert W(X)\rvert$}){$u\in[n]$}{%
            Run Algorithm~\ref{alg:ext:atomspaceoracle} with parameters:
            \begin{gather*}
              (n,m,\widetilde{m},\cO_G,\cO_U,(\sigma,(x_j)_{j=0}^{\lvert\sigma\rvert-1}),u)
              \df
              (n,m,\widetilde{m},\cO_G,\cO,(\sigma,(x_j)_{j=0}^{\lvert\sigma\rvert-1}),\place)
            \end{gather*}
            \;
            \lIf(\tcp*[f]{$u\notin W$}){The algorithm returned $0$}{\Continue}
            Run Algorithm~\ref{alg:runPRG} with parameters:
            \begin{gather*}
              (P,\widetilde{p},h,n,Z,u)\df(P_i,\widetilde{p}_i,h,n,Z^i,u).
            \end{gather*}
            \;
            \lIf(\tcp*[f]{$u\in W(X)$}){The algorithm returned $1$}{$\widetilde{s}^i\assign\widetilde{s}^i+1$}
          }
          \lIf{$\lvert\widetilde{s}^i - p\cdot s^i\rvert > c_{\adj}\cdot c_{\samp}\cdot\widetilde{\xi}\cdot p\cdot s^i$}{%
            \Continue
          }
          \uFor(\tcp*[f]{Check if $W(X)$ is a $(\widetilde{\delta},\widetilde{\xi},m)$-atom of $G$}){$z\in[n]^m$}{%
            \label{alg:equiPRG:goodspace:xFor}
            $k^0\assign 0$\tcp*[r]{Will count $\lvert\{u\in W(X) \mid t_G^{\widetilde{\delta}}(u,z) = 0\}\rvert$}
            $k^1\assign 0$\tcp*[r]{Will count $\lvert\{u\in W(X) \mid t_G^{\widetilde{\delta}}(u,z) = 1\}\rvert$}
            \uFor(\tcp*[f]{Compute $k^0$ and $k^1$}){$u\in [n]$}{%
              Run Algorithm~\ref{alg:ext:atomspaceoracle} with parameters:
              \begin{gather*}
                (n,m,\widetilde{m},\cO_G,\cO_U,(\sigma,(x_j)_{j=0}^{\lvert\sigma\rvert-1}),u)
                \df
                (n,m,\widetilde{m},\cO_G,\cO,(\sigma,(x_j)_{j=0}^{\lvert\sigma\rvert-1}),\place)
              \end{gather*}
              \;
              \lIf(\tcp*[f]{$u\notin W$}){The algorithm returned $0$}{\Continue}
              Run Algorithm~\ref{alg:runPRG} with parameters:
              \begin{gather*}
                (P,\widetilde{p},h,n,Z,u)\df(P_i,\widetilde{p}_i,h,n,Z^i,u).
              \end{gather*}
              \;
              \lIf(\tcp*[f]{$u\notin W(X)$}){The algorithm returned $0$}{\Continue}
              %
              \setcounter{algosplit}{\theAlgoLine}
            }
          }
        }
      }
    \end{algorithm}

    \begin{algorithm}[htbp]
      \caption*{(continued).}
      \DontPrintSemicolon
      \setcounter{AlgoLine}{\thealgosplit}
      \let\oldnl\nl
      \let\nl\relax
      \vspace{-\baselineskip}\InvisibleBegin{
        \vspace{-\baselineskip}\InvisibleBegin{
          \vspace{-\baselineskip}\InvisibleBegin{
            \vspace{-\baselineskip}\InvisibleBegin{
              \global\let\nl\oldnl 
              $a\assign 0$\tcp*[r]{Will count $\lvert\{j\in[m] \mid z_j\in N_G(u)\}\rvert$}
              \For{$j\in[m]$}{%
                \lIf{$\cO_G(u,z_j)=1$}{$a\assign a+1$}
              }
              \lIf(\tcp*[f]{$t_G^{\widetilde{\delta}}(u,z)=1$}){%
                $a\geq m - \widetilde{m}$
              }{%
                \label{alg:equiPRG:goodspace:k1update}
                $k^1\assign k^1+1$
              }
              \lElseIf(\tcp*[f]{$t_G^{\widetilde{\delta}}(u,z)=0$}){%
                $a\leq\widetilde{m}$
              }{%
                \label{alg:equiPRG:goodspace:k0update}
                $k^0\assign k^0+1$
              }
            }
            \If(\tcp*[f]{$W(X)$ is $(\widetilde{\delta},\widetilde{\xi}_i)$-split by $z$}){%
              $k^0 > \widetilde{\xi}\cdot\widetilde{s}^i$ and $k^1 > \widetilde{\xi}\cdot\widetilde{s}^i$}{%
              \Break\tcp*[r]{$W(X)$ is not a $(\widetilde{\delta},\widetilde{\xi}_i,m)$-atom}
            }
          }
          \Else(\tcp*[f]{Executes only if ``for'' loop of line~\ref{alg:equiPRG:goodspace:xFor} didn't break}){%
            \Break\tcp*[r]{$W(X)$ is a $(\widetilde{\delta},\widetilde{\xi}_i,m)$-atom}
          }
        }
        $\widetilde{n}\assign \widetilde{n} - r_k$\;
      }
      \Return{$(m,\widetilde{m},h,r,(\sigma^k,x^k,s^k,\widetilde{s}^k,r_k,P_k,\widetilde{p}_k,Z^k)_{k=1}^K)$}
    \end{algorithm}

    The complexity analysis of Algorithm~\ref{alg:equiPRG:goodspace} is also analogous to that of
    Algorithm~\ref{alg:equiPRG:atomspace} and yields the space-complexity bound
    \begin{multline*}
      O\bigggl(
      \epsilon^{-\ell-1}\cdot
      \left(\log\left(\frac{n}{\epsilon^\ell}\right)
      +
      \ell
      \cdot\bigl(\len(c_{\Root}) + \len(c_{\ZZ}) + \ell\cdot\len(\epsilon) + \len(c_{\samp}) + \len(c_{\adj})\bigr)
      \right)
      \\
      + \ell\cdot\len(c_{\adj})
      \bigggr)
    \end{multline*}
    and the time-complexity bound
    \begin{multline*}
      2^{O(\ell\cdot(\len(c_{\Root}) + \len(c_{\ZZ}) + \ell\cdot\len(\epsilon) + \len(c_{\samp}) + \len(c_{\adj})))}
      \cdot n^{O(\epsilon^{-\ell-1})}
      \\
      +
      \ell^4\cdot\log(\ell+1)\cdot\len(c_{\ZZ})^2\cdot\len(c_{\samp})^2\cdot\len(c_{\Root})^2\cdot\len(c_{\Sum})\cdot\len(\epsilon)^2
      \cdot\bigl(\log(n+1)\bigr)^2.
      \qedhere
    \end{multline*}
  \end{description}
\end{proof}

\begin{discussion}\label{dsc:dyadic}
  Since in Appendix~\ref{sec:PRGgen} we also provide a pseudorandom generator capable of handling any rational bias, one might
  wonder if that could be used in the algorithms of Theorem~\ref{thm:equiPRG} instead of ensuring that the $\gamma_i$ are
  dyadics. Unfortunately, doing so would hurt the complexity of the algorithms considerably: in the worst case $\gamma_i$ would
  be represented as $a/b$ with $b$ having $\Omega(\log(\epsilon^{-1}))$ prime factors and to use Theorem~\ref{thm:PRGgen}, we
  would need to take an exponent $e_i$ of each prime factor $q_i$ to be at least $\ceil{\log_{q_i}(n)}$, and more specifically,
  we would have $Q\geq\Omega(n^{\log(\epsilon^{-1})})$, which is much worse than what is obtained with Theorem~\ref{thm:equiPRG}
  (which is equivalent to $Q=\Theta(n)$).
\end{discussion}


\section{Lower bounds}

\begin{definition}\label{def:selector}
  Let $\ell,r\in\NN_+$. The \emph{selector graph} of height $\ell$ and width $r$ (see Figure~\ref{fig:selector}) is the graph
  $G_{\ell,r}$ given by
  \begin{align*}
    V(G_{\ell,r})
    & \df
    \bigl([r]^{<\ell+1}\setminus\{\varnothing\}\bigr)\times\{1\}
    \cup[r]^\ell\times\{2\},
    \\
    E(G_{\ell,r}) & \df \bigl\{\{(\tau,2),(\tau\rest_{[s]},1)\}\;\bigm\vert\; \tau\in[r]^\ell\land s\in[\ell]\bigr\}.
  \end{align*}
  (There is no vertex corresponding to the empty string.) Note that $G_{\ell,r}$ is bipartite with bipartition
  $(([r]^{<\ell+1}\setminus\{\varnothing\})\times\{1\},[r]^\ell\times\{2\})$, every vertex in $[r]^\ell\times\{2\}$ has degree
  exactly $\ell$ and for every $\sigma\in[r]^{<\ell+1}\setminus\{\varnothing\}$, the neighborhood of $(\sigma,1)$ is exactly the
  set of $(\tau,2)\in[r]^\ell\times\{2\}$ such that $\tau$ extends $\sigma$ (which has size $r^{\ell-\lvert\sigma\rvert}$).

  \begingroup
\def\ptsize{2pt}
\def\absep{2}
\def\scale{1}

\newcommand{\drawselector}[4]{
  \centering
  \def\height{#1}
  \def\width{#2}
  \def\horbasesep{#3}
  \def\layersep{#4}
  \begin{tikzpicture}[scale=\scale,rotate={-90}]
    \def\strings{v,}
    \foreach \t [%
      evaluate=\t as \y using -(\t-1)*\layersep,
      evaluate=\t as \basehorstep using \width^(\height-\t)%
    ] in {0,...,\height}{
      \global\let\oldstrings\strings
      \def\strings{}
      \foreach \s [%
        count=\c,
        evaluate=\s as \x using (\basehorstep*(\c-1) + (\basehorstep-1)/2)*\horbasesep%
      ]
      in \oldstrings {
        \ifx\s\empty\relax
        \else
        \coordinate (\s) at (\x,\y);
        \StrGobbleLeft{\s}{1}[\l]
        \IfEq{\s}{v}{%
        }{%
          \fill (\s) circle (\ptsize);
          \node[left] at (\s) {$(\l,1)$};%
        }
        \foreach \j in {1,...,\width}{
          \edef\newstrings{\strings \s\j,}
          \global\let\strings\newstrings
        }
        \fi
      }
    }
    \pgfmathsetmacro{\y}{\absep}
    \foreach \s [%
      count=\c,
      evaluate=\s as \x using (\c-1)*\horbasesep%
    ] in \oldstrings {
      \ifx\s\empty\relax
      \else
      \coordinate (w\s) at (\x,\y);
      \fill (w\s) circle (\ptsize);
      \StrGobbleLeft{\s}{1}[\l]
      \node[right] at (w\s) {$(\l,2)$};
      \foreach \t [%
        evaluate=\t as \cy using -(\t-1)*\layersep%
      ] in {1,...,\height} {
        \pgfmathtruncatemacro{\n}{\t+1}
        \StrLeft{\s}{\n}[\p]
        \draw (\p) -- (w\s);%
      }
      \fi
    }
  \end{tikzpicture}
  \caption*{\normalsize $G_{\height,\width}$}
}

\begin{figure}
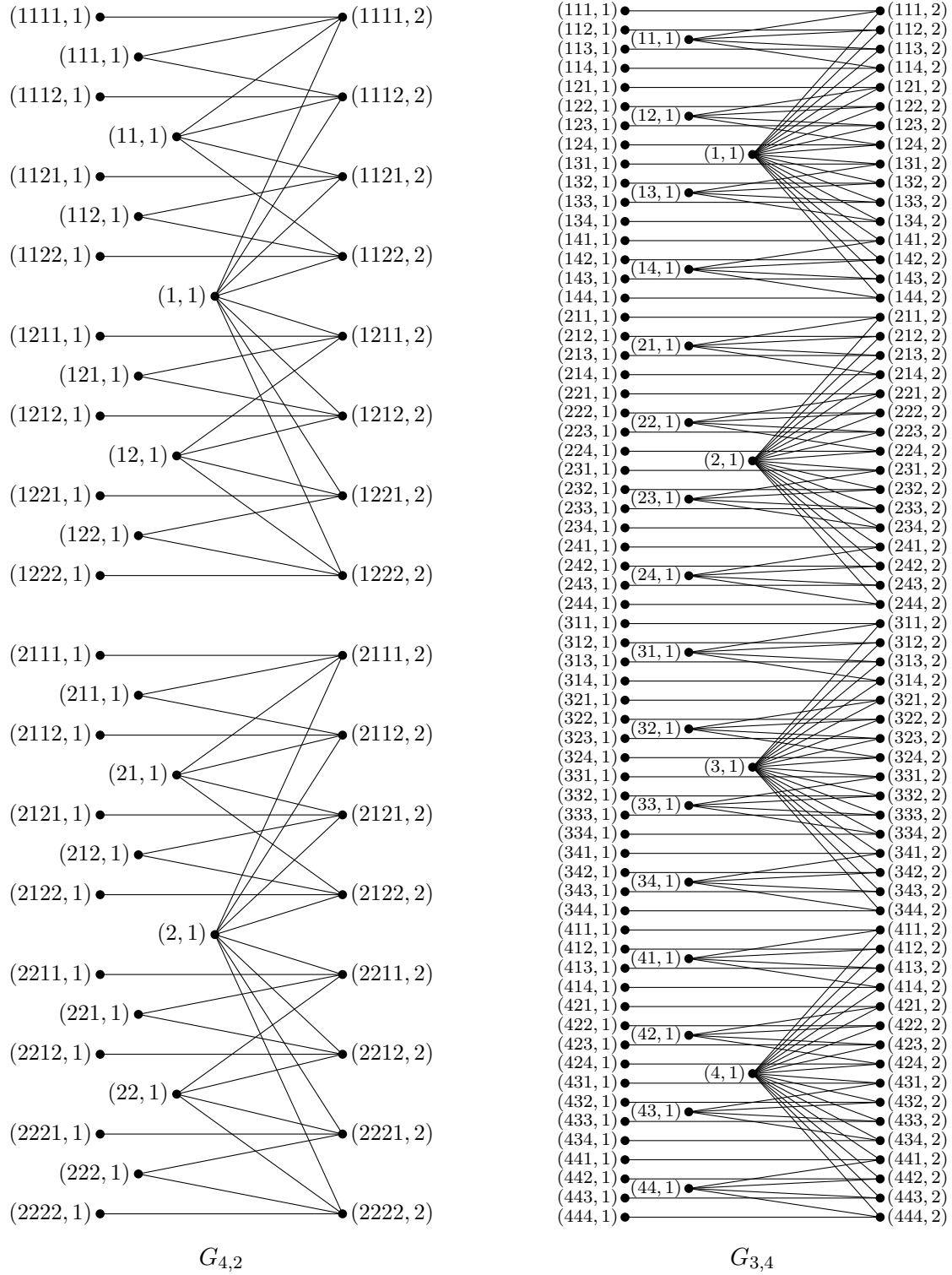

  \centering
  \begin{subfigure}{0.45\linewidth}
    \small
    \drawselector{4}{2}{1.25}{0.6}
  \end{subfigure}
  \qquad
  \begin{subfigure}{0.45\linewidth}
    \scriptsize
    \drawselector{3}{4}{0.3}{1}
  \end{subfigure}
  \caption{Selector graph $G_{\ell,r}$ of Definition~\ref{def:selector}. (Due to space constraints, height is represented
    horizontally in the picture and width is represented vertically.)}
  \label{fig:selector}
\end{figure}
\endgroup


  Given further $n\in\NN_+$ and $\lambda\in\RR_{\geq 0}$, we let $G_{\ell,r}(n,\lambda)$ be the independent blow-up
  $G_{\ell,r}^b$ of $G_{\ell,r}$ (see Definition~\ref{def:indepblowup}) corresponding to the function $b\colon
  V(G_{\ell,r})\to\NN$ given by
  \begin{gather*}
    b(\eta,i) \df
    \begin{dcases*}
      1, & if $i=1$,\\
      n, & if $i=2$ and $\eta\neq\One$,\\
      \floor{\lambda\cdot n}, & if $i=2$ and $\eta=\One$,
    \end{dcases*}
  \end{gather*}
  where $\One\df(1,1,\ldots,1)\in[r]^\ell$ is the all $1$s string. With a slight abuse of notation, since the vertices in
  $([r]^{<\ell+1}\setminus\{\varnothing\})\times\{1\}$ do not get duplicated, we will think of
  $([r]^{<\ell+1}\setminus\{\varnothing\})\times\{1\}$ as a subset of $V(G_{\ell,r}(n,\lambda))$ and for each $\tau\in[r]^\ell$,
  we let $U_\tau\subseteq V(G_{\ell,r}(n,\lambda))$ be the set of all copies generated from the vertex $(\tau,2)\in
  V(G_{\ell,r})$.

  For each $\sigma\in[r]^{<\ell+1}$, we also let
  \begin{gather*}
    U_\sigma \df \bigcup_{\substack{\tau\in[r]^\ell\\\tau\rest_{[\lvert\sigma\rvert]} = \sigma}} U_\tau
  \end{gather*}
  be the union of all $U_\tau$ indexed by strings that extend $\sigma$. In this notation, $G_{\ell,r}(n,\lambda)$ is bipartite
  with partition $(([r]^{<\ell+1}\setminus\{\varnothing\})\times\{1\}, U_\varnothing)$ and for
  $\sigma\in[r]^{<\ell+1}\setminus\{\varnothing\}$, the neighborhood of $(\sigma,1)$ in $G_{\ell,r}(n,\lambda)$ is exactly $U_\sigma$.
\end{definition}

\begin{lemma}\label{lem:GellrLit}
  For $\ell,r\in\NN_+$, we have $\Lit(G_{\ell,r})\leq\ell$. In particular, for every $n\in\NN_+$ and $\lambda\in\RR_{\geq 0}$,
  we have $\Lit(G_{\ell,r}(n,\lambda))\leq\ell$.
\end{lemma}

\begin{proof}
  The second assertion follows from the first along with Lemma~\ref{lem:indepblowup}. For the first assertion, suppose for a
  contradiction that $((u_\sigma)_{\sigma\in\{0,1\}^{<\ell+1}}, (N_{G_{\ell,r}}(w_\tau))_{\tau\in\{0,1\}^{\ell+1}})$ is an
  $(\ell+1)$-Littlestone tree in $G_{\ell,r}$. Since the vertex $w_{\One}$ corresponding to the all $1$s string $\One\in\{0,1\}^{\ell+1}$
  has at least $\ell+1$ neighbors, namely, the vertices
  \begin{gather}\label{eq:GellrLit:nhd}
    u_\varnothing, u_{(1)}, u_{(1,1)},\ldots, u_{\mathop{\underbrace{\scriptstyle (1,\ldots,1)}}\limits_{\mathclap{\text{$\ell$ times}}}}
  \end{gather}
  (which must be distinct), it follows that $w_{\One}\in([r]^{<\ell+1}\setminus\{\varnothing\})\times\{1\}$ (as all vertices in
  $[r]^\ell\times\{2\}$ have degree $\ell$). In turn, since $G_{\ell,r}$ is bipartite, all vertices in~\eqref{eq:GellrLit:nhd}
  must be in $[r]^\ell\times\{2\}$. If we inspect how $u_\varnothing\in[r]^\ell\times\{2\}$ behaves in the Littlestone tree, we
  see that it must have at least $2^\ell$ neighbors, namely, the vertices in
  \begin{gather*}
    \{w_\tau \mid \tau\in\{0,1\}^{\ell+1}\land\tau_1 = 1\}
  \end{gather*}
  (which must be distinct) and since $u_\varnothing\in[r]^\ell\times\{2\}$ implies that its degree is $\ell$, we must have
  $2^\ell\leq\ell$, which is a contradiction.
\end{proof}

\begin{lemma}\label{lem:Gellrgood}
  Let $\ell,r,n\in\NN_+$, let $\lambda\in\RR_{\geq 0}$, let $\epsilon\in(0,1/2)$ and suppose that $r < (1-\epsilon)/\epsilon$
  and $P\subseteq U_\varnothing$ is $\epsilon$-good in $G_{\ell,r}(n,\lambda)$. Then the following hold:
  \begin{enumerate}
  \item\label{lem:Gellrgood:maj} For every $\sigma\in[r]^{<\ell+1}$, $\lvert P\cap U_\sigma\rvert/\lvert
    P\rvert\in[0,\epsilon]\cup[1-\epsilon,1]$.
  \item\label{lem:Gellrgood:unique} For every $t\in\{0,\ldots,\ell\}$, there exists a unique $\sigma\in[r]^t$ such that $\lvert
    P\cap U_\sigma\rvert/\lvert P\rvert\in[1-\epsilon,1]$.
  \end{enumerate}
\end{lemma}

\begin{proof}
  Item~\ref{lem:Gellrgood:maj} is easy: for $\sigma=\varnothing$, it is obvious as $P\subseteq U_\varnothing$ and for
  $\sigma\neq\varnothing$, we note that the neighborhood of the vertex $(\sigma,1)$ in $G_{\ell,r}(n,\lambda)$ is exactly
  $U_\sigma$, so since $P$ is $\epsilon$-good, it must have an $\epsilon$-majority opinion, which implies
  \begin{gather*}
    \frac{\lvert P\cap U_\sigma\rvert}{\lvert P\rvert}
    =
    \frac{\lvert P\cap N_{G_{\ell,r}(n,\lambda)}((\sigma,1))\rvert}{\lvert P\rvert}
    \in
    [0,\epsilon]\cup[1-\epsilon,1].
  \end{gather*}

  \medskip

  The proof of item~\ref{lem:Gellrgood:unique} is by induction in $t$. For $t=0$, the result is obvious as $P\subseteq
  U_\varnothing$. For $t > 0$, by inductive hypothesis, there must exist a unique $\sigma'\in[r]^{t-1}$ such that $\lvert P\cap
  U_{\sigma'}\rvert/\lvert P\rvert\geq 1-\epsilon$.

  If $\sigma\in[r]^t$ does not extend $\sigma'$, then $U_\sigma$ is disjoint from $U_{\sigma'}$, from which we get $\lvert P\cap
  U_\sigma\rvert/\lvert P\rvert\leq\epsilon$.

  On the other hand, we know that
  \begin{gather*}
    U_{\sigma'} = \bigcup_{s\in[r]} U_{\sigma'\conc s},
  \end{gather*}
  so there must exist $s_*\in[r]$ such that
  \begin{gather*}
    \frac{\lvert U_{\sigma'\conc s_*}\cap P\rvert}{\lvert P\rvert}
    \geq
    \frac{1}{r}\cdot\frac{\lvert U_{\sigma'}\cap P\rvert}{\lvert P\rvert}
    \geq
    \frac{1-\epsilon}{r}
    >
    \epsilon,
  \end{gather*}
  where the last inequality follows since $r < (1-\epsilon)/\epsilon$. By item~\ref{lem:Gellrgood:maj}, the string
  $\sigma\df\sigma'\conc s_*$ must then satisfy
  \begin{gather*}
    \frac{\lvert U_\sigma\cap P\rvert}{\lvert P\rvert}\geq 1 - \epsilon.
  \end{gather*}

  Finally, uniqueness follows since the sets $U_\sigma$ for $\sigma\in[r]^t$ form a collection of pairwise
  disjoint sets and $\epsilon < 1/2$.
\end{proof}

\begin{theorem}\label{thm:Gellr:nonequi}
  Let $\ell,n\in\NN_+$, let $\epsilon\in(0,1/2)$, let $r\df\ceil{(1-\epsilon)/\epsilon}-1$ and $\lambda\df 1$. Then $r\geq 1$,
  so the graph $G_{\ell,r}(n,\lambda)$ is well-defined, every $\epsilon$-good in $G_{\ell,r}(n,\lambda)$ set $P\subseteq
  U_\varnothing$ satisfies $\lvert P\rvert\leq n/(1-\epsilon)$ and every partition $\cP$ of $U_\varnothing$ into $\epsilon$-good
  sets of $G_{\ell,r}(n,\lambda)$ satisfies
  \begin{gather*}
    \lvert\cP\rvert
    \geq
    (1-\epsilon)\cdot r^\ell
    =
    \bigl(1 + o_{\epsilon\to 0}(1)\bigr)\cdot\epsilon^{-\ell}.
  \end{gather*}
\end{theorem}

\begin{proof}
  Note that since $\epsilon < 1/2$, we have
  \begin{gather*}
    r \df \Ceil{\frac{1-\epsilon}{\epsilon}} - 1\geq 1.
  \end{gather*}
  Furthermore, we clearly have $r < (1-\epsilon)/\epsilon$.

  If $P\subseteq U_\varnothing$ is an $\epsilon$-good in $G_{\ell,r}(n,\lambda)$ set, then by Lemma~\ref{lem:Gellrgood}, there
  exists a unique $\tau_P\in[r]^\ell$ such that
  \begin{gather*}
    (1-\epsilon)\cdot\lvert P\rvert
    \leq
    \lvert P\cap U_{\tau_P}\rvert
    \leq
    n,
  \end{gather*}
  where the last inequality follows since $\lambda=1$ (so all $U_\tau$ have size $n$), hence $\lvert P\rvert\geq
  n/(1-\epsilon)$.

  Finally, if $\cP$ is a partition of $U_\varnothing$ into $\epsilon$-good in $G_{\ell,r}(n,\lambda)$ sets, then
  \begin{gather*}
    r^\ell\cdot n = \lvert U_\varnothing\rvert = \sum_{P\in\cP}\lvert P\rvert \leq \lvert\cP\rvert\cdot\frac{n}{1-\epsilon},
  \end{gather*}
  from which it follows that
  \begin{gather*}
    \lvert\cP\rvert
    \geq
    (1-\epsilon)\cdot r^\ell
    =
    \bigl(1 + o_{\epsilon\to 0}(1)\bigr)\cdot\epsilon^{-\ell}.
    \qedhere
  \end{gather*}
\end{proof}

\begin{corollary}\label{cor:Gellr:nonequi}
  For every $\ell\in\NN_+$, every $\epsilon\in(0,1/2)$ and every $n\in\NN_+$, there exists a graph $G$ with $\lvert G\rvert=n$,
  $\Lit(G)\leq\ell$ and such that every $\epsilon$-good partition of $G$ has size at least
  \begin{gather*}
    \bigl(1 + o_{n\to\infty, \epsilon\to 0,\ell}(1)\bigr)\cdot\epsilon^{-\ell}.
    \epsilon^{-\ell}.
  \end{gather*}
\end{corollary}

\begin{proof}
  We consider the graph $G_{\ell,r}(n,\lambda)$ as in Theorem~\ref{thm:Gellr:equi} using $\epsilon'$ in place of $\epsilon$ and $n$
  large enough. Lemma~\ref{lem:GellrLit} guarantees that $\Lit(G_{\ell,r}(n,\lambda))\leq\ell$ regardless of how we pick
  $\epsilon'$ and $n$.

  Our choice is as follows: given $\alpha\in(0,1)$ small enough so that $\epsilon\cdot(1+\alpha) < 1/2$, we set $\epsilon'\df
  (1+\alpha)\cdot\epsilon$.

  Let $\cP$ be an $\epsilon$-good equipartition $\cP$ of $G_{\ell,r}(n,\lambda)$ (but each part can have vertices in
  $([r]^{<\ell+1}\setminus\{\varnothing\})\times\{1\}$) and let
  \begin{gather*}
    \cP' \df \{P\in\cP \mid \lvert P\setminus U_\varnothing\rvert \leq \alpha\cdot\lvert P\cap U_\varnothing\rvert\}.
  \end{gather*}
  Since each $P\in\cP\setminus\cP'$ must necessarily intersect the set $V(G_{\ell,r}(n,\lambda))\setminus U_\varnothing$, whose
  size is $(r^{\ell+1}-r)/(r-1)$, it follows that $\lvert\cP\setminus\cP'\rvert \leq (r^{\ell+1}-r)/(r-1)$. In turn, we get
  \begin{gather*}
    \sum_{P\in\cP\setminus\cP'} \lvert P\cap U_\varnothing\rvert
    \leq
    \sum_{P\in\cP\setminus\cP'} \frac{\lvert P\setminus U_\varnothing\rvert}{\alpha}
    \leq
    \frac{\lvert V(G_{\ell,r}(n,\lambda))\setminus U_\varnothing\rvert}{\alpha}
    =
    \frac{r^{\ell+1}-r}{\alpha\cdot(r-1)}.
  \end{gather*}

  We now claim that for every $P\in\cP'$, the set $P\cap U_\varnothing$ is $\epsilon'$-good in $G_{\ell,r}(n,\lambda)$. Indeed,
  since $P$ is $\epsilon$-good in $G_{\ell,r}(n,\lambda)$, by Lemma~\ref{lem:downcont}\ref{lem:downcont:good}, we know that
  $P\cap U_\varnothing$ is $\epsilon_-$-good in $G_{\ell,r}(n,\lambda)$, where
  \begin{gather*}
    \epsilon_-
    \df
    \epsilon\cdot\frac{\lvert P\rvert}{\lvert P\cap U_\varnothing\rvert}
    \leq
    \epsilon + \frac{\lvert P\setminus U_\varnothing\rvert}{\lvert P\cap U_\varnothing\rvert}
    \leq
    \epsilon + \alpha
    =
    \epsilon',
  \end{gather*}
  where the second inequality follows since $P\in\cP'$.

  By Theorem~\ref{thm:Gellr:nonequi}, we then conclude that for every $P\in\cP'$, we have $\lvert P\cap U_\varnothing\rvert\leq
  n/(1-\epsilon)$. Thus, we have
  \begin{align*}
    r^\ell\cdot n
    & =
    \lvert U_\varnothing\rvert
    =
    \sum_{P\in\cP} \lvert P\cap U_\varnothing\rvert
    =
    \sum_{P\in\cP\setminus\cP'} \lvert P\cap U_\varnothing\rvert
    +
    \sum_{P\in\cP'} \lvert P\cap U_\varnothing\rvert
    \\
    & \leq
    \frac{r^{\ell+1}-r}{\alpha\cdot(r-1)}
    +
    \lvert\cP'\rvert\cdot\frac{n}{1-\epsilon}
    \leq
    \frac{r^{\ell+1}-r}{\alpha\cdot(r-1)}
    +
    \lvert\cP\rvert\cdot\frac{n}{1-\epsilon},
  \end{align*}
  from which we conclude that
  \begin{gather*}
    \lvert\cP\rvert
    \geq
    (1-\epsilon)\cdot\left(r^\ell - \frac{r^{\ell+1}-r}{\alpha\cdot(r-1)\cdot n}\right)
    =
    \bigl(1 + o_{n\to\infty,\epsilon\to 0,\alpha,\ell}(1)\bigr)\cdot\epsilon^{-\ell}.
  \end{gather*}
  Since $\alpha$ can be chosen arbitrarily in function of $\epsilon$ and then $n$ chosen large enough as a function of $\alpha$,
  our final bound is
  \begin{gather*}
    \bigl(1 + o_{n\to\infty,\epsilon\to 0,\ell}(1)\bigr)\cdot\epsilon^{-\ell}.
    \qedhere
  \end{gather*}
\end{proof}

\begin{theorem}\label{thm:Gellr:equi}
  Let $\ell,n\in\NN_+$, let $\epsilon\in(0,1/3)$, let $r\df\ceil{(1-\epsilon)/\epsilon}-1$ and $\lambda\df (r^\ell-1)\cdot r$.
  Then $r\geq 2$, $\lambda\in\NN_+$, hence the graph $G_{\ell,r}(n,\lambda)$ is well-defined, and every equipartition of
  $U_\varnothing$ into $\epsilon$-good sets of $G_{\ell,r}(n,\lambda)$ has size at least
  \begin{align*}
    \Floor{\frac{(1-\epsilon)\cdot(r^\ell-1+\lambda)}{1+(1-\epsilon)/n}} + 1
    & \geq
    \bigl(1 + o_{n\to\infty,\epsilon}(1)\bigr)\cdot\ceil{(1-\epsilon)\cdot(r^\ell-1+\lambda)}
    \\
    & =
    \bigl(1 + o_{n\to\infty,\epsilon\to 0}(1)\bigr)\cdot\epsilon^{-\ell-1}.
  \end{align*}
\end{theorem}

\begin{proof}
  Note that since $\epsilon<1/3$, we have
  \begin{gather*}
    r \df \Ceil{\frac{1-\epsilon}{\epsilon}} - 1\geq 2,
  \end{gather*}
  which in turn implies
  \begin{gather*}
    \lambda \df (r^\ell-1)\cdot r \geq (2^\ell-1)\cdot 2 \geq 2.
  \end{gather*}
  Note also that $\lambda$ is an integer, which means that in the graph $G_{\ell,r}(n,\lambda)$, we have $\lvert U_{\One}\rvert =
  \lambda\cdot n$, hence $\lvert U_\varnothing\rvert = (r^\ell - 1 + \lambda)\cdot n$. Finally, note that $r <
  (1-\epsilon)/\epsilon$.

  Let $\cP$ be an equipartition of $U_\varnothing$ into $\epsilon$-good sets of $G_{\ell,r}(n,\lambda)$ and suppose for a contradiction
  that
  \begin{gather*}
    \lvert\cP\rvert
    <
    \Floor{\frac{(1-\epsilon)\cdot(r^\ell-1+\lambda)}{1+(1-\epsilon)/n}} + 1,
  \end{gather*}
  which implies
  \begin{gather*}
    \lvert\cP\rvert
    \leq
    \frac{(1-\epsilon)\cdot(r^\ell-1+\lambda)}{1+(1-\epsilon)/n}.
  \end{gather*}
  Since $\cP$ is an equipartition, it follows that for every $P\in\cP$, we have
  \begin{align*}
    (1-\epsilon)\cdot\lvert P\rvert
    & \geq
    (1-\epsilon)\cdot\Floor{\frac{\lvert U_\varnothing\rvert}{\lvert\cP\rvert}}
    >
    (1-\epsilon)\cdot\left(\frac{\lvert U_\varnothing\rvert}{\lvert\cP\rvert} - 1\right)
    \\
    & =
    \frac{(1-\epsilon)\cdot(r^\ell-1+\lambda)\cdot n}{\lvert\cP\rvert} - (1 - \epsilon)
    \geq
    (n + 1-\epsilon) - (1 - \epsilon)
    =
    n.
  \end{align*}

  Since $r < (1-\epsilon)/\epsilon$, by Lemma~\ref{lem:Gellrgood}, for each part $P\in\cP$, there must exist a unique
  $\tau_P\in[r]^\ell$ such that
  \begin{gather*}
    \frac{\lvert P\cap U_{\tau_P}\rvert}{\lvert P\rvert}
    \geq
    1 - \epsilon
  \end{gather*}
  and since $(1-\epsilon)\cdot\lvert P\rvert > n$, we must have $\lvert U_{\tau_P}\rvert > n$, so it must be the case that
  $\tau_P = \One$. In turn, this implies $\lvert P\setminus U_{\One}\rvert\leq\epsilon\cdot\lvert P\rvert$

  But then we have
  \begin{align*}
    (r^\ell-1)\cdot n
    & =
    \sum_{\tau\in[r]^\ell\setminus\{\One\}} \lvert U_\tau\rvert
    =
    \sum_{\tau\in[r]^\ell\setminus\{\One\}} \sum_{P\in\cP} \lvert P\cap U_\tau\rvert
    =
    \sum_{P\in\cP} \lvert P\setminus U_{\One}\rvert
    \\
    & \leq
    \sum_{P\in\cP}
    \epsilon\cdot\lvert P\rvert
    =
    \epsilon\cdot\lvert U_\varnothing\rvert
    =
    \epsilon\cdot(r^\ell-1+\lambda)\cdot n,
  \end{align*}
  which by rearranging gives
  \begin{gather*}
    \lambda \geq \frac{(r^\ell-1)\cdot(1-\epsilon)}{\epsilon} > (r^\ell-1)\cdot r = \lambda,
  \end{gather*}
  where the last inequality follows since $r < (1-\epsilon)/\epsilon$. Thus, we get a contradiction: $\lambda > \lambda$.
\end{proof}

\begin{corollary}\label{cor:Gellr:equi}
  For every $\ell\in\NN_+$, every $n_0\in\NN_+$ and every $\epsilon\in(0,1/3)$, there exists a graph $G$ with $\lvert
  G\rvert=n\geq n_0$, $\Lit(G)\leq\ell$ and such that every $\epsilon$-good equipartition of $G$ has size at least
  \begin{gather*}
    \bigl(1 + o_{n\to\infty,\epsilon\to 0,\ell}(1)\bigr)\cdot
    \epsilon^{-\ell-1}.
  \end{gather*}
\end{corollary}

\begin{proof}
  The idea is similar to that of Corollary~\ref{cor:Gellr:nonequi}, except that we reduce to Theorem~\ref{thm:Gellr:equi}
  instead of Theorem~\ref{thm:Gellr:nonequi}; furthermore, in our reduction, we need to ensure that we construct an equitable
  partition.

  We consider the graph $G_{\ell,r}(n,\lambda)$ as in Theorem~\ref{thm:Gellr:equi} using $\epsilon'$ in place of $\epsilon$ and
  $n$ large enough. Lemma~\ref{lem:GellrLit} guarantees that $\Lit(G_{\ell,r}(n,\lambda))\leq\ell$ regardless of how we pick
  $\epsilon'$ and $n$.

  Our choice is as follows: given $\alpha\in(0,1)$ small enough so that $\epsilon\cdot(1+\alpha) < 1/3$, we set $\epsilon'\df
  (1+\alpha)\cdot\epsilon$ and recalling that $r\df\ceil{(1-\epsilon')/\epsilon'}-1$ and $\lambda\df(r^\ell-1)\cdot r$, let $n$
  be large enough so that $n\geq n_0$ and
  \begin{gather}\label{eq:Gellr:equi:nlarge}
    \epsilon + \frac{(r^{\ell+1}-r)\cdot(1-\epsilon')}{(r-1)\cdot n}
    \leq
    \epsilon'.
  \end{gather}

  Let $\cP$ be an $\epsilon$-good equipartition $\cP$ of $G_{\ell,r}(n,\lambda)$ (but each part can have vertices in
  $([r]^{<\ell+1}\setminus\{\varnothing\})\times\{1\}$) and suppose for a contradiction that
  \begin{gather*}
    \lvert\cP\rvert < \Floor{\frac{(1-\epsilon')\cdot(r^\ell-1+\lambda)}{1+(1-\epsilon')/n}} + 1.
  \end{gather*}
  Note that we can remove or add at most $(r^{\ell+1}-r)/(r-1)$ vertices of each part to turn $\cP$ into an equipartition of
  $U_\varnothing$. Let $\cQ$ be the resulting equipartition of $U_\varnothing$ and note that
  \begin{gather*}
    \lvert\cQ\rvert\leq\lvert\cP\rvert \leq \frac{(1-\epsilon')\cdot(r^\ell-1+\lambda)}{1+(1-\epsilon')/n}.
  \end{gather*}
  If we can prove that all parts of $\cQ$ are $\epsilon'$-good in $G_{\ell,r}(n,\lambda)$, then this will contradict
  Theorem~\ref{thm:Gellr:equi}.

  But indeed, note that for every $Q\in\cQ$, we have
  \begin{gather*}
    \lvert Q\rvert
    \geq
    \Floor{\frac{(r^\ell-1+\lambda)\cdot n}{\lvert\cQ\rvert}}
    >
    \frac{(r^\ell-1+\lambda)\cdot n}{\lvert\cQ\rvert} - 1
    \geq
    \frac{n + (1-\epsilon')}{1-\epsilon'} - 1
    =
    \frac{n}{1-\epsilon'}.
  \end{gather*}

  Now since each part $Q\in\cQ$ is obtained from a part $P_Q$ of $\cP$, which is $\epsilon$-good in $G_{\ell,r}(n,\lambda)$, by
  either adding or removing at most $(r^{\ell+1}-r)/(r-1)$ vertices, by Corollary~\ref{cor:cont}\ref{cor:cont:good}, it follows
  that $Q$ is $\epsilon_{\symdiff}$-good in $G_{\ell,r}(n,\lambda)$, where
  \begin{gather*}
    \epsilon_{\symdiff}
    \df
    \epsilon + \frac{\lvert Q\symdiff P_Q\rvert}{\lvert Q\rvert}
    <
    \epsilon + \frac{(r^{\ell+1}-r)\cdot(1-\epsilon')}{(r-1)\cdot n}
    \leq
    \epsilon',
  \end{gather*}
  where the last inequality follows from~\eqref{eq:Gellr:equi:nlarge}.

  Thus, we get our desired contradiction. Therefore, we must have had
  \begin{align*}
    \lvert\cP\rvert
    & \geq
    \Floor{\frac{(1-\epsilon')\cdot(r^\ell-1+\lambda)}{1+(1-\epsilon')/n}} + 1
    =
    \bigl(1 + o_{n\to\infty,\epsilon'\to 0}(1)\bigr)\cdot(\epsilon')^{-\ell-1}
    \\
    & =
    \bigl(1 + o_{n\to\infty,\epsilon\to 0,\alpha,\ell}(1)\bigr)\cdot(1+\alpha)^{-\ell-1}\cdot\epsilon^{-\ell-1}
    \\
    & \leq
    \bigl(1 + o_{n\to\infty,\epsilon\to 0,\alpha\to 0,\ell}(1)\bigr)\cdot\epsilon^{-\ell-1}.
  \end{align*}
  Since $n$ needs to be made large depending on $\alpha$ and $\epsilon$, our final bound is
  \begin{gather*}
    \bigl(1 + o_{n\to\infty,\epsilon\to 0,\ell}(1)\bigr)\cdot\epsilon^{-\ell-1}.
    \qedhere
  \end{gather*}
\end{proof}

\begin{lemma}\label{lem:GellrVC}
  For $\ell,r\in\NN_+$, we have $\VC(G_{\ell,r})\leq 2$. In particular, for every $n\in\NN_+$ and $\lambda\in\RR_{\geq 0}$,
  we have $\VC(G_{\ell,r}(n,\lambda))\leq 2$.
\end{lemma}

\begin{proof}
  The second assertion follows from the first along with Lemma~\ref{lem:indepblowup}. For the first assertion, suppose for a
  contradiction that $\{u_1,u_2,u_3\}$ is a set of size $3$ that is shattered by $\cH_G\df\{N_{G_{\ell,r}}(w)\mid w\in
  V(G_{\ell,r})\}$. Since $G_{\ell,r}$ is bipartite, all three vertices $u_1$, $u_2$ and $u_3$ must belong to the same side of
  the bipartition, so we consider two cases.

  In the first case, we suppose $u_1,u_2,u_3\in([r]^{<\ell+1}\setminus\{\varnothing\})\times\{1\}$. Then we can write
  $u_i=(\sigma_i,1)$ ($i\in[3]$) and since there must exist an $w\in V(G_{\ell,r})$ that is adjacent to all three of the
  vertices, it follows that $\sigma_1$, $\sigma_2$ and $\sigma_3$ must have a common extension $\tau\in[r]^\ell$ (so that
  $(\tau,2)$ is the vertex that is adjacent to all three vertices). By possibly reordering $u_1$, $u_2$ and $u_3$, we may
  suppose that $\lvert\tau_1\rvert < \lvert\tau_2\rvert < \lvert\tau_2\rvert$, but then the neighborhoods of the vertices must
  satisfy $N_{G_{\ell,r}}(u_1)\supseteq N_{G_{\ell,r}}(u_2)\supseteq N_{G_{\ell,r}}(u_3)$, contradicting the fact that
  $\{u_1,u_2,u_3\}$ is shattered (since there is no vertex adjacent to $u_1$ that is not adjacent to $u_2$).

  In the second case, we suppose $u_1,u_2,u_3\in[r]^\ell\times\{2\}$. Then we can write $u_i=(\tau_i,2)$ ($i\in[3]$) and we can
  let $j_*\in[\ell]$ be the first position in which at least two of the $\tau_i$ differ. By possibly reordering $u_1$, $u_2$ and
  $u_3$, we may suppose that $(\tau_1)_{j_*}\neq(\tau_2)_{j_*}\neq(\tau_3)_{j_*}$ (but the first and third could be equal or
  not). We claim that there is no vertex $w$ that is adjacent to $u_1$ and $u_2$ and not adjacent to $u_3$ (which then
  contradicts the shattering of $\{u_1,u_2,u_3\}$). Suppose not, which means that the vertex $w$ must be of the form
  $(\sigma,1)\in([r]^{<\ell+1}\setminus\{\varnothing\})\times\{1\}$. Since $w$ is adjacent to both $u_1$ and $u_2$, we know that
  $\sigma$ must be a prefix of both $\tau_1$ and $\tau_2$, and since the first position in which these differ is $j_*$, it
  follows that $\lvert\sigma\rvert\leq j_*$. However, since $j_*$ is also the first position in which $\tau_3$ differs from
  $\tau_1$, it follows that $\sigma$ is also a prefix of $\tau_3$, which implies that $w$ is adjacent to $u_3$, a contradiction.
\end{proof}

\begin{remark}
  Corollary~\ref{cor:Gellr:equi} and Lemma~\ref{lem:GellrVC} together imply that the $\VC$-regularity lemma of
  Fox--Pach--Suk~\cite{FPS19} cannot be shown to yield $\epsilon$-good equipartitions. This is because the $\VC$-regularity
  lemma is guaranteed to yield at most $O_{n\to\infty,\epsilon\to 0,d}(\epsilon^{-2\cdot d-1})$ parts where $d\df\VC(G)$.
  However, Corollary~\ref{cor:Gellr:equi} along with Lemma~\ref{lem:GellrVC} provide a family of graphs of increasing sizes with
  $\VC$-dimension $2$ (which means that the $\VC$-regularity lemma yields at most $O_{n\to\infty,\epsilon\to
    0,d}(\epsilon^{-5})$ parts), Littlestone dimension $\ell$ and whose smallest $\epsilon$-good equipartition has size $(1 +
  o_{n\to\infty,\epsilon\to 0,\ell}(1))\cdot \epsilon^{-\ell-1}$, which is larger than the $\VC$-regularity bound provided
  $\ell\geq 5$.
\end{remark}


\section*{Acknowledgments}

The first author would like to thank Maryanthe Malliaris for some discussions about a high-level version of the statements of
the results of the manuscript, Yoshiharu Kohayakawa for pointing out the reference~\cite{AN08} containing
Proposition~\ref{prop:highchebyshev} and Gabriel Conant for pointing out the fact that goodness implies homogeneity
(Lemma~\ref{lem:good->hom}).

\section*{AI usage disclosure}

All results regarding algorithms and upper bounds, i.e., almost the entirety of the paper, is the product of humans. ChatGPT~5.5
Pro was used as a research assistant to help explore possible lower-bound constructions. ChatGPT~5.5 Pro was also used as an
editorial aid and as an additional verification aid during the preparation of the manuscript.

Regarding the lower bound constructions, the authors independently checked the arguments and rewrote, revised, simplified and
improved the exposition and citations. The authors take full responsibility for the correctness, exposition, and attribution in
the final manuscript.

\printbibliography

\appendix

\section{Pseudorandom generator of arbitrary rational bias}
\label{sec:PRGgen}

\begin{theorem}\label{thm:PRGgen}
  Let $n,h,t\in\NN_+$ with $h\leq n$ and let $(q_i,e_i)_{i=1}^t\in(\NN_+\times\NN_+)^t$ be such that for every $i\in[t]$, $q_i$
  is prime and $q_i^{e_i}\geq n$. Let also $Q\df\prod_{i\in[t]} q_i^{e_i}$, let $\widetilde{p}\in\{0,\ldots,Q\}$ and let
  $p\df\widetilde{p}/Q\in\QQ\cap(0,1]$. Then Algorithms~\ref{alg:encodePRGgen} and~\ref{alg:runPRGgen} together form a
  pseudorandom generator of bitlength $n$, independence $h$, bias $p$ and seed space $[Q]^h$ in the sense that:
  \begin{itemize}
  \item Algorithm~\ref{alg:encodePRGgen} outputs a tuple $P = (P_i)_{i=1}^t$, where $P_i\in\FF_{q_i}[X]$ is an irreducible monic
    polynomial of degree $e_i$.
  \item If $f(Z)_u$ is the output of Algorithm~\ref{alg:runPRGgen} when given the output $P$ of
    Algorithm~\ref{alg:encodePRGgen}, the index $u\in[n]$ and $Z\in[Q]^h$, then $f\colon[Q]^h\to\{0,1\}^n$ is a pseudorandom
    generator of bitlength $n$, independence $h$, bias $p$ and seed space $[Q]^h$ (see Definition~\ref{def:PRGbiased}).
  \end{itemize}

  The space-complexities of Algorithms~\ref{alg:encodePRGgen} and~\ref{alg:runPRGgen} are
  \begin{align*}
    O\bigl(\log(t) + \log(Q)\bigr),
    & &
    O\bigl(h\cdot\log(Q)\bigr).
  \end{align*}
  respectively, and the time-complexities are
  \begin{gather*}
    \sum_{i=1}^t O\Bigl(q_i^{2\cdot e_i}\cdot e_i\cdot\log(q_i)\cdot\log\bigl(\log(q_i)\bigr)\Bigr)
    \leq
    O\Bigl(Q^2\cdot\log(Q)\cdot\log\bigl(\log(Q)\bigr)\Bigr),
    \\
    \sum_{i=1}^t O\Bigl(h\cdot e_i^2\cdot\log(q_i)\cdot\log\bigl(\log(q_i)\bigr)\Bigr)
    \leq
    O\Bigl(h\cdot\bigl(\log(Q)\bigr)^2\cdot\log\bigl(\log(Q)\bigr)\Bigr).
  \end{gather*}
  respectively.
\end{theorem}

\begin{proof}
  We start by arguing correctness of Algorithm~\ref{alg:encodePRGgen} and for that, it suffices to argue correctness for each
  iteration of the loop of line~\ref{alg:encodePRGgen:For}.

  \begin{algorithm}[htbp]
    \caption{Computation algorithm of query-into-oracle model with space- and time-complexities
      \begin{gather*}
        O\bigl(\log(t) + \log(Q)\bigr),
        \\
        \sum_{i=1}^t O\Bigl(q_i^{2\cdot e_i}\cdot e_i\cdot\log(q_i)\cdot\log\bigl(\log(q_i)\bigr)\Bigr)
        \leq
        O\Bigl(Q^2\cdot\log(Q)\cdot\log\bigl(\log(Q)\bigr)\Bigr).
      \end{gather*}
      respectively, that returns $P=(P_i)_{i=1}^t$, where $P_i\in\FF_{q_i}[X]$ is an irreducible monic polynomial of degree
      $e_i$, such that when given to Algorithm~\ref{alg:runPRGgen} produces an oracle for a pseudorandom generator of bitlength
      $n$, independence $h$, bias $p$ and seed space $[Q]^h$, where $Q\df\prod_{i=1}^t q_i^{e_i}$. This particular algorithm
      does not need the values of $n$, $h$ or $p$.}
    \label{alg:encodePRGgen}
    \DontPrintSemicolon
    \KwIn{A tuple $(q_i,e_i)_{i=1}^t\in(\NN_+\times\NN_+)^t$ such that $q_i$ is prime for every $i\in[t]$.}
    \KwOut{A tuple $P=(P_i)_{i=1}^t$, where $P_i(X)\in\FF_{q_i}[X]$ is an irreducible monic polynomial of degree $e_i$.}
    \For{$i\in[t]$}{%
      \label{alg:encodePRGgen:For}
      \If{$e_i=1$}{%
        $P_i(X)\df X$\;
        \Continue
      }
      Let $d$ be the largest divisor of $e_i$ smaller than $e_i$.\;
      $m\assign (q_i^{e_i}-1)/(q_i^d-1)$\;
      Find the lexicographically first monic polynomial $P_i(X)\in\FF_{q_i}[X]$ of degree $e_i$ that divides
      \begin{gather*}
        \widetilde{P}(X) \df \sum_{j=0}^m X^j \in \FF_{q_i}[X].
      \end{gather*}
      \;\label{alg:encodePRGgen:Pi}
    }
    \Return{$(P_i)_{i=1}^t$}
  \end{algorithm}

  \begin{algorithm}[htbp]
    \caption{Oracle algorithm of query-into-oracle model with space- and time-complexities
      \begin{gather*}
        O\bigl(h\cdot\log(Q)\bigr),
        \\
        \sum_{i=1}^t O\Bigl(h\cdot e_i^2\cdot\log(q_i)\cdot\log\bigl(\log(q_i)\bigr)\Bigr)
        \leq
        O\Bigl(h\cdot\bigl(\log(Q)\bigr)^2\cdot\log\bigl(\log(Q)\bigr)\Bigr).
      \end{gather*}
      respectively, that when given $P$ computed via Algorithm~\ref{alg:encodePRGgen}, a seed $Z\in[Q]^h$ and an index $u\in[n]$,
      outputs a bit $X_u$ such that the corresponding function $Z\mapsto (X_u)_{u\in[n]}$ is a pseudorandom generator of
      bitlength $n$, independence $h$, bias $p$ and seed space $[Q]^h$.}
    \label{alg:runPRGgen}
    \DontPrintSemicolon
    \KwIn{A tuple $(q_i)_{i=1}^t$ of prime numbers, a tuple $P=(P_i)_{i=1}^t$, where $P_i(X)\in\FF_{q_i}[X]$ is an irreducible
      monic polynomial, an integer $\widetilde{p}\in\{0,\ldots,Q\}$, where $Q\df\prod_{i=1}^t q_i^{e_i}$ and $e_i\df\deg(P_i)$, numbers
      $h,n\in\NN_+$ such that $q_i^{e_i}\geq n\geq h$ for every $i\in[t]$, a seed $Z\in[Q]^h$, and an index $u\in[n]$.}
    \KwOut{A bit $X_u\in\{0,1\}$ such that when $Z$ is random with uniform distribution over $[Q]^h$, then $(X_u)_{u=1}^n$ is a
      collection of $p$-Bernoulli bits that is $h$-wise independent, where $p\df\widetilde{p}/Q$.}
    \lIf{$\widetilde{p}=0$}{\Return{$0$}}
    $M\assign\prod_{i=1}^t q_i^{e_i}$\tcp*[r]{$M$ initialized as $Q$}
    \lIf{$\widetilde{p}=M$}{\Return{$1$}}
    $V\assign 0$\tcp*[r]{Accumulated value in $\{0,\ldots,Q-1\}$}
    $W\assign (Z_j - 1)_{j=1}^h$\tcp*[r]{Remainder seed in $\{0,\ldots,Q-1\}^h$}
    \uFor{$i\in[t]$}{%
      \label{alg:runPRGgen:For}
      $s\assign (W_j\bmod q_i^{e_i})_{j=1}^h$\tcp*[r]{Seed in $\{0,\ldots,q_i^{e_i}-1\}$}
      $W\assign (\floor{W_j/q_i^{e_i}})_{j=1}^h$\;
      Write each $s_j\in\{0,\ldots,q_i^{e_i}-1\}$ in base $q_i$ as
      \begin{gather*}
        s_j = c^j_{e_i-1}c^j_{e_i-2}\cdots c^j_0.
      \end{gather*}
      \;
      Write $u-1\in\{0,\ldots,n-1\}\subseteq\{0,\ldots,q_i^{e_i}-1\}$ in base $q_i$ as
      \begin{gather*}
        u-1 = b_{e_i-1}b_{e_i-2}\cdots b_0.
      \end{gather*}
      \;
      \setcounter{algosplit}{\theAlgoLine}
    }
  \end{algorithm}

  \begin{algorithm}[htbp]
    \caption*{(continued).}
    \DontPrintSemicolon
    \setcounter{AlgoLine}{\thealgosplit}
    \let\oldnl\nl
    \let\nl\relax
    \vspace{-\baselineskip}\InvisibleBegin{
      \global\let\nl\oldnl 
      %
      Interpreting the numbers $c^j_\ell$ and $b_\ell$ as elements of $\FF_{q_i}$, compute
      \begin{gather*}
        R(X)
        \df
        \left(
        \sum_{j\in[h]}
        \left(
        \left(\sum_{\ell=0}^{e_i-1} c^j_\ell\cdot X^\ell\right)
        \cdot\left(\sum_{\ell=0}^{e_i-1} b_\ell\cdot X^\ell\right)^{j-1}
        \right)
        \right)
        \bmod P_i(X).
      \end{gather*}
      \;\label{alg:runPRGgen:poly}
      Writing $R(X) = \sum_{\ell=0}^{e_i-1} r_\ell\cdot X^\ell$ and interpreting the coefficients $r_\ell\in\FF_{q_i}$ as
      elements of $\{0,\ldots,q_i-1\}$, let
      \begin{gather*}
        v\df\sum_{\ell=0}^{e_i-1} r_\ell\cdot q_i^\ell\in\{0,\ldots,q_i^{e_i}-1\}.
      \end{gather*}
      \;
      $M\assign M/q_i^{e_i}$\;
      $V\assign V + M\cdot v$\;
      \lIf{$V>\widetilde{p}$}{\label{alg:runPRGgen:return0}\Return{$0$}}
    }
    \Return{$1$}
  \end{algorithm}

  It is clear that if $e_i=1$ then $P_i(X)=X\in\FF_{q_n}[X]$ is an irreducible monic polynomial of degree $e_i=1$. Suppose then
  that $e_i\geq 2$. We know that for any sought irreducible monic polynomial $P_i(X)\in\FF_{q_i}[X]$ of degree $e_i$, the
  quotient $\FF_{q_i}[X]/(P_i)$ is isomorphic to the field $\FF_{q_i^{e_i}}$ of size $q_i^{e_i}$.

  In turn, we know that in the algebraic closure $\overline{\FF}_{q_i}$ of $\FF_{q_i}$, the field $\FF_{q_i^{e_i}}$ is exactly
  given by the $q_i^{e_i}$ distinct roots of the polynomial $R_{e_i}(X)\df X^{q_i^{e_i}} - X$ (note that the derivative of $R$ is
  $-1$, so all its roots have multiplicity $1$). Furthermore, the field $\FF_{q_i^{e_i}}$ contains as subfields exactly the
  fields of the form $\FF_{q_i^f}$, where $f\in\NN_+$ is a divisor of $e_i$; in fact, $\FF_{q_i^f}$ consists exactly of the $q_i^f$
  roots of $R_f(X)\df X^{q_i^f} - X$, which are themselves roots of $R_{e_i}(X)$. Now, the algorithm takes $d\in\NN_+$ as the
  largest divisor of $e_i$ smaller than $e_i$ and since
  \begin{gather*}
    R_{e_i}(X)
    =
    R_f(X)\cdot \widetilde{P}(X),
  \end{gather*}
  where
  \begin{align*}
    \widetilde{P}(X)
    & \df
    \sum_{j=0}^m X^j,
    &
    m
    & \df
    \frac{q_i^{e_i}-1}{q_i^d-1}
    =
    \sum_{j=0}^{e_i/d-1} q_i^{j\cdot d},
  \end{align*}
  it follows that the roots of $\widetilde{P}$ in the algebraic closure $\overline{\FF}_{q_i}$ are exactly the elements of
  $\FF_{q_i^{e_i}}\setminus\FF_{q_i^d}$ and since $\FF_{q_i^d}$ is the largest proper subfield\footnote{This is because the
  subfields of $\FF_{q_i^{e_i}}$ are exactly those of the form $\FF_{q_i^m}$ with $e_i$ divisible by $m$; one direction is
  because $\FF_{q_i^{e_i}}$ must be an $\FF_{q_i^m}$-vector space so $q_i^{e_i}$ has to be a power of $q_i^m$, hence $m$ must
  divide $e_i$; and the other direction is because when $m$ divides $e_i$, we can factor $X^{q_i^{e_i}}-X =
  (X^{2^m}-X)\cdot\sum_{j=0}^{(q_i^{e_i}-1)/(q_i^d-1)} X^j$.} of $\FF_{q_i^{e_i}}$, all irreducible factors of $\widetilde{P}$
  over $\FF_{q_i}[X]$ must be of degree $e_i$, which means that any monic polynomial in $\FF_{q_i}[X]$ of degree $e_i$ that
  divides $\widetilde{P}$ must be irreducible. Thus, Algorithm~\ref{alg:encodePRGgen} eventually finds one such monic
  irreducible $P_i(X)\in\FF_{q_i}[X]$ of degree $e_i$.

  \smallskip

  We now analyze the space-complexity of Algorithm~\ref{alg:encodePRGgen}. First, we need to store the index $i\in[t]$ and in
  each iteration of the loop of line~\ref{alg:encodePRGgen:For}, we need $O(\log(e_i+1))$ space to compute $d$ and we need at
  most $O(e_i\cdot\log(q_i))$ space to compute $m$.

  Regarding $P_i$ and $\widetilde{P}$, for the former, we need space at most $O(e_i\cdot\log(q_i))$ to store the current $P_i$
  being tested. For $\widetilde{P}$, we don't actually need to store it in memory: we know that all of its coefficients are $1$,
  so we only need to keep track of its degree $m$. Furthermore, when testing divisibility of $\widetilde{P}$ by $P_i$, when
  performing long division, we only need to keep track of the topmost $e_i$ coefficients of the remainder as the standard long
  division algorithm progresses (again, because all other coefficients of the remainder will be the ones of $\widetilde{P}$,
  which are all $1$). Thus, the space-complexity of line~\eqref{alg:encodePRGgen:Pi} is at most $O(e_i\cdot\log(q_i))$. Thus,
  the total space-complexity of Algorithm~\ref{alg:encodePRGgen} is at most
  \begin{gather*}
    O\bigl(\log(t+1)\bigr) + \sum_{i=1}^t O\bigl(e_i\cdot\log(q_i)\bigr) = O\bigl(\log(t) + \log(Q)\bigr)
  \end{gather*}
  (recalling that $Q\df\prod_{i\in[t]} q_i^{e_i}$).

  \smallskip

  For time-complexity of Algorithm~\ref{alg:encodePRGgen}, note that to compute the largest divisor $d$ that is less than $e_i$,
  we can instead compute the smallest non-trivial divisor $d'\in\NN$ since $d = e_i/d'$ (this still only takes space
  $O(\log(e_i+1))$). We only need to search for $d'$ in $\{2,\ldots,\floor{\sqrt{e_i}}\}$, so the time-complexity of computing
  $d$ is at most
  \begin{gather*}
    O\Bigl(\sqrt{e_i}\cdot\bigl(\log(e_i+1)\bigr)^2\Bigr)
    \leq
    O(e_i).
  \end{gather*}
  For the computation of $m$, the time-complexity is easily seen to be at most $O(e_i^2\cdot(\log(q_i))^2)$. Finally, the
  time-complexity of line~\ref{alg:encodePRGgen:Pi} amounts to the time it takes to enumerate all the potential $P_i$ and the
  time it takes to test divisibility of $\widetilde{P}$ by it. There are a total of $q_i^{e_i}$ many $P_i$ and to check
  divisibility of $\widetilde{P}$ by one such $P_i$, we can perform long polynomial division and since $P_i$ is monic, this
  never involves computing inverses in $\FF_{q_i}$, which means that the polynomial division takes time at most the product of
  the degrees times the time it takes to multiply numbers in $\FF_{q_i}$ (the time-complexity of addition in $\FF_{q_i}$ is
  completely dominated by that), which with the Fast-Fourier-Transform Multiplication Algorithm is
  $O(\log(q_i)\cdot\log(\log(q_i)))$, that is, the time-complexity of checking divisibility of $\widetilde{P}$ by $P_i$ is at
  most
  \begin{gather*}
    O\bigl(
    m\cdot e_i\cdot\log(q_i)\cdot\log\bigl(\log(q_i)\bigr)
    \bigr)
    \leq
    O\Bigl(
    q_i^{e_i}\cdot e_i\cdot\log(q_i)\cdot\log\bigl(\log(q_i)\bigr)
    \Bigr).
  \end{gather*}

  Thus, the total time-complexity of Algorithm~\ref{alg:encodePRGgen} is at most
  \begin{gather*}
    \sum_{i=1}^t O\Bigl(q_i^{2\cdot e_i}\cdot e_i\cdot\log(q_i)\cdot\log\bigl(\log(q_i)\bigr)\Bigr)
    \leq
    O\Bigl(Q^2\cdot\log(Q)\cdot\log\bigl(\log(Q)\bigr)\Bigr).
  \end{gather*}

  \medskip

  We now argue correctness of Algorithm~\ref{alg:runPRGgen}. The result is completely obvious if $\widetilde{p}=0$ (in which
  case $p=0$) or if $\widetilde{p}=Q$ (in which case $p=1$), so we assume that $\widetilde{p}\in[Q-1]$. We fix $P$, $h$ and $n$
  of the input of the algorithm, we let $\rn{Z}$ be picked uniformly at random in $[Q]^h$ and we make the following definitions:
  \begin{itemize}
  \item $\rn{X}_u$ as the bit that Algorithm~\ref{alg:runPRGgen} outputs when further given $\rn{Z}$ and $u\in[n]$.
  \item For a variable of Algorithm~\ref{alg:runPRGgen}, we will write its name followed by $[i]$ to mean its value at the end
    of the $i$th iteration of the loop of line~\ref{alg:runPRGgen:For} and we will write its name followed by $[0]$ to mean its
    value right as the loop of line~\ref{alg:runPRGgen:For} starts. We will write it in bold face if the variable depends on
    $\rn{Z}$ (which means it will be a random variable that is $\rn{Z}$-measurable) and we will write a subscript $u$ if it
    depends on the value of $u$. Under these conventions, it is easy to check the following dependencies:
    \begin{align*}
      M[i], & &
      \rn{V}[i]_u, & &
      \rn{W}[i], & &
      \rn{s}[i], & &
      \rn{c}^j_\ell[i], & &
      b_\ell[i]_u, & &
      \rn{R}(X)[i]_u, & &
      \rn{v}[i]_u.
    \end{align*}
  \item If the algorithm terminates before all $t$ iterations of the loop of line~\ref{alg:runPRGgen:For} executed (due to the
    return statement of line~\ref{alg:runPRGgen:return0}), we define the values of later iterations as if the algorithm would
    have continued executing for them (this will be convenient for the proof).
  \end{itemize}

  We will now list several claims, culminating in the final claim of correctness, and provide a short proof of these claims
  later:
  \begin{enumerate}[wide, label={\textbf{Claim~\thetheorem.\arabic*.}}, ref={\thetheorem.\arabic*}]
  \item\label{clm:PRGgen:W} $\rn{W}[i]$ is uniformly distributed in $\{0,\ldots,\prod_{j=i+1}^t q_j^{e_j}-1\}^h$.
  \item\label{clm:PRGgen:s} For $i > 0$, $\rn{s}[i]$ is uniformly distributed in $\{0,\ldots,q_i^{e_i}-1\}$. Furthermore,
    $(\rn{s}[i])_{i\in[t]}$ is mutually independent.
  \item\label{clm:PRGgen:c} For $i > 0$, $(\rn{c}^j_\ell[i])_{(j,\ell)\in[h]\times\{0,\ldots,e_i-1\}}$ is uniformly distributed
    in the set $\{0,\ldots,q_i-1\}^{[h]\times\{0,\ldots,e_i-1\}}$. Furthermore,
    $((\rn{c}^j_\ell[i])_{(j,\ell)\in[h]\times\{0,\ldots,e_i-1\}})_{i\in[t]}$ is mutually independent in the index $i$.
  \item\label{clm:PRGgen:wR} For $i > 0$, if we define the (random) polynomial
    \begin{gather*}
      \rn{\widetilde{R}}_i(Y)\in \frac{\FF_{q_i}[X]}{(P_i)}[Y]
    \end{gather*}
    as the polynomial in the variable $Y$ with coefficients in the field $\FF_{q_i}[X]/(P_i)\cong\FF_{q_i^{e_i}}$ given by
    \begin{gather*}
      \rn{\widetilde{R}}_i(Y)
      \df
      \sum_{j\in[h]}
      \left(
      \left(\sum_{\ell=0}^{e_i-1} \rn{c}^j_\ell[i]\cdot X^\ell\right)\bmod P_i(X)
      \right)\cdot Y^{j-1},
    \end{gather*}
    then $\rn{\widetilde{R}}_i(Y)$ is a uniformly at random polynomial in $(\FF_{q_i}[X]/(P_i))[Y]$ of degree at most $h-1$.
    Furthermore, $(\rn{\widetilde{R}}_i)_{i\in[t]}$ is mutually independent.
  \item\label{clm:PRGgen:R} For $i > 0$ and $\rn{\widetilde{R}}_i(Y)$ as in Claim~\ref{clm:PRGgen:wR}, we have
    \begin{gather*}
      \rn{R}(X)[i]_u
      =
      \rn{\widetilde{R}}_i\left(\sum_{\ell=0}^{e_i-1} b_\ell[i]_u\cdot X^\ell \bmod P_i(X)\right).
    \end{gather*}
  \item\label{clm:PRGgen:Runif} For $i > 0$, $(\rn{R}(X)[i]_u)_{u\in[n]}$ is a sequence of $h$-wise independent random elements
    of $\FF_{q_i}[X]/(P_i)\cong\FF_{q_i^{e_i}}$ such that each $\rn{R}(X)[i]_u$ is uniformly distributed in
    $\FF_{q_i}[X]/(P_i)$. Furthermore, $((\rn{R}(X)[i]_u)_{u\in[n]})_{i\in[n]}$ is mutually independent in the index $i$.
  \item\label{clm:PRGgen:v} For $i > 0$, $(\rn{v}[i]_u)_{u\in[n]}$ is a sequence of $h$-wise independent random elements of
    $\{0,\ldots,q_i^{e_i}-1\}$ such that each $\rn{v}[i]_u$ is uniformly distributed in $\{0,\ldots,q_i^{e_i}-1\}$. Furthermore,
    $((\rn{v}[i]_u)_{u\in[n]})_{i\in[t]}$ is mutually independent in the index $i$.
  \item\label{clm:PRGgen:M} $M[i] = \prod_{j=i+1}^t q_j^{e_j}$.
  \item\label{clm:PRGgen:V} $(\rn{V}[i]_u)_{u\in[n]}$ is a sequence of $h$-wise independent random variables such that each
    $\rn{V}[i]_u$ is uniformly distributed in
    \begin{gather*}
      \{0, M[i], 2\cdot M[i],\ldots, Q - M[i]\}.
    \end{gather*}
  \item\label{clm:PRGgen:X} Algorithm~\ref{alg:runPRGgen} is correct, i.e., $(\rn{X}_u)_{u\in[n]}$ is a sequence of $h$-wise
    independent random variables such that each $\rn{X}_u$ is $p$-Bernoulli.
  \end{enumerate}

  Claims~\ref{clm:PRGgen:W} and~\ref{clm:PRGgen:s} follow straightforwardly by induction in $i\in\{0,\ldots,t\}$.
  Claim~\ref{clm:PRGgen:c} is an immediate corollary of Claim~\ref{clm:PRGgen:s} and the way Algorithm~\ref{alg:runPRGgen}
  computes $\rn{c}^j_\ell[i]$ from $\rn{s}[i]$. In turn, Claim~\ref{clm:PRGgen:wR} follows directly from
  Claim~\ref{clm:PRGgen:c} and the definition of $\rn{\widetilde{R}}_i$.

  Claim~\ref{clm:PRGgen:R} follows since (polynomial) modulus satisfies
  \begin{align*}
    (A+B)\bmod P_i & = (A\bmod P_i + B\bmod P_i)\bmod P_i,
    \\
    (A\cdot B)\bmod P_i & = \bigl((A\bmod P_i)\cdot (B\bmod P_i)\bigr)\bmod P_i,
  \end{align*}

  Claim~\ref{clm:PRGgen:Runif} follows from Claims~\ref{clm:PRGgen:wR} and~\ref{clm:PRGgen:R} along with
  Proposition~\ref{prop:PRG}: Claim~\ref{clm:PRGgen:wR} says the $\rn{\widetilde{R}_i}(Y)$ are uniformly at random
  polynomials of degree at most $h-1$ with coefficients in the field $\FF_{q_i}[X]/(P_i)\cong\FF_{q_i^{e_i}}$ and
  Claim~\ref{clm:PRGgen:R} says that $\rn{R}(X)[i]_u$ is the evaluation of this polynomial at the point
  \begin{gather*}
    \sum_{\ell=0}^{e_i-1} b_\ell[i]_u\cdot X^\ell \bmod P_i(X)\in\frac{\FF_{q_i}[X]}{(P_i)}\cong\FF_{q_i^{e_i}}
  \end{gather*}
  and since different values of $u$ lead to different values of $(b_\ell[i]_u)_{\ell\in\{0,\ldots,e_i-1\}}$, changing $u$
  changes the point of the field $\FF_{q_i}[X]/(P_i)\cong\FF_{q_i^{e_i}}$ on which we are evaluating the polynomial, so
  Proposition~\ref{prop:PRG} gives us the $h$-wise independence of $(\rn{R}(X)[i])_{u\in[n]}$ and the fact that each
  $\rn{R}(X)[i]_u$ is uniformly distributed in $\FF_{q_i}[X]/(P_i)\cong\FF_{q_i^{e_i}}$. The mutual independence in the index
  $i\in[t]$ is inherited from Claim~\ref{clm:PRGgen:wR} via Claim~\ref{clm:PRGgen:R}.

  Claim~\ref{clm:PRGgen:v} follows directly from Claim~\ref{clm:PRGgen:Runif} and the definition of the $\rn{v}[i]_u$.

  Claim~\ref{clm:PRGgen:M} follows with a simple induction in $i\in\{0,\ldots,t\}$.

  We prove Claim~\ref{clm:PRGgen:V} by induction in $i\in\{0,\ldots,t\}$. For $i=0$, this is trivial since $M[0]=Q$ and
  $\rn{V}[0] = 0$. For $i > 0$, we know that
  \begin{gather*}
    \rn{V}[i]_u = \rn{V}[i-1] + M[i]\cdot\rn{v}[i].
  \end{gather*}
  Since $\rn{V}[i-1]_u$ depends only on the values of $\rn{v}[1]_u,\ldots,\rn{v}[i-1]_u$, $\rn{v}[i]_u$ is independent from
  these and uniformly distributed in $\{0,\ldots,q_i^{e_i}-1\}$ (and $h$-wise independent in the index $u$), and $M[i]\cdot
  q_i^{e_i} = M[i-1]$, using the inductive hypothesis, it follows that $\rn{V}[i]_u$ is uniformly distributed in
  \begin{gather*}
    \{0, M[i], 2\cdot M[i],\ldots, Q-M[i]\}
  \end{gather*}
  and is $h$-wise independent in the index $u$.

  Finally, for Claim~\ref{clm:PRGgen:X}, we note that $\rn{X}_u=1$ if and only if there exists $i\in[t]$ such that
  $\rn{V}[i]_u>\widetilde{p}$. Since clearly $\rn{V}[1]_u\leq\rn{V}[2]_u\leq\cdots\leq\rn{V}[t]_u$, we conclude that $\rn{X}_u =
  \One[\rn{V}[t]_u > \widetilde{p}]$. By Claim~\ref{clm:PRGgen:V}, we know that $(\rn{V}[t]_u)_{u\in[n]}$ is $h$-wise
  independent in the index $u$ with each $\rn{V}[t]_u$ uniformly distributed in
  \begin{gather*}
    \{0, M[t], 2\cdot M[t],\ldots, Q-M[t]\} = \{0,1,2,\ldots,Q-1\},
  \end{gather*}
  from which it follows that
  \begin{gather*}
    \PP[\rn{X}_u] = \frac{\widetilde{p}}{Q} = p,
  \end{gather*}
  as desired.

  \smallskip

  We now analyze the space-complexity of Algorithm~\ref{alg:runPRGgen}. Storing the variables
  \begin{align*}
    M, & &
    V, & &
    W, & &
    s,
    \\
    \bigl((c^j_\ell)_{\ell=0}^{e_i-1}\bigr)_{j\in[h]}, & &
    (b_\ell)_{\ell=0}^{e_i-1}, & &
    R, & &
    v
  \end{align*}
  takes space
  \begin{align*}
    O\bigl(\log(Q)\bigr), & &
    O\bigl(\log(Q)\bigr),
    \\
    O\bigl(h\cdot\log(Q)\bigr), & &
    \max_{i=1}^t O\bigl(e_i\cdot\log(q_i)\bigr)\leq O\bigl(\log(Q)\bigr),
    \\
    \max_{i=1}^t O\bigl(h\cdot e_i\cdot\log(q_i)\bigr)\leq O\bigl(h\cdot\log(Q)\bigr), & &
    \max_{i=1}^t O\bigl(e_i\cdot\log(q_i)\bigr)\leq O\bigl(\log(Q)\bigr),
    \\
    \max_{i=1}^t O\bigl(e_i\cdot\log(q_i)\bigr)\leq O\bigl(\log(Q)\bigr), & &
    \max_{i=1}^t O\bigl(e_i\cdot\log(q_i)\bigr)\leq O\bigl(\log(Q)\bigr),
  \end{align*}
  respectively. Furthermore, it is straightforward to check that all computations require extra space at most proportional to
  the size of the variables being computed (see more details about line~\ref{alg:runPRGgen:poly} below), so the space-complexity
  of Algorithm~\ref{alg:runPRGgen} is at most
  \begin{gather*}
    O\bigl(h\cdot\log(Q)\bigr).
  \end{gather*}

  \smallskip

  Finally, we analyze the time-complexity of Algorithm~\ref{alg:runPRGgen}. It is straightforward to check that the most
  time-expensive step is the computation of line~\ref{alg:runPRGgen:poly}. This can be done efficiently by inductively computing
  for each $\widetilde{h}\in\{0,\ldots,h\}$ the values:
  \begin{gather*}
    \begin{aligned}
      A_{-1}(X) & \df 0, &
      B_0(X) & \df 1,
    \end{aligned}
    \\
    B_{\widetilde{h}}(X)
    \df
    \left(\sum_{\ell=0}^{e_i-1} b_\ell\cdot X^\ell\right)\cdot B_{\widetilde{h}-1}(X) \bmod P_i(X),
    \\
    A_{\widetilde{h}}(X)
    \df
    \left(
    A_{\widetilde{h}-1}(X) + \left(\sum_{\ell=0}^{e_i-1} c^j_\ell\cdot X^\ell\right)\cdot B_{\widetilde{h}}(X)
    \right) \bmod P_i(X).
  \end{gather*}
  At each step we need to compute two polynomial moduli by $P_i(X)$ of polynomials in $\FF_{q_i}[X]$ of degree at most $2\cdot
  e_i$, two polynomial multiplications in $\FF_{q_i}[X]$ of polynomials of degree at most $e_i$ each and a polynomial addition
  in $\FF_{q_i}[X]$ of two polynomials of degree at most $2\cdot e_i$. Using the Fast-Fourier-Multiplication Transform for
  multiplication over $\FF_{q_i}$ (and recalling that $P_i(X)$ is monic, so long polynomial division by it does \emph{not}
  involve computing inverses in $\FF_{q_i}$), the time-complexity of line~\ref{alg:runPRGgen:poly} is at most
  \begin{gather*}
    O\Bigl(h\cdot e_i^2\cdot\log(q_i)\cdot\log\bigl(\log(q_i)\bigr)\Bigr).
  \end{gather*}

  Thus, the time-complexity of Algorithm~\ref{alg:runPRGgen} is at most
  \begin{gather*}
    \sum_{i=1}^t O\Bigl(h\cdot e_i^2\cdot\log(q_i)\cdot\log\bigl(\log(q_i)\bigr)\Bigr)
    \leq
    O\Bigl(h\cdot\bigl(\log(Q)\bigr)^2\cdot\log\bigl(\log(Q)\bigr)\Bigr).
    \qedhere
  \end{gather*}
\end{proof}


\end{document}